\documentclass{amsart}
 \usepackage[top=1in, bottom=1in, left=1.375in, right=1.375in]{geometry}
 \usepackage{amsmath,amssymb}
 \usepackage{algpseudocode}
 \usepackage{algorithm}
 \usepackage{algorithmicx}
 \usepackage{caption}
 \usepackage{subcaption}
 \usepackage{amsthm}
 \numberwithin{equation}{section}
 \usepackage{amsfonts}
 \usepackage{graphicx}
\usepackage{hyperref}

\usepackage{amsmath}
\usepackage{amssymb}
\usepackage{amsthm}
\usepackage{amsmath}
\usepackage{setspace}
\usepackage{xcolor}
\usepackage{fancyhdr}
\usepackage{amssymb}
\usepackage{amsthm}
\usepackage{listings}
    \usepackage{cite}
    \usepackage{tabularx}
    \usepackage{graphicx}
    \usepackage{epstopdf}    
    \usepackage{epsfig}
    \usepackage{float}
    \usepackage{ listings}
    \usepackage{appendix}
    \usepackage{subcaption}
 \theoremstyle{plain}
  \newtheorem{thm}{Theorem}
 \newtheorem{prop}{Proposition}[section]
 \newtheorem{lem}[prop]{Lemma}
 
 \theoremstyle{definition}

 \theoremstyle{remark}
 \newtheorem{remark}[prop]{Remark}

 \let\pa=\partial
 \let\al=\alpha
 \let\b=\beta
 \let\d=\delta
 \let\g=\gamma
 \let\e=\varepsilon
  
 \let \kp = \kappa
 \let\lam=\lambda

 \let\f=\frac
 \let\inf = \infty
 \let \les = \lesssim
 \let\om=\omega
 
 \let \th = \theta
 \let \vp = \varphi
 
\let\B = \Big
 \let\D=\Delta
 \let\Lam=\Lambda
 
 \let\Om=\Omega
 \let\td = \tilde

 \let\teq \triangleq
 
 \let\pa=\partial

 \def\cK{{\mathcal K}}
 \def\cL{{\mathcal L}}
 \def\cM{{\mathcal M}}

 \def\cR{{\mathcal R}}
 
 \def\cT{{\mathcal T}}

 \def\cM{{\mathcal M}}

 \def\na{\nabla}
 \def\la{\langle}
 \def\ra{\rangle}
\def\lt{\left}
\def\rt{\right}
\def\one{\mathbf{1}}

 \newcommand{\beq}{\begin{equation}}
 \newcommand{\eeq}{\end{equation}}
  \newcommand{\bal}{\begin{aligned} }
  \newcommand{\eal}{\end{aligned}}
 \newcommand{\ben}{\begin{eqnarray}}
 \newcommand{\een}{\end{eqnarray}}
 \newcommand{\beno}{\begin{eqnarray*}}
 \newcommand{\eeno}{\end{eqnarray*}}

 \newcommand{\ee}{\mathbf{e}}

 \newcommand{\uu}{\mathbf{u}}

 \newcommand{\R}{\mathbb{R}}

 \newcommand{\tr}{\mathrm{Tr}}
\newcommand{\sgn}{\mathrm{sgn}}

 \author{Jiajie Chen, Thomas Y. Hou, De Huang}
 \address{Applied and Computational Mathematics, Caltech, Pasadena, CA 91125. Emails: jchen@caltech.edu, hou@cms.caltech.edu, dhuang@caltech.edu}
  \date{\today}

\title[Self-similar blowup of the HL model]{Asymptotically self-similar blowup of the Hou-Luo model for the 3D Euler equations}
 \begin{document}
 \begin{abstract}
Inspired by the numerical evidence of a potential 3D Euler singularity \cite{luo2014potentially,luo2013potentially-2}, we prove finite time singularity from smooth initial data for the HL model introduced by Hou-Luo in \cite{luo2014potentially,luo2013potentially-2} for the 3D Euler equations with boundary. Our finite time blowup solution for the HL model and the singular solution considered in \cite{luo2014potentially,luo2013potentially-2} share some essential features, including similar blowup exponents, symmetry properties of the solution, and the sign of the solution. 
 We use a dynamical rescaling formulation and the strategy proposed in our recent work in \cite{chen2019finite} to establish the nonlinear stability of an approximate self-similar profile.
 The nonlinear stability enables us to prove that the solution of the HL model with smooth initial data and finite energy will develop a focusing asymptotically self-similar singularity in finite time. 
 Moreover the self-similar profile is unique within a small energy ball and the 
$C^\gamma$ norm of the density $\th$ with $\gamma\approx 1/3$ is uniformly bounded  up to the singularity time. 
  \end{abstract}
 \maketitle

\vspace{-0.2in}

\section{ Introduction}

The three-dimensional (3D) incompressible Euler equations are one of the most fundamental equations in fluid dynamics. Despite their wide range of applications, the global well-posedness of the 3D incompressible Euler equations is one of the most outstanding open questions in the theory of nonlinear partial differential equations. The interested readers may consult the excellent surveys \cite{constantin2007euler,gibbon2008three,hou2009blow,kiselev2018,
majda2002vorticity} and the references therein. 
The difficulty associated with the global regularity of the 3D Euler equations can be described by the vorticity equation:
\begin{equation}
  \omega_{t} + u \cdot \nabla \omega = \omega \cdot \nabla u,
  \label{eqn_eu_w}
\end{equation}
where $\omega = \nabla \times u$ is the \emph{vorticity vector} of the fluid, and $u$ is related to $\omega$ via the \emph{Biot-Savart law}. Formally, $\nabla u$ has the same scaling as $\omega$, which implies that the vortex stretching term $\omega \cdot \nabla u$ formally scales  like $\omega^{2}$. However, $\nabla u$ is related to $\omega$ through the Riesz transform. Various previous studies indicate that the nonlocal nature of the vortex stretching term and the local geometric regularity of the vorticity vector may lead to dynamic depletion of the nonlinear vortex stretching (see e.g. \cite{constantin1996geometric,deng2005geometric,hou2006dynamic}), which may prevent singularity formation in finite time. 

In \cite{luo2014potentially,luo2013potentially-2}, Luo and Hou investigated the 3D axisymmetric Euler equations with a solid boundary and presented some convincing numerical evidence that the 3D Euler equations develop a potential finite time singularity. They considered a class of smooth initial data with finite energy that satisfy certain symmetry properties. The potential singularity occurs at a stagnation point of the flow along the boundary. 
The presence of the boundary and the hyperbolic flow structure near the singularity play an important role in the singularity formation. 
To understand the mechanism for this potential 3D Euler singularity, Hou and Luo \cite{luo2013potentially-2} proposed the following one-dimensional model along the boundary at $r=1$: 
\beq\label{eq:HL}
\bal
\om_t + u \om_x &=  \th_x , \\
\th_t + u \th_{x} &= 0,  \quad  u_x = H \omega.
\eal
\eeq
Here $u=u^z$, $\om = \om^\phi$, and $\theta = (u^\phi)^2$, with $ u^\phi$ and $\om^\phi$ being the angular velocity and angular vorticity, respectively. 
Numerical study presented in \cite{luo2013potentially-2} shows that the HL model develops a finite time singularity from smooth initial data with blowup scaling properties surprisingly similar to those observed for the 3D Euler equations. By exploiting the symmetry properties of the solution and some monotonicity property of the velocity kernel, Choi et al 
have been able to prove that the HL model develops a finite time singularity in \cite{choi2014on} 
using a Lyapunov functional argument. Part of our analysis to be presented is inspired by the sign property of a quadratic interaction term between $u$ and $\omega$ obtained in \cite{choi2014on}.
However, there seems to be some essential difficulties in extending the method 
in \cite{choi2014on} to the 3D Euler equations.

There has been a number of subsequent developments inspired by the singularity scenario reported in\cite{luo2014potentially,luo2013potentially-2}, see e.g. \cite{kiselev2013small,choi2015finite,choi2014on,kryz2016} and the excellent survey article \cite{kiselev2018}.
Although various simplified models have been proposed to study the singularity scenario reported in\cite{luo2014potentially,luo2013potentially-2}, currently there is no rigorous proof of the Luo-Hou blowup scenario with smooth data.  Recently, Elgindi \cite{elgindi2019finite} (see also \cite{elgindi2019stability}) proved an important result that the 3D axisymmetric Euler equations without swirl can develop a finite time singularity for $C^{1,\alpha}$ initial velocity. In a setting similar to the Luo-Hou scenario, singularity formation of the 2D Boussinesq and the 3D axisymmetric Euler equations with $C^{1,\alpha}$ velocity and boundary has been established by  the first two authors \cite{chen2019finite2}.

The main result of this paper is stated by the informal theorem below. The more precise 
and stronger statement will be given by Theorem \ref{thm2} in Section \ref{sec:dyn}. 

\begin{thm}\label{thm1}

There is a family of initial data $(\th_0, \om_0)$ with $\th_{0, x}, \om_0 \in C_c^{\inf}$, such that the solution of the HL model \eqref{eq:HL} will develop a focusing asymptotically self-similar singularity in finite time.
 The self-similar blowup profile $(\th_{\inf},\om_{\inf})$ is unique within a small energy ball and its associated scaling exponents $c_{l,\inf}, c_{\om,\inf}$ satisfy $|\lam - 2.99870|\leq 6 \cdot 10^{-5}$
with $\lam = c_{l, \inf}|c_{\om,\inf}|^{-1}$.
  Moreover, the $C^{\g}$ norm of $\th$ is uniformly bounded up to the blowup time $T$, and the $C^\beta$ norm of $\th$ blows up at $T$ for any $ \b \in (\g , 1]$ with $ \g = \f{ \lam - 2 }{ \lam }$.
 \end{thm}

 Using the self-similar profile $(\th_{\infty},\om_{\inf})$, we can construct the self-similar blowup solution 
\beq\label{eq:self_similar}
\om_*(x, t) = \f{1}{ (1- t )|c_{\om,\inf}| }   \om_{\inf}( \f{x}{  (1-t)^{ \lam } }  ), \quad 
\th_*(x, t) = \f{ 1 }{ (1- t)^{ 2- \lam} |c_{\om,\inf}| } \th_{\inf}( \f{x}{(1-t)^{ \lam } } ) , 
\eeq
that blows up at $T=1$.
The blowup exponent $ \lam \approx 2.99870$ in the HL model is surprisingly close to the blowup exponent $ \lam \approx 2.9215$ of the 3D Euler equations considered by Luo-Hou \cite{luo2014potentially,luo2013potentially-2}. 
An important property that characterizes the stable nature of the blowup in
the HL model is that $\bar{c}_l x + \bar{u} \geq 0.49 x, \bar c_l =3, \bar u <0$ for any $x \geq 0$,
here $\bar{u}, \bar c_l$ are the velocity and the scaling exponents of an approximate self-similar profile.
We use this property to extract the main damping effect from the linearized operator in the near field using some carefully designed singular weights.

As we will show later, $\bar{c}_l x + \bar{u}$ is the velocity field for the linearized equation in the dynamic rescaling formulation. The inequality  $\bar{c}_l x + \bar{u} \geq 0.49 x, x \geq 0$ implies that the perturbation is transported from the near field to the far field and then damped by the damping term $\bar{c}_\om \om$ in the $\om$ equation and by $2 \bar{c}_\om \th_x$ in the $\th_x$ equation. This is the main physical mechanism that generates the dynamic stability of the self-similar blowup in the HL model. We believe that this also captures the dynamic stability of the blowup scenario considered by Luo-Hou along the boundary \cite{luo2013potentially-2,luo2014potentially}, whose numerical evidence has been reported in \cite{liu2017spatial}.

There are four important components of our analysis for the HL model. 
The first one is to construct the approximate steady state with sufficiently small residual error by decomposing it into a semi-analytic part that captures the far field behavior of the solution and a numerically computed part that has compact support. See more discussion in Section \ref{sec:on_steady_state}. The second one is that we extract the damping effect from the local terms in the linearized equations by using carefully designed singular weights. The third one is that the contributions from the advection terms are relatively weak compared with those coming from the vortex stretching terms. As a result, we can treat those terms coming from advection as perturbation to those from vortex stretching. The last one is to apply some sharp functional inequalities to control the nonlocal terms and take into account cancellation among various nonlocal terms. This enables us to show that the contributions from the nonlocal terms are relatively small compared with those from the local terms and can be controlled by the damping terms. 
We refer to Section \ref{sec:dyn} for more detailed discussion of the main ingredients in our stability analysis.

We believe that the analysis of the 2D Boussinesq equations and 3D Euler equations with smooth initial data and boundary would benefit from the four important components mentioned above.
The stability analysis of the HL model is established based on some weighted $L^2$ space. For the 2D Boussinesq equations and 3D Euler equations, 
a wider class of functional spaces, e.g. weighted $L^p$ or weighted $C^{\alpha}$ spaces, can be explored to derive larger damping effect from the linearized equations and to further establish stability analysis.

There is an interesting implication of our blowup results for the self-similar solution $(\om_*, \th_* )$ defined in \eqref{eq:self_similar}. In Section \ref{sec:decay_profile}, we show that the profile satisfies $\lim_{x\to \infty} \th_{\infty}(x)  |x|^{ -\g}  = C $ for some $C >  0$ (see \eqref{eq:self_similar}). Thus, we have $\lim_{ t \to 1}\th_*(x,t) \to C |x|^{ \g} $ for any $x \neq 0$. Since $0 < \g < 1$, the self-similar solution forms a cusp singularity at $x = 0$ as $t \to 1$. Moreover, from Theorem \ref{thm1}, for a class of initial data $\th_0$, the $C^{\g}$ norm of the singular solution $\th$ is uniformly bounded up to the blowup time. Note that from Theorem \ref{thm1}, we have $|\g - 0.33304| < 2 \cdot 10^{-5}$, thus $\g \approx \f{1}{3}$ and $\lim_{ t \to 1}\th_*(x,t) = C |x|^{ \g} \approx C |x|^{1/3} $. Similarly, we can generalize the method of analysis to prove $\lim_{ t \to 1}\om_*(x,t) = C_2 |x|^{ (\g-1)/2 } \approx C_2 |x|^{ -1/3} $. Interestingly, the limiting behavior is closely related to a family of explicit solutions of \eqref{eq:HL} discovered by Hoang and Radosz in \cite{hoang2020singular}
\beq\label{eq:th_13}
\om(x, t)  = k |x|^{-1/3} \sgn(x), \quad \th(x, t) = c_1 k^2 |x|^{1/3} + c_2 k^3 t ,
\eeq
where $c_1, c_2 > 0$ are suitable constants and $k > 0$ is arbitrary. We remark that from Theorem \ref{thm1}, the $C^{1/3}$ norm of $\th$ from a class of smooth initial data that we consider  blows up at the singularity time since $\f{1}{3} > \g$, while the non-smooth $\th$ in \eqref{eq:th_13} remains in $C^{1/3}$ for all time.


The cusp formation and the H\"older regularity on $\th$ are related to the $C^{1/2}$ conjecture by Silvestre and Vicol in \cite{silvestre2014transport} and the cusp formation on the Cordoba-Cordoba-Fontelos (CCF) model \cite{cordoba2005formation,kiselev2010regularity,li2008blow,chen2019finite}, which is the $\theta-$equation in \eqref{eq:HL} coupled with $u = H \th$.
The cusp formation of a closely related model was established in \cite{hoang2017cusp}, and the $C^{1/2}$ conjecture was studied in \cite{Elg17,Elg19} for a class of $C^{1,\al}$ initial data with small $\al$. Using the same method for the HL model, we have obtained an approximate self-similar profile for the CCF model with residual $O(10^{-8})$ and $\gamma = 0.5414465$, which is accurate up to six digits. This blowup exponent $\gamma$ is qualitatively similar to that obtained in 
\cite{lushnikov2020collapse} for the generalized Constantin-Lax-Majda model (gCLM) (see\cite{OSW08}) with $a=-1$. In a follow-up work, we will generalize our method of analysis to study the cusp formation of the CCF model, and rigorously prove that $\theta \in C^\gamma$ up to the singularity time with $\gamma > 1/2$. Moreover, the $C^\beta$ norm of $\th$ will blow up at the singularity time for any $\beta > \gamma$.

There has been a lot of effort in studying potential singularity of the 3D Euler equations using various simplified models. In 
\cite{choi2015finite,kiselev2018finite,hoang2018blowup,hoang2020singular}, the authors proposed several simplified models to study the Hou-Luo blowup scenario \cite{luo2014potentially,luo2013potentially-2} and established finite time blowup of these models. In these works, the velocity is determined by a simplified Biot-Savart law in a form similar to the key lemma in the seminal work of Kiselev-Sverak \cite{kiselev2013small}. 
In \cite{hou2015self}, Hou and Liu established the self-similar singularity of the CKY model \cite{choi2015finite} using the 
property that the CKY model can be reformulated as a local ODE system. 
The HL model does not enjoy a similar local property, and our method to prove  self-similar singularity is completely different from that in \cite{hou2015self}. In \cite{elgindi2017finite,elgindi2018finite} , Elgindi and Jeong proved finite time singularity formation for the 2D Boussinesq and 3D axisymmetric Euler equations in a domain with a corner using $\mathring{C}^{0,\al}$ data.

Several other 1D models, including the Constantin-Lax-Majda (CLM) model \cite{CLM85}, the De Gregorio (DG) model \cite{DG90,DG96}, and the gCLM
model \cite{OSW08}, have been introduced to study the effect of advection and vortex stretching in the 3D Euler. 
Singularity formation from smooth initial data has been established for the CLM model in \cite{CLM85},  for the DG model in \cite{chen2019finite}, and for the gCLM model with various parameters in \cite{Elg17,chen2020singularity,chen2020slightly,chen2019finite,Elg19,Cor10}. 
In the viscous case, singularity formation of the gCLM model with some parameters has been established in \cite{chen2020singularity,schochet1986explicit}.

The rest of the paper is organized as follows. In Section \ref{sec:dyn}, we outline some main ingredients in our stability analysis by using the dynamic rescaling formulation. 
Section \ref{sec:lin} is devoted to linear stability analysis. 
In Section \ref{sec:on_steady_state}, we discuss some technical difficulty in obtaining an approximate steady state with a residual error of order $10^{-10}$. In Section \ref{sec:non}, we perform nonlinear stability analysis and establish the finite time blowup result. 
In Section \ref{sec:hol}, we estimate the H\"older regularity of the singular solution.
In Section \ref{sec:con}, we give a formal derivation to demonstrate that both the HL model and the 2D Boussinesq equations with $C^{1,\alpha}$ initial data for velocity and $ \th$ and with boundary have the same leading system for small $\alpha$. We make some concluding remarks in Section \ref{sec:conclude}. Some technical estimates and derivations are deferred to the Appendix.

\section{Outline of the main ingredients in the stability analysis}\label{sec:dyn}

In this section, we will outline the main ingredients in our stability analysis by using the dynamic rescaling formulation for the HL model. The most essential part of our analysis lies in the linear stability. We need to use a number of techniques to extract the damping effect from the linearized operator around the approximate steady state of the dynamic rescaling equations and obtain sharp estimates of various nonlocal terms. Since the damping coefficient we obtain is relatively small (about $0.03$), we need to construct an approximate steady state with a very small residual error of order $10^{-10}$. This is extremely challenging since the solution is supported on the whole real line with a slowly decaying tail
 in the far field. We use analytic estimates and numerical analysis with rigorous error control to verify that the residual error is small in the energy norm. See detailed discussions in
Section \ref{sec:on_steady_state} and Section 10 of the Supplementary Material \cite{chen2021HLsupp}.

Passing from linear stability to nonlinear stability is relatively easier since the perturbation is quite small due to the small residual error. Yet we need to verify various inequalities involving the approximate steady state using the interval arithmetic \cite{moore2009introduction,rump2010verification,gomez2019computer} and numerical analysis with computer assistance. The most essential part of the linear stability analysis can be established based on the grid point values of the approximate steady state constructed on a relatively coarse grid, which does not involve the lengthy rigorous verification. See more discussion in Section \ref{sec:lin_verif}. The reader who is not interested in the rigorous verification can skip the lengthy verification process presented in the Supplementary Material \cite{chen2021HLsupp}.

\subsection{Dynamic rescaling formulation} An essential tool in our analysis is the dynamic rescaling formulation. Let $ \om_{phy}(x, t), \th_{phy}(x,t)$ be the solutions of the physical equations \eqref{eq:HL}, then it is easy to show that 
\[
  \om(x, \tau) = C_{\om}(\tau) \om_{phy}(   C_l(\tau) x,  t(\tau) ), \quad   \th(x, \tau) = C_{\th}(\tau) \th_{phy}(C_l(\tau) x, t(\tau))
  \]
are the solutions to the dynamic rescaling equations
\beq\label{eq:HL_dyn}
\om_{\tau} + (c_l x + u) \om_x = c_{\om} \om + \th_x , \quad \th_{\tau} + (c_l x + u) \th_x = c_{\th} \th, \quad u_x = H \om,
\eeq
where $t(\tau) = \int_0^{\tau} C_{\om}(s ) d s$ and
\[
  C_{\om}(\tau) = \exp\lt( \int_0^{\tau} c_{\om} (s)  d s \rt), \ C_l(\tau) = \exp\lt( \int_0^{\tau} -c_l(s) ds    \rt) , \ 
  C_{\th}(\tau) = \exp\B( \int_0^{\tau} c_{\th}(s) d s \B).
\]
In order for the dynamic rescaling formulation to be equivalent to the original HL model, we must enforce a relationship among the three scaling parameters, $c_l$, $c_\om$ and $c_\th$, i.e. $c_{\th} = c_l + 2 c_{\om}$.

The dynamic rescaling formulation was introduced in \cite{mclaughlin1986focusing,  landman1988rate} to study the self-similar blowup of the nonlinear Schr\"odinger equations. This formulation
is closely related the modulation technique, which has been developed by Merle, Raphael, Martel, Zaag and others, see e.g. \cite{merle1997stability,kenig2006global,merle2005blow,
martel2014blow,merle2015stability}.
  The dynamic rescaling formulation and modulation technique have been very effective in analyzing singularity formation for many nonlinear PDEs including the nonlinear Schr\"odinger equation \cite{kenig2006global,merle2005blow}, the nonlinear wave equation \cite{merle2015stability}, the nonlinear heat equation \cite{merle1997stability}, the generalized KdV equation \cite{martel2014blow}, the De Gregorio model and the generalized Constantin-Lax-Majda model \cite{chen2019finite,chen2020slightly,chen2020singularity}, and singularity formation in 3D Euler equations \cite{elgindi2019finite,chen2019finite2}.

To simplify our presentation, we still use $t$ to denote the rescaled time in \eqref{eq:HL_dyn}. Taking the $x$ derivative on the $\th$ equation in \eqref{eq:HL_dyn} yields 
\beq\label{eq:HLdyn}
\bal
\om_t + (c_l x + u  ) \om_x &= c_{\om} \om + \th_x , \\
(\th_x)_t + (c_l x + u) \th_{xx} &= (c_{\th}  - c_l - u_x) \th_x = (2 c_{\om} - u_x) \th_x, 
\quad  u_x = H \om,
\eal
\eeq
where $c_{\th} =c_l + 2 c_{\om}$. We still have two degrees of freedom in choosing $c_l, c_{\om}$ to uniquely determine the dynamic rescaled solution. We impose the following normalization conditions on $c_{\om}, c_l$
\beq\label{eq:normal0}
c_l = 2 \f{\th_{xx}(0) }{\om_x(0) }, \quad c_{\om} = \f{1}{2} c_l + u_x(0). 
\eeq
These two normalization conditions play the role of forcing 
\beq\label{eq:normal1}
\theta_{xx}(t,0)=\theta_{xx}(0,0), \quad \omega_x(t,0)=\omega_x(0,0)
\eeq
for all time. Our study shows that enforcing $\theta_{xx}(t,0)$ to be independent of time is essential for stability by eliminating a dynamically unstable mode in the dynamic rescaling formulation.

\subsection{Main result}\label{sec:thm2}
Throughout this paper, we will consider solution of \eqref{eq:HL_dyn} with odd $\om, \th_x$ and $\th(t, 0)=0$.
Under this setting, it is not difficult to show that the odd symmetries of $\th_x, \om, u$ and the condition $\th(t, 0) = 0$
are preserved by the equations. 

Due to the symmetry, we restrict 
the inner product and $L^2$ norm to $\R_+$
\beq\label{nota:inner}
\la f, g \ra  \teq \int_0^{\infty} fg dx  , \quad  || f||_2^2  = \int_0^{\infty} f^2 dx .
\eeq

Let $\psi, \vp$ be the singular weights defined in \eqref{eq:mode}, and $\lam_i$ be the parameter given in \eqref{eq:para3}. We use the following energy in our energy estimates
\beq\label{energy:thm}
\bal
E^2( f, g ) \teq &||  f \psi^{1/2}||_2^2 + \lam_1 ||  g \psi^{1/2}||_2^2
+ \lam_2 \f{\pi}{2} ( H g(0) )^2 + \lam_3 \la f , x^{-1}\ra^2+  \lam_4 ( ||  D_x f \psi^{1/2}||_2^2 + \lam_1 || D_x g \vp^{1/2}||_2^2  ) ,
\eal
\eeq
where $H g(0) = -\f{1}{\pi}\int_{\R} g x^{-1} dx$ is related to $c_{\om}$ in \eqref{eq:normal0}. Our main result is the following.

\begin{thm}\label{thm2}
Let $(\bar{\th},\bar{\om}, \bar c_l, \bar c_{\om})$ be the approximate self-similar profile constructed in Section \ref{sec:on_steady_state}, and $E_* = 2.5 \cdot 10^{-5}$.
For odd initial data $\th_{0,x}, \om_0$ of \eqref{eq:HL_dyn} with $\th_0(0) = 0$ and a small perturbation to $(\bar \th_x, \bar \om ), E( \th_{0,x} - \bar \th_x , \om_0 - \bar \om) \leq E_*$, we have (a) $E( \th_x - \bar \th_x, \om - \bar \om) \leq E_*$ for all time.

(b) The solution $(\th, \om,c_l, c_{\om})$ converges to a steady state of \eqref{eq:HL_dyn} $(\th_{\inf}, \om_{\inf}, c_{l,\inf}, c_{\om,\inf})$ 
\[
 || (\th_x(t) - \th_{\infty,x}) \psi^{1/2}||_2 + || (\om(t) - \om_{\infty} ) \vp^{1/2} ||_2 + || c_l (t)- c_{l,\infty}||_2 + || c_{\om} (t)- c_{\om, \infty}||_2 \leq C  e^{-\kp_2 t} 
\]
exponentially fast, for some $ \kp_2 > 0, C >0$. Moreover, $(\th_{\inf}, \om_{\inf}, c_{l,\inf}, c_{\om,\inf})$ enjoys the regularity $E( \th_{x,\inf} - \bar \th_x , \om_{\inf} - \bar \om) \leq E_*$, and is the unique steady state in the class 
$E( \th_x- \bar \th_x, \om - \bar \om) \leq E_*$ with normalization conditions \eqref{eq:normal0}and $\th(0) = 0$, and odd assumption on $\th_x, \om$.

(c) Let $\g = \f{ c_{\th,\infty}}{c_{l,\infty}} = \f{ c_{l,\infty} + 2 c_{\om,\inf} }{c_{l,\infty}}$. We have $| \frac{ c_{\om,\inf} }{c_{l,\infty}} - 2.99870| \leq 6 \cdot 10^{-5}$. Moreover, the solution enjoys the H\"older estimates $ \th_{\infty} \in C^{\g}$ and 
 $\sup_{t \geq 0 } || \th ||_{C^{\g}} \les 1.$

 (d) For the physical equations \eqref{eq:HL} with the above initial data, the solution blows up in finite time $T$ with the following blowup estimates for any $\g < \b \leq 1$
\[
|| \th_{phy}(t)||_{C^{\b}}  \gtrsim (T-t)^{-\d},
  \quad \d = \f{ 2( \b- \g)}{1-\g}  >0.
 \]
If in addition $\th_{0,x} |x|^{1-\g} \in L^{\inf}$, the $C^{\g}$ norm is uniformly bounded 
up to the blowup time: $ \sup_{t \in [0, T) } || \th_{phy}(t) ||_{C^{\g}} \les 1 $. 

\end{thm}

The assumption $\th_{0,x} |x|^{1-\g } \in L^{\inf}$ in (d) is to ensure the decay $|\th_{0,x}|\leq C |x|^{\g-1}$, which is consistent with $\th_0 \in C^{\g}$. In fact, if $\th_0 \in C^{\g}$, we get $|\th_0(x) | \les 1 + |x|^{\g}$. Then, formally, $\th_{0,x}$ has a decay rate $|x|^{\g-1}$.

\subsection{Main ingredients in our stability analysis}\label{sec:main_idea}
The key step to prove Theorem \ref{thm2} is the stability analysis. We will outline several important ingredients to establish it in this subsection

\vspace{-0.05in}
\subsubsection{ The stability of the linearized operator }

The most essential part of our analysis is the linear stability of the linearized operator around the approximate steady state $(\bar \th, \bar \om, \bar c_l,  \bar{c}_\omega)$. To simplify our notation, we still use $\omega$, $u$, $\theta, \ c_l$,  and $c_\omega$ to denote the perturbation. The linearized system for the perturbation is given below by neglecting the nonlinear and error terms:
\beq\label{eq:intro_lin}
\bal
\pa_t \th_x + ( \bar{c_l} x + \bar{u}) \th_{xx}  &=  (2 \bar c_{\om} - \bar u_x) \th_x
+ ( 2c_{\om} - u_x) \bar \th_x - u \bar \th_{xx}, \\
\om_t +  ( \bar{c_l} x + \bar{u}) \om_x & =  \bar{c}_{\om} \om  +\th_x +  c_{\om} \bar{\om} -   u \bar{\om}_x , \quad c_{\om} = u_x(t,0), \quad c_l = 0.
\eal
\eeq
The condition $c_{\om} = u_x(t,0), c_l = 0 $ is a consequence of the normalization conditions \eqref{eq:normal1}. There are two groups of terms in the above system, one representing the local terms and the other representing the nonlocal terms. 
Among the nonlocal terms, we can further group them into three subgroups, one from the vortex stretching term, one from the advection term, and the remaining from the rescaling factor $c_{\om}$.

As in our previous works \cite{chen2019finite,chen2019finite2}, we design some singular weights to extract the damping effect from the local terms. As we mentioned before, we will use
$\bar{c_l} x + \bar{u} \geq 0.49 x$ to extract an $O(1)$ damping effect. 
Since the damping coefficient that we can extract from the local terms is relatively small and the linearized operator is not a normal operator,
we typically expect to have a transient growth for a standard energy norm 
of the solution to \eqref{eq:intro_lin}. This will present considerable difficulty for us to obtain nonlinear stability since the approximate steady state also introduces a residual error. To overcome this difficulty, we need to design a weighted energy norm carefully so that the energy of the solution to the linearized equations decreases monotonically in time. We remark that weighted energy estimates with singular weights have also been used in \cite{elgindi2019finite,Sve19,chen2020slightly,chen2020singularity} for nonlinear stability analysis.

\subsubsection{Control of nonlocal terms}\label{sec:models}

The most challenging part of the linear stability analysis is how to control several nonlocal terms that are of $O(1)$. It is essential to obtain sharp estimates of these nonlocal terms by 
applying sharp weighted functional inequalities, e.g. Lemma \ref{lem:wg}, and
taking into account the cancellation among different nonlocal terms and the structure of the coupled system. 
We have also used the $L^2$ isometry property and several other properties of the Hilbert transform in an essential way. 
We remark that some of these properties of the Hilbert transform have been used in the previous works, see, e.g. \cite{Cor06,Elg17,chen2019finite,Cor10,chen2020singularity}.
Based on our observation that the blowup is driven by vortex stretching and the advection is relatively weak compared with vortex stretching, we will treat the nonlocal terms that are generated by the advection terms, e.g. $u \bar \th_{xx}$ in \eqref{eq:intro_lin}, as perturbation to the linearized vortex stretching terms, e.g. $u_x \bar \th_x$ in \eqref{eq:intro_lin}. 
We will use the following five strategies in our analysis.

\vspace{0.1in}
\noindent
{\bf (1) The decomposition of the velocity field.}
We first denote $\tilde{u} \teq u - u_x(0) x$ and choose a constant $c=1/(2p-1)$ where $p$ is related to the order of the singular weight $|x|^{-p}$ being used. We further decompose $\tilde{u}$ into a main term and a remainder term as follows:
\beq\label{eq:decomp_u}
\tilde{u} = c x \tilde{u}_x + (\tilde{u} - c x \tilde{u}_x) \teq \tilde{u}_M + \tilde{u}_R,
\eeq
where $\tilde{u}_M = c x \tilde{u}_x$ and $ \tilde{u}_R = (\tilde{u} - c x \tilde{u}_x)$. The contribution from the remainder term $\tilde{u}_R$ is smaller than $x \td u_x$ due to an identity (see Appendix \ref{app:IBPu})
\beq\label{eq:IBPu}
 || (\td{u} - \f{1}{ 2p  - 1} \td{u}_x x ) x^{- p}||_2^2 = \f{1}{(2p-1)^2}  \int_{\R_+} \f{\td{u}^2_x}{ x^{ 2p -2}}  dx.
\eeq
We can choose $p=3$ in the near field, which enables us to gain a small factor of $1/5$ in estimating the $\tilde{u}_R$ term in terms of the weighted norm of $\tilde{u}_x$.

\vspace{0.05in}
\noindent
{\bf (2) Exploiting the nonlocal cancellation between $\tilde{u}_x$ and $\omega$.}
For the main term $\td u_{M}= c x \td u_x$ and the vortex stretching term $ -u_x \bar \th_{x}$, 
we use an orthogonality between $\tilde{u}_x$ and $\omega$ 
\beq\label{eq:model2_cancel}
\la \tilde{u}_x, \om x^{-3} \ra
= \la H \om - H\om(0), \om x^{-3} \ra \;= \;0
\eeq
(see Lemma \ref{lem:cancel}). 
We will use similar orthogonal properties to exploit the cancellation between $- \tilde{u}_x \bar\th_x$ in the $\theta_x$ equation and $\theta_x$ in the $\omega$ equation in \eqref{eq:intro_lin} by performing the weighted $L^2$ estimates for $\theta_x$ and $\omega$ together. To illustrate this idea, we consider the following model:
 
\vspace{0.05in}
\noindent
{\bf Model 1 for nonlocal interaction} 
\beq\label{eq:model2}
\pa_t \th_x = - ( u_x - u_x(0)) \bar \th_x, \quad \om_t  = \th_x.
\eeq
The above system is derived by dropping other terms in \eqref{eq:intro_lin}. The profile $\bar \th_x$ satisfies $\bar \th_x(0) = 0$ and $\bar \th_x > 0$ for $x >0$. 

By performing $L^2(\rho_1)$ estimate on $\th_x$ and $L^2(\rho_2)$ estimate on $\om$, we get 
\beq\label{eq:est_coupled}
\f{1}{2} \f{d}{dt} ( \la \th_x, \th_x \rho_1  \ra  + \la \om, \om \rho_2 \ra  )
= - \la (u_x - u_x(0) ) \bar \th_x \rho_1, \th_x \ra + \la \om \rho_2, \th_x \ra \teq I.
\eeq

From \eqref{eq:model2_cancel}, we know that $( u_x - u_x(0)) x^{-2}$ and $\om x^{-1}$ are orthogonal. Formally, $I$ is the sum of the projections of $\th_x$ onto two directions that are orthogonal. To exploit this orthogonality, we choose $ \rho_1 =(\mu x \bar \th_x)^{-1} \rho_2$ with any $\mu >0$. We can rewrite $I$ as follows 
\[
I =  \la -(u_x - u_x(0) )x^{-2}, \th_x   \bar \th_x \rho_1 x^2 \ra 
+ \la \mu \om x^{-1}, \th_x  \bar \th_x \rho_1 x^2  \ra \teq \la  A + B,  \th_x   \bar \th_x \rho_1 x^2 \ra ,
\]
where $A = -(u_x - u_x(0) )x^{-2}$ and $B =\mu \om x^{-1} $. Applying the Cauchy-Schwarz inequality yields 
\[
I \leq || A+ B||_2  ||  \th_x   \bar \th_x \rho_1 x^2 ||_2.
\]

The equality can be achieved if $ \th_x  \bar \th_x \rho_1 x^2 = c (A+B)$ for some $c$. Expanding $|| A+B||_2$ and using Lemma \ref{lem:cancel} with $f = \om$ and $ g=u$, we get 
\beq\label{eq:model2_impro}
|| A+ B||_2^2 =  || A||_2^2 + || B||_2^2 + 2 \la A, B \ra = || A ||_2^2 + ||B||_2^2,
\eeq
which is sharper than the trivial estimate $|| A + B||_2 \leq  || A||_2 + ||B||_2$. The $||A||_2^2$ term can be further bounded by $|| \om x^{-1}||_2^2$ using the $L^2$ isometry of the Hilbert transform in Lemma \ref{lem:iso}. The $||B||_2^2$ term can be bounded by the weighted $L^2$ norm of $\om$ directly.

\vspace{0.05in}
\noindent
{\bf (3) Additional damping effect from $c_\omega$.}
Another nonlocal term in \eqref{eq:intro_lin} is $c_\omega=u_x(t,0)=H(\omega)(t,0)$. Physically, the role of $c_\omega$ is to rescale the amplitude of the blowup profile $\omega$ in the original physical variable so that the magnitude of the dynamic rescaled profile remains $O(1)$ for all time. Thus, we expect that the dynamic rescaling parameter $c_\omega$ should also offer some 
stabilizing effect to the blowup profile and the linearized system \eqref{eq:intro_lin}. 
Indeed, by deriving an ODE for $c_\omega$, we can extract an additional damping term, which will be used to control other nonlocal terms associated with $c_\omega$. 
To illustrate this idea, we consider the following model: 

\vspace{0.05in}
\noindent
{\bf Model 2 for the $c_{\om}$ term}
\beq\label{eq:model3}
\pa_t \th_x = c_{\om} \bar f, \quad \pa_t \om = \th_x + c_{\om} \bar g,
\eeq
where $\bar f, \bar g$ are odd and $\bar f, \bar g >0$ for $x > 0$ with $ \bar f x^{-1}, \bar g x^{-1} \in L^1$. Note that the profile satisfies that $\bar \th_x - x \bar \th_{xx}, \bar \om - x \bar \om_x$ are odd and positive for $x>0$.
This system models the $c_{\om}$ terms in \eqref{eq:intro_lin} with coupling $\th_x$ in $\om$ equation by dropping other terms. Recall 
\[
c_{\om} = -\f{1}{\pi} \int_{\R} \f{\om}{x} dx = -\f{2}{\pi} \la \om , x^{-1} \ra.
\]

Obviously, it can be bounded by some weighted $L^2$ norm of $\om$ using the Cauchy-Schwarz inequality. Yet, the constant in this estimate is large. Denote $A = \la \bar f , x^{-1}\ra,  B = \la \bar g , x^{-1} \ra$. By definition, $A, B > 0$.
We derive an ODE for $c_{\om}$ using the $\om$ equation
\[
  \pa_t \la \om, x^{-1} \ra = c_{\om} \la \bar g, x^{-1} \ra + \la \th_x , x^{-1} \ra
  = - \f{2}{\pi} B \la \om, x^{-1} \ra  + \la \th_x , x^{-1} \ra.
\]

We see that the $c_{\om}$ term in the $\om $ equation in \eqref{eq:model3} has a damping effect, which is not captured by the weighted $L^2$ estimates. To handle the coupled term, we also derive an ODE for $\la \th_x, x^{-1} \ra$ using the $\th_x$ equation 
\[
\pa_t \la \th_x , x^{-1} \ra = c_{\om} \la \bar f, x^{-1} \ra
= -\f{2}{\pi} A \la \om, x^{-1} \ra.
\]

Multiplying both sides of these ODEs by $\la \om, x^{-1} \ra$ or $\la \th_x, x^{-1}\ra$, we yield 
\beq\label{eq:model3_ode}
\bal
\f{1}{2} \f{d}{dt} \la \om, x^{-1} \ra^2 
&= -\f{2}{\pi} B \la \om, x^{-1} \ra^2 + \la \th_x , x^{-1} \ra \la \om, x^{-1} \ra,  \\
\f{1}{2} \pa_t \la \th_x , x^{-1} \ra^2 
&= -\f{2}{\pi} A  \la \th_x, x^{-1} \ra \la \om, x^{-1} \ra  .
\eal
\eeq

The $\la \th_x, x^{-1} \ra \la \om, x^{-1} \ra$ terms in the above ODEs have cancellation. This implies that the $c_{\om}$ term in the $\th_x$ equation and $\th_x$ term in the $\om$ equation have cancellation, which is not captured by the weighted $L^2$ estimate. We will derive similar ODEs in the analysis of \eqref{eq:intro_lin} and obtain damping term similar to $-\f{2}{\pi} B \la \om, x^{-1} \ra^2$ in the above ODEs, which enables us to control the $c_{\om}$ terms in \eqref{eq:intro_lin} effectively.

\vspace{0.05in}
\noindent
{\bf (4) Estimating the ${\bf u}$ term in \eqref{eq:intro_lin}.} To estimate the $u$ terms in \eqref{eq:intro_lin} effectively, we have two approaches. The first approach is to exploit the cancellation between $u$ and $\om$ similar to that in Model 1. See Lemma \ref{lem:cancel}. The second approach is to decompose $\td  u$ into the main term $ \td u_M = c x \td u_x$ and an error term $\td u_R$ as \eqref{eq:decomp_u1}. For $\td u_M$, we employ the estimates on $u_x$ discussed previously. The error term $ \td u_R$ enjoys better estimate \eqref{eq:IBPu} and is 
treated as a perturbation.

\vspace{0.05in}
\noindent
{\bf (5) Obtaining sharp estimates for other interaction terms.}
To obtain sharper estimates for a number of quadratic interaction terms, we introduce a number of parameters in various intermediate steps and optimize these parameters later by 
solving a constrained optimization problem. In the ODE for $c_{\om}$ and the weighted $L^2$ estimates,
we need to control a number of quadratic interaction terms, e.g. $\la \om, x^{-1}\ra \cdot \la \th_x , x^{-1}\ra$. 
We treat these interaction terms as the products of projection of $\th_x$ and $\om$ onto some low dimensional subspaces and reduce them to some quadratic forms in a finite dimensional space. This connection enables us to reduce the problem of obtaining sharp estimates of these terms to computing the largest eigenvalue $\lam_{\max}$ of a matrix. 
We then compute $\lam_{\max}$ as part of the constrained optimization problem to determine these parameters and obtain a sharper upper bound in the energy estimate.

\section{Linear stability}\label{sec:lin}

In this section, we establish the linear stability of \eqref{eq:lin0} in some weighted $L^2$ spaces.

\subsection{ Linearized operators around approximate steady state}\label{sec:linop}  \

The approximate steady state of \eqref{eq:HLdyn} $( \bar \th_x, \bar \om)$ we construct are odd with scaling factors 
\[
 \bar c_l = 3, \quad |\bar c_{\om}+ 1.00043212| < 10^{-8}, \quad \bar c_{\om} \approx -1.
 \]
It  has regularity $ \bar \om, \bar \th_x \in C^{3}$ and decay rates $ \pa_x^i \bar \om \sim x^{ \al-i}, \pa_x^i \bar \th_x \sim x^{2\al-i},i=0,1,2$ with $\al$ slightly smaller than $-1/3$. One can find plots of $(\bar \om, \bar \th_x)$ (with particular rescaling) in Figure \ref{fig:compare_initial_data} in Section \ref{sec:numerical_uniqueness}. See detailed discussion in Section \ref{sec:on_steady_state}. Note that we do not require a $C^{\inf}$ approximate steady state in our analysis, since the $C^3$ approximate steady state 
is regular enough for us to perform weighted $H^1$ estimates and establish its nonlinear stability.

Linearizing around the approximate steady state $(\bar \th_x, \bar \om)$, we obtain the equations for the perturbation 
\beq\label{eq:lin00}
\bal
\pa_t \th_x + ( \bar{c_l} x + \bar{u}) \th_{xx}  &=  (2 \bar c_{\om} - \bar u_x) \th_x
+ (2 c_{\om} - u_x) \bar \th_x - u \bar \th_{xx} + F_{\th} + N(\th) , \\
\om_t +  ( \bar{c_l} x + \bar{u}) \om_x & =  \bar{c}_{\om} \om  +\th_x +  c_{\om} \bar{\om} -   u \bar{\om}_x +  F_{\om} + N(\om) ,
\eal
\eeq
where the error terms $F_{\th}, F_{\om}$ and the nonlinear terms $N(\th), N(\om)$ read
\beq\label{eq:NF}
\bal
F_{\th} &= ( 2 \bar{c}_{\om} - \bar u_x) \bar\th_x - (\bar{c_l} x + \bar{u}) \cdot \bar{\th}_{xx} , \quad
F_{\om} =\bar{\th_x} + \bar{c_{\om} } \bar{\om } -  ( \bar{c_l} x + \bar{u} ) \cdot \bar{\om}_x , \\
N(\th) & = (2 c_{\om} - u_x) \th_x - u \th_{xx},  \quad N(\om)  = c_{\om} \om -  u \om_x . \\
 \eal
 \eeq

We consider odd initial perturbation $\om_0, \th_{0,x}$ with $\om_{0,x}(0) = 0, \th_{0, xx}(0)=0$.
Note that the normalization conditions \eqref{eq:normal0},\eqref{eq:normal1} implies 
\beq\label{eq:normal}
c_{\om} = u_x(0), \quad c_l = 0,  \quad \th_{xx}(t,0) = \th_{0,xx}(0) = 0, \quad \om_{x}(t,0) = \om_{0x}(0) = 0,
\eeq
for the perturbation. Since $\om, \th_{x}$ are odd, these normalization conditions imply that 
near $x=0$, $\om = O(x^3),  \th_{x} = O(x^3)$ for sufficient smooth solution. This important property enables us to use a more singular weight in our stability analysis to extract a larger damping coefficient.

We rewrite the $c_{\om}$ and $u$ terms as follows
\beq\label{eq:lin_rewrite}
\bal
 (2 c_{\om} - u_x) \bar \th_x - u \bar \th_{xx} 
 &= - (u_x - u_x(0)) \bar \th_x - (u - u_x(0) x) \bar \th_{xx} + c_{\om} ( \bar \th_x - x \bar \th_{xx} ), \\
c_{\om} \bar \om - u \bar \om_x & = -( u -u_x(0) x) \bar \om_x + c_{\om} (\bar \om - x \bar \om_x) \\
 \eal
\eeq

Denote $\Lambda = (-\D)^{1/2}$. From $ \pa_x u =  H \om$ and $\Lambda = \pa_x H$, we have $ u (x)= - \Lambda^{-1} \om (x) =\f{1}{\pi} \int \log| x-y| \om(y) dy$. Using this notation, we get $u - u_x(0)x = -\Lam^{-1} \om -H \om(0)x $. We introduce the following linearized operators 
\beq\label{eq:linop}
\bal
\cL_{\th 1}(f, g) & =  - ( \bar{c_l} x + \bar{u}) f_{x} + (2 \bar c_{\om} - \bar u_x)  f - (  H g - Hg(0)) \bar \th_x - ( - \Lam^{-1} g - Hg(0) x) \bar \th_{xx} , \\
\cL_{\om 1}(f, g)  & =  - ( \bar{c_l} x + \bar{u}) g_x  +  \bar{c}_{\om} g  + f  - (- \Lam^{-1} g - Hg(0) x) \bar \om_x ,  \\
\cL_{\th}(f, g) &= \cL_{\th 1}(f, g) + Hg(0)(\bar \th_x - x\bar \th_{xx}), \quad \cL_{\om}(f, g) = \cL_{\om 1}(f, g) + Hg(0)  (\bar \om - x \bar \om_x).
\eal
\eeq

Using these operators, we can rewrite \eqref{eq:lin00} as follows 
\beq\label{eq:lin0}
\bal
\pa_t \th_x & = \cL_{\th 1}(\th_x, \om) + c_{\om}(\bar \th_x - x\bar \th_{xx})  + F_{\th} + N(\th), \\
\pa_t \om & =  \cL_{\om 1}(\th_x, \om) + c_{\om}(\bar \om - x\bar \om_{x})  + F_{\om} + N(\om).
\eal
\eeq

Clearly, $\cL_{\th}, \cL_{\om}$ are the linearized operators associated to \eqref{eq:lin00}. The motivation of introducing $\cL_{\th 1}, \cL_{\om 1}$ is that the estimates of these operators will be used importantly in both the weighted $L^2$ and weighted $H^1$ estimates.

\subsection{Singular weights}

  For some $e_1, e_2, e_3 > 0$ determined by the profile $(\bar \om, \bar \th)$, we introduce 
\beq\label{eq:ext}
\xi_1 = e_1 x^{-2/3}  - ( \bar \th_x + \f{1}{5} x \bar \th_{xx}  ) , \quad \xi_2 = 
 e_2 x^{-2/3}  - ( \bar \th_x + \f{3}{7} x \bar \th_{xx}  ) , \quad \xi_3 = -\f{e_3}{3} x^{-4/3} - \bar \om_x.
\eeq

Following the guideline of the construction of the singular weight in \cite{chen2019finite}, we design different parts of the singular weight that have different decays as follows 
\beq\label{eq:mode}
\bal
\psi_n &= \f{1}{ \bar \th_x + \f{1}{5} x \bar \th_{xx} + \chi \xi_1 } ( \al_1 x^{-4} + \al_2 x^{-3}),
\quad \psi_f =  \f{1}{\bar \th_x + \f{3}{7} x \bar \th_{xx} + \chi \xi_2} \al_3 x^{-4/3} , \\
\vp_s & = \al_4 x^{-4} , \quad \vp_n  = \al_5 (  \al_1 x^{-3} + \al_2 x^{-2} ) , \quad \vp_f = \al_6 x^{-2/3} ,
\eal
\eeq
where the parameters are positive and chosen in \eqref{eq:para2}, and the cutoff function $\chi$ defined in Appendix \ref{app:chi} is supported in $|x| \geq \rho_2$ for $\rho_2 > 10^{8}$. The subscripts $f, n, s$ are short for \textit{far, near, singular}. We use the following weights in the weighted Sobolev estimate 
\beq\label{eq:wg}
\psi = \psi_n + \psi_f, \quad \vp = \vp_s + \vp_n + \vp_f.
\eeq

We introduce $\chi, \xi_1, \xi_2$ and add them in the definition of $\psi_n, \psi_f$ for the following purposes. Firstly, recall from the beginning of Section \ref{sec:linop} that $\bar \th_x + c x \bar \th_{xx}$ with $c = \f{1}{5}$ or $c = \f{3}{7}$ has decay $ x^{ 2 \al}$ which is close to $x^{-2/3}$. In particular, for sufficient large $x$, it can be well approximated by $e x^{-2/3}$ for some constant $e$. The parameters $e_1, e_2$ in \eqref{eq:ext} are determined in this way. Secondly, in the far field, where $\chi(x)=1$, the weights $\psi_n, \psi_f$ reduce to $c_1 x^{-7/3}, c_2 x^{-2/3}$ for some $c_1, c_2$, respectively. These explicit powers are much simpler than the weights in the near field and have forms similar to those in $\vp$. They will be useful for the 
analytic estimates (see Section \ref{sec:slow}) and simplify the computer-assisted verification of the estimates in the far field. We introduce $\xi_3$ similar to $\xi_1,\xi_2$ and it will be used later.

\begin{remark}\label{rem:chi}

Since $\chi$ is supported in $|x| >  10^8$ and the profile $(\bar \om, \bar \th_x)$ decays for large $|x|$, we gain a small factor in the estimates of the terms involving $\chi$. Thus, the upper bound in these estimates are very small.
 The reader can safely skip the technicalities due to the $\chi$ terms.
\end{remark}

\subsubsection{ The form of the singular weights}\label{sec:wg_form}

We add $\bar \th_x, \bar \th_{xx}$ terms in the denominators in $\psi_n,\psi_f$ to cancel the variable coefficients in our energy estimates. In Model 1 in Section \ref{sec:models}, we have chosen $ \rho_1 = ( \mu x \bar \th_x )^{-1}  \rho_2$ so that we can combine the estimates of two interactions in \eqref{eq:est_coupled}. Here, we design $\psi_n, \vp_n$ with a similar relation $ \psi_n = \f{1}{f} x^{-1} \vp_n, f = \bar \th_x + \f{1}{5} x \bar \th_{xx} + \chi \xi_1$ for the same purpose. Similar consideration applies to $\psi_f, \vp_f$. See also estimates \eqref{eq:varcoe_near}, \eqref{eq:varcoe_far}. This idea has been used in \cite{chen2019finite,chen2019finite2} for stability analysis. 

The profile satisfies $\bar \th_x + \f{1}{5} x \bar \th_{xx} , \bar \th_x + \f{3}{7} x \bar \th_{xx} >0 $ for $x > 0$. The weight $\psi$ is of order $x^{-5}$ for $x$ close to $0$, while it is of order $x^{-2/3}$ for large $x$. We choose $\vp$ of order $x^{-4}$ near $0$ so that we can apply the sharp weighted estimates in Lemma \ref{lem:wg} to control $ u_x$ and $u$.

We will use the following notations repeatedly 
\beq\label{eq:ut}
\td{u} \teq u - u_x(0)x, \quad \td{u}_x = u_x - u_x(0) .
\eeq

\subsection{Weighted $L^2$ estimates}
Performing weighted $L^2$ estimates on \eqref{eq:lin0} with weights $\psi, \vp$, we obtain 
\beq\label{eq:l2h1}
\bal
\f{1}{2} \f{d}{dt} \la \th_x , \th_x  \psi \ra = &  
\la \cL_{\th 1} \th_x, \th_x \psi \ra + c_{\om} \la \bar \th_x - x \bar \th_{xx} , \th_x \psi \ra+ \la  N(\th)  , \th_x  \psi  \ra   +  \la F_{\th} , \th_x \psi \ra  \\
= &  \B( - \la ( \bar c_l x + \bar u) \th_{xx} , \th_x \psi \ra  + \la (2 \bar c_{\om} - \bar u_x) \th_x, \th_x \psi \ra \B)   \\
& +\B(-  \la ( u_x -u_x(0) ) \bar \th_x + (u - u_x(0) x) \bar \th_{xx}, \th_x \psi \ra   + c_{\om} \la \bar \th_x - x \bar \th_{xx} , \th_x \psi \ra  \B) \\
&+ \la  N(\th)  , \th_x  \psi  \ra   +  \la F_{\th} , \th_x \psi \ra
\teq D_1 + Q_1 + N_1 + F_1 ,  \\
\f{1}{2} \f{d}{dt} \la \om , \om  \vp \ra = &  \la \cL_{\om 1} \om, \om \vp \ra + c_{\om} \la \bar \om_x - x \bar \om_{xx} , \om_x \vp \ra  + \la N(\om) , \om  \vp  \ra   +  \la F_{\om} , \om \vp  \ra \\
 = &\B(- \la  (   \bar{c_l} x + \bar{u}) \om_x,  \om  \vp  \ra   + \la  \bar{c}_{\om} \om,  \om  \vp \ra     \B) + \B( \la \th_x ,\om \vp\ra - \la (u -u_x(0) x ) \bar \om_x, \om \vp \ra \\
&+  c_{\om}\la \bar{\om}  -   x \bar{\om}_x, \om \vp \ra  \B)+ \la N(\om) , \om  \vp  \ra   +  \la F_{\om} , \om \vp  \ra 
\teq D_2 + Q_2 + N_2 + F_2 . \\
\eal
\eeq

Our goal in the remaining part of this Section is to establish an estimate similar to 
\beq\label{eq:goal}
D_1 + \lam_1 D_2 + Q_1 + \lam_1 Q_2 \leq -c ( || \th_x \psi^{1/2}||_2^2 + \lam_1 || \om \vp^{1/2}||_2^2 ),
\eeq
for some $\lam_1 > 0$ with $c>0$ as large as possible. This implies the linear stability of \eqref{eq:lin0} with $N_i , F_i = 0$ in the energy norm $ || \th_x \psi^{1/2}||_2^2 + \lam_1 || \om \vp^{1/2}||_2^2$. The actual estimate is slightly more complicated and we will add $c_{\om}^2, \la \th_x, x^{-1}\ra^2$ to the energy. We ignore the term $c_\omega$ and $\la \theta_x,x^{-1} \ra^2$ for now to illustrate the main ideas. See \eqref{eq:L2}.

The $D_1, D_2$ terms only involve the local terms about $\th_x, \om$ and we treat them as damping terms. The $Q_i$ term denotes the quadratic terms other than $D_i$ in the weighted $L^2$ estimates; The $N_i$ and $F_i$ terms represent the nonlinear terms and error terms in \eqref{eq:lin0}. 

For $D_1, D_2$, performing integration by parts on the transport term, we obtain 
\beq\label{eq:damp}
D_1 = \la D_{\th} , \th_x^2 \psi \ra, \quad  D_2  = \la D_{\om},  \om^2 \vp \ra,
\eeq
where $D_{\th}, D_{\om}$ are given by 
\[
\bal
D_{\th } & = \f{1}{2\psi} ( (\bar c_l x + \bar u ) \psi )_x + 2 \bar c_{\om} - \bar u_x ,\quad D_{\om}   = \f{1}{ 2 \vp} ( (\bar c_l x + \bar u ) \vp )_x + \bar c_{\om} .
\eal
\]
We will verify that $D_{\th}, D_{\om} \leq - c< 0$ for some constant $c>0$ in \eqref{ver:damp}, Appendix \ref{app:ver}. The weight $\psi$ in \eqref{eq:wg} involves three parameters $\al_1,\al_2,\al_3$. We choose the approximate values of $\al_i$ with $\al_i>0$ so that $D_{\th} \leq -c$ with $c$ as large as possible and varies slowly. This enables us to obtain a large damping coefficient. After we choose $\al_1, \al_2, \al_3$, 
we choose positive $\alpha_4, \alpha_5$ and $\alpha_6$ in the weight $\vp$ in \eqref{eq:wg} so that $D_\omega \leq -c_1$ with $c_1$ as large as possible and varies slowly.
The final values are given in \eqref{eq:para2}. See also Figure \ref{fig:linear_stability} for plots of the grid point values of $D_{\th}, D_{\om}$.

Using the notations in \eqref{eq:ut}, we can rewrite $Q_1 + \lam_1 Q_2$ as follows 
\beq\label{eq:bad}
\bal
Q_1 + \lam_1 Q_2 
&= - \la \td{u}_x \bar \th_x + \td{u} \bar \th_{xx} , \th_x \psi \ra 
+ \lam_1 \la \om, \th_x \vp \ra
 -\lam_1 \la \td u \bar \om_x, \om \vp \ra \\
 & + c_{\om} \la ( \bar{\th}_x - x\bar{\th}_{xx}) ,\th_x \psi \ra
 + \lam_1 c_{\om} \la  ( \bar{\om} - x \bar{\om}_x ) , \om \vp \ra.
 \eal
\eeq

The terms in $Q_1 +\lam_1 Q_2$ are the interactions among $u, \om, \th_x$ and do not have a favorable sign. Our goal is to prove that they are perturbation to the damping terms $D_1, D_2$ and establish \eqref{eq:goal}. This is challenging since the coefficients of the quadratic terms in $Q_1 + \lam_1 Q_2$ and in $D_i$ are comparable.

\subsubsection{ Decompositions on $Q_i$}\label{sec:decomp_Q}

Recall different parts of the weights in \eqref{eq:mode}. They provide a natural decomposition of the global interaction among $u, \om, \th_x$ into the near field and the far field interaction. We have a straightforward partition of unity 
\beq\label{eq:unity}
 \psi_n \psi^{-1} + \psi_f \psi^{-1} = 1, \quad \vp_n \vp^{-1}+ \vp_f \vp^{-1} + \vp_s \vp^{-1} = 1.
\eeq
According to different singular orders and decay rates of the weights in \eqref{eq:mode}, $ \psi_f \psi^{-1}, \vp_f \vp^{-1} $ are mainly supported in the far field, $\psi_n \vp^{-1}$ in the near field, $\vp_n \vp^{-1} $ near $|x| \approx 1$, and $\vp_s \vp^{-1}$ near $0$. Next, we decompose the interaction using these weights. Using $\psi = \psi_f + \psi_n$, we get
\[
- \la \td{u}_x \bar \th_x , \th_x \psi \ra
= - \la \td{u}_x ( \bar \th_x + \chi \xi_1 ) , \th_x \psi_n \ra
 - \la \td{u}_x ( \bar \th_x + \chi \xi_2 ) , \th_x \psi_f \ra
 + \la \td{u}_x  \chi (\xi_1  \psi_n + \xi_2 \psi_f ), \th_x \ra.
\]

We decompose the first two terms on the right hand side of \eqref{eq:bad} as follows 
\beq\label{eq:decomp}
\bal
- \la \td{u}_x \bar \th_x  + \td{u} \bar \th_{xx} , \th_x \psi \ra & + \lam_1 \la \om, \th_x \vp \ra
= \B( - \la \td{u}_x (\bar \th_x + \chi \xi_2 ) + \td{u} \bar \th_{xx} , \th_x \psi_f \ra
+ \lam_1 \la \th_x , \om \vp_f \ra \B) \\
& +  \B( - \la \td{u}_x ( \bar \th_x + \chi \xi_2 ) + \td{u} \bar \th_{xx} , \th_x \psi_n \ra
+ \lam_1 \la \th_x , \om \vp_n \ra \B) 
+ \lam_1 \la \th_x , \om \vp_s \ra  \\
&+ \la \td{u}_x  \chi (\xi_1  \psi_n + \xi_2 \psi_f ), \th_x \ra
\teq I_f + I_n + I_s + I_{r1} .
\eal
\eeq
The subscripts $f, n, s, r$ are short for \textit{far, near, singular, remainder}. Denote 
$I_{u\om} = -\lam_1 \la \td u \bar \om_x, \om \vp \ra$ in \eqref{eq:bad}. The main terms in \eqref{eq:bad} are $I_f, I_n$ and $I_s$. From the above discussion on \eqref{eq:unity}, the interactions in $I_n, I_f, I_s$ are mainly supported in different regions. Since $u$ depends on $\om$ linearly, $I_{u\om}$ can be seen as the interaction between $\om$ and itself. 
This type of interaction is different from $I_n, I_f, I_s$. Since $c_{\om} = u_x(0) = -\f{1}{\pi}\int_{\R} \om dx $, the terms $c_{\om} \la ( \bar{\th}_x - x\bar{\th}_{xx}) ,\th_x \psi \ra , \lam_1 c_{\om} \la  ( \bar{\om} - x \bar{\om}_x ) , \om \vp \ra$ in \eqref{eq:bad} are the projections of $\om, \th_x$ onto some rank-1 space. The estimate of the $c_{\om}$ terms is smaller than that of $I_n, I_f, I_s, I_{u\om}$. The term $I_{r1}$ is very small compared to other terms and will be estimated directly.

We will exploit the structure of the interactions in \eqref{eq:bad} using the above important decompositions.

\subsection{Outline of the estimates}


In order to establish the weighted $L^2$ estimates similar to \eqref{eq:goal}, we first develop sharp estimates on each term in the above decomposition. In these estimates, we introduce several parameters, when we apply the Cauchy-Schwarz or Young's inequality. These parameters are important in our estimates. Since the coefficients in the damping term $D_1, D_2$ \eqref{eq:damp} are relative small, we can treat the interaction term as perturbation to the damping term using the energy estimates, only for certain range of parameters. See more discussion in Remark \ref{rem:ineq_para}. Thus, the upper bound in these estimates depend on several parameters. Then, using these estimates, we reduce the estimate similar to \eqref{eq:goal} to some inequality constraints on the parameters with explicit coefficients. See \eqref{eq:goal_para2} for an example.
Finally, to obtain an overall sharp energy estimate, e.g. \eqref{eq:goal} with $c> 0$ as large as possible, we determine these parameters guided by solving a constrained optimization problem.

In our energy estimates, to obtain the sharp weighted estimates of $x u_x, u$ with singular weight $x^{-2p}$ by applying Lemma \ref{lem:wg}, we can only use a few exponents $p = 3,2 ,\f{5}{3}$. Thus, we need to perform the energy estimates very carefully. The linear combinations of different powers in Lemma \ref{lem:wg}, e.g. $\al x^{-4} + \b x^{-2}$, plays a role similar to that of interpolating different singular weights, e.g. $x^{-4}, x^{-2}$. It enables us to obtain sharp weighted estimates with singular weight $x^{-2q}$ and intermediate exponent $q$. In our weighted estimates of $u_x, u$, we choose some weights with a few parameters, see e.g. \eqref{eq:Su1}. Moreover, to generalize the cancellations and estimates in the Model 1 in Section \ref{sec:models} to the more complicated linearized system \eqref{eq:intro_lin}, we also need to perform the energy estimates carefully so that we can apply the cancellation in Lemma \ref{lem:cancel}.

\subsection{Estimates of the interaction in the near field $I_n$}\label{sec:fast}

We use ideas in Model 1 in Section \ref{sec:models} to estimate the main term introduced below and ideas in Section \ref{sec:models} to estimate $u$. 

Firstly, we choose $c = \f{1}{5}$ in the decomposition \ref{eq:decomp_u} $\td{u} = \f{1}{5} x \td{u}_x + \td{u} -\f{1}{5} x\td{u}_x $, and decompose $\td{u}_x (\bar{\th}_x + \chi \xi_1 ) + \td u \bar \th_{xx}$ into the main term $\cM$ and the remainder $\cR$ as follows
\beq\label{eq:decomp_u1}
\td{u}_x ( \bar{\th}_x + \chi \xi_1) + \td u \bar \th_{xx} = \td{u}_x ( \bar \th_x + \f{1}{5} \bar \th_{xx} x  + \chi \xi_1 ) + ( \td u - \f{1}{5} \td u_x x ) \bar \th_{xx} \teq \cM + \cR.
\eeq

This term also appears in $I_f$ and we will use another decomposition in Section \ref{sec:slow}.

Recall $I_n$ in \eqref{eq:decomp}. Using the above decomposition, we yield
\[
\bal
I_n  & =  - \la \td{u}_x ( \bar{\th}_x + \chi \xi_1) + \td u \bar \th_{xx}, \th_x \psi_n \ra 
+ \lam_1 \la \om, \th_x \vp_n \ra   \\
&= \B( - \la \cM, \th_x \psi_n \ra + \lam_1 \la \om, \th_x \vp_n \ra    \B) 
 + \la - \cR, \th_x \psi_n \ra \teq I_{\cM} + I_{\cR}. \\
\eal
\]

The estimates of $I_{\cM}$ are similar to that in Model 1 in Section \ref{sec:models}.
Recall the formulas of $\psi_n, \vp_n$ in \eqref{eq:mode}.  Using Young's inequality $ab \leq t_2 a^2 + \f{1}{4t_2} b^2$ for $t_2>0$, we obtain 
\beq\label{eq:varcoe_near}
\bal
I_{\cM} & = -\la  \td u_x , \th_x  ( \al_2 x^{-3} + \al_1 x^{-4} ) \ra
+ \la \om, \th_x \lam_1 \al_5  ( \al_2 x^{-2} + \al_1 x^{-3} ) \\
&= \la - \td u_x x^{-2} + \lam_1 \al_5 \om x^{-1}, \th_x ( \al_2 x^{-1} + \al_1 x^{-2}) \ra \\
& \leq t_2 || - \td u_x x^{-2} + \lam_1 \al_5 \om x^{-1}||_2^2 
+ \f{1}{4t_2} || \th_x ( \al_2 x^{-1} + \al_1 x^{-2})||_2^2 .
\eal
\eeq

\begin{remark}\label{rem:moti_wg}
We design the special form $ \psi_n$ in \eqref{eq:mode} so that the denominator in $\psi_n$ and the coefficient $ \bar \th_x + \f{1}{5} \bar \th_{xx} x  + \chi \xi_1 $ in $\cM$ cancel each other. This allows us to obtain a desirable term of the form $ J \teq - \td u_x x^{-2} + \lam_1 \al_5 \om x^{-1}$. The term $t_2 ||J||_2^2$ in \eqref{eq:varcoe_near} is a quadratic form in $\om$, where we can exploit the cancellation between $\td u_x$ and $\om$ to obtain a sharp estimate. 
See Model 1 for the motivation.
\end{remark}

Using the weighted estimate in Lemma \ref{lem:wg} and Lemma \ref{lem:cancel} with $f= \om$ and $g = u$,  we get 
\beq\label{eq:can_near}
\bal
&t_2 || - \td u_x x^{-2} + \lam_1 \al_5 \om x^{-1}||_2^2 
= t_2 \B( || \td u_x x^{-2}||_2^2 - 2 \lam_1 \al_5 \la \td u_x , \om x^{-3} \ra
+ (\lam_1 \al_5)^2 || \om x^{-1}||_2^2 \B) \\
 =& t_2 \B( || \om x^{-2}||_2^2 + (\lam_1 \al_5)^2 || \om x^{-1}||_2^2 \B)  = t_2 \B\la \om^2, x^{-4} + (\lam_1 \al_5)^2 x^{-2} \B\ra.
\eal
\eeq

The cancellation is exactly the same as \eqref{eq:model2_impro} in Model 1. For $I_{\cR}$, using  Young's inequality $ab \leq t_{22} a^2 + \f{1}{4 t_{22}} b^2$, 
\eqref{eq:IBPu} with $p=3$ and the weighted estimate in Lemma \ref{lem:wg}, we obtain 
\beq\label{eq:est_fast_IR}
\bal
I_{\cR}  & = \la ( \td u - \f{1}{5} \td u_x x ) \bar \th_{xx}, \th_x \psi_n \ra   \leq  t_{22} || (\td{u} -  \f{1}{5} \td{u}_x x ) x^{-3 } ||_2^2  + \f{1}{4 t_{22} } || x^{3} \bar \th_{xx}  \psi_n  \th_x ||_2^2 \\
&  = \f{ t_{22}}{25} || \td{u}_x x^{-2}||_2^2 + \f{1}{4 t_{22} } || x^{3} \bar \th_{xx}  \psi_n  \th_x ||_2^2   = \f{ t_{22}}{25} || \om x^{-2}||_2^2 + \f{1}{4 t_{22} } || x^{3} \bar \th_{xx}  \psi_n  \th_x ||_2^2 .
\eal
\eeq

The remainder $I_{\cR}$ is much smaller than $I_{\cM}$ since we get a small factor $\f{1}{2p-1}=\f{1}{5}$ from \eqref{eq:IBPu}. Combining the above estimates, we establish the estimate for $I_n = I_{\cM} + I_{\cR}$
\beq\label{eq:est_fast}
\bal
 I_n& \leq \B\la \om^2, t_2  x^{-4} + \f{t_{22}}{25} x^{-4} + t_2 (\lam_1 \al_5)^2 x^{-2}  \B\ra 
 + \B\la \th_x^2, \f{1}{4 t_2} (\al_2 x^{-1} + \al_1 x^{-2})^2 + \f{1}{4 t_{22}} (x^{3} \bar \th_{xx}  \psi_n )^2 \B\ra  .
 \eal
\eeq

\begin{remark}\label{rem:ineq_para}
If we neglect other terms in \eqref{eq:bad} except $I_n$, a necessary condition for \eqref{eq:goal} is
\beq\label{eq:goal_para1}
I_n + D_1 +  \lam_1 D_2 \leq -c ( || \th_x \psi^{1/2}||_2 + \lam_1 || \om \vp^{1/2}||_2^2 )
\eeq
with $c >0$, where $D_1, D_2$ are the damping terms in \eqref{eq:damp}.
 We cannot determine the ratio $\lam_1$ between two norms and $t_i$ in Young's inequality without using the profile $(\bar \th, \bar \om)$. For example, if we use equal weights $\lam_1 = 1, t_2 = t_{22} = \f{1}{2}$, we cannot apply estimate \eqref{eq:est_fast} to establish \eqref{eq:goal_para1} even with $c= 0 $. Therefore, we introduce several parameters, especially when we apply Young's inequality.  
At this step, we do not fix $\lam_1, t_{ij}$ such that the subproblem \eqref{eq:goal_para1} holds with $c>0$ as large as possible. In fact, such parameters may not be ideal for \eqref{eq:goal} since the final energy estimate involves other terms in \eqref{eq:goal},\eqref{eq:bad} to be estimated later on.
Instead, we identify the ranges $\lam_1 \in [0.31, 0.33] ,t_2 \in [5,5.8], t_{22} \in [13, 14]$, such that a weaker version of \eqref{eq:goal_para1} with $c= 0.01$ holds with the estimate \eqref{eq:est_fast} on $I_n$.
See Appendix \ref{app:range} for rigorous verification. Similarly, we will obtain the ranges of other parameters $t_i$ introduced in later estimates. We will determine the values of $\lam_1, t_{ij}$ in these ranges by combining the estimates of $I_f, I_n$ and other terms in \eqref{eq:bad}.

The estimates \eqref{eq:varcoe_near}, \eqref{eq:can_near} on the main term is crucial. If we estimate two inner products separately without using the cancellation between $\td u_x, \om$ in Lemma \ref{lem:cancel} with $f = \om$ and $g=u$, we would fail to establish \eqref{eq:goal_para1} even with $c=0$ since the damping term $D_i$ is relatively small. \end{remark}

\begin{remark}
Several key ideas in the above estimates will be used repeatedly later. Firstly, we will perform decompositions on $\td u$ into the main term and the remainder similar to \eqref{eq:decomp_u1}. Secondly, we will use Lemmas \ref{lem:cancel}, \ref{lem:rota} to estimate the inner product between $\td u$ and $\om$ similar to \eqref{eq:can_near}. Thirdly, we will use Lemma \ref{lem:wg} to estimate weighted norms of $\td u_x, \td u$ similar to \eqref{eq:can_near}. 
\end{remark}

\subsection{Estimates of the interaction in the far field $I_f$}\label{sec:slow}

We use ideas and estimates similar to that of $I_n$ to estimate $I_f$. The main difference is that to estimate the inner product between $\td u_x$ and $\om$, instead of using Lemma \ref{lem:cancel}, we will use Lemma \ref{lem:rota}. See estimates \eqref{eq:varcoe_near} and \eqref{eq:varcoe_far}.

Firstly, we choose $c = \f{3}{7}$ in \eqref{eq:decomp_u} and decompose $ \td u_x (\bar \th_x + \chi \xi_2) + \td u \bar \th_{xx}$ into the main term $\cM$ and the remainder $\cR$ as follows
\[
\td u_x (\bar \th_x + \chi \xi_2)+  \td{u} \bar \th_{xx} 
= \td u_x ( \bar \th_x + \f{3}{7} x \bar \th_{xx} + \chi \xi_2) 
+ (\td u - \f{3}{7} x \td u_x ) \bar \th_{xx} \teq \cM + \cR.
\]
We choose $c=\f{3}{7}$, which is different from that in \eqref{eq:decomp_u}, since we will apply \eqref{eq:IBPu} with a different power $p$ later. Recall $I_f$ in \eqref{eq:decomp}. The above formula implies
\[
\bal
I_f & =  - \la \td{u}_x (\bar \th_x + \chi \xi_2 ) + \td{u} \bar \th_{xx} , \th_x \psi_f \ra + \lam_1 \la \th_x , \om \vp_f \ra \\
& = \B( - \la \cM , \th_x \psi_f \ra + \lam_1 \la \om, \th_x \vp_f \ra   \B)
+  \la  -\cR,  \th_x \psi_f \ra  \teq I_{\cM} + I_{\cR}. 
\eal
\]

Recall the weights $\psi_f, \vp_f$ in \eqref{eq:mode}. Using Young's inequality 
$a \cdot b \leq t_1 a^2 + \f{1}{4 t_1} b^2$ for some $t_1>0$ to be determined,
we obtain 
\beq\label{eq:varcoe_far}
\bal
I_{\cM} &=  \la - \al_3 \td{u}_x x^{-4/3} + \lam_1 \al_6  \om x^{-2/3}, \th_x \ra 
= \la - \al_3 \td{u}_x x^{-1} + \lam_1 \al_6 \om x^{-1/3}, \th_x x^{-1/3} \ra  \\
& \leq t_1 || - \al_3 \td{u}_x x^{-1} + \lam_1 \al_6 \om x^{-1/3}||_2^2 
+ \f{1}{4t_1} || \th_x x^{-1/3} ||_2^2 \teq I_{\cM,1} + I_{\cM,2}. 
\eal
\eeq

We design the special form $\psi_f$ in \eqref{eq:mode} to obtain a desirable term of the form $- \al_3 \td{u}_x x^{-1} + \lam_1 \al_6 \om x^{-1/3}$. See also Remark \ref{rem:moti_wg}. We further estimate $I_{\cM,1}$. Applying  Lemma \ref{lem:wg} and Lemma \ref{lem:rota}, we derive 
\beq\label{eq:can_far}
\bal
I_{\cM,1} = & t_1 ( || \al_3 \td{u}_x x^{-1} ||_2^2 
 - 2 \al_3 \lam_1 \al_6 \la \td{u}_x , \om x^{-4/3} \ra + (\lam_1 \al_6)^2 || \om x^{-1/3}||_2^2 ) \\
 = & t_1 ( \al_3^2 || \om x^{-1}||_2^2 
-  \f{2 \al_3 \lam_1 \al_6}{2 \sqrt{3}} ( || \td{u}_x x^{-2/3}||_2^2 - || \om x^{-2/3} ||_2^2 )
+ (\lam_1 \al_6)^2 || \om x^{-1/3}||_2^2  ) \\
=&  t_1 \la \om^2, \al_3^2 x^{-2} + \f{\al_3 \lam_1 \al_6}{\sqrt{3}} x^{-4/3} 
+ ( \lam_1 \al_6)^2 x^{-2/3} \ra - \f{t_1 \al_3 \lam_1 \al_6}{\sqrt{3}} || \td{u}_x x^{-2/3}||_2^2.
 \eal
\eeq

\begin{remark}
The negative sign in $- t_1 2 \al_3 \lam_1 \al_6 \la \td{u}_x , \om x^{-4/3} \ra $ in \eqref{eq:can_far} is crucial. Firstly, we can bound the positive term $   \f{\al_3 \lam_1 \al_6 t_1}{\sqrt{3}} || \om x^{-2/3}||_2$ derived from the identity in Lemma \ref{lem:rota} directly without an overestimate. Secondly, $- \f{t_1 \al_3 \lam_1 \al_6}{\sqrt{3}} || \td{u}_x x^{-2/3}||_2^2$ from the same identity provides a good quantity that allows us to control the weighted norm of $\td{u}, \td{u}_x$ with a slowly decaying weight using Lemma \ref{lem:wg}.
\end{remark}

We introduce $D_u$ to denote the parameter in \eqref{eq:can_far} 
\beq\label{eq:Du}
D_{u} = \f{t_1 \al_3 \lam_1 \al_6}{\sqrt{3}}, 
\eeq
We use Young's inequality $a b \leq t_{12} a^2 + \f{1}{4 t_{12}} b^2$ for some $t_{12}>0$ and \eqref{eq:IBPu} with $p = \f{5}{3}$ to estimate $I_{\cR}$ directly 
\beq\label{eq:est_slow_I3}
\bal
I_{\cR}  & = - \la  ( \td u - \f{3}{7} x \td{u}_x),   \bar \th_{xx} \psi_f \th_x  \ra \leq  t_{12} || (\td{u} - \f{3}{7} \td{u}_x x ) x^{-5/3}||_2^2 + \f{1}{ 4 t_{12} }  || \th_x \psi_f  \bar \th_{xx} x^{5/3}||_2^2 \\
& = t_{12} \cdot \f{9}{49} || \td{u}_x x^{-2/3}||_2^2 + \f{1}{ 4 t_{12} }  || \th_x \psi_f  \bar \th_{xx} x^{5/3}||_2^2 .\\
\eal
\eeq

The remainder $I_{\cR}$ is smaller since we get a factor $\f{1}{2p-1} = \f{3}{7}$ from \eqref{eq:IBPu}.

Combining the above estimates, we obtain the estimate of $I_f = I_{\cM,1} + I_{\cM,2} + I_{\cR}$
\beq\label{eq:est_slow}
\bal
 I_f \leq & t_1 \B\la \om^2, \al_3^2 x^{-2} + \f{\al_3 \lam_1 \al_6}{\sqrt{3}} x^{-4/3}  + ( \lam_1 \al_6)^2 x^{-2/3} \B\ra \\
& + \B\la \th_x^2, \f{1}{4 t_1} x^{-2/3} + \f{1}{4 t_{12}} (\psi_f \bar \th_{xx} x^{5/3})^2 \B\ra
- \B( D_u - \f{9}{49} t_{12} \B) || \td{u}_x x^{-2/3}||_2^2 \ .
\eal
\eeq

Similar to the discussion in Remark \ref{rem:ineq_para}, in order for $I_f + D_1 + D_2 \leq 
-c ( || \th_x \psi^{1/2}||_2 + \lam_1 || \om \vp^{1/2}||_2^2 )$ with $c = 0.01$, we can choose $t_1 \in [1.2,1.4], t_{12} \in [0.55, 0.65]$. See Appendix \ref{app:range} for the verification.

\subsection{Estimates of the interaction with the most singular weight $I_s$}\label{sec:singu}

Recall $I_s $ in \eqref{eq:decomp} and $\psi_s = \al_4 x^{-4}$ 
in \eqref{eq:mode}. Using Young's inequality $ab \leq t_4 a^2 + \f{1}{4t_4} b^2$ for $t_4 > 0$, we yield
\beq\label{eq:est_singu}
I_s =  \lam_1 \la \om, \th_x \vp_s \ra = \lam_1 \al_4 \la \om, \th_x x^{-4} \ra
\leq t_4 \la \om^2, x^{-3}\ra 
+ \f{ (\lam_1 \al_4)^2}{4 t_4} \la \th_x^2,  x^{-5} \ra .
\eeq

In order for $I_s + D_1 + D_2 \leq - c ( || \th_x \psi^{1/2}||_2 + \lam_1 || \om \vp^{1/2}||_2^2 )$ with $c = 0.01$, we can choose $t_4 \in [3,5]$. See Appendix \ref{app:range} for the verification. We do not combine estimates of $I_s$ with the estimates for the interaction between $\td u$ and $\th_x$ in Section \ref{sec:fast} since the weight $x^{-4}$ is too singular. In fact, to apply estimate similar to that in \eqref{eq:est_coupled} in Model 1, the weight for $\th_x$  near $0$ is $2$ order more singular than that of $\om$. In this case, it is of order $x^{-6}$ near $0$ and more singular than $\psi$. 
\subsection{Estimates of the interaction $u$ and $\om$}\label{sec:est_uw}

Firstly, we rewrite $ - \lam_1 \la \td{u} , \bar \om_x  \om \vp \ra $ in \eqref{eq:bad} by decomposing $\bar \om_x$ into the main term $ \bar \om_x+ \chi \xi_3$ and the error term $\xi_3$
\beq\label{eq:chi_term2}
I_{u \om}  \teq - \lam_1 \la \td{u} , \bar \om_x  \om \vp \ra 
= - \lam_1 \la \td{u} , ( \bar \om_x + \chi \xi_3) \om \vp \ra + \lam_1 \la \td{u} , \chi \xi_3 \om \vp \ra \teq  J + I_{r2}.
\eeq

We will estimate $I_{r2}$ in Section \ref{sec:chi} and show that it is very small. See also Remark \ref{rem:chi}.

\vspace{0.1in}
\paragraph{\bf{Difficulty}}
The main difficulty in establishing a sharp estimate for $J$ is the slow decay of the coefficient $( \bar \om_x + \chi \xi_3)  \vp$. A straightforward estimate similar to \eqref{eq:est_slow_I3} yields $ |J| \leq \lam_1 || \td u x^{-p}||_2 || \bar \om_x \om \rho x^p ||_2 $. 
In view of the weighted estimate in Lemma \ref{lem:wg}, we have effective estimates of $ (u - u_x(0) x) x^{-p}$ for exponent $p= 3,2$ or $\f{5}{3}$. In order to further control $ || \bar \om_x \om \vp x^p ||_2$ by the weighted $L^2$ norm $|| \om \vp^{1/2}||_2$, we cannot choose $ p =3$ or $p=2$ due to the fast growth of $x^p$ for large $|x|$. On the other hand, if we choose $p=\f{5}{3}$, the resulting constant $\f{36}{49}$ in Lemma \ref{lem:wg} is much larger than the constant $\f{4}{25}, \f{4}{9}$ corresponding to $p=3, p=2$. 

To overcome this difficulty, we exploit the cancellations between $u$ and $\om$ in both the near field and the far field, which is similar to that in the estimate of $I_n, I_f$.
We decompose the coefficient $(\bar \om_x + \chi \xi_3)  \vp$ in $J$ into the main terms $\cM_i$ and the remainder $\cK_{u\om}$
\beq\label{eq:Kuw}
\bal
( \bar \om_x + \chi \xi_3)  \vp &= 
- \f{1}{3} e_3 \al_6 x^{-2}   +\tau_1 x^{-4} + 
\B( ( \bar \om_x  + \chi \xi_3 )\vp + \f{1}{3} e_3 \al_6 x^{-2} - \tau_1 x^{-4}  \B) \\
&\teq \cM_1 + \cM_2 + \cK_{u\om},
\eal
\eeq
where $e_3, \al_6$ are defined in \eqref{eq:ext} and \eqref{eq:mode} and $\tau_1 >0$ is some parameter.

Let us motivate the above decomposition. From the definitions of $\xi_3, \vp$ in \eqref{eq:ext}-\eqref{eq:wg} and the discussion therein, we have $\bar \om_x  + \chi \xi_3 \approx - \f{1}{3} e_3 x^{-4/3} , \vp \approx \al_6 x^{-2/3}$ for some $e_3,\al_6>0$ and large $|x|$. 
Thus, $( \bar \om_x  + \chi \xi_3 )\vp$ can be approximated by $-\f{1}{3} e_3 \al_6 x^{-2}$ for large $|x|$. Since $\vp \approx \al_4 x^{-4}$ and $\bar \om_x \approx \bar \om_x(0) > 0$ near $0$, $ ( \bar \om_x  + \chi \xi_3 )\vp  $ is approximated by $\tau_1 x^{-4}$ for some $\tau_1>0$ in the near field.

Using the above formula, we can decompose $J$ as follows 
\[
J = - \lam_1 \la \td{u} , ( \bar \om_x + \chi \xi_3) \om \vp \ra 
= -  \lam_1 \la  \td u, \cM_1  \om \ra - \lam_1 \la \td u, \cM_2 \om \ra 
- \lam_1 \la \td u , \cK_{u \om} \om \ra \teq I_{\cM1} + I_{\cM2} + I_{\cR}.
\]
 To estimate the main terms, we use cancellations in Lemma \ref{lem:cancel}. Using $\td u = u - u_x(0) x$ defined in \eqref{eq:ut}, 
$ -u_x(0) \la x, \om x^{-2}\ra = \f{\pi}{2} u_x^2(0)$ and Lemma \ref{lem:cancel} with $f= \om$ and $g = u$, we get 
\[
I_{\cM_1} =  \f{ \lam_1 e_3 \al_6 }{3}  \la \td{u}, \om x^{-2} \ra =\f{ \lam_1 e_3 \al_6 }{3}  \f{\pi}{4} u_x(0)^2  -  \f{\lam_1 e_3 \al_6 }{3}   \la \Lam \f{u}{x},  \f{u}{x} \ra ,
\]
where $\Lam =  (-\pa_x^2)^{1/2}$. We denote by $A(u)$ the right hand side of the above equation 
\beq\label{eq:Au}
A(u) \teq \f{ \lam_1 e_3 \al_6 }{3}  \f{\pi}{4} u_x(0)^2  -  \f{\lam_1 e_3 \al_6 }{3}   \la \Lam \f{u}{x},  \f{u}{x} \ra  .
\eeq

Since $e_3 \al_6 \lam_1 >0$ and $ \la \Lambda \f{u}{x}, \f{u}{x} \ra \geq 0$, the second term in $A(u)$ is a good term and we will use it in the weighted $H^1$ estimate.

Although $I_{\cM 2}$ is a quadratic form on $\om$, it does not have a good sign similar to the identities in Lemma \ref{lem:cancel}. Yet, we can approximate $\td u$ by $\td u_x$ using \eqref{eq:decomp_u} and then use the cancellation between $\td u_x$ and $\om$.
Choosing $c = \f{1}{5}$ in \ref{eq:decomp_u} and using the cancellation in Lemma \ref{lem:cancel}
with $f = \om$ and $g = u$, we obtain 
\[
I_{\cM2} =-\lam_1 \tau_1 \la \td{u}, x^{-4} \om \ra = -\lam_1 \tau_1 \la \td{u} - \f{1}{5} \td{u}_x x, \om x^{-4} \ra
- \lam_1 \tau_1  \la \f{\td{u}_x}{5} , \om x^{-3} \ra 
= -\lam_1 \tau_1 \la \td{u} - \f{1}{5} \td{u}_x x, \om x^{-4} \ra .
\]

The form $\td{u} - \f{1}{5} \td{u}_x x$ allows us to gain a small factor $\f{1}{5}$ using \eqref{eq:IBPu} with $p=3$. 
Using Young's inequality $ab \leq ca^2 + \f{1}{4c} b^2$, \eqref{eq:IBPu} with $p=3$
and Lemma \ref{lem:wg}, we obtain 
\beq\label{eq:est_Iuw}
\bal
I_{\cM2}  & \leq \lam_1 \tau_1 \B( t_{34} || ( \td{u} - \f{1}{5} \td{u}_x ) x^{-3}||_2^2 
+ \f{1}{4 t_{34}} || \om x^{-1 } ||_2^2    \B)  =  \lam_1 \tau_1 \B(  \f{t_{34}}{25} || \td{u}_x x^{-2} ||_2^2 
  + \f{1}{4 t_{34}} || \om x^{-1 } ||_2^2\B)  \\
  &= \lam_1 \tau_1 ( \f{ t_{34}}{25 } || \om x^{-2} ||_2^2 
+ \f{1}{4 t_{34}} || \om x^{-1 } ||_2^2    ) 
 = \lam_1 \tau_1 \B\la \om^2, \f{t_{34}}{25} x^{-4} + \f{1}{4 t_{34}} x^{-2} \B\ra,
\eal
\eeq
for some $t_{34} >0$. For $t_{31}, t_{32} >0$ to be defined, we denote 
\beq\label{eq:Su1}
S_{u1} = t_{31} x^{-6} + t_{32} x^{-4} +  2 \cdot 10^{-5} x^{-10/3} , 
\eeq

 We estimate $I_{\cR}$ directly using Young's inequality and the weighted estimate in Lemma \ref{lem:wg}
\beq\label{eq:est_kuw}
\bal
I_{\cR} & = -\lam_1 \la \td{u}, \cK_{u \om} \om \ra
\leq \lam_1 ( || \td{u} S_{u1}^{1/2} ||_2^2  + \f{1}{4} || S_{u1}^{-1/2} \cK_{u\om} \om ||_2^2 )  \\
& \leq \lam_1 \B\la \om^2,  \f{4 t_{31}}{25} x^{-4} + \f{4t_{32} }{9} x^{-2}  
+ \f{1}{4} \cK_{u\om}^2 S_{u1}^{-1}  \B\ra + \f{36 \lam_1  }{49} \cdot 2 \cdot 10^{-5}|| \td{u}_x x^{-2/3}||_2^2,
\eal
\eeq
where $\cK_{u\om}$ in defined \eqref{eq:Kuw}.

\begin{remark}\label{rem:est_uw_S1}
From \eqref{eq:ext},\eqref{eq:mode} and \eqref{eq:wg}, we have asymptotically $\cK_{u\om} \sim C x^{-4}$ for $x$ close to $0$. The slowly decaying part in $\cK_{u\om}$ is given by $ f =
( \bar \om_x  + \chi \xi_3 )\al_6 x^{-2/3} + \f{1}{3} e_3 \al_6 x^{-2}  =  (1 - \chi) (\bar \om_x + \f{1}{3} e_3x^{-4/3}) \al_6 x^{-2/3}$. In the support of $1 - \chi$, since $-\f{1}{3} e_3x^{-4/3}$ approximates $\bar \om_x$, $f$ can be approximated by $c x^{-2}$ with some very small constant $c$. We add  $x^{-6}$,  $  2 \cdot 10^{-5} x^{-10/3}$ in $S_{u1}$ so that $\cK_{u\om}^2 S_{u1}^{-1}$ can be bounded by $\vp$. 
We also add the power $x^{-4}$ in $S_{u1}$ to obtain a sharper estimate.
See also Remark \ref{rem:ineq_para} for the discussion on the parameters.

\end{remark}

Combining the above estimates on $I_{\cM i}, I_{\cR}$, we prove 
\beq\label{eq:est_uw}
\bal
I_{u \om} & = J + I_{r2}\leq 
\lam_1 \B\la  \om^2,   \f{4 t_{31}}{25} x^{-4} + \f{4t_{32} }{9} x^{-2}  
+ \f{1}{4} \cK_{u\om}^2 S_{u1}^{-1}
+ \tau_1 ( \f{t_{34}}{25} x^{-4} + \f{1}{4 t_{34}} x^{-2})\B\ra  \\
& \quad + A(u)  +  \f{72 \lam_1  }{49} \cdot 10^{-5} || \td{u}_x x^{-2/3}||_2^2 + I_{r2}.
\eal
\eeq

The term $I_{r2}$ was not estimated and we keep it on both sides. We can determine the ranges of parameters $t_{31}, t_{32}, t_{34}, \tau_1$ so that $ J + D_1 + D_2 \leq -0.01 ( || \th_x \psi^{1/2}||_2 + \lam_1 || \om \vp^{1/2}||_2^2 )$. 

\subsection{Estimates of the $I_{r1}, I_{r2}$}\label{sec:chi} Recall $I_{r1}, I_{r2}$ in \eqref{eq:decomp} and \eqref{eq:chi_term2}. Since $\chi$ is supported in the far field $|x|\geq \rho_2 > 10^{8} $ and the profile $(\bar \om, \bar \th_x)$ decays, we can get a small factor in the estimate of these terms. We establish the following estimate in Appendix \ref{app:chi}
 \beq\label{eq:est_chi} 
|I_{r1} | + |I_{r2} |\leq  \la G_{\th}, \th_x^2 \ra + \la G_{\om} , \om^2 \ra 
+  G_{c} c_{\om}^2,
\eeq
where $G_{\th}, G_{\om}, G_{c }$ are given by 
\beq\label{eq:est_chi2}
\bal
G_{\th} &= 10^{10} \cdot \f{ (2+\sqrt{3})^2}{4 }  \chi^2 (\xi_1  \psi_n + \xi_2 \psi_f )^2 ,\quad 
G_{c} = \f{ \lam_1^2 || x \xi_3 \chi^{1/2} \vp^{1/2}||_2^2}{ 4 } \cdot 10^2, \\
G_{\om} &= 10^{-10} x^{-4/3} + 10^{-5} x^{-2/3} + 
\f{10^5}{4  }  ( \f{ 6\lam_1 (2 + \sqrt{3})}{5}  )^2 ( x^{4/3} \chi \xi_3 \vp)^2 
+ 10^{-2} \chi \vp  .
\eal
\eeq
These functions are very small compared to the weights $\vp, \psi$ \eqref{eq:mode}-\eqref{eq:wg}. We focus on a typical term $G_{\th}$ to illustrate the smallness. From \eqref{eq:ext}-\eqref{eq:mode}, for large $|x|$, $\xi_i, \psi_f, \psi$ have decay rate $x^{-2/3}$, $\psi_n$ has a decay rate at least $x^{-2}$. For $ \pa_x^i  \bar \th_x$, we recall from the beginning of Section \ref{sec:linop} that it has decay rate $x^{ 2\al - i}$ with $\al$ slightly smaller than $-\f{1}{3}$. Thus $|\chi (\xi_1 \psi_n + \xi_2 \psi_f) |^2 /\psi$ has decay rate $x^{-2}$. Since $\chi$ is supported in $|x| \geq \rho_2 > 10^8$, we get a small factor $x^{-2} \one_{|x| > \rho_2}<10^{-16}$, which absorbs the large constant $10^{10}$ in $G_{\th}$. Therefore, $G_{\th}$ is very small compared to $\psi$.

\subsection{ Summary: estimates of $\cL_{\th1}, \cL_{\om 1}$}

Recall $\cL_{\th1}, \cL_{\om1}$ in \eqref{eq:l2h1}, the quadratic terms in \eqref{eq:l2h1}, \eqref{eq:bad}. Combining \eqref{eq:damp},\eqref{eq:est_slow},\eqref{eq:est_fast}, \eqref{eq:est_singu}, \eqref{eq:est_uw} and \eqref{eq:est_chi}, we obtain the estimate on $\cL_{\th1}, \cL_{\om1}$
\beq\label{eq:coer1}
\bal
\la \cL_{\th 1} \th_x, \th_x  \psi \ra
&+ \lam_1 \la \cL_{\om 1} \om, \om \vp \ra 
\leq   \la  D_{\th} +   A_{\th} \psi^{-1}, \th_x^2 \psi \ra  + \la  \lam_1 D_{\om} +  A_{\om}  \vp^{-1} ,\om^2 \vp \ra  \\
&- \B( D_u  - \f{9}{49} t_{12} - \f{72 \lam_1  }{49} \cdot 10^{-5}\B) || \td{u}_x x^{-2/3}||_2^2   +  A(u) + G_{c } c_{\om}^2 ,
\eal
\eeq
where $A(u)$ is defined in \eqref{eq:Au}, $D_u, t_{12}$ are given in \eqref{eq:para2},
and the $A_{\th}, A_{\om}$ terms are the sum of the integrals of $\om^2, \th_x^2$ in the upper bounds in \eqref{eq:est_slow},\eqref{eq:est_fast}, \eqref{eq:est_singu}, \eqref{eq:est_chi} given by
\beq\label{eq:coer1_cost}
\bal
A_{\th} &\teq 
 \lt( \f{1}{4 t_1} x^{-2/3} + \f{1}{4 t_{12}} (\psi_f \bar \th_{xx} x^{5/3})^2 \rt)
+ \lt( \f{1}{4 t_2} (\al_2 x^{-1} + \al_1 x^{-2})^2 + \f{1}{4 t_{22}} (x^{3} \bar \th_{xx}  \psi_n )^2 \rt)   + \f{ (\lam_1 \al_4 )^2}{4 t_4} x^{-5}  + G_{\th}, \\
A_{\om} & \teq t_1 \lt( \al_3^2 x^{-2} + \f{\al_3 \lam_1 \al_6}{\sqrt{3}} x^{-4/3}  + ( \lam_1 \al_6)^2 x^{-2/3}  \rt) + \lt( t_2  x^{-4} + \f{t_{22}}{25} x^{-4} + t_2 (\lam_1 \al_5)^2 x^{-2}  \rt) \\
&+ t_4 x^{-3}  + \lam_1 \lt(  \f{4 t_{31}}{25} x^{-4} + \f{4t_{32} }{9} x^{-2}  
+ \f{1}{4} \cK_{u\om}^2 S_{u1}^{-1}
+ \tau_1 ( \f{t_{34}}{25} x^{-4} + \f{1}{4 t_{34}} x^{-2}) \rt) + G_{\om}.
\eal
\eeq

In the previous estimates, we have obtained the ranges of $t_{ij}$ such that $ I + D_1 + D_2 \leq 
-0.01 ( || \th_x \psi^{1/2}||_2^2 + \lam_1 || \om \vp^{1/2}||_2^2)$ for several terms $I$ in \eqref{eq:bad}, e.g. $I = I_f, I_n$. We further determine the approximate values of  $\lam_1, t_{ij}$ so that 
\beq\label{eq:goal_para2}
D_{\th} + A_{\th} \psi^{-1} \leq -c , \quad \lam_1 D_{\om} + A_{\om} \vp^{-1} \leq -\lam_1 c
\eeq
with $c>0$ as large as possible. The functions in \eqref{eq:goal_para2} depend on the parameters and other explicit functions. The above task is equivalent to solving a constrainted optimization problem by maximizing $c$, subject to the constraints \eqref{eq:goal_para2} and $\lambda_1, t_{ij} > 0$ within an interval that we have determined.

Estimates \eqref{eq:coer1}, \eqref{eq:goal_para2} imply the linear stability estimate \eqref{eq:goal} up to  
the $c_{\om}$ terms in \eqref{eq:bad}, $A(u)$ \eqref{eq:Au} and $ G_c c_{\om}^2$.
In Section \ref{sec:est_cw}, we further control these $c_{\om}$ terms.
The estimate of these $c_{\om}$ terms are small. We will perturb $\lam_1,t_{ij}$ around their approximate values and finalize the choices of $\lam_1,t_{ij}$. The final values of these parameters are given in \eqref{eq:para2}, \eqref{eq:para3}.

The main reasons that we can establish \eqref{eq:goal_para2} are the followings. Firstly, we exploit several cancellations using Lemma \ref{lem:wg} and apply sharp weighted estimates in Lemma \ref{lem:wg} to estimate the nonlocal terms. Secondly, we have $ I + D_1 + D_2 \leq - c( || \th_x \psi^{1/2}||_2^2 + \lam_1 || \om \vp^{1/2}||_2^2)$ for $I$ being the main terms in \eqref{eq:bad}, i.e. $I = I_f, I_n$ or $I_s$. Thirdly, the interactions in $I_f, I_n, I_s$ are mainly supported in different regions. See the discussion after \eqref{eq:decomp}. 
Finally, the main term in $I_{u\om}$ is estimated using several cancellations. 

\begin{figure}[t]
\centering
    \begin{subfigure}[b]{0.48\textwidth}
    	\includegraphics[width=1\textwidth]{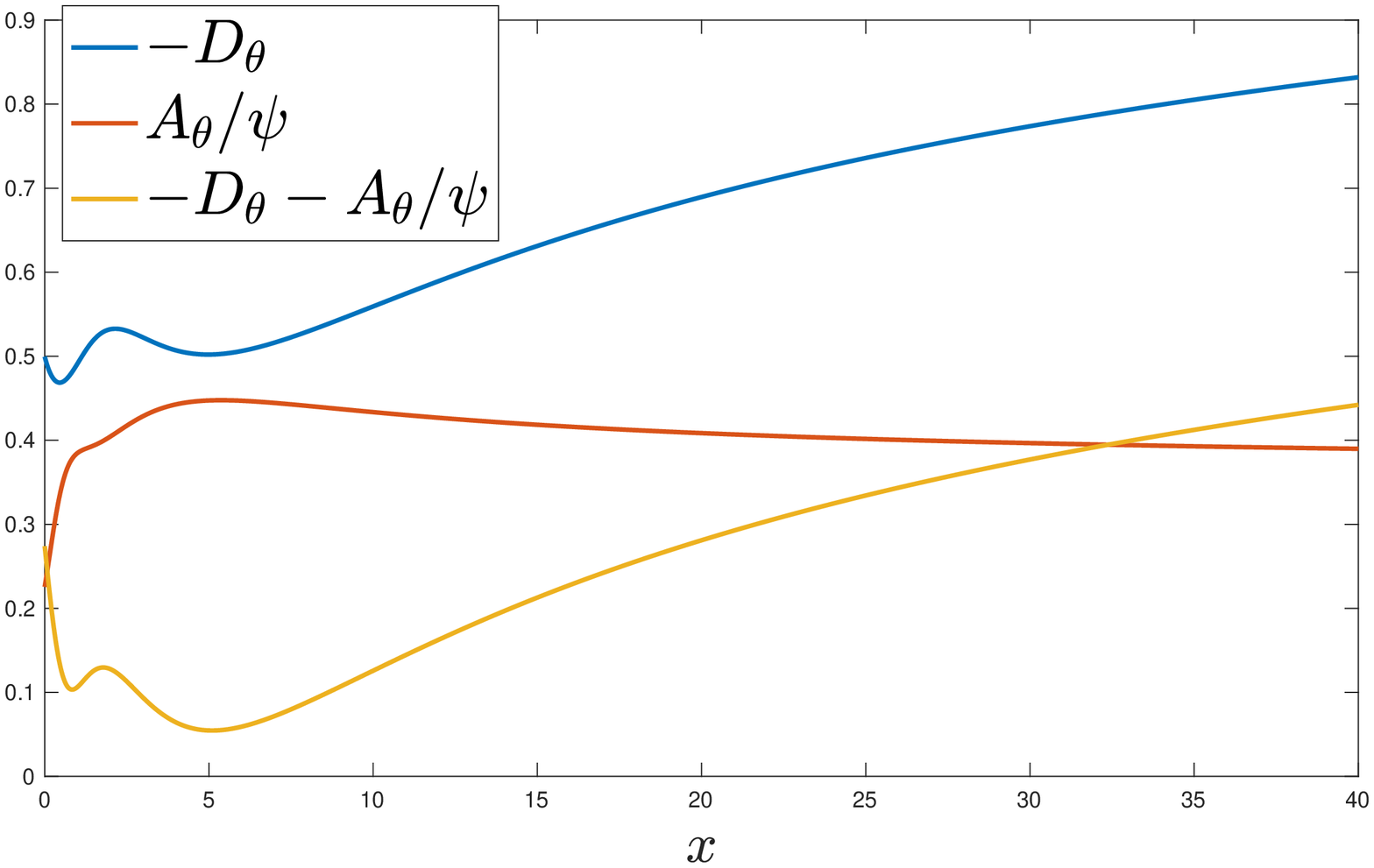}
    \end{subfigure}
    \begin{subfigure}[b]{0.48\textwidth}
    	\includegraphics[width=1\textwidth]{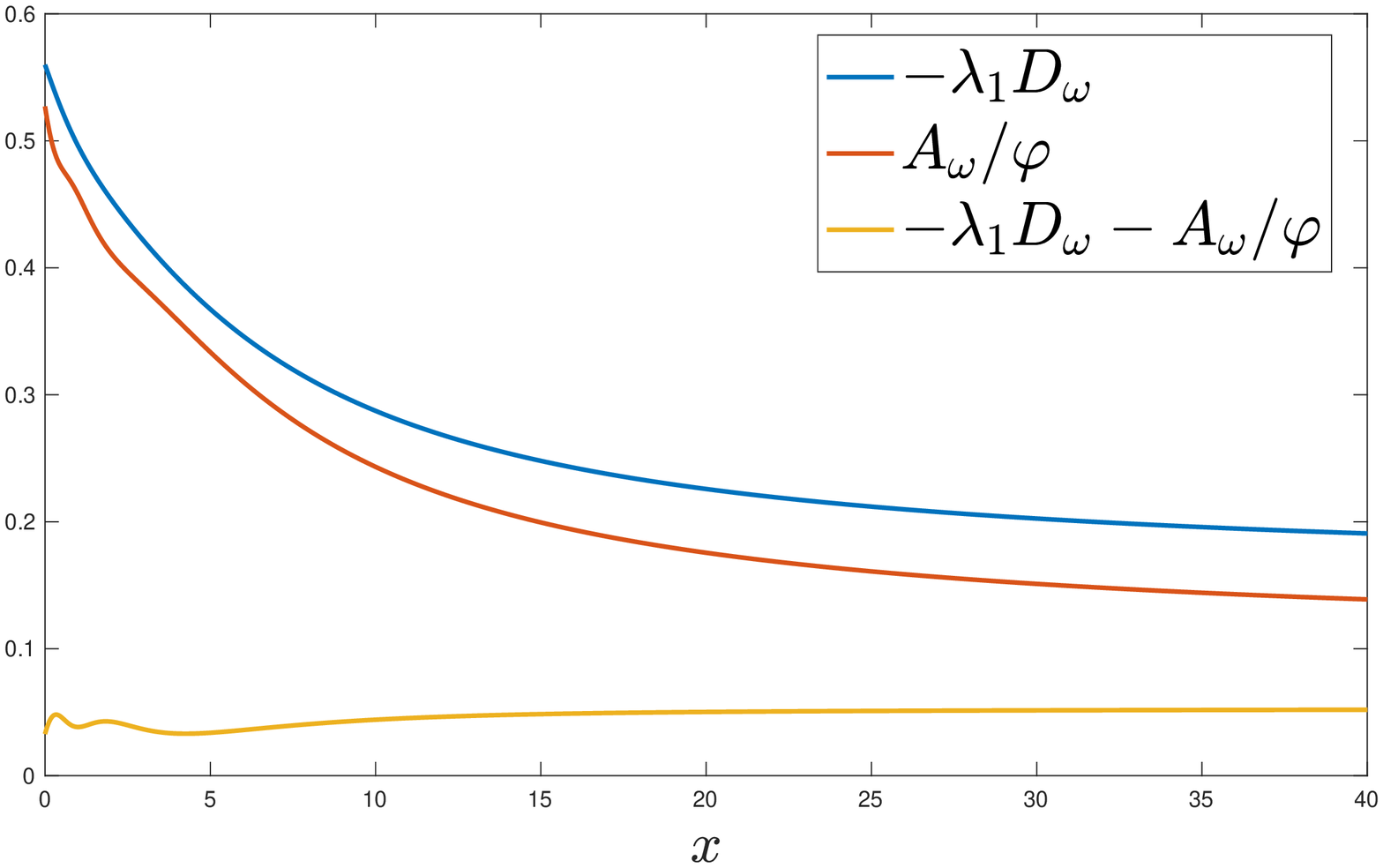}
    \end{subfigure}
    \caption[Illustration of linear stability estimates]
    { Left : Grid point values of  $-D_\th$, $A_\th\psi^{-1}$ and $-D_\th-A_\th\psi^{-1}$ for $x\in [0,40]$. Right: Those of $-\lambda_1D_\om$, $A_\om\phi^{-1}$ and $-\lambda_1D_\om-A_\om\phi^{-1}$ (with $\lambda_1 = 0.32$). }
         \label{fig:linear_stability}
    \vspace{-0.1in}
\end{figure}

To illustrate that the inequalities in \eqref{eq:goal_para2} can actually be achieved, we plot in Figure \ref{fig:linear_stability} the grid point values of the functions $-D_{\th} - A_{\th} \psi^{-1}$ and $-\lam_1 D_{\om} - A_{\om} \vp^{-1}$ for $x\in[0,40]$ with the parameters $\lambda_1, t_{ij}$ given in \eqref{eq:para2}, \eqref{eq:para3}. 
It is shown that their grid point values are all positive and uniformly bounded away from $0$. In fact, the minimum of the grid point values of $-D_\th-A_\th\psi^{-1}$ is above $0.032$ and that of $-\lambda_1D_\om-A_\om\phi^{-1}$ is above $0.054$.

Estimate \eqref{eq:coer1} on $\cL_{\th1}, \cL_{\om1}$ is important and we will use it in Section  
\ref{sec:non} to establish the weighted $H^1$ estimates.

\subsection{Estimate of the $c_{\om}$ term}\label{sec:est_cw}

We use the idea in Model 2 in Section \ref{sec:models} to obtain the damping term for $c_{\om}$ by deriving the ODEs for $c_{\om}^2$ and $\la \th_x, x^{-1}\ra^2$. We introduce some notations
\beq\label{eq:dth0}
d_{\th} \teq \la \th_x, x^{-1} \ra, \quad  \bar d_{\th} \teq \la \bar \th_x , x^{-1} \ra,
\quad \bar u_{\th ,x} \teq  H \bar \th_x , \quad u_{\D} = \td{u} - \f{1}{5} \td{u}_x x .
\eeq

\subsubsection{Derivation of the ODEs} 
Recall $c_{\om} = u_x(0)  = -\f{2}{\pi}  \int_0^{+\infty} \f{\om}{x} dx$ from \eqref{eq:normal}. 
Using a derivation similar to that in Model 2 in Section \ref{sec:models}, , we derive the following ODE in Appendix \ref{app:ode}
\beq\label{eq:cw2}
\bal
\f{1}{2} \f{d}{dt} \f{\pi}{2} c^2_{\om} &=   \f{\pi}{2} (\bar{c}_{\om} + \bar{u}_x(0)) c^2_{\om} +
 c_{\om}  \int_0^{\infty}  \f{  \bar{u} \om_x + u \bar \om_x } {x} dx
- c_{\om} d_{\th}   - c_{\om}  \int_0^{\infty} \f{ F_{\om} + N(\om)}{x} dx.
\eal
\eeq

The ODE for $d_{\th}^2$ \eqref{eq:dth2} is derived similarly in Appendix \ref{app:ODE_dth}. There is a cancellation between these two ODEs, which is captured by Model 2 in Section \ref{sec:models}. To exploit this cancellation, we combine two ODEs and derive the following ODE in Appendix \ref{app:ode} with $\lam_2, \lam_3 >0$ to be chosen
\beq\label{eq:ODE}
\bal
\f{1}{2} \f{d}{dt} ( \f{\lam_2 \pi }{2} c_{\om}^2 + \lam_3 d_{\th}^2 )
& = \f{\pi \lam_2 }{2} (\bar c_{\om} + \bar u_x(0)) c^2_{\om} +  2 \bar c_{\om} \lam_3 d_{\th}^2
+ \cT_0 + \cR_{ODE},
\eal
\eeq
where $\cT_0$ is the sum of the quadratic terms that do not have a fixed sign 
\beq\label{eq:cT0}
\bal
\cT_0  = &  - (\lam_2 - \lam_3 \bar d_{\th}) c_{\om} d_{\th}   + \lam_2 c_{\om} \la \om, f_2\ra
- \lam_3 d_{\th} \la \th_x , f_6 \ra  + \lam_3 d_{\th} \la \om, f_4 \ra \\
&+ \lam_2 c_{\om} \la u_{\D} x^{-1}, f_8 \ra - \lam_3 d_{\th} \la u_{\D} x^{-1}, f_9 \ra ,
\eal
\eeq
$u_{\D}$ is defined in \eqref{eq:dth0}, $f_i$ defined in \eqref{eq:func} are some functions depending on the profile $(\bar \om, \bar \th)$, and  $\cR_{ODE}$ is the sum of the remaining terms in the ODEs given by 
\beq\label{eq:Rode}
\cR_{ODE} \teq  - \lam_2 c_{\om} \la F_{\om} + N(\om), x^{-1} \ra + \lam_3 d_{\th} \la F_{\th} + N(\th), x^{-1} \ra. 
\eeq

Since the approximate steady state satisfies $\bar{c}_{\om} <0, \bar{u}_x(0) <0$, $\f{\pi \lam_2 }{2} (\bar c_{\om} + \bar u_x(0)) c^2_{\om}$ and  $2 \bar c_{\om} \lam_3 d_{\th}^2 $ in \eqref{eq:ODE} are damping terms. 

\subsubsection{Derivation of the $\cT_0$ term}\label{sec:model3_nonlocal}

Let us explain how we obtain \eqref{eq:ODE}. The ODEs of $c_{\om}^2, d_{\th}^2$ (\eqref{eq:cw2} and \eqref{eq:dth2}) involves the integrals of the nonlocal terms $u, u_x$ in the form of $\la \td{u}, f \ra$ or $\la \td{u}_x,  g \ra$ for some functions $f,g$. To estimate these terms effectively, we use the antisymmetry property of the Hilbert transform in Lemma \ref{lem:anti} to transform these terms into integrals of $\om$. We first consider $\la \td{u}_x,  g\ra$ and $\la u_x ,  g\ra$. Using $ u_x = H \om$, $\f{u_x - u_x(0)}{x} = H( \f{\om}{x})$ and Lemma \ref{lem:anti}, we get 
\beq\label{eq:IBP_int}
\bal
\la u_x ,  g \ra &= \la H \om , g \ra = - \la \om,  H   g \ra , \quad
\la \td{u}_x ,  g \ra  = \B\la H \lt( \f{\om}{x} \rt), x   g \B\ra=  - \B\la \f{\om}{x}, H( x  g) \B\ra. 
\eal
\eeq

For $\la \td{u},  f \ra$, we first approximate $ f$ by $ p_x $ for some function $p$ and then perform a decomposition $\td{u} = c x \td{u}_x + (\td{u} - c x \td{u}_x)$. We obtain 
\[
\la \td{u}, f \ra = \la \td{u}, p_x \ra + \la \td{u}, f - p_x \ra 
= \la \td{u}, p_x \ra + \la c x \td{u}_x, f - p_x \ra + \la \td{u} - c x \td{u}_x, f - p_x \ra
\teq I_1 + I_2 + I_3.
\] 

The last term enjoys much better estimate than $\la \td{u}, f\ra$ due to \eqref{eq:IBPu} and the fact that $f - p_x$ is much smaller than $f$. For $I_1, I_2$, using integration by parts, we get 
\[
I_1 + I_2 = \la \td{u}_x, -p + cx(f -p_x) \ra .
\]

Using \eqref{eq:IBP_int}, we can further rewrite the above term as an integral of $\om$. 

In addition, we introduce the function $f_i$ to simplify the integrals of $\om, \th_x$. These derivations lead to the $\cT_0$ term. We refer the details to Appendix \ref{app:ode}.

We remain to estimate the $c_{\om}$ terms in \eqref{eq:l2h1} in the weighted $L^2$ estimates and $f_3, f_7$ that are defined in \eqref{eq:func}.  We combine $\cT_0$ and these $c_{\om}$ terms, and define  
\beq\label{eq:cT}
\cT \teq  \cT_0 + c_{\om} \la  \bar{\th}_x - x\bar{\th}_{xx}  , \  \th_x \psi \ra + \lam_1 c_{\om} \la  \bar{\om} - x \bar{\om}_x  , \om \vp \ra  
= \cT_0 + c_{\om} \la \om, f_3 \ra + c_{\om} \la \th_x , f_7 \ra.
\eeq

In the weighted $L^2$ estimates, it remains to estimate $\cT$.
Though each term in $\cT$ can be estimated by the weighted $L^2$ norms of $\om, \th_x$ and $ c^2_{\om}, \la \th_x, x^{-1}\ra^2$ using the Cauchy-Schwarz inequality, these straightforward estimates do not lead to sharp estimates since these Cauchy-Schwarz inequalities do not achieve (or are close to) equalities for the same functions. We use the optimal-constant argument in \cite{chen2019finite} to obtain a sharp estimate on $\cT$.

\subsubsection{ Sharp estimates on $\cT$}\label{sec:cT}

For positive $T_1 , T_2, T_3 \in C(\R_+) $ and positive parameter $s_1, s_2 > 0$ to be determined, we consider the following inequality with sharp constant $C_{opt}$
\beq\label{eq:Copt0}
\cT \leq C_{opt} ( || \om T_1^{1/2} ||_2^2 + || \th_x T_2^{1/2} ||_2^2 + || \f{u_{\D}}{x} T_3^{1/2} ||_2^2 + s_1 c^2_{\om} + s_2 d^2_{\th} ),
\eeq
where $u_{\D}$ is defined in \eqref{eq:dth0}. We define several functions 
\beq\label{eq:func2}
\bal
X  &= \om T_1^{1/2}, \quad Y = \th_x T_2^{1/2}, \quad Z  =  u_{\D} x^{-1} T_3^{1/2} , \\
g_1 & = -\f{2}{\pi} x^{-1 } T_1^{-1/2}, \ g_2 =  f_2  T_1^{-1/2}, \  g_3 =  f_3 T_1^{-1/2}, \ g_4 = f_4 T_1^{-1/2} ,  \\
g_5 & = x^{-1} T_2^{-1 / 2},  \ g_6 = f_6 T_2^{-1/2} , \ g_7 = f_7 T_2^{-1/2} , \  g_8  = f_8 T_3^{-1/2} ,  \ g_9 = f_9 T_3^{-1/2}. 
\eal
\eeq

Notice that each term in \eqref{eq:cT0} and \eqref{eq:cT} can be seen as the projection of $X, Y, Z$ onto some function $g_i$. For example, $c_{\om}, d_{\th}$ can be written as follows 
\[
c_{\om} = u_x(0) =  -\f{2}{ \pi} \int_0^{\infty} \f{\om}{x} dx = \la X , g_1 \ra, 
\quad d_{\th} = \int_0^{\infty} \f{\th_x}{x} dx = \la Y , g_5 \ra.
\]

Using the definition of $\cT$ in \eqref{eq:cT0}, \eqref{eq:cT} and the functions in \eqref{eq:func2}, we rewrite \eqref{eq:Copt0} as
 \beq\label{eq:Copt1}
\bal
\cT = &\la X, g_1 \ra \la  X, g_3 \ra
+ \la X, g_1 \ra \la Y, g_7 \ra
- ( \lam_2 - \lam_3 \bar d_{\th}) \la X, g_1 \ra \la Y, g_5 \ra 
+ \lam_2 \la X, g_1 \ra \la X, g_2 \ra   \\
&-  \lam_3 \la Y, g_5 \ra \la Y, g_6 \ra  + \lam_3 \la Y, g_5 \ra \la X, g_4 \ra + 
\lam_2 \la X, g_1 \ra \la Z, g_8 \ra  - \lam_3 \la Y, g_5 \ra \la Z, g_9 \ra  \\
\leq  & C_{opt} ( ||X||_2^2 + || Y||_2^2 + || Z||_2^2  
+ s_1 \la X, g_1 \ra^2 + s_2 \la Y, g_5 \ra^2 )  .
\eal
 \eeq

  We project $X, Y, Z$ onto the following finite dimensional spaces 
 \beq\label{eq:Copt_space}
X \in \mathrm{ span }\{ g_1, g_2, g_3, g_4 \} \teq \Sigma_1,
\quad Y  \in \mathrm{ span }\{ g_5, g_6, g_7 \} \teq \Sigma_2, 
\quad Z \in \mathrm{ span } \{ g_8, g_9 \} \teq \Sigma_3,
 \eeq
which only makes the right hand side of \eqref{eq:Copt1} smaller. Then \eqref{eq:Copt1} completely reduces to an optimization problem on the finite dimensional space.  
Using the optimal-constant argument in \cite{chen2019finite}, we obtain 
\[
C_{opt} = \lam_{\max}( D^{-1/2} M_{s} D^{-1/2}),
\]
where $D, M_s$ defined in \eqref{eq:Copt_M} are symmetric matrices with entries determined by the inner products among $g_i$. In particular, $C_{opt}$ can be computed rigorously and we present the details in Appendix \ref{app:Copt}.

\subsection{Summary of the estimates}\label{sec:est_sum}

Recall the $c_{\om}$ terms in \eqref{eq:l2h1}, the operators in \eqref{eq:linop}. Combining \eqref{eq:coer1}, \eqref{eq:cT0} and \eqref{eq:cT}, we yield 
\beq\label{eq:comb1}
\bal
 & \la \cL_{\th} \th_x, \th_x \psi \ra + \lam_1 \la \cL_{\om} \om, \om \vp \ra + \cT_0
= \la \cL_{\th1} \th_x, \th_x \psi \ra + \lam_1 \la \cL_{\om1} \om, \om \vp \ra + \cT \\
\leq &  \la  D_{\th} +   A_{\th} \psi^{-1}, \th_x^2 \psi \ra  + \la  \lam_1 D_{\om} +  A_{\om}  \vp^{-1} ,\om^2 \vp \ra  \\
-& \B( D_u - \f{9}{49} t_{12} - \f{72 \lam_1  }{49}\cdot 10^{-5} \B) || \td{u}_x x^{-2/3}||_2^2  + \cT  + A(u)  + G_c c_{\om}^2 \teq J.\\
\eal
\eeq

We use the remaining damping of $\om, \th_x, \td{u}_x$ and the argument in Section \ref{sec:cT} to control $\cT$. In \eqref{eq:para1}, Appendix \ref{app:cT}, we define $T_i > 0, s_i > 0$ that are used to compute the upper bound of $C_{opt} < 1$ in \eqref{eq:Copt0}. These functions and scalars are essentially determined by four parameters $\lam_2, \lam_3, \kp, t_{61} >0$.
Using the estimate \eqref{eq:Copt0}, we obtain 
\beq\label{eq:est_cT}
 \cT \leq  || \om T_1^{1/2} ||_2^2  + || \th_x T_2^{1/2} ||_2^2  + || \f{ u_{\D}}{x} T_3^{1/2}||_2^2  + s_1 c_{\om}^2 + s_2 d_{\th}^2.
\eeq

The $u_{\D}$ term can be further bounded by $ || \td{u}_x x^{-2/3}||_2$ and $|| \om x^{-2}||_2$ 
similar to \eqref{eq:est_kuw}, which is established in \eqref{eq:udel} in Appendix \ref{app:cT}.
Plugging \eqref{eq:udel} and \eqref{eq:est_cT} in \eqref{eq:comb1}, we obtain 
\beq\label{eq:comb2}
J \leq 
- \kp  ||   \th_x \psi^{1/2}||_2^2 
- \kp \lam_1 || \om \vp^{1/2} ||_2^2  + (s_1 + G_{c }) c_{\om}^2 + s_2 d_{\th}^2 
- 10^{-6}   || \td{u}_x x^{-2/3}||_2^2  + A(u) ,
\eeq
for $ \kp > 0$ determined in Appendix \ref{app:para}. The details are elementary and presented in Appendix \ref{app:cT}. For $\lam_2, \lam_3 >0$ given in \eqref{eq:para3}, we define the weighted
 $L^2$ energy 
\beq\label{energy:l2}
E_1^2( \th_x, \om) = ||\th_x \psi^{1/2}||_2^2 + \lam_1 || \om \psi^{1/2}||_2^2 
+ \lam_2 \f{\pi}{2} \cdot \f{4}{\pi^2} \la \om, x^{-1}\ra^2  + \lam_3 \la \th_x, x^{-1} \ra^2.
\eeq
Note that $\f{2}{\pi}\la \om, x^{-1} \ra  = - u_x(0) = - c_{\om}$ \eqref{eq:normal}. 
Recall the relations of different operators in \eqref{eq:linop}. Combining the equations \eqref{eq:l2h1}, \eqref{eq:ODE} and using the estimates \eqref{eq:comb1} and \eqref{eq:comb2},  we establish
\[
\bal
\f{1}{2}\f{d}{dt} E_1^2(\th, \om)
& =  \la \cL_{\th} \th_x, \th_x \psi \ra + \lam_1 \la \cL_{\om} \om, \om \vp \ra + \cT_0
+ \f{\pi \lam_2 }{2} (\bar c_{\om} + \bar u_x(0)) c^2_{\om} +  2 \bar c_{\om} \lam_3 d_{\th}^2 + \cR_{L^2}\\
&\leq  - \kp  || \th_x  \psi^{1/2} ||_2^2
- \kp \lam_1 || \om \vp^{1/2} ||_2^2
- 10^{-6}   || \td{u}_x x^{-2/3}||_2^2 
+ A (u)  \\
&\quad + \B( \f{\pi \lam_2 }{2} (\bar c_{\om} + \bar u_x(0)) + s_1+ G_{c } \B) c^2_{\om}  
 + ( 2 \bar c_{\om} \lam_3  + s_2) d_{\th}^2 + \cR_{L^2} , \\
\eal
\]
where $\cR_{L^2}$ is given by 
\beq\label{eq:Rl2}
\cR_{L^2} \teq N_1 + F_1 + \lam_1 N_2 + \lam_1 F_2 + \cR_{ODE} 
\eeq
and $N_i, F_i$ are defined in \eqref{eq:l2h1} and $\cR_{ODE}$ in \eqref{eq:Rode}. Recall $A(u)$ in \eqref{eq:Au}, $c_{\om} = u_x(0)$. Using the definitions of $s_i$ in \eqref{eq:para1}, we get
\[
\f{\pi \lam_2 }{2} (\bar c_{\om} + \bar u_x(0)) + s_1+ G_{c }
+ \f{ \pi\lam_1 e_3 \al_6}{12}  = -r_{c_{\om}} ,\quad s_2 + 2 \bar c_{\om} \lam_3  = - \kp \lam_3,
\]
for $r_{c_{\om}}, \kp>0$ determined in Appendix \ref{app:para}.
Hence, we obtain
\[
\bal
 &A (u) 
+ \B( \f{\pi \lam_2 }{2} (\bar c_{\om} + \bar u_x(0)) + s_1 + G_{c } \B) c^2_{\om} 
+ ( 2 \bar c_{\om} \lam_3  + s_2) d_{\th}^2 
=  - r_{c_{\om}} c^2_{\om} 
 -  \kp \lam_3 d_{\th}^2 -\f{\lam_1 e_3 \al_6 }{3}   \la \Lam \f{u}{x},  \f{u}{x} \ra.
 \eal
\]

Therefore, we obtain 
\beq\label{eq:L2}
\bal
\f{1}{2}\f{d}{dt} E_1^2(\th, \om) & \leq - \kp || \th_x \psi^{1/2} ||_2^2
- \kp \lam_1  || \om \vp^{1/2} ||_2^2 - 10^{-6}   || \td{u}_x x^{-2/3}||_2^2  -\f{\lam_1 e_3 \al_6 }{3}   \la \Lam \f{u}{x},  \f{u}{x} \ra \\
&- r_{c_{\om}} c^2_{\om}  - \kp \lam_3 d_{\th}^2 + \cR_{L^2} \teq Q + \cR_{L^2},
\eal
\eeq
These parameters satisfy $r_{c_{\om}} \geq   \f{\pi}{2} \kp \lam_2$.
Thus \eqref{eq:L2} implies 
\beq\label{eq:L2_2}
\f{1}{2}\f{d}{dt} E_1^2(\th, \om)  \leq - \kp E_1^2(\th, \om) + \cR_{L^2}, 
\eeq
and we establish the linear stability. See also \eqref{eq:goal}. 
Compared to \eqref{eq:L2_2}, \eqref{eq:L2} contains extra damping terms $- || \td{u}_x x^{-2/3}||_2^2$, $- \la \Lam \f{u}{x},  \f{u}{x} \ra$ and $-( r_{c_{\om}} -  \f{\pi}{2} \kp\lam_2 ) c_{\om}^2$. We  choose $r_{c_{\om}} > \kp \f{\pi}{2}\lam_2$ and keep these terms in \eqref{eq:L2} mainly to obtain sharper constants in our later weighted $H^1$ estimates.

\subsection{From linear stability to nonlinear stability with rigorous verification}\label{sec:lin_verif} 

In this subsection, we describe some main ideas how to go from linear stability to nonlinear stability with computer-assisted proof.

 (1) As we discuss at the beginning of Section \ref{sec:dyn}, the most challenging and essential part in the proof is the weighted $L^2$ linear stability analysis established in Section \ref{sec:lin}, since there is \textit{no} small parameter and the linearized equations \eqref{eq:intro_lin} are complicated. 

(2) The weighted $L^2$ linear stability estimates can be seen as a-priori estimates on the perturbation, and 
     we proceed to perform higher order energy estimates in a similar manner and establish the nonlinear energy estimate for some energy $E(t)$ of the perturbation
 \beq\label{eq:intro_boost}
       \f{d}{dt} E^2 \leq C E^3 - \lam E^2 + \e E  .
   \eeq
   Here, $-\lam E^2$ with $\lam >0$ comes from the linear stability, $C E^3$ with some constant $C( \bar \om, \bar \th) > 0$ controls the nonlinear terms, and $\e$ is the weighted norm of the residual error of the approximate steady state. See more details in Section \ref{sec:non}. To close the bootstrap argument $E(t) < E^*$ with some threshold $E^* > 0$, a sufficient condition is that 
              $ \e <\e^* = \lam^2 / ( 4 C)$,
 which provides an upper bound  on the required accuracy of the approximate steady state. 
 
The essential parts of the estimates in (1), (2) are established based on the grid point values of $(\bar \om, \bar \th)$ constructed using a moderate fine mesh. These parts do not involve the lengthy rigorous verification in the Supplementary Material \cite{chen2021HLsupp}. These estimates already provide a strong evidence of nonlinear stability.

 A significant difference from this step and step (1) is that we have a small parameter $\e$. As long as $\e$ is sufficiently small, thanks to the damping term $-\lam E^2$ established in step (1),  we can afford a large constant $C(\bar \om, \bar \th)$ in the estimate of the nonlinear terms and close the nonlinear estimates. We can complete all the nonlinear estimates in this step.

 (3) We follow the general approach in \cite{chen2019finite} to construct an approximate steady state with residual error below a required level $\e^*$. To achieve the desired accuracy, the construction is typically done by solving \eqref{eq:HLdyn} for a sufficiently long time using a fine mesh.  The difficulty of the construction depends on the target accuracy $\e^*$, and we refer to Section \ref{sec:on_steady_state} for more discussion on the new difficulty and the construction 
 of the approximate steady state for the HL model.
 Here, the mesh size plays a role similar to a small parameter that we can use. 
In practice, the profile $(\bar \om_1, \bar \th_1)$ constructed using a moderate fine mesh $\Om_1$ is close to the one $(\om_2, \th_2)$ constructed using a finer mesh $\Om_2$ with higher accuracy. As a result, the constants $C(\bar \om, \bar \th)$ and $\lam$ that we estimate in \eqref{eq:intro_boost} using different approximate steady states $(\om_i, \th_i)$ are nearly the same. This refinement procedure allows us to obtain an approximate steady state, based on which we close the nonlinear estimates \eqref{eq:intro_boost}. 
We refer more discussion of this philosophy to \cite{chen2019finite}.

 (4) Finally, we follow the standard procedure to perform rigorous verification on the estimates to pass from the grid point value to its continuous counterpart. 
Estimates that require rigorous verification with computer assistance are recorded in Appendix \ref{app:ver}. In the verification step, we can evaluate the approximate steady state on a much finer mesh $\Om_3$ with many more grid points so that they almost capture the whole behavior of the solution. Then, we use the regularity of the solution to pass from finite grid points to the whole real line. In this procedure, the mesh size in $\Om_3$ plays a role similar to a small parameter that we can exploit. In practice, to perform the rigorous verification, we evaluate the solution computed in a mesh with about $5000$ grid points using a much denser mesh with more than $5 \cdot 10^5$ grid points.

 In summary, in steps (2)-(4), we can take advantage of a small parameter which can be either the small error or the small mesh size, while there is no small parameter in step (1). Though these three steps could be technical, they are relatively standard from the viewpoint of analysis.

We remark that the approach of computer-assisted proof has played an important role in the analysis of many PDE problems, especially in computing explicit tight bounds of complicated (singular) integrals \cite{castro2014remarks,cordoba2017note,gomez2014turning} or bounding the norms of linear operators \cite{castelli2018rigorous,enciso2018convexity}. We refer to \cite{gomez2019computer} for an excellent survey on computer-assisted proofs in 
establishing rigorous analysis for PDEs, which also explains the use of interval arithmetics that guarantees rigorous computer-assisted verifications. Examples of highly nontrivial results established by the use of interval arithmetics can be found in, for example, \cite{lanford2017computer,hales2005proof,gabai2003homotopy,schwartz2009obtuse}. Our approach to establish stability analysis with computer assistance is different from existing computer-assisted approach, e.g. \cite{castro2020global}, where the stability is established by numerically tracking the spectrum of a given operator and quantifying the spectral gap. The key difference between their approach and ours is that we \textit{do not} use direct computation to quantify the spectral gap of the linearized operator. 
One of the main reasons is that the linearized operator in our case is not compact due to the Hilbert transform, and the non-compact component cannot be treated as a small perturbation. Thus we cannot approximate the linearized operator by a finite rank operator that can be further analyzed using matrix computation.


\section{On the approximate steady state}\label{sec:on_steady_state}

The proof of the main Theorem \ref{thm2} heavily relies on an approximate steady state solution $(\bar{\th},\bar{\om}, \bar c_l, \bar c_{\om})$ to the dynamic rescaling equations \eqref{eq:HLdyn},
which is smooth enough, e.g. $\bar \om, \bar \th_x \in C^3$. Moreover, as discussed in Section \ref{sec:lin_verif}, the residual error of the approximate steady state must be small enough 
in order to close the nonlinear estimates. In particular, the residual error $\e$ requirement depends on the stability gap $\lambda$ via the inequality $\e < \lam^2 / (4 C)$.

For comparison, we refer the reader to our previous work on proving the finite-time, approximate self-similar blowup of the 1D De Gregorio model via a similar computer-aided strategy \cite{chen2019finite}, where the corresponding approximate steady state  is constructed numerically on a compact domain $[-10,10]$. The stability gap that the authors proved in that work is relatively large (around $0.3$), and thus the point-wise error requirement on the residual can be relaxed to $10^{-6}$.

For the HL model, however, the stability gap $ \lam \approx \kp = 0.03$ (see \eqref{eq:para3}) 
that we can prove in the linear stability analysis \eqref{eq:L2_2} is much smaller, which leads to a much stronger requirement on the residual error. More precisely, we need to bound the residual in a weighted norm by $5.5\times 10^{-7}$ with weights \eqref{eq:wg}  that are singular of order $x^{-k}, k \geq 4 $ near $0$ and decay slowly for large $x$. This effectively requires the point-wise values of the residual to be as small as $10^{-10}$. To achieve this goal, it is not sufficient to simply follow the method in \cite{chen2019finite}, mainly due to the following reasons:

\begin{enumerate}
\item The steady state solution to \eqref{eq:HLdyn} is supported on the whole real line and has a slowly decaying tail in the far field (see below). If we approximate the steady state on a finite domain $[-L,L]$, we need to use an unreasonably large $L$ (roughly $L \ge 10^{30}$)
for the tail part beyond $[-L,L]$ to be considered as a negligible error, since truncating the tail leads to an error of order $L^{c_\om/c_l}\approx L^{-1/3}$. It is then impractical to achieve a uniformly small residual by only using mesh-based algorithms such as spline interpolations. 

\item Numerically computing the Hilbert transform of a function supported on the whole real line $\mathbb{R}$ is sensitively subject to round-off errors. For example, when computing $u$ from 
an odd function $\om$ via the Hilbert transform, we need to evaluate the convolution kernel $\log(|y-x|/|y+x|)$, which will be mistaken as $0$ by a computer program using double-precision if $|x/y|< 10^{-16}$. Such round-off errors, when accumulated over the whole mesh, are unacceptable in our case since we have a very high accuracy requirement for the computation of the approximate steady state solution. 
\end{enumerate}

To design a practical method of obtaining a sufficiently accurate construction, we must have some a priori knowledge 
on the behavior of a steady state $ (\om_{\inf}, \th_{\inf}, c_{l,\inf}, c_{\om,\inf})$. Assume that the velocity $u_{\inf}$ grows (if it grows) only sub-linearly in the far field, i.e. $ u_{\inf}(x) / x ,  u_{\inf,x}(x) \to 0$ as $x \to \infty$. Substituting this ansatz into the steady state equation of $\theta_x$ in \eqref{eq:HLdyn} yields 
\[
\frac{{\th}_{\inf, xx}}{ \th_{\inf,x}} \sim \frac{2 {c}_{ \om, \inf} }{ {c}_{l, \inf} } \cdot x^{-1}, \quad \text{which implies} \quad  \th_{\inf, x} \sim x^{2 {c}_{\om, \inf} / {c_{ l, \inf} } }.
\]
Furthermore, using these results to the steady state equation of $\om$ in \eqref{eq:HLdyn} yields 
\[
\frac{ {\om}_{\inf,x}} { {\om_{\inf}} } \sim \frac{  {c}_{\om,\inf}  }{ {c}_{l, \inf}  } \cdot x^{-1}, \quad \text{which implies} \quad \om_{\inf} \sim x^{ {c}_{\om, \inf} / {c_{l,\inf} }}.
\]
From our preliminary numerical simulation, we have $ {c}_{\om, \inf} / {c_{ l, \inf}} $ close to  $-1/3$. This straightforward argument implies that ${\om}_{\inf}$ and ${\th}_{\inf, x} $ should behave asymptotically like $x^{ {c}_{\om, \inf} / {c_{ l, \inf} } }$, $x^{2 {c}_{\om, \inf} / {c_{ l, \inf} } }$ as $x\rightarrow +\infty$, respectively, which in turn justifies the sub-linear growth of $u_{\inf}$. 

Guided by these observations, we will construct our approximate steady state as the combination of two parts:
\beq\label{eq:solu_decomp}
\bar{\om} = \om_b + \om_p,\quad \bar{\th} = \th_b + \th_p.
\eeq
We will call $(\om_b, \th_b)$ the explicit part and $(\om_p, \th_p)$ the perturbation part. The explicit part $(\om_b, \th_b)$ is constructed analytically to approximate the asymptotic tail behavior of the steady state for $x\geq L$, and satisfies $\om_b, \th_{b, x} \in C^5 $ and $\om_b \sim x^{\al}, \th_{x, b} \sim x^{2\al}$ with $\al \approx \bar c_{\om}/ \bar c_l < -\f{1}{3}$. 
The construction of $\om_b$ and its Hilbert transform relies on the following crucial identity 
\beq\label{eq:Hil_alpha}
H  ( \sgn(x) |x|^{-a } )= -\cot \f{\pi a}{2} \cdot |x|^{-a}, \quad a \in (0, 1),
\eeq
which is proved in the proof of Lemma \ref{lem:rota} in the Appendix.
It indicates that the leading order behavior of $H f$ for large $x$ is given by $ -\cot \f{\pi a}{2} \cdot |x|^{-a}$, if $f$ is odd with a decay rate $ |x|^{-a} $. By perturbing $\sgn(x) |x|^{-a }$ and \eqref{eq:Hil_alpha}, we construct $\om_b \in C^5$ and obtain the leading order behavior of $H \om_b$ for large $x$. This is one of the main reasons that we can compute the Hilbert transform of a function with slow decay accurately and overcome large round-off error in its computation. The perturbation part $(\om_p, \th_p)$ is constructed numerically using a quintic spline interpolation and methods similar to those in \cite{chen2019finite} in the domain $[-L,L]$ for some reasonably large $L$ (around $10^{16}$). 
By our construction, they satisfy that $ \om_p, \th_{p, x} \in C_c^3$. 
Since achieving a small residual error is critical to our proof,
a large portion of the Supplementary Material \cite{chen2021HLsupp} is devoted to the construction (Section 10) and error estimates of the approximate steady state (Section 11-15) with the above decomposition, especially the $\om_b$ part.

\subsection{Connection to the approximate steady state of the 2D Boussinesq in $\R^2_+$}\label{sec:2D_steady_state}

To generalize the current framework to the 2D Boussinesq equations, an important step is to construct an approximate steady state with a sufficiently residual error. The construction of the approximate steady state of the HL model provides important guidelines on this. The steady state equations of the dynamic rescaling formulation of the 2D Boussinesq, see e.g. \cite{chen2019finite2}, read
\[
\bal
(c_l x + \uu) \cdot \na \om &= c_{\om} \om + \th_x,  \\
 (c_l x + \uu) \cdot \na \th &= (c_l + 2 c_{\om}) \th,  \quad  \uu = \na^{\perp}(-\D)^{-1}\om.
\eal
\]
Denote $r = |x|$. Assume that the velocity $u$ grows sub-linearly in the far field : $\f{u(x)}{r} \to 0 $ as $r \to \infty$ and the scaling factors satisfy $c_l > 0, c_{\om} < 0$. Note that $ x \cdot \na = r \pa_r $. Passing to the polar coordinate $(r, \b), r = |x|, \b = \arctan \f{x_2}{x_1}$ and dropping the lower order terms, we yield 
\[
c_l r \pa_r \om(r, \b) = c_{\om} \om + \th_x + l.o.t., \quad c_l r \pa_r \th(r, \b) = (2 c_{\om} + c_l) \th + l.o.t.
\] 

Using an argument similar to the above argument for the HL model, we obtain
\[
 \om(r, \b) \sim p(\b) r^{ \al}, \quad  \th(r, \b) \sim q(\b) r^{1 + 2 \al},  \quad \al = \f{c_{\om}}{c_l} < 0 .
\]
We remark that $\th_x$ has a decay rate $r^{2\al}$ faster than that of $\om$. The computation in \cite{luo2014potentially} suggests that $ \al  \approx -\f{1}{3}$. Thus, the 
 profile (if it exists) for the 2D Boussinesq does not have a fast decay, and we also encounter the difficulties similar to (1) and (2).
  In particular, the 2D analog of difficulty (2) is to obtain the stream function $ \psi = (-\D)^{-1} \om$ accurately in $\R^2_+$. To design a practical method that overcomes these difficulties, it is important to perform a decomposition similar to \eqref{eq:solu_decomp}, 
 where $\om_p, \na \th_p$ have compact support and $\om_b, \th_b$ capture the tail behavior of the steady state. For the 2D Boussinesq, $\om_b, \th_b$ become semi-analytic since the angular part $ p(\b), q(\b)$ cannot be determined a-priori. To overcome the difficulty of solving the stream function in the far field, we seek a generalization of \eqref{eq:Hil_alpha}. We consider the ansatz $\psi = r^{2 + \al} f(\b)$ and solve 
 \[
  - \D (  r^{2 + \al} f(\b)) = r^{\al} p(\b)
 \]
 with boundary condition $f(0) = f(\pi/2) = 0$ due to the Dirichlet boundary condition and the odd symmetry for the solution $\om$. In the polar coordinate, the above equation is equivalent to 
\[
 (-\pa_{\b}^2 - (2 + \al)^2 ) f(\b) = p(\b) , \quad f(0) = f(\pi/2) = 0.
\]

Solving the above equation is significantly simpler than solving $ -\D \psi = \om$ in $\R^2_+$ since it is one-dimensional and in a compact domain. The above two formulas are a generalization of \eqref{eq:Hil_alpha} that connects the leading order far field behavior of $\om$ with that of the velocity. We believe that the above decomposition is crucial to construct the approximate steady state with sufficiently small residual error  for the 2D Boussinesq equations. The supplementary material on the analysis of the decomposition 
\eqref{eq:solu_decomp} for the HL model can be seen as a preparation for the more complicated decomposition in the 2D Boussineq equations.

\section{ Nonlinear stability and finite time blowup}\label{sec:non}

In this section, we further establish nonlinear stability analysis of \eqref{eq:lin0}.

\subsection{Weighted $H^1$ estimate}

In order to obtain nonlinear stability, we first establish the weighted $H^1$ estimate similar to 
\beq\label{eq:goal_H1}
\bal
\f{1}{2} \f{d}{dt} ( || D_x \th_x   \psi^{1/2} ||_2^2  +  \lam_1 || D_x \om  \vp^{1/2} ||_2^2 )
\leq   - c ( || D_x \th_x   \psi^{1/2} ||_2^2  +  \lam_1 || D_x \om  \vp^{1/2} ||_2^2 )
+ C E_1^2(\th, \om) + \cR_{H^1}
\eal
\eeq
for some $c , C>0$, where $D_x = x \pa_x$, $E_1$ is defined in \eqref{energy:l2} and $\cR_{H^1}$ are the error terms and nonlinear terms in the weighted $H^1$ estimate to be introduced.

In the work of Elgindi-Ghoul-Masmoudi \cite{Elg19}, they made a good observation that
the weighted $H^1$ estimates of the equation studied in \cite{Elg19} can be established by performing weighted $L^2$ estimates of $ x \pa_x f$ with the same weight as that in the weighted $L^2$ estimate, since the commutator between the linearized operator and $x \pa_x$ is of lower order. Inspired by this observation, we perform weighted $L^2$ estimates on $x \pa_x \th_x $ and $x \pa_x \om$. However, one important difference between our problem and that considered in \cite{Elg19} is that the commutator between the linearized operator in \eqref{eq:lin0} and $x\pa_x$ is not of lower order.

Denote $D_x = x \pa_x $. Similar weighted derivatives have been used in \cite{Elg19,elgindi2019finite,chen2019finite2} for stability analysis. We derive the equations for $D_x \th_x, D_x \om$. Taking $D_x$ on both side of \eqref{eq:lin0}, we get 
\beq\label{eq:est_H1}
\bal
\pa_t D_x \th_x & = \cL_{\th 1} (D_x \th_x, D_x \om) 
 + c_{\om} D_x ( \bar \th_x - x \bar \th_{xx} ) 
+  [ D_x , \cL_{\th1}] (\th_x, \om) + D_x F_{\th} + D_x N(\th) , \\
\pa_t D_x \om &= \cL_{\om 1} (D_x \th_x, D_x \om) 
 + c_{\om} D_x ( \bar \om - x \bar \om_{x} ) +  [D_x , \cL_{\om 1}] (\th_x, \om)+ D_x F_{\th} + D_x N(\th) ,
\eal
\eeq
where $[D_x, \cL](f, g) \teq D_x \cL(f,g) - \cL(D_x f, D_x g)$. In the Appendix \ref{app:commu}, we obtain the following formulas for the commutators 
\beq\label{eq:commu}
\bal
{ [ D_x , \cL_{\th 1} ] (\th_x, \om) } &= - (\bar u_x - \f{\bar u}{x}) D_x \th_x- D_x \bar u_x \th_x -D_x \td{u} \bar \th_{xx}  - \td{u} ( \bar \th_{xx} + D_x \bar \th_{xx}) , \\
[D_x , \cL_{\om 1}] (\th_x, \om) &= - (\bar u_x - \f{\bar u}{x}) D_x \om - \td{u} (\bar \om_x + D_x \bar \om_x)  ,
\eal
\eeq
where $\td{u}, \td{u}_x $ are defined in \eqref{eq:ut}. 

Performing the weighted $H^1$ estimates, we get 
\beq\label{eq:linH1}
\bal
  \f{1}{2} \f{d}{dt}\B( \la D_x \th_x, D_x \th_x \psi \ra  & + \lam_1 \la D_x \om, D_x \om \vp \ra \B)
  =  \B( \la \cL_{\th 1} ( D_x \th, D_x \om) , D_x \th_x \psi \ra
+ \lam_1 \la \cL_{\om 1} ( D_x \th, D_x \om) , D_x \om \vp \ra \B) \\
&  + \B( \la  [ D_x , \cL_{\th 1} ] (\th_x, \om) ,  D_x \th_x \psi \ra + 
\lam_1 \la [D_x , \cL_{\om 1}] (\th_x, \om),  D_x \om \vp \ra \B) \\
& + \B( \la c_{\om} D_x ( \bar \th_x - x \bar \th_{xx} ) , D_x \th_x \psi \ra
+ \lam_1  \la c_{\om} D_x ( \bar \om - x \bar \om_{x} ), D_x \om \vp \ra   \B) + \cR_{H^1} \\
&\teq  Q_1 + Q_2 + Q_3 +  \cR_{H^1},
\eal
\eeq
where $\cR_{H^1}$ is the remaining term in the weighted $H^1$ estimate 
\beq\label{eq:Rh1}
\cR_{H^1} = \la D_x N(\th), D_x \th_x \psi \ra + \lam_1 \la D_x N(\om), D_x \om \vp \ra
+ \la D_x F_{\th}, D_x \th_x \psi \ra + \lam_1 \la D_x F_{\om}, D_x \om \vp \ra.
\eeq

\subsubsection{Estimate of $Q_1$}\label{sec:H1_Q1}
Applying the estimate of $\cL_{\th,1}, \cL_{\om 1}$ in \eqref{eq:coer1} to $(D_x \th_x, D_x \om)$, we obtain 
\beq\label{eq:est_H11}
\bal
 Q_1 &\leq    \la  D_{\th} +   A_{\th} \psi^{-1}, ( D_x \th_x)^2 \psi \ra  + \la  \lam_1 D_{\om} +  A_{\om}  \vp^{-1} , ( D_x \om)^2 \vp \ra   \\
 & \quad +  A(-\Lam^{-1}(D_x\om)) + G_c \cdot (H D_x\om(0))^2,
\eal
\eeq
where $G_c$ is defined in \eqref{eq:est_chi2}, and we have dropped the term related to $|| \td{u}_x x^{-2/3}||_2^2$ in \eqref{eq:coer1} since $D_u -\f{9}{49} t_{12} - \f{72 \lam_1 }{49}\cdot 10^{-5} > 0$. In addition, we have replaced $ u = -\Lam^{-1}\om$ in  $A(u)$ in \eqref{eq:coer1} by $-\Lam^{-1} (D_x\om) $ and replaced $ c_{\om} = H \om(0)$ by $H D_x\om (0)$. Recall the definition of $A(u)$ in \eqref{eq:Au}. Since $\Lam = H \pa_x$ and $H \circ H = - Id$, we yield
\[
\pa_x ( - \Lam^{-1} D_x \om)(0) =  H D_x \om(0) = - \f{1}{\pi} \int_{\R} \om_x dx = 0,
\]
which implies 
\beq\label{eq:est_H112}
G_c \cdot (H (D_x\om)(0))^2 = 0 , \quad  A(- \Lam^{-1} D_x \om ) \leq 0.
\eeq
We treat $Q_1$ as the damping terms in the weighted $H^1$ estimate since from \eqref{ver:para1}, we have 
\beq\label{eq:damp_H1}
D_{\th} + A_{\th} \psi^{-1}\leq - \kp,  \quad \lam_1 D_{\om} + A_{\om} \vp^{-1} \leq -\lam_1 \kp , \quad \kp >0.
 \eeq

\subsubsection{Estimate of $Q_2$}\label{sec:H1_Q2}
Recall the commutators in \eqref{eq:commu}. The profile satisfies $ \bar u_x - \f{\bar u}{x} >0$ and thus $ - (\bar u_x - \f{\bar u}{x}) f$ with $f= D_x \th_x, D_x \om$ is a damping term in the $D_x \th_x$ or $D_x \om$ equation. We do not estimate these terms. 

For the term $D_x \bar u_x \th_x$ in \eqref{eq:commu}, using integration by parts, we get
\[
- \la D_x \bar u_x \th_x, D_x \th_x \psi \ra
= - \la x^2 \bar u_{xx} \psi, \f{1}{2} \pa_x (\th_x)^2 \ra
= \f{1}{2} \la  (x^2 \bar u_{xx} \psi)_x , \th_x^2 \psi\ra.
\]

The approximate steady state satisfies the following inequality 
\beq\label{eq:ver_H1_1}
 (x^2 \bar u_{xx} \psi)_x \leq  0.02 \psi , 
\eeq
which will be verified rigorously by the methods in the Supplementary Material \cite{chen2021HLsupp}. We record it in \eqref{ver:ver_H1_1}, Appendix \ref{app:ver}. Using \eqref{eq:ver_H1_1}, we obtain
\beq\label{eq:est_H113}
- \la D_x \bar u_x \th_x, D_x \th_x \psi \ra \leq \e_1  || \th_x \psi^{1/2} ||_2^2, \quad \e_1 = 0.01.
\eeq

The nonlocal terms in \eqref{eq:commu}  are of lower order than $D_x \om$ and we estimate then directly. We introduce some weights 
\beq\label{eq:Su23}
S_{u2} = t_{71} x^{-6} + t_{72} x^{-4} + 2 \cdot 10^{-6} x^{-10/3}, 
\quad S_{u3} = t_{81} x^{-6} + t_{82} x^{-4}  + 2 \cdot 10^{-6} x^{-10/3},
\eeq
for some parameters $t_{ij}>0$ to be determined. Using Young's inequality, we get 
\beq\label{eq:Young_ep1}
\bal
 &| \la D_x \td{u} \bar \th_{xx} , D_x \th_x  \psi \ra| +  | \la  \td{u} (\bar \th_{xx} + D_x \bar \th_{xx}) , D_x \th_x  \psi \ra| \\
\leq &   ||   D_x \td{u} S_{u2}^{1/2} ||_2^2 +\f{1}{4}  || S_{u2}^{-1/2}\bar \th_{xx}  D_x \th_x  \psi ||_2^2 
+ ||    \td{u} S_{u3}^{1/2} ||_2^2  +  \f{1}{4} || S_{u3}^{-1/2} ( \bar \th_{xx} + D_x \bar \th_{xx} )  D_x \th_x  \psi ||_2^2 .
 \eal
\eeq

We introduce the weights $S_{u2}, S_{u3}$ for a reason similar to that of $S_{u1}$ in Remark \ref{rem:est_uw_S1}.
Recall $D_x \td{u} = \td{u}_x $ and $\td{u}, \td{u}_x$ in \eqref{eq:ut}. Using the weighted estimates in Lemma \ref{lem:wg} yields
\beq\label{eq:Young_ep12}
\bal
 &||   D_x \td{u} S_{u2}^{1/2} ||_2^2
 +  ||    \td{u} S_{u3}^{1/2} ||_2^2    \\
 \leq  &    \B\la \om^2, ( t_{71} +\f{ 4 t_{81}}{25})   x^{-4} + ( t_{72}   +  \f{4 t_{82}}{9} )x^{-2} \B\ra + (1 +  \f{ 36  }{49}  )\cdot 2 \cdot 10^{-6} || \td{u}_x x^{-2/3}||_2^2.
    \eal
\eeq
In \eqref{eq:Young_ep12}, we do not estimate $ || \td u_x x^{-2/3}||_2^2$ in $ || D_x \td u S_{u2}^{1/2}||_2$ and keep it on both sides.

\begin{remark}\label{rem:Young_ep}
We will choose large enough parameters $t_{ij}$ in $S_{u2}, S_{u3}$ \eqref{eq:Su23} so that the weighted $L^2$ norm of $D_x \th_x$ terms in \eqref{eq:Young_ep1} are relative small compared to the damping term of $D_x \th_x$ in the weighted $H^1$ estimate \eqref{eq:linH1}, e.g. $Q_1$ in \eqref{eq:est_H11}. See also \eqref{eq:damp_H1}. The weighted $L^2$ norm of $\om$ and $|| \td u_x x^{-2/3}||_2$ in \eqref{eq:Young_ep12} will be bounded using the damping terms in the weighted $L^2$ estimate \eqref{eq:L2}.  
The same argument applies to controlling the weighted $L^2$ norm of $ D_x \om$ term in \eqref{eq:Young_ep2}.
\end{remark}

Next, we estimate the $\td{u}( (\bar \om_x + D_x \bar \om)$ term in \eqref{eq:commu}. The idea is similar to that in Section \ref{sec:est_uw}. We perform the following decomposition 
\beq\label{eq:chi_term3}
\bal
 - \lam_1 \la \td{u} (\bar \om_x + D_x \bar \om_x ) , D_x \om \vp \ra
 &= - \lam_1 \la \td{u} (\bar \om_x + D_x \bar \om_x - \f{1}{3} \chi \xi_3 ) , D_x \om \vp \ra
 - \f{1}{3} \lam_1 \la \td u \chi \xi_3 , D_x \om \vp \ra \\
 &\teq J + I_{r3}.
 \eal
\eeq

The estimate of $I_{r3}$ is similar to \eqref{eq:est_chi} and we obtain the following estimate in Appendix \ref{app:chi}
\beq\label{eq:est_chi3}
|I_{r3}| \leq \la  G_{\om 2} , \om^2 \ra +  \la G_{\om 3}, (D_x \om)^2 \ra + G_{c 2} c_{\om}^2, 
\eeq
where $G_{\om2}, G_{\om3}$ and $G_{c2}$ are given by 
\beq\label{eq:est_chi4}
\bal
 G_{\om 2}  & =  \f{1}{4 \cdot 10^6 }( \f{ 2\lam_1 (2 + \sqrt{3})}{5}  )^2 x^{-2/3}   , \quad 
G_{c2} = \f{ \lam_1^2|| x \xi_3 \chi^{1/2} \vp^{1/2}||_2^2}{ 36  } \cdot 10^3, \\
 G_{\om 3} & =  10^6 ( x^{4/3} \chi \xi_3 \vp)^2 +  10^{-3} \chi \vp .
\eal
\eeq
These functions are small due to the same reason that we describe in Section \ref{sec:chi}.

For $J$, we perform a decomposition 
\[
  J  = - \lam_1 \B\la \td{u} , \B( (\bar \om_x + D_x \bar \om_x - \f{1}{3} \chi \xi_3 )\vp - \f{e_3 \al_6}{9} x^{-2} \B) D_x\om \B\ra - \f{\lam_1 e_3\al_6}{9} \la \td{u}, D_x \om x^{-2} \ra \teq I_1 + I_2
\]

Note that $\td{u} = u - u_x(0) x$ and $\int_0^{\infty} x D_x \om x^{-2} dx =\int_0^{\infty} \om_x  dx =0 $. Using Lemma \ref{lem:cancel} with $f = \om$ and $g =u$, we get 
\[
I_2 = - \f{\lam_1 e_3 \al_6}{9} \B( \la u, D_x \om x^{-2} \ra  - u_x(0) \int_0^{\infty} x D_x \om x^{-2} dx   \B) =-\f{\lam_1 e_3 \al_6}{9} \la u, \om_x x^{-1} \ra =\f{\lam_1 e_3 \al_6}{ 9} \la \Lam \f{u}{x} , \f{u}{x} \ra,
\]
which can be controlled using the damping term in \eqref{eq:L2}. Denote 
\beq\label{eq:Kuw2}
S_{u4} = t_{91}x^{-6 } + t_{92} x^{-4} + 5 \cdot 10^{-4} x^{-10/3},  \quad \cK_{u \om 2} = 
(\bar \om_x + D_x \bar \om_x - \f{1}{3} \chi \xi_3 )\vp - \f{e_3 \al_6}{9} x^{-2}  .
\eeq

For $I_1$, using Young's inequality and the weighted estimate in Lemma \ref{lem:wg}, we get 
\beq\label{eq:Young_ep2}
\bal
&|I_1 |  \leq \lam_1 \la S_{u4}, \td{u}^2 \ra + \f{ \lam_1}{4} \la \cK_{u\om2}^2 S_{u4}^{-1}, (D_x\om)^2 \ra \\
\leq & \lam_1 \la \om^2, \f{4 t_{91}}{25} x^{-4} + \f{4 t_{92} }{9} x^{-2} \ra
+  \f{36 \lam_1 }{49} \cdot 5 \cdot 10^{-4} || \td{u}_x x^{-2/3}||_2^2 + \f{ \lam_1}{4} \la \cK_{u\om2}^2 S_{u4}^{-1}, (D_x\om)^2 \ra.
\eal
\eeq
We introduce the weight $S_{u4}$ for a reason similar to that of $S_{u1}$ in Remark \ref{rem:est_uw_S1}. The $|| \td{u}_x x^{-2/3}||^2_2$ term is further controlled by the corresponding damping term in \eqref{eq:L2}.

Combining the above estimates on the commutators in \eqref{eq:commu}, we obtain 
\beq\label{eq:est_H12}
\bal
Q_2\leq & \la - (\bar u_x - \f{ \bar u}{x}) + B_{\th} \psi^{-1} , (D_x \th_x)^2 \psi \ra
+ \la -\lam_1  (\bar u_x - \f{\bar u}{x}) + B_{\om} \vp^{-1}, (D_x \om)^2 \vp \ra + \e_1 || \th_x \psi^{1/2} ||_2^2 \\
&+\la A_{\om 2} , \om^2  \ra +  \f{\lam_1 e_3 \al_6}{ 9} \la \Lam \f{u}{x} , \f{u}{x} \ra + \B( 
(1 + \f{36}{49})\cdot 2 \cdot 10^{-6} + \f{36 \lam_1  }{49}  \cdot 5 \cdot 10^{-4} \B) || \td{u}_x x^{-2/3}||_2^2
+ G_{c2} c_{\om}^2,
\eal
\eeq
where $G_{c2}$ is defined in \eqref{eq:est_chi4}. The term $(\bar u_x - \f{ \bar u}{x})$ comes from the commutators \eqref{eq:commu} and we do not estimate them in $Q_2$ in \eqref{eq:linH1}. 
The terms $B_{\th}, B_{\om}, A_{\om 2}$ are the sum of the coefficients in the integrals of $(D_x \th_x)^2, (D_w \om)^2, \om^2$ in the above estimates 
\beq\label{eq:coer2_cost}
\bal
B_{\th} &\teq \f{1}{4} S_{u2}^{-1} (\bar \th_{xx} \psi)^2 + \f{1}{4} S_{u3}^{-1} ( \bar \th_{xx} + D_x \bar \th_{xx} )^2 \psi^2 , \quad B_{\om}  \teq   \f{ \lam_1}{4} \cK_{u\om2}^2 S_{u4}^{-1}  + G_{\om 3}, \\
 A_{\om 2} & \teq ( t_{71} +\f{ 4 t_{81}}{25})   x^{-4} + ( t_{72}   +  \f{4 t_{82}}{9} )x^{-2}
+ \lam_1 ( \f{4 t_{91}}{25} x^{-4} + \f{4 t_{92} }{9} x^{-2} ) + G_{\om 2}. \\
\eal
\eeq

\subsubsection{Estimate of $Q_3$}\label{sec:H1_Q3}
Recall the $c_{\om}$ terms in \eqref{eq:est_H1}. Denote by $K_1, K_2$ the following $L^2$ norms 
\beq\label{eq:K12}
K_1 \teq  || \pa_x ( x^3 \bar \th_{xxx} \psi ) \psi^{-1/2} ||_2, \quad K_2\teq || \pa_x ( x^3 \bar \om_{xx} \vp ) \vp^{-1/2} ||_2.
\eeq

Using integration by parts and the Cauchy-Schwarz inequality, we obtain
\[
\bal
|c_{\om} \la D_x (\bar \th_x - x \bar \th_{xx}), D_x \th_x \psi \ra |
= & |c_{\om}     \la  - x^2 \bar \th_{xxx} \cdot ( x \psi ), \ \pa_x \th_x \ra |  
= |c_{\om} \la  \pa_x( x^3 \bar \th_{xxx} \psi) , \th_x \ra  | \\
\leq & |c_{\om}| \cdot || \pa_x ( x^3 \bar \th_{xxx} \psi ) \psi^{-1/2} ||_2 || \th_x \psi^{1/2}||_2
= K_1 |c_{\om} | \cdot|| \th_x \psi^{1/2}||_2,
\eal
\]
where we have used $ x \pa_x ( f - x \pa_x f) = -x^2 f_{xx}, f = \bar \th_x$ in the first equality. Similarly, we have 
\[
\lam_1 |c_{\om} \la D_x (\bar \om - x \bar \om_{x}), D_x \om \vp \ra |
\leq \lam_1 |c_{\om}| \cdot || \pa_x ( x^3 \bar \om_{xx} \vp ) \vp^{-1/2} ||_2 || \om \vp^{1/2}||_2
= \lam_1 K_2 |c_{\om}| \cdot|| \om \vp^{1/2}||_2.
\]

Using Young's inequality, we obtain 
\beq\label{eq:est_H13}
\bal
 Q_3 \leq & K_1 |c_{\om}|  \cdot || \th_x \psi^{1/2}||_2 +  \lam_1 K_2 |c_{\om}| \cdot || \om \vp^{1/2}||_2 \\
\leq &   \g_1 || \th_x \psi^{1/2}||_2^2 + \g_2 || \om \vp^{1/2}||_2^2 
+  c_{\om}^2 \B(  \f{ K_1^2}{4 \g_1} +\f{ (\lam_1  K_2)^2}{4 \g_2}  \B),
\eal
\eeq
where $\g_1, \g_2>0$ are chosen in \eqref{eq:para_H1}.

\subsubsection{Summary of the estimates}
We determine the parameters $t_{ij}$ in the estimates in Sections \ref{sec:H1_Q1}-\ref{sec:H1_Q3} 
and choose $\kp_2$ so that  
\beq\label{eq:est_H14}
\bal
&D_{\th_2 } + B_{\th} \psi^{-1} \leq -\kp_2, \quad 
  D_{\om 2} + B_{\om} \vp^{-1} \leq - \kp_2 \lam_1 ,  \\
&D_{\th 2} \teq D_{\th} + A_{\th} \psi^{-1}  - (\bar u_x - \f{ \bar u}{x}) ,
  \quad D_{\om 2} \teq  \lam_1 D_{\om} + A_{\om} \vp^{-1} -\lam_1  (\bar u_x - \f{\bar u}{x}). \\
\eal
\eeq
 The terms $D_{\th 2},  D_{\om2}$ are the coefficients of the damping terms in the weighted $H^1$ estimate \eqref{eq:linH1} and are already determined in the weighted $L^2$ estimates. 
The terms $ B_{\th} \psi^{-1},  B_{\om} \vp^{-1} $ defined in \eqref{eq:coer2_cost} are the coefficients in the weighted $L^2$ norm of $D_x \th_x, D_x \om$ in \eqref{eq:Young_ep1}, \eqref{eq:Young_ep2}. The motivation of \eqref{eq:est_H14} is that we use the damping terms to control the weighted $L^2$ norms of $D_x \th_x, D_x \om$ in the estimates of $Q_i$. The idea is the same as that in Remark \ref{rem:Young_ep}.

We first choose $\kp_2 < \kp= 0.03$ in Appendix \ref{app:para}, where $\kp$ is related to \eqref{eq:L2_2}. This choice is motivated by our estimate \eqref{eq:goal_H1_2}. 
The dependences of $A_{\om 2}, B_{\th}, B_{\om}$ on $t_{ij}$ are given in \eqref{eq:coer2_cost}, \eqref{eq:Su23}, \eqref{eq:Kuw2}. Inequalities in \eqref{eq:est_H14} can be seen as constraints on $t_{ij}$. We choose $t_{ij}$ subject to the constraints \eqref{eq:est_H14} such that $|| A_{\om 2} \vp^{-1} ||_{\infty}$ is as small as possible. This enables us to obtain sharper constant $a_{H^1}$ in the weighted $H^1$ estimate \eqref{eq:H1}. 
After $t_{ij}$ are determined, we verify \eqref{eq:est_H14} and 
\beq\label{eq:est_H15}
|| A_{\om 2} \vp^{-1} ||_{\infty}\leq a_{H^1},
\eeq
using the methods in the Supplementary Material \cite{chen2021HLsupp}, and record them in \eqref{ver:est_H14}, Appendix \ref{app:ver}, where $a_{H^1}$ is given in \eqref{eq:para_H1}.

Combining \eqref{eq:est_H11}, \eqref{eq:est_H112}, \eqref{eq:est_H113}, \eqref{eq:est_H12}, \eqref{eq:est_H13} and \eqref{eq:est_H14}, we prove 
\beq\label{eq:H1}
\bal
&\f{1}{2} \f{d}{dt} ( || ( D_x \th_x)  \psi^{1/2} ||_2^2  +  \lam_1 || (D_x \om ) \vp^{1/2} ||_2^2 )
\leq   - \kp_2 || (D_x \th_x)  \psi^{1/2} ||_2^2  - \kp_2 \lam_1 || (D_x \om)  \vp^{1/2}||_2^2 \\ 
&+  (\e_1 + \g_1) || \th_x \psi^{1/2} ||_2^2  + ( a_{H^1} + \g_2) || \om   \vp^{1/2} ||_2^2  +  
  \B( G_{c2}   +  \f{ K_1^2}{4 \g_1} +\f{ (\lam_1  K_2)^2}{4 \g_2}      \B) c_{\om}^2 
 \\
   &+ \f{\lam_1 e_3 \al_6}{ 9} \la \Lam \f{u}{x} , \f{u}{x} \ra + 
\B( (1 + \f{36}{49})\cdot 2 \cdot 10^{-6} + \f{36 \lam_1  }{49}  \cdot 5 \cdot 10^{-4} \B) || \td{u}_x x^{-2/3}||_2^2 + \cR_{H^1}.
   \eal
\eeq

Recall the weighted $L^2$ energy $E_1$ in \eqref{energy:l2}. For some $\lam_4 > 0$, we construct the energy 
\beq\label{energy:h1}
\bal
&E^2(\th_x, \om) =  E_1^2(\th_x, \om) +  \lam_4 ( ||  D_x \th_x \psi^{1/2}||_2^2 + \lam_1 || D_x \om \vp^{1/2}||_2^2  )  \\
=&||  \th_x \psi^{1/2}||_2^2 + \lam_1 ||  \om \psi^{1/2}||_2^2
+ \lam_2 \f{\pi}{2} c_{\om}^2 + \lam_3 d_{\th}^2 +  \lam_4 ( ||  D_x \th_x \psi^{1/2}||_2^2 + \lam_1 || D_x \om \vp^{1/2}||_2^2  ) .
\eal
\eeq

Note that $c_{\om}, || \th_x \psi^{1/2}||_2, || \om \vp^{1/2}||_2$ in \eqref{eq:H1} can be bounded by the energy $E_1$ in \eqref{energy:l2}. The terms $\la \Lambda \f{u}{x}, \f{u}{x} \ra$ and $|| \td u_x x^{-2/3}||_2^2$ can be bounded by their damping terms in \eqref{eq:L2}. To motivate later estimates and the choice of several parameters, we neglect these two terms. Then \eqref{eq:H1} implies \eqref{eq:goal_H1} with $c= \kp_2$ and some $C>0$. Combining \eqref{eq:L2_2} and \eqref{eq:goal_H1}, we get 
\beq\label{eq:goal_H1_1}
\f{1}{2}\f{d}{dt} E^2(\th_x, \om)
\leq - ( \kp - \lam_4 C ) E_1^2 
- \kp_2 (|| D_x \th_x   \psi^{1/2}||_2^2  +  \lam_1 || D_x \om  \vp^{1/2}||_2^2) + \cR_{L^2} +  \lam_4 \cR_{H^1},
\eeq
where $\kp = 0.03$. We first choose $\kp_2< \kp$ and then $\lam_4$ small enough, such that 
\beq\label{eq:goal_H1_2}
\kp  - \lam_4 C \geq \kp_2.
\eeq
Then we obtain the linear stability of \eqref{eq:lin0} in the energy norm $E$.

\subsection{Nonlinear stability}\label{sec:non_l2h1}

Combining \eqref{eq:L2} and \eqref{eq:H1}, we derive 
\[
\bal
&\f{1}{2}\f{d}{dt} E^2(\th_x, \om) 
 \leq  
- \kp || \th_x  \psi^{1/2} ||_2^2
- \kp \lam_1 || \om \vp^{1/2} ||_2^2  - r_{c_{\om}} c^2_{\om}  -  \kp \lam_3 d_{\th}^2  \\
&- \lam_4  \kp_2 ||  (D_x \th_x) \psi^{1/2} ||_2^2  - \lam_4  \kp_2 \lam_1  || (D_x \om)  \vp^{1/2} ||_2^2
+ \lam_4 (\e_1 + \g_1) || \th_x \psi^{1/2} ||_2^2  \\
& + \lam_4 ( a_{H^1} + \g_2) || \om  \vp^{1/2} ||_2^2  
   + \lam_4   \B( G_{c2}   +  \f{ K_1^2}{4 \g_1} +\f{ (\lam_1  K_2)^2}{4 \g_2}      \B) c_{\om}^2
       + ( \f{\lam_1 e_3 \al_6}{ 9} \lam_4 -   \f{\lam_1 e_3 \al_6}{ 3}  ) \la \Lam \f{u}{x} , \f{u}{x} \ra  \\
    &+ 
\B( \B( (1 + \f{36}{49})\cdot 2 \cdot 10^{-6} + \f{36 \lam_1  }{49}  \cdot 5 \cdot 10^{-4}\B) \lam_4 -10^{-6} \B) || \td{u}_x x^{-2/3}||_2^2     + \cR_{L^2} + \lam_4 \cR_{H^1}    .
   \eal
\]

Since $\kp_2 < \kp$, we choose small $\lam_4 > 0$ in Appendix \ref{app:para} so that 
\beq\label{eq:para_ineq}
\bal
&\lam_4  \cdot \f{\lam_1 e_3 \al_6}{9} < \f{\lam_1 e_3 \al_6}{3} , \quad 
\B( (1 + \f{36}{49})\cdot 2 \cdot 10^{-6} + \f{36 \lam_1  }{49}  \cdot 5 \cdot 10^{-4}\B) \lam_4 <10^{-6}  ,\\
& r_{c_{\om}} -  \lam_4  \B( \f{ K_1^2}{4 \g_1} +\f{ (\lam_1  K_2)^2}{4 \g_2} \B)
-  \lam_4 G_{c2} >  \kp_2 \cdot \f{\pi \lam_2}{2},    \\
& \kp  - \lam_4 \g_1 -\lam_4 \e_1 \geq \kp_2 , \quad  \kp \lam_1 - \lam_4 \g_2 - \lam_4 a_{H^1} \geq \kp_2 \lam_1, \quad \kp \lam_3 \geq \kp_2 \lam_3,
\eal
\eeq
where $K_1, K_2$ are defined in \eqref{eq:K12}. The above inequalities will be verified rigorously by the methods in the Supplementary Material \cite{chen2021HLsupp}.
Note that $r_{c_{\om}} > \f{\pi}{2} \lam_2 \kp$ and $\kp_2 < \kp$. 
The above conditions are essentially the same as \eqref{eq:goal_H1_2}.
We keep the damping term $\la \Lam \f{u}{x}, \f{u}{x} \ra$ and $|| \td u_x x^{-2/3 }||_2^2$ in \eqref{eq:L2} to control the corresponding terms in \eqref{eq:H1}. Plugging the above estimates and  \eqref{eq:para_ineq} into the differential inequality, we yield 
\beq\label{eq:non1}
\bal
\f{1}{2}\f{d}{dt} E^2(\th_x, \om) 
 \leq  & 
 - \kp_2 || \th_x \psi^{1/2}||_2^2 - \kp_2 \lam_1 || \om \vp^{1/2}||_2^2
 - \kp_2 \f{ \pi \lam_2}{2} c^2_{\om}  -  \kp_2 \lam_3 d_{\th}^2  \\
&- \lam_4 \kp_2 || (D_x \th_x)  \psi^{1/2} ||_2^2  - \lam_4 \kp_2 \lam_1 || (D_x \om)  \vp^{1/2}||_2^2 +
 \cR_{L^2} + \lam_4 \cR_{H^1}  \\
 \leq & -\kp_2 E^2(\th_x, \om) +  \cR_{L^2} + \lam_4 \cR_{H^1} \teq  -\kp_2 E^2(\th_x, \om) + \cR,
 \eal
\eeq
where $\cR =  \cR_{L^2} + \lam_4 \cR_{H^1} $ and $\kp_2 = 0.024$ is given in \eqref{eq:para_H1}.

\subsubsection{Outline of the estimates of the nonlinear and error terms}
Recall the definitions of $\cR_{L^2}$ and $\cR_{H^1}$ in \eqref{eq:Rl2} and \eqref{eq:Rh1}. The nonlinear terms in $\cR_{L^2}, \cR_{H^1}$, e.g. $\la D_x N(\th), D_x \th_x \psi \ra$, depend cubically on $\th_x, \om$. In the Supplementary Material \cite{chen2021HLsupp}, we use the energy $E(\th_x, \om)$ and interpolation to control $|| u_x ||_{\infty}$ and $||\th_x||_{\infty}$.  Using these $L^{\infty}$ estimates, we further estimate the nonlinear terms in $\cR$. For example, a typical nonlinear term in $\cR$ can be estimated as follows 
\[
\bal
 |\la u \th_{xx} , \th_x \psi \ra | &= \f{1}{2} |\la u \psi, \pa_x (\th_x)^2 \ra|
= \f{1}{2} | \la (u\psi)_x , \th_x^2 \ra  | = \f{1}{2} |\la u_x \psi + u \psi_x, \th_x^2  \ra | \\
&\leq \f{1}{2} ( || u_x||_{L^{\infty} } + ||\f{u}{x}||_{\infty} || \f{x\psi_x}{\psi} ||_{L^{\infty}}) || \th_x \psi^{1/2}||_2^2\leq \f{1}{2} || u_x||_{L^{\infty}} (1 + || \f{x\psi_x}{\psi} ||_{L^{\infty}}) || \th_x \psi^{1/2}||_2^2,
\eal
\]
where we have used $ |\f{u}{x}| \leq ||u_x||_{\infty}$ in the last inequality since $u(0) = 0$. The above upper bound can be further bounded by $E^3(\th_x, \om)$.

The error terms in $\cR_{L^2}, \cR_{H^1}$, e.g. $F_1 = \la F_{\th}, \th_x \psi \ra$, depend linearly on $\th_x, \om$. We estimate these terms using the Cauchy-Schwarz inequality. 
A typical term $F_1$ can be estimated as follows 
\[
|F_1| \leq  || F_{\th} \psi^{1/2}||_2 || \th_x \psi^{1/2}||_2.
\]
The error term $|| F_{\th} \psi^{1/2}||_2$ is small and $|| \th_x \psi^{1/2}||_2$ can be further bounded by $E(\th_x, \om)$. 

In the Supplementary Material \cite{chen2021HLsupp}, we work out the constants in these estimates and establish the following estimates 
\[
\cR = \cR_{L^2} + \lam_4 \cR_{L^2} \leq 36 E^3 + \e  E, \quad \e = 5.5 \cdot 10^{-7}.
\]

\subsubsection{Nonlinear stability and finite time blowup}
Plugging the above estimate on $ \cR$ in \eqref{eq:non1}, we establish the nonlinear estimate 
\[
\f{1}{2} \f{d}{dt} E^2(\th_x, \om) \leq -\kp_2 E(\th_x, \om)^2 + 36 E(\th_x, \om)^3 + \e E(\th_x, \om),
\]
where $\kp_2=0.024$ is given in \eqref{eq:para_H1}. We choose the threshold $E_* = 2.5 \cdot 10^{-5}$ in the Bootstrap argument. Since 
\[
- \kp_2 E_*^2 + 36 E_*^3 + \e E_* < 0,
\] 
the above differential inequality implies that if $E(0) < E_*$, the bootstrap assumption 
\beq\label{eq:boot}
E( \th_x(t), \om(t)) < E_* 
\eeq
holds for all $t>0$. Consequently, we can choose odd initial perturbations $\th_x, \om$ which satisfy $\om_x(0) = \th_{xx}(0) = 0, E(\th_x, \om) < E_*$ and modify the far field of $ \bar \th, \bar \om$ so that $ \bar \om + \om , \th_x + \bar \th_x \in C_c^{\infty}$. The bootstrap result implies that for all time $t>0$, the solution $\om(t)+\bar{\om}, \th_x(t) + \bar \th_x, c_l(t) +\bar{c}_l, c_{\om}(t) + \bar{c}_{\om}$ remain close to $\bar{\om}, \bar\th_x, \bar c_l,   \bar{c}_{\om} $, respectively. Using the rescaling argument in Section \ref{sec:dyn}, we obtain finite time blowup of the HL model.

\subsection{ Convergence to the steady state}\label{sec:conv}

We use the time-differentiation argument in \cite{chen2019finite} to establish convergence. 
The initial perturbations $(\th_x, \om)$ satisfy the properties in the previous Section. Since the linearized operators and the error terms in \eqref{eq:lin0} are time-independent, differentiating \eqref{eq:lin0} in $t$, we get 
\[
\pa_t (\th_x)_t = \cL_{\th}( (\th_x)_t, \om_t ) + \pa_t N(\th), \quad 
\pa_t ( \om )_t = \cL_{\om}( (\th_x)_t, \om_t ) + \pa_t N(\om).
\]

Applying the estimates of $\cL_{\th}, \cL_{\om}$ in Section \ref{sec:lin} and \eqref{eq:L2_2} to $(\th_x)_t, \om_t$, we obtain 
\[
\f{1}{2} \f{d}{dt} E_1( (\th_x)_t, \om_t)^2 
\leq -\kp E_1( (\th_x)_t, \om_t)^2 + \cR_2,
\]
where the energy notation $E_1$ is defined in \eqref{energy:l2} and $\cR_2$ is given by 
\[
\bal
&E_1( (\th_x)_t, \om_t )= || (\th_x)_t  \psi^{1/2}||_2^2 + \lam_1 || \om_t  \psi^{1/2} ||_2^2 
+ \lam_2 \f{\pi}{2} (\pa_t c_{\om})^2 + \lam_3 (\pa_t d_{\th})^2,\\
&\cR_2 = \la \pa_t N(\th), (\th_x)_t \psi \ra + \lam_1 \la \pa_t N(\om), \om_t \vp \ra 
 - \lam_2 \pa_t c_{\om} \la \pa_t N(\om), x^{-1} \ra + \lam_3  \pa_t d_{\th} \la \pa_t N(\th), x^{-1} \ra. 
 \eal
\] 

The term $\pa_t c_{\om}$ in the above estimates is from 
\[
 H \om_t (0) = \pa_t H \om(0) = \pa_t c_{\om}.
\]
Similarly, we obtain the term $\pa_t d_{\th} $. Using the a-priori estimate $E(\th_x, \om) < E_*$ 
in \eqref{eq:boot} and the energy $E_1( (\th_x)_t, \om_t )$, we can further estimate $\cR_2$. In the Supplementary Material \cite{chen2021HLsupp}, we prove 
\beq\label{eq:EE_time}
\f{1}{2} \f{d}{dt} E_1( (\th_x)_t, \om_t)^2 
\leq - 0.02 E_1( (\th_x)_t, \om_t)^2 .
\eeq

Using this estimate and the argument in \cite{chen2019finite}, we prove that the solution $\om + \bar \om, \th_x + \bar \th_x$ converge to the steady state $\om_{\infty}, \th_{\infty, x}$ in $L^2(\vp), L^2(\psi)$ and $c_l(t), c_{\om}(t) $ converge to $c_{l,\infty}, c_{\om,\infty}$ exponentially fast. Moreover, the steady state admits regularity $ (D_x)^i ( \om_{\infty} - \bar \om ) \in L^2(\vp) , (D_x)^i (\th_{x,\infty}- \bar \th_x )\in L^2(\psi) $ for $i=0,1$. We obtain $\th_{\inf}$ from $\th_{\inf, x}$ by imposing $\th_{\inf}(0)=0$ and integration.

Recall the energy $E$ in \eqref{energy:h1}. Since 
\[ 
\lam_2 \pi /2 > 3 > 1.5^2, \quad E > ( \lam_2 \pi / 2)^{1/2} |c_{\om}| \geq 1.5 |c_{\om}|
\]
(see \eqref{eq:para3}), using the convergence result, the a-priori estimate \eqref{eq:boot} and 
\eqref{eq:normal}, we obtain 
\beq\label{eq:profile_inf}
 E( \th_{x,\infty} - \bar \th_x,  \om_{\infty } - \bar \om ) \leq E_*, \quad c_{l,\infty} = \bar c_l = 3, \quad  | c_{\om, \infty} - \bar c_{\om} | \leq \f{2}{3} E_* = \f{5}{3}  \cdot 10^{-5}.
\eeq

Recall $\bar c_{\om} < -1.0004$ from the beginning of Section \ref{sec:linop}.
Thus, $c_{\om,\inf} < -1$ and we conclude that the blowup is focusing and asymptotically self-similar with blowup scaling $\lam = \f{ c_{l,\infty}}{ -c_{\om,\infty}}$ satisfying 
\[
| \lam - \bar \lam |
\leq | \f{ c_{l,\infty}}{ c_{\om,\infty}} - \f{\bar c_l}{\bar c_{\om}} |
+ | \f{\bar c_l}{\bar c_{\om}} + \bar \lam | < 3 | \bar c_{\om} -c_{\om,\inf} |
+   10^{-5}  < 6 \cdot 10^{-5}, \quad \bar \lam = 2.99870,
\]
where $\bar \lam$ is the determined by the first $6$ digits of $ - \bar c_{l} / \bar c_{\om}$.

\subsection{ Uniqueness of the self-similar profiles}

Suppose that  $(\om_1, \th_1)$ and $(\om_2, \th_2)$ are two initial perturbations which are small in the energy norm $ E(\om_i, \th_{i,x}) < E_*$. The associated solution $(\om_i, \th_{i,x})$  solves \eqref{eq:lin0}
\[
\pa_t \th_{i,x} = \cL_{\th}( \th_{i,x} ,\om_i) + F_{\th} + N(\th_i),
\quad \pa_t \om_{i} = \cL_{\om}( \th_{i,x} ,\om_i) + F_{\om} + N(\om_i).
\]

Denote 
\beq\label{eq:dN}
\d \om  \teq  \om_1 - \om_2, \quad \d \th \teq  \th_{1} - \th_{2}, \quad 
  \d N_{\th} = N(\th_1) - N(\th_2), \quad \d N_{\om} =  N(\om_1) - N(\om_2).
\eeq

A key observation is that the forcing terms $F_{\th}, F_{\om}$ do not depend on $(\om_i, \th_i)$. Thus, we derive
\[
\pa_t \d \th_x = \cL_{\th}( \d \th_x, \d \om) + \d N_{\th}, 
\quad \pa_t \d \om = \cL_{\om}( \d \th_x, \d \om) + \d N_{\om}.
\]

Applying the estimates of $\cL_{\th}, \cL_{\om}$ in Section \ref{sec:lin} and \eqref{eq:L2_2}, we get 
\[
\f{1}{2} \f{d}{dt} E_1(  \d \th_x, \d \om )^2 \leq - \kp  E_1( \d \th_x, \d \om)^2 + \cR_3
\]
where the energy notation $E_1$ is defined in \eqref{energy:l2} and $\cR_3$ is given by 
\[
\cR_3 = \la \d N_{\th}, \d \th_x \psi \ra 
+ \lam_1 \la \d N_{\om}, \d \om \vp \ra 
- \lam_2 c_{\om}( \d \om) \cdot \la \d N_{\om}, x^{-1} \ra
+  \lam_3 d_{\th}( \d \th_x) \cdot \la \d N_{\th}, x^{-1} \ra.
\]

The above formulations are very similar to that in Section \ref{sec:conv}. Formally, the difference operator $\d $ is similar to the time differentiation $\pa_t$. In the Supplementary Material \cite{chen2021HLsupp}, we show that $(\d \th_x, \d \om)$ enjoys the same estimates as that of $( \pa_t \th_x, \om_t)$ in \eqref{eq:EE_time}	
\beq\label{eq:EE_uni}
\f{1}{2} \f{d}{dt} E_1(\d \th_x, \d \om)^2 \leq - 0.02 E_1(\d \th_x, \d \om)^2 .
\eeq
As a result, $E_1(\d \th_x, \d \om)$ converges to $0$ exponentially fast and the two solutions $(\om_i + \bar \om_i, \th_i + \bar \th_i), i=1,2$ converge to steady states $(\om_{i,\inf}, \th_{i,\inf})$ with the same $\om_{\inf}$ and $\th_{\inf, x}$. Since $\th_{i,\inf}(0)= 0 $,  two steady states are the same.

\subsection{Numerical evidence of stronger uniqueness}\label{sec:numerical_uniqueness} The above discussion argues that the steady state is unique at least within a small energy norm ball. However, our numerical computation suggests that the steady state of the dynamical rescaling equations \eqref{eq:HLdyn},\eqref{eq:normal0} is unique (up to rescaling) for a much larger class of smooth initial data $\om,\theta$ with $\th(0) = 0$ that satisfy the following conditions:
\begin{enumerate}
\item odd symmetry: $\om(x)$ and $\theta_x(x)$ are odd functions of $x$;
\item non-degeneracy condition: $\om_x(0)>0$ and $\theta_{xx}(0)>0$;
\item sign condition: $\om(x) , \th_x(x) > 0$ for $x > 0$.
\end{enumerate}
In fact, these conditions are consistent with the initial data considered by Luo-Hou in \cite{luo2014potentially,luo2013potentially-2} restricted on the boundary. They are preserved by the equations as long as the solution exists.
Moreover, this class of initial data leads to finite time blowup of the HL model \cite{choi2014on}.

\begin{figure}[t]
\centering
    \begin{subfigure}[b]{0.32\textwidth}
    	\includegraphics[width=1\textwidth]{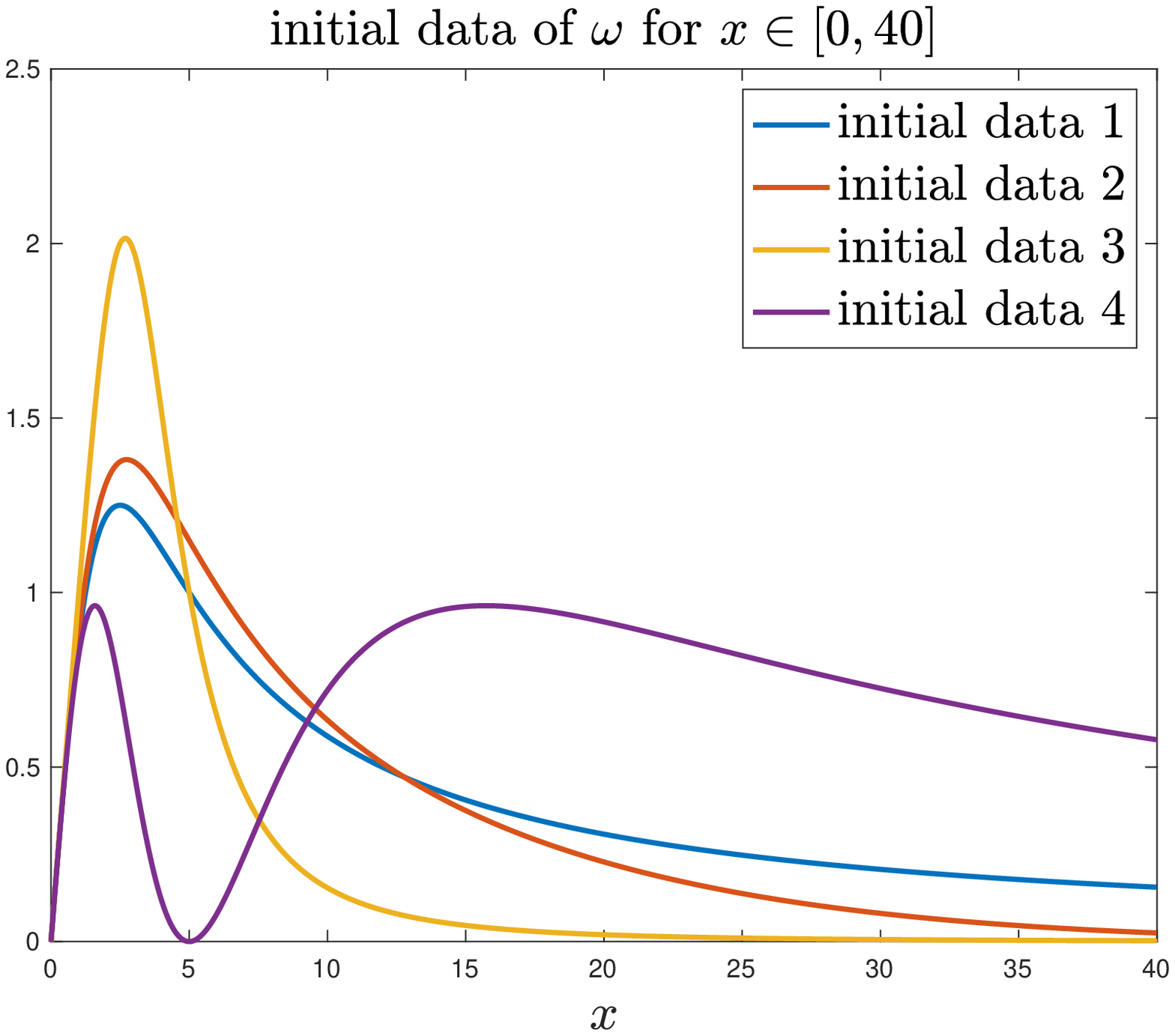}
    	\caption{}
    \end{subfigure}
    \begin{subfigure}[b]{0.32\textwidth}
    	\includegraphics[width=1\textwidth]{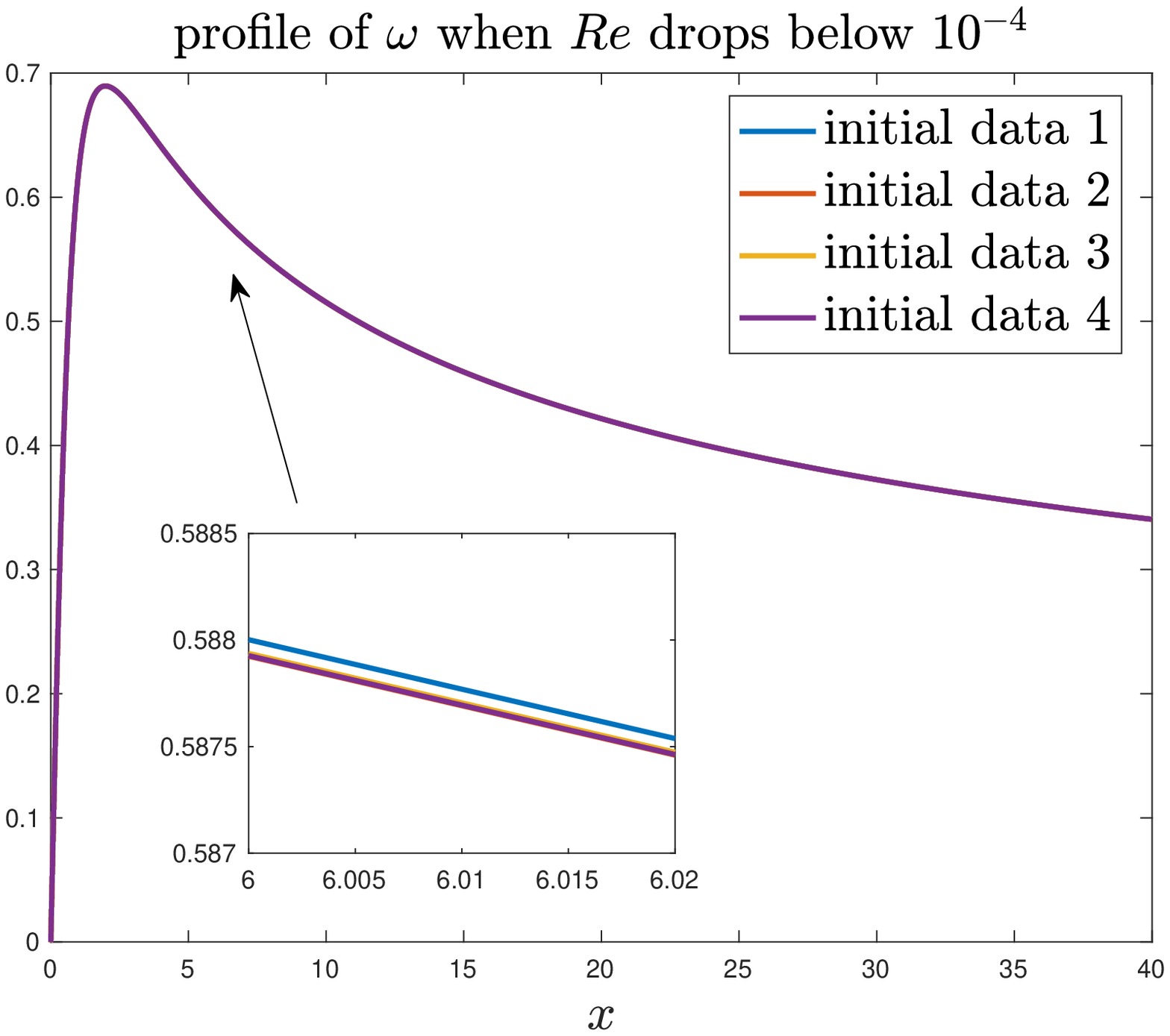}
    	\caption{}
    \end{subfigure}
    \begin{subfigure}[b]{0.32\textwidth}
    	\includegraphics[width=1\textwidth]{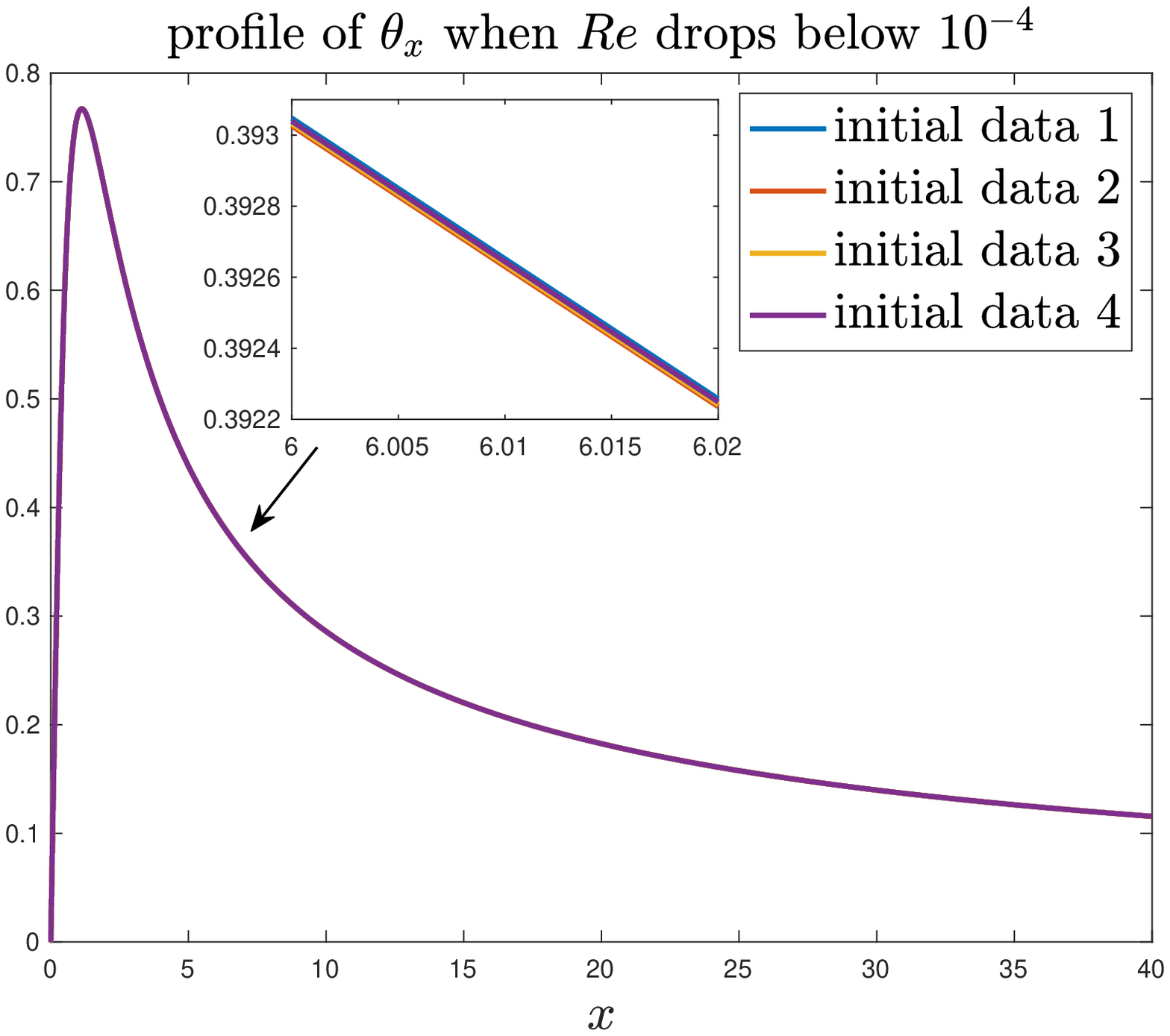}
    	\caption{}
    \end{subfigure}
    \begin{subfigure}[b]{0.32\textwidth}
    	\includegraphics[width=1\textwidth]{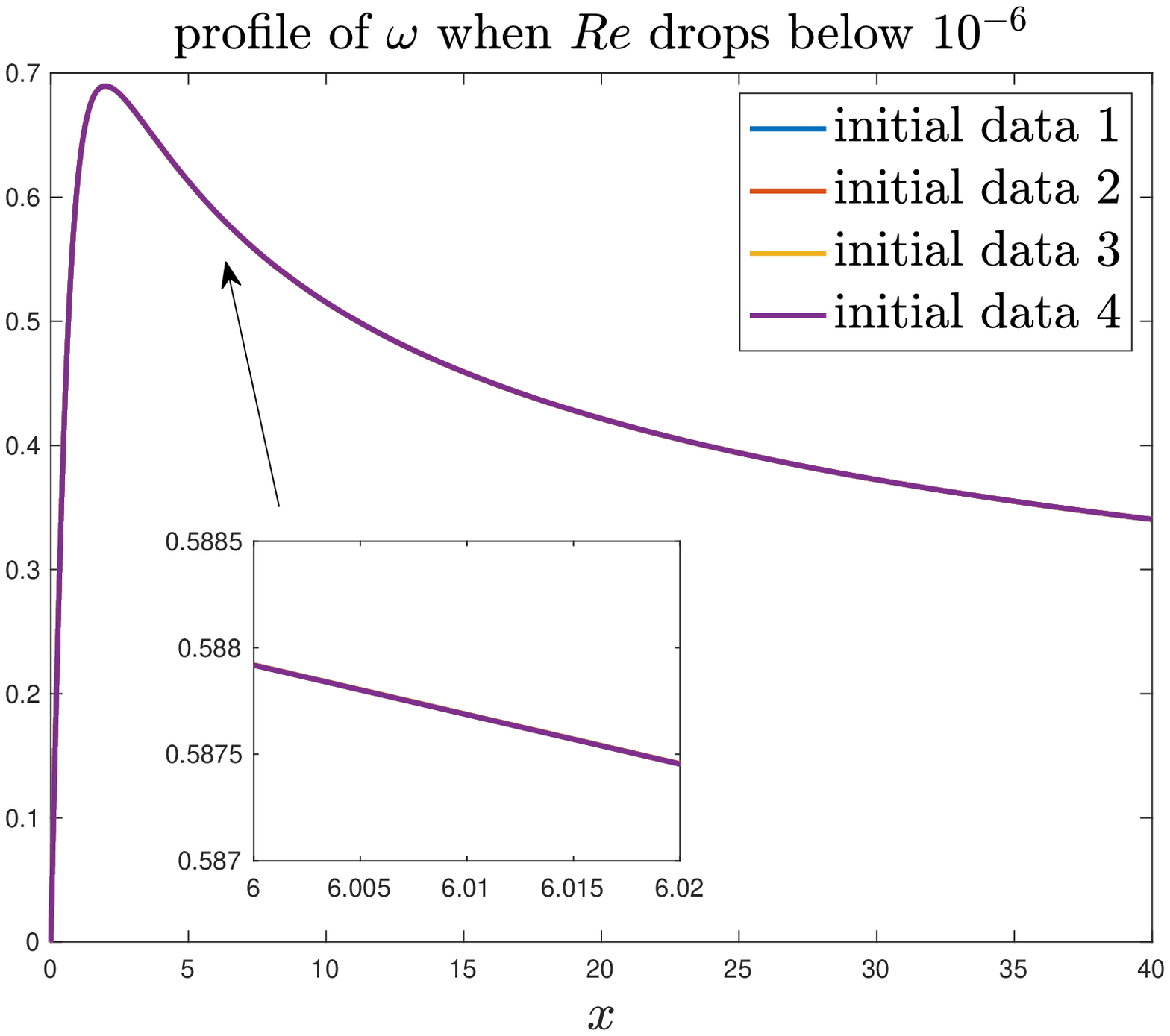}
    	\caption{}
    \end{subfigure}
    \begin{subfigure}[b]{0.32\textwidth}
    	\includegraphics[width=1\textwidth]{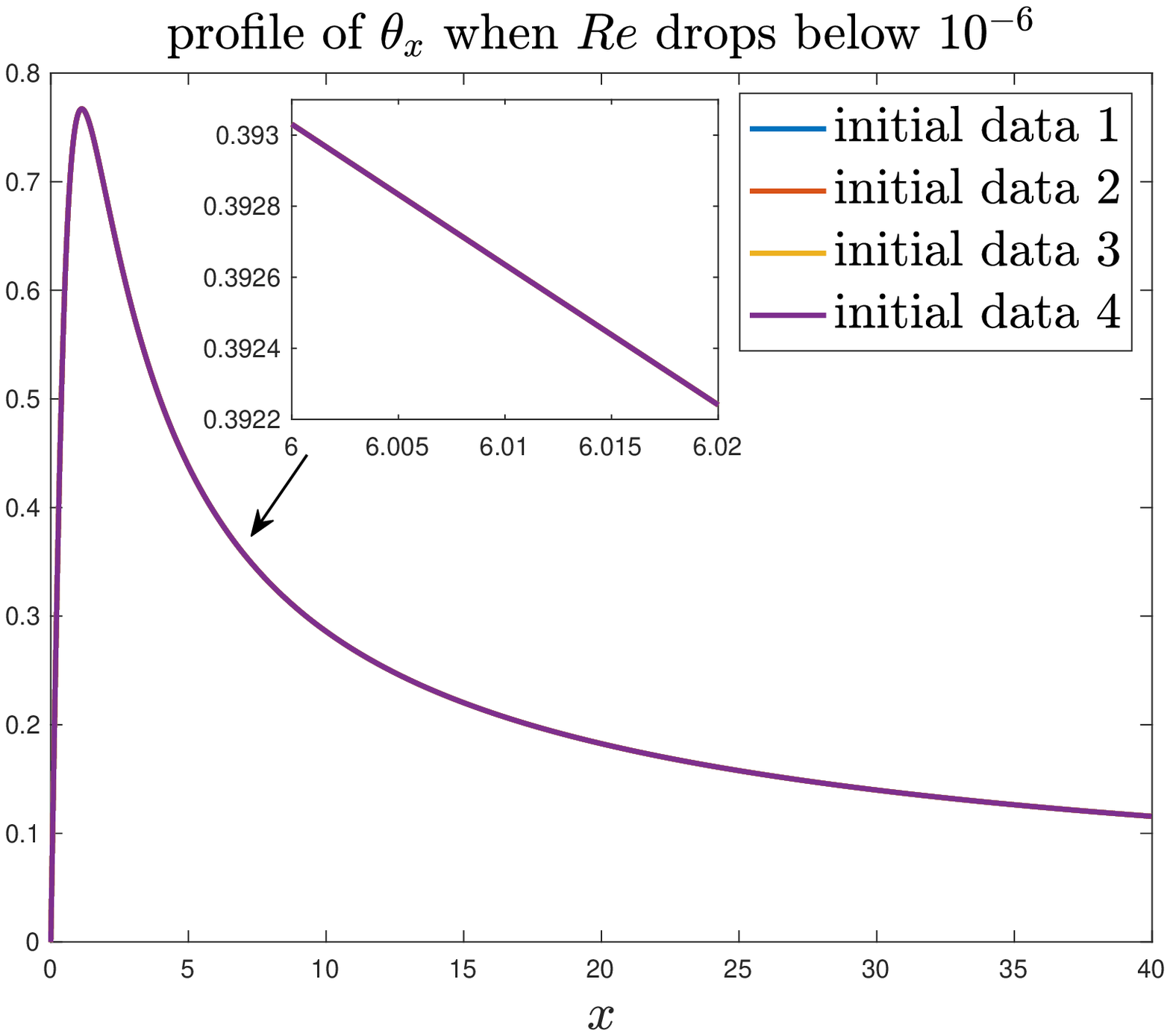} 
    	\caption{}
    \end{subfigure}
    \caption[Profiles of $\omega$]{(A) Four different initial data of $\omega$; (B)(C) Profiles of $\om$ and $\th_x$ when $Re$ drops below $10^{-4}$ the first time. (D)(E) Profiles of $\om$ and $\th_x$ when $Re$ drops below $10^{-6}$ the first time. }
    \label{fig:compare_initial_data}
\end{figure}

Here we present the convergence study for the dynamic rescaling equations for four sets of initial data that belong to the function class described above. The four initial data of $\om$ are given by $\om^{(i)}(x) = a_if_i(b_i x), i=1,2,3,4$, where 
\begin{align*}
f_1(x) = \frac{x}{1+x^2}, \quad f_2(x) = \frac{xe^{-(x/10)^2}}{1+x^2},\quad f_3(x) = \frac{x}{1+x^4},\quad f_4 = \frac{x(1-x^2)^2}{(1+x^2)^3},
\end{align*}
and the parameters $a_i,b_i$ are chosen to normalize the initial data such that they satisfy the same normalization conditions:
\[\om^{(i)}_x(0) = 1\quad \text{and}\quad u^{(i)}_x(0) = -2.5,\quad i = 1,2,3,4.\]
The initial data of $\th_x$ are chosen correspondingly as
\[\th_x^{(i)} = (c_lx + u^{(i)})\om^{(i)}_x - c_\om\om^{(i)},\quad i = 1,2,3,4,\]
so that the initial residual of the $\om$ equation is everywhere $0$. The initial value of the scaling parameters are set to be $c_l = 3$ and $c_\om = -1$, respecting our preliminary numerical result that $\bar{c}_l/\bar{c_\om} \approx -3$. Note that all these initial data of $\om,\th_x$ are far away from the approximate steady (with proper rescaling) with $O(1)$ distance in the the energy norm that is used in our analysis. In particular, we have $\om^{(1)}(x) = O(x^{-1})$, $\om^{(2)}(x) = O(x^{-1}e^{-(x/10)^2})$, $\om^{(3)}(x) = O(x^{-3})$ for $x\rightarrow +\infty$, while the approximate steady should satisfy $\bar{\om}(x) = O(x^{\bar{c}_\om/\bar{c}_l})$ where $\bar{c}_\om/\bar{c}_l$ is approximately $-1/3$ according to our numerical results. Moreover, $\om^{(4)}(x)$ has two peaks, while $\bar{\om}(x)$ only has one. Figure \ref{fig:compare_initial_data}(a) plots the four initial data of $\om$ for $x\in[0,40]$.

With each set of these initial data, we numerically solve the dynamic rescaling equations \eqref{eq:HLdyn} subject to the normalization conditions \eqref{eq:normal0} using the algorithm described in Section 10 of the Supplementary Material \cite{chen2021HLsupp} (by modifying the initial values of the part $\om_p$ and $(\th_x)_p$). We verify the uniqueness of the steady state by comparing the profiles of $\om$ at the first time the maximum grid-point residual $Re := \max_{i}\{|F_\om(x_i)|,|F_{\theta_x}(x_i)|\}$ drops below some small number $\epsilon$. Here the residuals $F_\om$ and $F_{\theta_x}$ are defined as
\beq\label{eq:residual}
F_\om = -(c_l x + u)\om_x + c_\om + v, \quad F_{\th_x} = -(c_l x + u)\th_{xx} + (2c_\om -u_x)\th_x.
\eeq
Figure \ref{fig:compare_initial_data} (b) and (c) plot the solutions of $\om$ when $Re\leq 10^{-4}$ and when $Re\leq10^{-6}$, respectively. We can see that the profiles of $\om$ from different initial data are barely distinguishable when the residual is smaller than $10^{-4}$; they become even closer to each other when the residual is even smaller. This implies that the solutions in the four cases of computation should converge to the same steady state. 

\section{ H\"older regularity of the blowup solution}\label{sec:hol}

To estimate the $C^{\g}$ norm with $\g= \f{c_{\th,\inf} }{c_{l,\inf}}$ of the solution $\th$, we will use the following estimate
\beq\label{eq:hol1}
\bal
\f{ |f(y) - f(x)| }{ |x-y|^{\g}} 
&= |x-y|^{-\g} |\int_x^y f_x(z) dz |
\les |x-y|^{-\g} \int_x^y z^{\g - 1} dz  \cdot || f_x x^{1-\g}||_{\infty} \\
&\les |x-y|^{-\g} ( y^{\g} - x^{\g}) \cdot || f_x x^{1-\g}||_{\infty} \les || f_x x^{1-\g}||_{\infty} .
\eal
\eeq
for any $0 \leq x < y$. The difficulty lies in the decay estimate of $\th_x$ since the previous a-priori estimates only imply that $\th_x$ decays with rate slower than $x^{\g-1}$.
The decay rate $x^{\g-1}$ is sharp since it is exactly the decay rate of the self-similar profile $\th_{\infty, x}$, which will be established in Section \ref{sec:decay_profile}. In Section \ref{sec:decay_pert}, we establish the decay estimates of the perturbation. In Section \ref{subsec:hol}, we estimate the H\"older norm of the solution.

\textbf{ Notations } In this Section, we use the notation $A \les B$ if there exists some finite constant $C > 0$, such that $A \leq C B $. The constant $C$ \textit{can} depend on the norms of the approximate steady state $( \bar \th, \bar \om)$ and the self-similar profile $(\th_{\inf}, \om_{\inf})$ constructed in Section \ref{sec:conv}, e.g. $|| \th_x||_{\inf}, || \bar \th_x||_{\inf}$, as long as these norms are finite. These constants do not play an important role in characterizing several exponents and thus we do not need to track them.

\subsection{ Decay estimates of the self-similar profile}\label{sec:decay_profile}
Recall that we have constructed the self-similar profile $(\th_{\inf}, \om_{\inf})$ in Section \ref{sec:conv}. Using the estimate \eqref{eq:profile_inf},we obtain 
\beq\label{eq:u_tail}
 |u_{\infty}(x) | \les |x|^{5/6},  \quad   | c_{l,\inf} x + u_{\infty}(x) | \geq 0.3 |x| , \quad u_{\inf, x}  \in L^{\inf} ,  \quad \th_{\infty}(1) \neq 0, \quad \th_{x, \inf} \in L^{\inf},
\eeq
whose proofs are referred to Section 10 in the Supplementary Material \cite{chen2021HLsupp}. Recall that the profile $( \th_{\infty}, \om_{\infty})$ solves 
\beq\label{eq:steady_th}
 ( c_{l,\infty} x + u_{\infty}   ) \th_{\infty, x} = c_{\th, \infty} \th_{\infty}, \quad  u_{\infty, x } = H \om_{\infty}.
\eeq

Solving the ODE on $\th_{\infty}$, we obtain 
\beq\label{eq:tail_1}
\th_{\infty}(x) = \th_{\infty}(1) \exp(J(x)), \quad J(x) \teq \int_1^x \f{c_{\th,\infty} }{ c_{l, \infty} y + u_{\infty}(y) } dy, 
\quad \th_{\infty, x} = \f{c_{\th,\infty} \th_{\infty}(x) }{ c_{l, \infty} x + u_{\infty}(x) }.
\eeq

 Denote $\g = \f{c_{\th,\infty}}{ c_{l,\infty}}$. Using the estimates on $u_{\inf}$ in \eqref{eq:u_tail}, we obtain $ | J(x) - \g \log(x)| \les 1$. Thus, for some constant $C_1 >0$ depending on the profile, we get 
 \[
\lim_{ x \to \infty }  \th_{\infty }(x) x^{-\g} = C_1 \th_{\infty}(1)  \neq 0.
\]

Plugging the above limit and \eqref{eq:u_tail} in the formula of $\th_{\infty,x}$ in \eqref{eq:tail_1}, we yield 
\beq\label{eq:hol20}
\lim_{x \to \infty} \th_{x,\infty} x^{1-\g} = 
\lim_{x \to \infty} 
 \f{c_{\th,\infty} x }{ c_{l, \infty} x + u_{\infty}(x) }
 \cdot  \th_{\infty}(x)  x^{-\g} =  C_1 \g \th_{\infty}(1) .
  \eeq
Combining the above estimate and $\th_{\inf, x} \in L^{\inf}$ from \eqref{eq:u_tail}, we prove 
\beq\label{eq:hol2}
|| \th_{\inf, x} x^{1- \g} ||_{\inf} \les 1.
\eeq

Differentiating \eqref{eq:steady_th} and using $c_{\th,\inf} = c_{l,\inf} + 2 c_{\om,\inf}$, we get 
\beq\label{eq:hol22}
(c_{l,\inf} x + u_{\inf}) \th_{\inf, xx} = (c_{\th, \inf} - c_{l,\inf} - u_{\inf,x }) \th_{\inf, x}
=( 2 c_{\om, \inf}  - u_{\inf,x }) \th_{\inf, x}.
\eeq

Using \eqref{eq:u_tail}, we further obtain 
\beq\label{eq:hol23}
\B| \f{x \th_{\inf, xx}}{\th_{\inf, x}} \B| =\B| \f{ ( 2 c_{\om, \inf}  - u_{\inf,x } ) x }{ c_{l,\inf} x+ u_{\inf}} \B|
\les \B| \f{x}{ c_{l,\inf} x+ u_{\inf}} \B| \les 1.
\eeq

\subsection{ Decay estimates of the perturbation}\label{sec:decay_pert}

Note that we have constructed $(\th_{\inf},\om_{\inf})$ in Section \ref{sec:conv} with estimate \eqref{eq:profile_inf}. We treat them as known functions. Similar to \eqref{eq:lin00}, \eqref{eq:NF}, linearizing the $\th_x$ equation around the self-similar profile, we get 
\[
\pa_t \th_x + ( c_{l,\inf} x + u_{\infty}  + u  ) \th_{xx} 
= (2 c_{\om,\inf}  -u_{\inf,x }) \th_x  + (2 c_{\om} - u_x ) \th_{\inf, x} - u \th_{\inf, xx}
+ (2 c_{\om} - u_x) \th_x,
\]
with normalization conditions 
\beq\label{eq:hol_normal}
c_{\om} = u_x(0) , \quad c_l = 0, \quad c_{\th} = c_l + 2 c_{\om} .
\eeq
Here, the nonlinear terms are given by $u \th_{xx}, (2 c_{\om} - u_x) \th_x$, and the error term is $0$ since we linearize the equation around the exact steady state. 
To obtain the decay estimates of $\th_x$ with a decay rate $O( |x|^{\g-1})$, we choose $\rho$ with a growth rate $O(|x|^{1-\g})$ and perform $L^{\inf}$ estimate on $\th_x \rho$, which will imply $|\th_x| \leq  |\rho^{-1}| \les |x|^{\g-1}$ for large $x$.
We derive the equation for $\th_x \rho$ as follows
\beq\label{eq:tail_2}
\bal
 \pa_t( \th_x \rho)
&+ ( c_{l,\inf} x + u_{\infty}  + u  ) ( \th_{x} \rho)_x 
= I(\rho) \th_x \rho  +  J , \\
 I(\rho) &\teq  2 c_{\om,\inf}  -u_{\inf,x } + (c_{l,\inf} x + u_{\infty}) \rho_x \rho^{-1} , \\
  J  &\teq  (2 c_{\om} - u_x ) \th_{\inf, x} \rho - u \th_{\inf, xx} \rho +
 u \th_x \rho_x+(2 c_{\om} - u_x) \th_x \rho.
\eal
\eeq
For a typical function $\rho$ with a growth rate $O(|x|^{\g-1})$, e.g. $\rho = \sgn(x) |x|^{\g -1}$,
since $u_{\infty}$ has sublinear growth \eqref{eq:u_tail}, for large $x>0$, we get 
\[
I(\rho) = 2 c_{\om, \infty} + c_{l,\infty} x  (x^{1-\g} )_x x^{\g-1} + l.o.t.
=2 c_{\om, \infty} + c_{l,\infty}( 1-\g) + l.o.t. = l.o.t.,
\]
where we have $ c_{l,\infty}( 1-\g)  = c_{l,\infty} - c_{\th, \infty} = -2 c_{\om,\infty} $ to obtain the last equality. Thus, we expect that $I(\rho)$ is not uniformly negative, i.e. $I(\rho) \leq -c$ for some $c>0$, and we do not obtain 
a damping term in the $L^{\infty}$ estimate of $\th_x \rho$, which is different from the weighted $L^2$ and $H^1$ estimates in Sections \ref{sec:lin}, \ref{sec:non}. In some sense, the decay rate $O( |x|^{\g-1})$ is critical. An ideal choice of $\rho$ with the desired growth rate is $ \th_{\infty,x}^{-1}$, since we have \eqref{eq:hol20} and $I(\rho)$ term in \eqref{eq:tail_2} vanishes : 
\[
I( \rho)= 2 c_{\om,\inf}  -u_{\inf,x } - \f{ (c_{l,\inf} x + u_{\infty}) \th_{\inf, xx}  }{ \th_{\inf, x}} = 
\f{ (2 c_{\om,\inf}  -u_{\inf,x } ) \th_{\inf, x} - (c_{l,\inf} x + u_{\infty}) \th_{\inf, xx}  }{\th_{\inf, x}} = 0,
\]
where we have used \eqref{eq:hol22} to obtain the last equality.

Recall $c_{\om} = u_x(t,0)$. Using $\rho = \th_{\inf, x}^{-1}$, $ | \f{ x \th_{\inf,xx} }{\th_{\inf, x}} | \les 1$ in \eqref{eq:hol23} and $ | \f{u}{x}| \les || u_x ||_{\inf} $,  we get 
\[
\bal
|J | &= \B| (2 c_{\om} - u_x ) \f{ \th_{\inf, x} }{ \th_{\inf, x}} - u  \f{ \th_{\inf, xx} }{ \th_{\inf,x}} - u \th_x  \f{\th_{\inf, xx}}{\th_{\inf,x}^2 } +(2 c_{\om} - u_x) \th_x \rho \B| 
\les  || u_x||_{\inf} (1 + || \th_x \rho ||_{\inf}).
\\
\eal
\]

For $\th_x(\cdot ,0) \rho \in L^{\infty}$, performing $L^{\infty}$ estimates in \eqref{eq:tail_2}, we yield 
\beq\label{eq:tail_3}
\f{d}{dt} || \th_x \rho ||_{\infty}
\les || u_x||_{\infty} (1 + || \th_x \rho||_{\inf} ).
\eeq

Next, we control $|| u_x ||_{\infty} $. Recall the energy $E$ in \eqref{energy:h1} and the a-priori estimates in \eqref{eq:boot},\eqref{eq:profile_inf} 
\[
 E( \th_{x, \infty} + \th_x - \bar \th_x, \om_{\inf} + \om - \bar \om )  \leq E_*,
 \quad  E( \th_{x, \infty}  - \bar \th_x, \om_{\inf}  - \bar \om ) \leq E_*.
\]

Using the triangle inequality, for any $ t\geq 0$, we get 
\beq\label{eq:hol3}
|| \th_x \psi^{1/2} ||_2 + || D_x \th_x \psi^{1/2} ||_2 + || \om \vp^{1/2} ||_2
+ || D_x \om \vp^{1/2} ||_2 + | c_{\om}(\om) | + | d_{\th}(\th_x)| \les 1.
\eeq

Denote $\kp_3 = 0.02$. Applying  \eqref{eq:EE_uni} to two solutions $( \th_{\infty}, \om_{\infty})$ and $ (\th_{\infty} + \th, \om_{\infty} + \om)$, we get 
\beq\label{eq:hol32}
|| \th_x(t) \psi^{1/2} ||_2 + || \om(t) \vp^{1/2}||_2  + |c_{\om}(t)|
\les E_1( \th_x(t), \om(t) ) \leq e^{-\kp_3 t} E_1( \th_x(0), \om(0)  \les e^{-\kp_3 t},
\eeq
where we have used \eqref{eq:hol3} to obtain the last inequality. Since $H(D_x \om)(0) = 0$, using Lemma \ref{lem:com}, we get $ H (D_x \om)= D_x H \om =x u_{xx}$. From \eqref{eq:mode} and \eqref{eq:wg}, we have $x^{-4/3} + x^{-2/3} \les \vp$. Applying Lemma \ref{lem:wg23} to $f = D_x \om $ and $f = \om$ (note that $H(D_x\om)(0) = 0$), we obtain 
\beq\label{eq:hol33}
\bal
|| u_x||_{\infty}^2 &\les \int_{\R}  |u_{xx} u_x| dx 
=  \int_{\R}  |  H (D_x \om )  \cdot H \om  x^{-1}| dx 
\les || H(D_x \om) x^{-2/3} ||_2 || H \om x^{-1/3}||_2 \\
&\les || D_x \om x^{-2/3}||_2 || \om x^{-1/3}||_2 
\les || D_x \om \vp^{1/2} ||_2 || \om \vp^{1/2}||_2 
\les e^{-\kp_3 t/ 2} .
\eal 
\eeq

Plugging the above estimate in \eqref{eq:tail_3}, we yield 
\[
\f{d}{dt} || \th_x \rho||_{\inf} \les e^{-\kp_3 t / 4} (1 + || \th_x \rho||_{\inf} ).
\]
Since $\kp_3>0$, solving the differential inequality and using $ | x^{1-\g}| \les | \th_{x,\inf}^{-1}|$ from \eqref{eq:hol2}, we prove 
\[
\sup_{t \geq 0} || \th_x (t) \rho ||_{\inf} \les 1, \quad 
\sup_{t \geq 0} || \th_x(t) x^{1-\g} ||_{\inf}
\les \sup_{t \geq 0} || \th_x(t) \th_{x,\inf}^{-1} ||_{\inf} \les 1,
\]
Since $\th$ is even, using \eqref{eq:hol1}, \eqref{eq:hol2} and the above estimate, we prove 
\beq\label{eq:hol4}
\sup_{t \geq 0} || \th_{\inf} + \th(t) ||_{C^{\g }} \les 1.
\eeq

\begin{remark}
Since we do not have a damping term in the $L^{\inf}$ estimate \eqref{eq:tail_3}, the exponential convergence estimates in \eqref{eq:hol32}, \eqref{eq:hol33} play a crucial role in obtaining \eqref{eq:hol4}.

\end{remark}

\subsection{ H\"older regularity}\label{subsec:hol}
Denote $\hat \th = \th_{\inf} + \th$ and by $\th_{phy}$ the solution with initial data $\hat \th(0, \cdot)$ in the physical space. Recall the rescaling relation and the normalization conditions \eqref{eq:hol_normal}
\beq\label{eq:hol50}
\bal
C_{\om}(\tau) &= \exp( \int_0^{\tau} c_{\om}(s) + c_{\om,\inf} ds), \quad  t(\tau) =\int_0^{\tau} C_{\om}(s) ds , \\
  C_{\th}(\tau) &= \exp(\int_0^{\tau}  c_{\th}(s) + c_{\th,\inf} ds ),
\quad   C_l(\tau) = \exp(- \int_0^{\tau} (c_l(s) + c_{l,\inf}) ds ) , \\
\hat \th(x, \tau) &= C_{\th}( \tau) \th_{phy}( C_l(\tau) x, t( \tau)), \quad  c_{\th} = 2 c_{\om}, \quad c_l = 0.
\eal
\eeq
From assumptions $\hat \th_x(0) |x|^{ 1- \g } \in L^{\inf}$ in (d) in Theorem \ref{thm2}, $E( \hat \th_x(0) - \bar \th_x, \hat \om(0) - \bar \om) \les 1$, and estimates \eqref{eq:hol20} and $E( \th_{\infty,x} - \bar \th_x, \om_{\infty} - \bar \om) \les 1$, it is not difficult to obtain that $\th_x \rho \in L^{\inf} $. Thus, $\th, \hat \th$ enjoys the energy estimates in Section \ref{sec:decay_pert}.  
Using \eqref{eq:hol32}, \eqref{eq:hol4}, \eqref{eq:hol50} and  $\g c_{l,\inf} = c_{\th,\inf}$, we prove 
\[
\sup_{ \tau \geq 0}
||\th_{phy}(t(\tau))||_{C^{\g}} = \sup_{\tau \geq 0} || \hat \th(\tau)||_{C^{\g}}  C_{\th}^{-1} C_l^{-\g}
=  \sup_{\tau \geq 0}  || \hat \th(\tau)||_{C^{\g}} \exp( \int_0^{\tau} -2 c_{\om} d \tau)
\les 1.
\]

\subsubsection{Blowup in higher H\"older norm}

We show that for any $\b > \g$, the $C^{\b }$ norm of the solution blows up. Since $
1 \les  \psi(x) $ for $x \in [0,1]$, using \eqref{eq:hol32} and Cauchy-Schwarz inequality, we get 
\beq\label{eq:hol5}
| \th(1) - \th(0)| = | \int_0^1  \th_x(y) dy | 
\les \B( \int_0^1 \th_x(y)^2 dy \B)^{1/2} \les || \th_x \psi^{1/2}||_2  \les e^{-\kp_3 \tau }.
\eeq

Recall the formulas in \eqref{eq:hol50}. Denote $T = t(\inf)$. Since $|c_{\om}(\tau)|$ decays exponentially \eqref{eq:hol32} and $c_{\om, \inf} < -\f{1}{2}$, 
we obtain 
\[
\bal
&C_{\om}(\tau) \gtrsim e^{ c_{\om,\inf} \tau},  \quad C_{\th}(\tau)^{-1} \gtrsim e^{ - c_{\th,\inf} \tau},\quad C_l(\tau)^{-1} \gtrsim e^{ c_{l,\inf} \tau } , \\
& T - t(\tau) = \int_{\tau}^{\inf} C_{\om}(s) ds \gtrsim \int_{\tau}^{\inf} e^{ c_{\om, \inf} s} ds 
\gtrsim e^{ c_{\om, \inf} \tau} .
\eal
\]
Recall $\g c_{l,\inf} =c_{\th, \inf } = c_{l,\inf} + 2 c_{\om,\inf}$. Denote $\d =- \f{ \b c_{l,\inf} - c_{\th,\inf}} {c_{\om,\inf}}  = \f{ 2( \b-  \g)}{1 - \g} > 0$. We have
\[
S \teq \liminf_{ \tau \to \inf} ||  \th_{phy}(x, \tau)||_{C^{\b}}
( T- t(\tau))^{\d} 
\gtrsim \liminf_{ \tau \to \inf} || \hat \th(x, \tau)||_{C^{\b}}
C_{\th}^{-1} C_l^{-\b} \exp( \d c_{\om, \inf} \tau).
\]

Note that $\th_{\inf}(0) = 0$. Using \eqref{eq:hol5}, we have  $ || \hat \th(\tau) ||_{C^{\b}} \geq |\hat \th(\tau, 1) - \hat \th(\tau , 0)|
\geq | \th_{\inf}(1)| - C \exp( -\kp_3 \tau)$. Using this estimate, $\d =- \f{ \b c_{l,\inf} - c_{\th,\inf}} {c_{\om,\inf}}$ and \eqref{eq:u_tail}, we establish 
\[
S \gtrsim \liminf_{\tau \to \inf} |\th_{\inf}(1)| 
\exp( ( -c_{\th,\inf}  + \b c_{l,\inf} + \d c_{\om,\inf}) \tau ) 
\gtrsim |\th_{\inf}(1)| >0.
\]
We conclude the proof of result (d) in Theorem \ref{thm2}.

\begin{remark}
The exponential convergence in \eqref{eq:hol32} is crucial for us to obtain the unique H\"older exponent $\g$ that characterizes the regularity of the singular solution and the sharp blowup rate. It enables us to essentially treat the perturbation as $0$.
\end{remark}

\section{ Connection between the HL model and the 2D Boussinesq equations in $\R_2^+$}\label{sec:con}

In this section, we discuss the connection between the leading order system of the HL model and that of the 2D Boussinesq equations in $\R_2^+$ with low regularity initial data. 

\subsection{The leading order system for the 2D Boussinesq equations }

The 2D Boussinesq equations in $\R_2^+$ read
\beq\label{eq:bous}
\bal
\om_t +  \uu \cdot \na \om  &= \th_{x},  \\
\th_t + \uu \cdot  \na \th & =  0 , \\
\eal
\eeq
where the velocity field $\uu = (u , v)^T : \R_+^2 \times [0, T) \to \R^2_+$ is determined via the Biot-Savart law
\[
 - \D \psi = \om , \quad  u =  - \psi_y , \quad v  = \psi_x,
 \]
with no flow boundary condition $\psi(x, 0 ) = 0   \quad x \in \R$.

Consider the polar coordinate $( r, \b)$ in $\R_2^+$ : $r = (x^2 + y^2)^{1/2}, \b = \arctan(y/x)$. For $\al > 0$, denote 
\[
R = r^{\al}, \quad  \Om(R, \b) = \om( x, y), \quad \eta(R, \b) = \th_x( x, y),  \quad \xi(R, \b) = \th_y(x, y).
\]

In \cite{chen2019finite2}, the following leading order system of \eqref{eq:bous} is derived based on the framework developed in \cite{elgindi2019finite} under the assumption that $\om, \na \th$ are in some  H\"older space $C^{\al}$ with sufficient small $\al$
\beq\label{eq:lead_bous}
\Om_t = \eta, \quad \eta_t = \f{2}{\pi \al} L_{12}(\Om) \eta, 
\quad L_{12}(\Om) =  \int_R^{\infty} \int_0^{\pi/2} \f{\Om( s ,\b) \sin(2\b) }{s} ds d \b.
\eeq

An important observation made in \cite{chen2019finite2} is that for certain class of $C^{\al}$ data, $\th$ is anisotropic in the sense that $|\th_y| \les \al |\th_x| $. Moreover, this property is preserved dynamically. Therefore, the $\th_y$ variable does not appear in the leading order system. Define the following operators 
\beq\label{eq:oper_PS}
P f(R) = \int_0^{\pi/2} f(R, \b) \sin( 2\b) d\b, \quad 
S f(R) = \f{2}{\pi \al} \int_R^{\infty} f( S) \f{dS}{S}.
\eeq

By definition, we have 
\beq\label{eq:oper_LS}
\f{2}{\pi \al} L_{12}(\Om) = \f{2}{\pi \al} \int_R^{\infty} P \Om(s) \f{ds}{s} = S( P \Om ).
\eeq

Since $L_{12}(\Om)$ does not depend on $\b$, we apply the operator $P$ to both sides of \eqref{eq:lead_bous} to obtain 
\beq\label{eq:lead_bous2}
\pa_P \Om = P \eta, \quad \pa_t P \eta = \f{2}{\pi \al} L_{12}(\Om) P \eta
= S (P\Om) \cdot P \eta.
\eeq

The above system is an 1D coupled system on $P\Om, P\eta$. Once $P\Om, P\eta$ are determined, we can obtain an explicit solution of \eqref{eq:lead_bous}.

\subsection{The leading order system for the HL model }\label{sec:lead_HL}

We use the observation made in \cite{Elg17} that the advection can be substantially weakened by choosing $C^{\alpha}$ data with sufficiently small $\al$. Suppose that $\om, \th_x \in C^{\al}$ with small $\al$. Then the advection terms in the system of $(\om,\th_x)$ in the HL model become
lower order terms  
\beq\label{eq:lead_HL1}
\om_t  = \th_x + l.o.t., \quad (\th_x)_t  = -u_x \th_x + l.o.t. ,\quad u_x = H \om.
\eeq

The above system is already very similar to \eqref{eq:lead_bous} by taking $\Om =\om, \eta =\th_x$. We further perform a simplification for the Hilbert transform. We impose extra assumptions that $\om,\th_x$ are odd, which are preserved dynamically. Due to these symmetries, it suffices to consider the HL model on $\R_+$. For $x >0$, symmetrizing the kernel, we get 
\[
\bal
H \om (x)  &= \f{1}{ \pi} \int_{\R_+} \om(y) ( \f{1}{x- y} -\f{1}{x+y}) dy 
= \f{1}{\pi} \int_{\R_+} \om(y) \f{2y}{x^2 - y^2} dy  = \f{1}{ \pi} \int_{\R_+}  \om(y) \f{ 2 }{ (x/y)^2 - 1} \f{dy}{y} .
\eal
\]

We learn the following formal derivation of the leading order part of general singular integral operator from Dr. Elgindi. \footnote{ Similar derivation was presented in the One World PDE Seminar "Singularity formation in incompressible fluids" by Dr. Elgindi. \url{https://www.youtube.com/watch?v=29zUjm7xFlI&feature=youtu.be}} Denote 
\beq\label{eq:change_HL}
X = x^{\al}, \quad Y = y^{\al}, \quad \Om( X) = \om( x ), \quad \eta(X) = \th_x( x).
\eeq
Using the above change of variables and $\f{dy}{y} = \f{1}{\al} \f{dY}{Y}$, we get 
\[
H \om(x) = \f{1}{\al \pi} \int_{\R_+} \om( Y^{1/\al}) \f{2}{ (\f{X}{Y})^{1/\al}  - 1} \f{dY}{ Y}
= \f{1}{\al \pi} \int_{\R_+} \Om(Y) K_{\al}(X, Y) \f{dY}{ Y}, 
\]
where $K_{\al}(X,Y) =  \f{2}{ (\f{X}{Y})^{1/\al}  - 1}$. Next, we consider the leading order part of $K_{\al}(X, Y)$ as $\al \to 0^+$. Note that 
\[
\lim_{\al \to 0^+} ( \f{X}{Y})^{1/\al} = 0 ,\ \mathrm{ for }  \ X < Y , \quad 
\lim_{\al \to 0^+} ( \f{X}{Y})^{1/\al} = \infty, \ \mathrm{ for } \ X > Y.
\]
Hence, for $X \neq Y$ and $X, Y >0$, we get 
\[
\lim_{\al \to 0^+} K_{\al}(X, Y) = -2 \cdot \one_{Y > X}.
\]
Therefore, formally, we get 
\beq\label{eq:Hilbert_simp}
H \om(x) = - \f{2}{\al \pi} \int_X^{\infty} \om(Y) \f{dY}{Y} + l.o.t. 
= -S \Om (X) + l.o.t.,
\eeq
where the operator $S$ is defined in \eqref{eq:oper_PS}. Now, plugging the above formula in \eqref{eq:lead_HL1}, dropping the lower order terms in \eqref{eq:lead_HL1} and applying the notations \eqref{eq:change_HL}, we derive another leading order system for the HL model 
\beq\label{eq:lead_HL2}
\pa_t \Om(X) = \eta(X) , \quad \pa_t  \eta (X)  = S \Om(X)  \cdot \eta (X).
\eeq

The above system is exactly the same as that in \eqref{eq:lead_bous2}. We remark that the lower order term in the simplification \eqref{eq:Hilbert_simp} needs to be estimated rigorously. In general, the system \eqref{eq:lead_HL1} is more complicated than \eqref{eq:lead_HL2} since the Hilbert transform is nonlocal and is a singular operator, while we can obtain a local relation between $ Sf$ and $f$ by taking derivative $\pa_X ( Sf)(X) = -\f{2}{\pi \al} \f{f(X)}{X}  $.

Note that $\one_{X < Y} = \one_{x < y}$. Undoing the change of variables in \eqref{eq:change_HL}, we get
\beq\label{eq:oper_SH}
S \Om(X) =  \f{2}{\pi \al} \int_{\R_+} \one_{x < y} \Om( Y ) \f{dY}{Y}
= \f{2}{\al\pi} \int_{\R_+} \one_{x< y} \om(y) \cdot \al \f{dy}{ y}
= \f{2}{\pi}\int_x^{\infty} \om(y) \f{dy}{y}.
\eeq

The operator on the right hand side is closely related to the Choi-Kiselev-Yao (CKY) simplification of the Hilbert transform \cite{choi2015finite}. Therefore, the leading order system \eqref{eq:lead_HL2} can be seen as the CKY's simplification of \eqref{eq:lead_HL1} without the lower order terms.

\section{Concluding remarks}\label{sec:conclude}

In this paper, we proved that the HL model develops a finite time focusing asymptotically self-similar blowup from smooth initial data with compact support and finite energy. Moreover, we showed that the solution of the dynamic rescaling equations converges to an exact steady state exponentially fast in time and the self-similar blowup profile is unique within a small energy ball. We also presented strong numerical evidence to demonstrate the uniqueness of the self-similar profile for a much larger class of initial data that satisfy certain symmetry and sign conditions consistent with the initial data considered by Luo-Hou in \cite{luo2014potentially,luo2013potentially-2}. The possibility of having a unique self-similar profile for a large class of initial data is very interesting and quite surprising if it can be justified rigorously.

One of the main difficulties in our stability analysis is to control a number of nonlocal terms with a relatively small damping coefficient. This is also the essential difficulty in generalizing the method of analysis presented in this work to prove the finite time blowup of the 2D Boussinesq equations or 3D axisymmetric Euler equations with smooth initial data and boundary. To establish linear stability, we designed singular weight functions carefully, applied several sharp weighted functional inequalities to control the nonlocal terms, and took into account cancellation among various nonlocal terms.

Our ultimate goal is to prove rigorously the Luo-Hou blowup scenario for the 2D Boussinesq equations and 3D Euler equations with smooth initial data and boundary.
Our numerical study suggested that the real parts of the eigenvalues of the discrete linearized operator for the 2D Boussinesq equations with smooth initial data and boundary are all negative and bounded away from 0 by a finite spectral gap. See also Section 3.4 in Dr. Pengfei Liu's Ph.D. thesis \cite{liu2017spatial} for an illustration of the eigenvalue distribution of the discretized linearized operator. Moreover, our numerical study shows that $|\th_y|$ is an order of magnitude smaller than $|\th_x|$. This seems to imply that the main driving mechanism for singularity formation is due to the coupling between $\om$ and $\th_x$, which is captured by our analysis for the HL model.

The framework of analysis that we established for the HL model provides a promising approach 
to studying the singularity formation of the 2D Boussinesq equations and 3D axisymmetric Euler equations with smooth initial data and boundary.
We can follow the general strategy developed
in this paper by (1) extracting the damping effect from the local terms, (2) treating the advection terms as perturbation to vortex stretching, and (3) controlling the nonlocal terms 
by developing sharp functional inequalities on the Biot-Savart law and  exploiting cancellation among them to control the nonlocal terms by using the damping effects from the local terms. Compared with the HL model, we will encounter some additional difficulties associated with the advection away from the boundary, and need to estimate more complicated Biot-Savart law in 2D Boussinesq and 3D Euler equations.
We will explore a more effective functional space, e.g. weighted $L^p$ or weighted $C^{\alpha}$ space, to establish the stability analysis. Such space offers the advantage of weakening the effect of the advection in the stability analysis and extracting larger damping effect from the local terms in the linearized equations. Moreover, it still allows us to estimate the Biot-Savart law effectively.

Guided by the singularity analysis presented in this paper, we have recently made some encouraging progress towards the ultimate goal of proving finite time self-similar blowup of the 2D Boussinesq equations and 3D Euler equations with smooth initial data and boundary. We will report our results in our future work.

\vspace{0.2in}
{\bf Acknowledgments.} The research was in part supported by NSF Grants DMS-1907977 and DMS-1912654. D. Huang would like to acknowledge the generous support from the Choi Family Postdoc Gift Fund.

 \appendix

\section{Properties of the Hilbert transform}
Throughout this section, we assume that $\om$ is smooth and decays sufficiently fast. The general case can be obtained by approximation. The properties of the Hilbert transform in  Lemmas \ref{lem:com}-\ref{lem:anti} are well known, see e.g. \cite{chen2019finite,duoandikoetxea2001fourier,chen2020singularity}.

\begin{lem}\label{lem:com}
Assume that $\om$ is odd. We have 
\[
 H\om(x) - H \om(0) = x H( \f{\om}{x}).
\]
\end{lem}

\begin{lem}\label{lem:iso}

Assume that $\om$ is odd and $\om_x(0) = 0$. For $p=1,2$, we have
\beq\label{eq:hil}
\bal
( u_x  - u_x(0)) x^{-p} = H( \om x^{-p}).
\eal
\eeq
Consequently, the $L^2$ isometry property of the Hilbert transform implies 
\[
||( u_x  - u_x(0)) x^{-p} ||_2^2 = || \om x^{-p} ||_2^2.
\]
\end{lem}

Recall the inner product $\la f, g \ra = \int_0^{\infty} fg dx $ (see \eqref{nota:inner} )and $\Lam = (-D)^{1/2} = H\pa_x$.

\begin{lem}\label{lem:anti}
For $f \in L^p, g \in L^q$ with $\f{1}{p} + \f{1}{q} = 1$ and $1< p < \infty$, we have
\beq\label{eq:anti}
\la H f , g\ra = -\la f, H  g\ra.
\eeq
\end{lem}


\begin{lem}\label{lem:cancel}
Denote $\Lam = (-\pa_x^2)^{1/2}$. Assume that $f$ is odd and $g_x = H f, g(0) = 0$. We have 
\[
\bal
& \la H f - H f (0), f x^{-3} \ra  = 0, 
 \quad \la g, f_x x^{-1} \ra = - \B\la \Lambda \f{ g}{x} , \f{ g }{x}  \B\ra ,
\quad \la g, f x^{-2} \ra = - \B\la \Lambda \f{ g}{x} , \f{ g }{x}  \B\ra - \f{\pi}{4} g_x(0)^2.
\eal
\]

\end{lem}

Identities similar to those in Lemma \ref{lem:cancel} have been used in \cite{chen2019finite,Elg19,Cor10,chen2020singularity}. 

\begin{proof}
The proof of the first identity can be found, e.g. in \cite{chen2019finite,elgindi2019stability}. We focus on the last two identities. Recall the Cotlar identity from, e.g. \cite{chen2019finite,duoandikoetxea2001fourier},
\beq\label{eq:colt}
(H F)^2 = F^2 + 2 H ( F\cdot H F)  .
\eeq

Since $f$ is odd, we get that $g, H(\f{g}{x})$ are odd and $H(\f{g}{x})(0) =0$. Applying Lemma \ref{lem:com} with $\om = f$ and \eqref{eq:colt}, we have the following identities 
\[
\bal
 \f{ f}{x} &= - H H  (\f{f}{x}) = - H ( \f{g_x - g_x(0)}{x} ), 
 \\
\int_{\R} H( \f{g}{x}) \f{g }{x^2}  dx
&= - \pi  H( \f{g}{x} H(\f{g}{x}))(0)
= - \f{\pi}{2}  \B(  ( H( \f{g}{x} )(0))^2 -( \f{g}{x}(0))^2 \B)
= \f{\pi}{2} g_x^2(0),
 \\
\int_{\R} H( \f{g}{x} ) \f{g_x(0) x}{x^2} dx 
&=  g_x(0)  \int_{\R} H(\f{g}{x}) \f{1}{x} dx
= -  \pi g_x(0)   H( H ( \f{g}{x}) )(0)
= \pi g_x(0) \f{g}{x} \B|_{x=0} 
= \pi g_x^2(0).
\eal
\]

For the third identity, using these identities and Lemma \ref{lem:anti}, we obtain 
\[
\bal
 \int_{\R} \f{ g f }{x^2}  dx
 &= - \int_{\R} \f{ g  }{x} H \f{ g_x -g_x(0)}{x} dx 
 = \int_{\R} H (\f{g}{x} ) \f{ g_x - g_x(0)}{x} dx 
 = \int_{\R} H (\f{g}{x} ) \B( (\f{g}{x})_x + \f{g - g_x(0)x}{x^2} \B)dx \\
& = \int_{\R} - \pa_x H (\f{g}{x} ) \cdot \f{g}{x} 
+ H (\f{g}{x} )  \f{g - g_x(0)x}{x^2})dx 
= -\int_{\R} \Lam (\f{g}{x}) \f{g}{x} dx- \f{\pi}{2} g_x^2(0),
\eal
\]
where we have used $\Lambda = (- \pa_{xx})^{1/2} =  H \pa_x$ in the last identity. Restricting the integrals to $\R_+$, we prove the third identity. For the second identity, using integration by parts, we yield 
\[
\int_{\R} \f{g f_x}{x} dx 
= -\int_{\R} \f{g_x f }{x} + \f{g f}{x^2} dx 
= \pi H( g_x f)(0) -\int_{\R} \Lam (\f{g}{x}) \f{g}{x} dx- \f{\pi}{2} g_x^2(0).
\]
 Using \eqref{eq:colt}, we yield 
\[
\pi H( g_x f)(0) = \pi H( f \cdot Hf)(0) = \f{\pi}{2} ( (H f(0) )^2 - (f(0))^2 =  \f{\pi}{2} 
(Hf(0))^2 = \f{\pi}{2} g_x^2(0).  \\
\]
Combining the above identities, we  complete the proof of the Lemma.
\end{proof}

\begin{lem}\label{lem:rota}
Assume that $\om \in L^2( |x|^{-4/3} + |x|^{-2/3})$ is odd and $u_x = H \om$. We have 
\[
\bal
\int_{\R} \f{ (u_x (x) - u_x(0))^2} {  |x|^{4/3} } =  \int_{\R} \B( \f{ w^2} { | x|^{4/3} }  +
2  \sqrt{3} \cdot \mathrm{sgn}(x) \f{  \om  ( u_x(x) - u_x(0) ) }{|x|^{4/3} } \B) dx .
\eal
\]
\end{lem}

It seems that the identity \eqref{eq:Hil_alpha} $H (  |x|^{-\al} )    = \tan\lt( \f{  \al \pi }{2}   \rt) sgn(x)  |x|^{-\al}$, which will be used in the proof of Lemma \ref{lem:rota}, is difficult to locate in the literature. We thus give a proof.

\begin{proof}
Firstly, we compute $H (|x|^{-\al})$.
For $\al \in (0, 1)$, we have $H (  |x|^{-\al} )    = C_{\al} sgn(x)  |x|^{-\al}$, 
for some constant $C_{\al}$. We determine $C_{\al}$ by applying Lemma \ref{lem:anti} to 
\[
\bal
f = |x|^{-\al}, \quad H f = C_{\al} sgn(x)  |x|^{-\al} , \quad g = - \f{x}{1 + x^2}, \quad Hg = \f{1}{1 + x^2},
\eal
\]
which implies 
\[
C_{\al} \int_0^{\inf } \f{x^{1-\al}}{ 1 +x^2} d x = \int_0^{\inf} \f{ 1}{x^{\al}(1 + x^2)} dx.
\]
The integrals can be evaluated using the Beta function $\mathrm{B}(x, y)$ and $\mathrm{B}( \b, 1-\b) =\f{\pi}{\sin(\b\pi)}$ for $\b \in (0,1)$. In particular, we get 
\[
C_{\al} = \f{ \mathrm{B}(\f{\al + 1}{2} , \f{1 - \al}{2})} {        \mathrm{B}(\f{2 - \al}{2} , \f{ \al}{2})    }  =  \f{ \pi /  \sin(  (\al + 1)  \pi/ 2     )         }  {         \pi /  \sin(  (2 - \al )  \pi / 2     )           }
=\tan\lt( \f{  \al \pi }{2}   \rt) .
\]
Choosing $\al = 1/3$, we get
\beq\label{eq:pow}
\bal
H( |x|^{-1/3} ) =  \f{1}{\sqrt{3}} sgn(x) |x|^{-1/3} , 
\quad H (sgn(x) |x|^{-1/3}  )  = -\sqrt{3} |x|^{-1/3}.
\eal
\eeq

Recall that $\om$ is odd. We assume that $\om$ is in the Schwartz space. Applying \eqref{eq:colt}, we yield 
\[
\bal
I \teq \int_{\R} &\f{ (u_x (x) - u_x(0))^2} {  |x|^{4/3} } =  \int_{\R} |x|^{2/3}  \lt( H\lt( \f{\om}{x} \rt)  \rt)^2 dx =  \int_{\R} \B\{ |x|^{2/3}  \lt(\f{ \om }{x} \rt)^2 +
2    |x|^{2/3} H \lt(  \f{ \om }{x}  H \lt(       \f{ \om }{x}   \rt)    \rt) \B\} dx . \\
\eal
\]

Since the Hilbert transform is antisymmetric ( Lemma \ref{lem:anti}), we get $H( \om H(\f{\om}{x})) = - \f{1}{\pi}\int_{\R} \f{\om}{x} H(\f{\om}{x}) dx =0$. Using Lemma \ref{lem:com}, we obtain
\[
|x|^{2/3} H \B( \f{\om}{x} H( \f{\om}{x}) \B) = |x|^{2/3} \f{1}{x} H\B( \om H( \f{\om}{x})  \B)
= \sgn(x) |x|^{-1/3} H\B( \om H( \f{\om}{x})  \B).
\]
Thus, applying Lemma \ref{lem:anti}, then \eqref{eq:pow} and $H (\f{\om}{x} )= \f{u_x -u_x(0)}{x}$ in Lemma \ref{lem:com}, we prove 
\[
\bal
I &= \int_{\R} \B\{ \f{ \om^2} { | x|^{4/3} }  - 2    H \lt( sgn ( x)  |x|^{-1/3}  \rt)    \om  H \lt(       \f{ \om }{x}   \rt) \B\}   dx 
 =  \int_{\R} \B\{ \f{ \om^2} { | x|^{4/3} }  +  2  \sqrt{3}      |x|^{-1/3}     \om  H \lt(       \f{ \om }{x}   \rt)  \B\} dx \\
 &    =  \int_{\R} \B\{ \f{ \om^2} { | x|^{4/3} }  +  2  \sqrt{3}      |x|^{-1/3}     \om         \f{ u_x -u_x(0) }{x}    \B\} dx  =
  \int_{\R} \B( \f{ \om^2} { | x|^{4/3} }  +
2  \sqrt{3}   sgn(x) \f{  \om  ( u_x(x) - u_x(0) ) }{|x|^{4/3} } \B) dx.
\eal
\]

To prove the Lemma for general odd $ \om \in L^2( |x|^{-4/3} + |x|^{-2/3})$, or equivalently $\f{\om}{x} \in L^2( |x|^{2/3} + |x|^{4/3})$, we approximate $ \f{\om}{x}$ by the Schwartz function and use the fact that $|x|^{2/3}$ is an $A_2$ weight \cite{duoandikoetxea2001fourier}.
\end{proof}

The  weighted estimates in Lemma \ref{lem:wg23} were established in \cite{Cor06}.

\begin{lem}\label{lem:wg23}
For $f  \in L^2( x^{-4/3} + x^{-2/3})$, we have 
\[
\bal
|| (H f - Hf(0)) x^{-2/3}||_2 &\leq  \cot \f{\pi}{12} || f x^{-2/3}||_2 = (2 + \sqrt{3})|| f x^{-2/3}||_2 , \\
||  Hf x^{-1/3}||_2 &\leq  \cot \f{\pi}{12} || f x^{-1/3}||_2 = (2 + \sqrt{3})|| f x^{-1/3}||_2 .
\eal
\]
\end{lem}

The estimate in the following Lemma  is the Hardy inequality \cite{hardy1952inequalities}.
\begin{lem} \label{lem:hardy}
Assume that $u$ is odd. Then for $p > \f{3}{2}$, we have 
\[
\int_0^{+\infty} \f{ (u(x) -u_x(0) x)^2 }{x^{ 2p}}  \leq \f{4}{(2p-1)^2}  \int_0^{+\infty}    \f{ (u_x(x) - u_x(0))^2 }{x^ {2p-2}} .
\]
\end{lem}

\begin{lem}\label{lem:wg}
Assume that $\om$ is odd and $ \om \in L^2( x^{-4} + x^{-2/3})$. Let $u_x = H \om$. For any $\al, \b, \g \geq 0$, we have 
\[
\bal
|| (u_x - u_x(0) ) ( \al x^{-4} + \b x^{-2})^{1/2} ||_2^2 &= || \om ( \al x^{-4} + \b x^{-2})^{1/2}||_2^2 \\
\B|\B| (u - u_x(0)  x ) \B( \f{\al}{x^6} + \f{\b}{x^4}  +  \f{\g}{x^{10/3}} \B)^{1/2} \B|\B|_2^2 & \leq \B|\B| \om \B( \f{4\al}{25 x^4}  + \f{ 4\b}{9 x^2} \B)^{1/2} \B|\B|_2^2  +  \f{36\g}{49}  || (u_x - u_x(0)) x^{-2/3}||_2^2. 
\eal
\]
\end{lem}

The first identity follows from Lemma \ref{lem:iso}. Applying Lemma \ref{lem:hardy} with $p=3,2, \f{5}{3}$ and then Lemma \ref{lem:iso} to the power $x^{-4}, x^{-2}$ yield the second inequality. The constants $\f{4}{25}, \f{4}{9}, \f{36}{49}$ are determined by $\f{4}{(2p-1)^2}$ with $p=3,2, \f{5}{3}$.

\section{Derivations and estimates in the linear stability analysis}

\subsection{Derivation of \eqref{eq:IBPu}}\label{app:IBPu}
For $ p \in [1,3]$, using integration by parts yields
\[
\bal
 &|| (\td{u} - \f{1}{ 2p  - 1} \td{u}_x x ) x^{- p}||_2^2 
 = \int_{\R_+} \B( \f{1}{ (2 p  - 1)^2 } \f{\td{u}^2_x}{ x^{2p-2} } - \f{2}{ 2p-1} \f{ \td{u} \td{u}_x }{ x^{2p-1}} + \f{ \td{u}^2}{ x^{2p}} \B) dx  \\
 = &  \int_{\R_+} \B( \f{1}{ (2p-1)^2 } \f{\td{u}^2_x}{ x^{2p-2 }}  
+ \f{1}{2p-1} (\pa_x x^{-(2p-1)} ) \td{u}^2 + \f{ \td{u}^2}{ x^{ 2p}} \B) dx 
 = \f{1}{(2p-1)^2}  \int_{\R_+} \f{\td{u}^2_x}{ x^{ 2p -2}}  dx.
 \eal
\]

\subsection{ Estimate of $I_{r1}, I_{r2}, I_{r3}$ }\label{app:chi} 
We construct the cutoff function $\chi$ in \eqref{eq:mode} as follows 
\[
 \chi(x) = \f{2}{\pi}\arctan( ( \f{ x -  l_1  }{ l_2} )^3 ) \one_{x \geq l_1} \ , \quad l_1 = 5 \cdot 10^8, \quad l_2 = 10 l_1.
\]

Recall $I_{r1}, I_{r2}$ in \eqref{eq:decomp}, \eqref{eq:chi_term2} and \eqref{eq:chi_term3}
\beq\label{eq:Ir}
I_{r1}  =  \la \td{u}_x  \chi (\xi_1  \psi_n + \xi_2 \psi_f ) , \th_x \ra, \quad 
I_{r2} = \lam_1 \la \td{u} , \chi \xi_3 \om \vp \ra, \quad 
I_{r3} = - \f{1}{3} \lam_1 \la \td u \chi \xi_3 , D_x \om \vp \ra.
\eeq

Recall from the beginning of Section \ref{sec:linop} that $\bar \om, \bar \th_x, \bar \om_x, \bar \th_{xx}$ have decay rates $ x^{ \al}, x^{ 2 \al}, x^{\al-1}, x^{2\al-1}$, respectively, with $\al$ slightly smaller than $-\f{1}{3}$. Using the formulas of $\xi_i$ in \eqref{eq:ext} and $\vp_f, \vp_n, \psi$ in \eqref{eq:mode}, \eqref{eq:wg}, we obtain the decay rates $\chi (\xi_1  \psi_n + \xi_2 \psi_f ) \sim C_1 x^{-4/3}, \ \chi \xi_3 \vp \sim C_2 x^{-2}$ for sufficiently large $x$, where $C_1, C_2$ are some constants.

Recall $\td{u}_x = u_x -u_x(0)$. Using the Cauchy-Schwarz inequality and Lemmas \ref{lem:wg23}, we obtain
\[
\bal
|I_{r1} | \leq || \td{u}_x x^{-2/3}||_2 
|| \chi (\xi_1  \psi_n + \xi_2 \psi_f )  \th_x||_2
\leq (2 + \sqrt{3}) || \om x^{-2/3}||_2 || \chi (\xi_1  \psi_n + \xi_2 \psi_f )  \th_x||_2 .
\eal
\]

For $I_{r2}$, we first decompose it as follows using $\td u = u - u_x(0) x$
\[
I_{r2} = \lam_1 \la u , \chi \xi_3 \om \vp \ra - u_x(0) \lam_1 \la x , \chi \xi_3 \om \vp \ra 
\teq J_1 + J_2.
\]

Using the Cauchy-Schwarz inequality, Lemma \ref{lem:hardy} with $p = \f{4}{3}$ and Lemma \ref{lem:wg23}, we get 
\[
\bal
| J_1 | &\leq \lam_1  || u x^{-4/3} ||_2 
||  x^{4/3} \chi \xi_3 \om \vp ||_2 
\leq  \f{ 6\lam_1}{5}   || u_x x^{-1/3} ||_2 ||  x^{4/3} \chi \xi_3 \om \vp ||_2  \\
&\leq \f{ 6\lam_1 (2 + \sqrt{3})}{5}   || \om x^{-1/3} ||_2 ||  x^{4/3} \chi \xi_3 \om \vp ||_2 .
\eal
\]

Recall $c_{\om} = u_x(0)$. For $J_2$, using Cauchy-Schwarz inequality, we yield
\[
|J_2| \leq \lam_1 |c_{\om}| \cdot ||  \chi^{1/2} \om \vp^{1/2}||_2  || x \xi_3 \chi^{1/2} \vp^{1/2}||_2.
\]

In the above estimates of $I_{r1}$, if we further bound $|| \chi (\xi_1  \psi_n + \xi_2 \psi_f )  \th_x||_2$ by the weighted $L^2$ norm $|| \th_x \psi^{1/2}||_2$, we obtain a small factor $\rho_2^{-1/3}$ since $\chi$ is supported in $|x| \geq \rho_2$ and the profile has decay. See also the above discussion on the decay rates. Similarly, we get a small factor in the estimates of $J_1, J_2$ from $||  x^{4/3} \chi \xi_3 \om \vp ||_2, ||  x^{4/3} \chi \xi_3 \om \vp ||_2$, respectively.

Using Young's inequality $ab \leq ta^2 + \f{1}{4t} b^2$, we obtain 
\[
\bal
|I_{r1}| + |I_{r2}|
&\leq t_{51} || \om x^{-2/3}||_2^2
+ \f{ (2+\sqrt{3})^2}{4 t_{51}}  || \chi (\xi_1  \psi_n + \xi_2 \psi_f )  \th_x||_2^2
+ t_{52} || \om x^{-1/3} ||_2^2 \\ 
& \quad + \f{1}{4 t_{52}}  ( \f{ 6\lam_1 (2 + \sqrt{3})}{5}  )^2
||  x^{4/3} \chi \xi_3 \om \vp ||_2^2 + t_{53} ||  \chi^{1/2} \om \vp^{1/2}||_2^2 
+  \f{  \lam_1^2 || x \xi_3 \chi^{1/2} \vp^{1/2}||_2^2}{ 4 t_{53}} c_{\om}^2,
\eal
\]
where $t_{51} = 10^{-10} , t_{52} = 10^{-5}, t_{53} =10^{-2}$. We choose these weights $t_{5i}$ so that the terms $ta^2 , \f{1}{4t}b^2$ in Young's inequality are comparable. It follows the estimate \eqref{eq:est_chi}.

Note that replacing $\om$ in $I_{r2}$ in \eqref{eq:Ir} by $-\f{1}{3} D_x \om$, we obtain $I_{r3}$. Therefore, applying the same estimate as that of $I_{r2}$ to $I_{r3}$, we yield
\[
|I_{r3}| \leq \f{2 \lam_1(2 + \sqrt{3})}{5} || \om x^{-1/3}||_2 || x^{4/3} \chi \xi_3 D_x \om \vp ||_2 +  \f{\lam_1}{3}|c_{\om}| \cdot ||  \chi^{1/2} D_x\om \vp^{1/2}||_2  || x \xi_3 \chi^{1/2} \vp^{1/2}||_2.
\]
Using Young's inequality $ab \leq ta^2 + \f{1}{4t} b^2$, we establish 
\[
|I_{r3}| \leq 
 t_{94} ||  x^{4/3} \chi \xi_3 D_x \om \vp ||_2^2   + \f{1}{4 t_{94}}  ( \f{ 2\lam_1 (2 + \sqrt{3})}{5}  )^2  || \om x^{-1/3} ||_2^2 + t_{95} ||  \chi^{1/2} D_x \om \vp^{1/2}||_2^2 
+  \f{ \lam_1^2 || x \xi_3 \chi^{1/2} \vp^{1/2}||_2^2}{ 36 t_{95}} c_{\om}^2 .
\]
where $t_{94} = 10^6, t_{95}  = 10^{-3}$. We choose these weights $t_{94}, t_{95}$ so that the terms $ta^2 , \f{1}{4t}b^2$ in Young's inequality are comparable. It follows \eqref{eq:est_chi3}.

\subsection{ Derivations of the ODE \eqref{eq:ODE} in Section \ref{sec:est_cw}}\label{app:ode}
We use the following functions in the derivations\beq\label{eq:func}
\bal
f_2 & \teq \f{1}{4} \f{\bar u_x}{x} - \f{1}{5} ( \f{3}{4} \bar u_{xx} + \f{1}{4} \f{\bar u_x}{x}) 
- \f{\bar u_x}{x} + \f{\bar u}{x^2}  , \quad  f_3  \teq \lam_1 ( \bar \om - x \bar \om_x) \vp  ,\\
f_4  & \teq \f{3}{5} \f{ \bar u_{\th,x}}{x} + \f{1}{5} 
\lt( \f{3}{5} \bar u_{\th, xx} + \f{2}{5} \f{\bar u_{\th, x}}{x} \rt)  , \quad f_6 \teq \f{\bar u}{x^2}  , \\
 f_7 & \teq (\bar \th_x - x \bar \th_{xx}) \psi, \quad f_8 \teq  \f{3}{4} \bar \om_x + \f{1}{4} \f{\bar \om}{x} , \quad f_9 \teq \f{3}{5} \bar \th_{xx} + \f{2}{5} \f{\bar \th_x}{x} . \\
\eal
\eeq

\subsubsection{Derivations of the ODE for $c_{\om}^2, d_{\th}^2$ and \eqref{eq:ODE}}\label{app:ODE_dth}

Recall $c_{\om} = u_x(0)  = -\f{2}{\pi}  \int_0^{+\infty} \f{\om}{x} dx$ from \eqref{eq:normal}. Multiplying the equation of $\om$ in \eqref{eq:lin00} by $- \f{1}{x}$ and then taking the integral from $0, +\infty$ yield
\[
\bal
 \f{d}{dt} \f{\pi}{2} c_{\om} = &  \f{d}{dt}  \int_0^{+\infty} \f{\om}{-x} dx  = \int_0^{\infty}  \f{ ( \bar c_l x + \bar{u} ) \om_x + u \bar \om_x } {x} dx  - \int_0^{\infty} \f{\th_x} {x}  dx
 \\
& 
+ \int_0^{\infty} \f{ \bar c_{\om} \om + c_{\om} \bar \om}{- x} dx
 - \int_0^{\infty} \f{ F_{\om} + N(\om)}{x} dx \\
 = &\int_0^{\infty}  \f{  \bar{u} \om_x + u \bar \om_x } {x} d x - d_{\th} + \f{\pi}{2} (\bar{c}_{\om} + \bar{u}_x(0)) c_{\om}
- \int_0^{\infty} \f{ F_{\om} + N(\om)}{x} dx ,
 \eal
\]
where we have used the notation $d_{\th}$ in \eqref{eq:dth0} and $\int_0^{\infty} \f{f}{-x} = \f{\pi}{2} H f(0)$ with $f = \om, \bar \om$ in the last identity. Multiplying $c_{\om}$ on both sides, we yield
\beq\label{eq:cw20}
\bal
\f{1}{2} \f{d}{dt} \f{\pi}{2} c^2_{\om} &=   \f{\pi}{2} (\bar{c}_{\om} + \bar{u}_x(0)) c^2_{\om} +
 c_{\om}  \int_0^{\infty}  \f{  \bar{u} \om_x + u \bar \om_x } {x} dx
- c_{\om} d_{\th}   - c_{\om}  \int_0^{\infty} \f{ F_{\om} + N(\om)}{x} dx.
\eal
\eeq
which is exactly \eqref{eq:cw2}.

We derive the ODE for $d_{\th}$ using the $\th$ equation in \eqref{eq:lin00}. Since $ \int_{\R_+} \f{\bar c_l x\th_{xx}}{x} dx = 0 $, we get 
\beq\label{eq:dth1}
\bal
\f{d}{dt} d_{\th} & = 2 \bar c_{\om} \int_{\R_+} \f{\th_x}{x} + 2 c_{\om} \int_{\R_+} \f{ \bar \th_x  }{x} - \int_0^{\infty} \f{ \bar u_x \th_x + \bar u \th_{xx}}{x} dx - \int_0^{\infty} \f{ u \bar {\th}_{xx} +  u_x  \bar \th_{x} } {x} dx  \\
 & + \int_0^{\infty} \f{ F_{\th } + N(\th) }{x} dx \teq I_1 + I_2 + I_3 + I_4 + I_5.
\eal
\eeq

We use the notation $\la \cdot , \cdot \ra$ in \eqref{nota:inner} to simplify the integral.
For $I_3$, using integration by parts, we obtain 
\[
I_3 = -\la ( \bar u \th_x)_x, x^{-1} \ra= \la \bar u \th_x , \pa_x x^{-1} \ra  = -\la \bar u \th_x, x^{-2} \ra.
 \]

Similarly, for $I_4$, we get 
\[
I_4 = -\la u \bar \th_x, x^{-2}\ra.
\]

Recall $c_{\om} = u_x(0)$. We rewrite the above term using the decomposition $ u = \td{u} + u_x(0) x$ \eqref{eq:ut}
\[
I_4 = - \la (\td{u} +  u_x(0) x ) \bar \th_x, x^{-2} \ra = - \la \td{u} \bar \th_x, x^{-2} \ra - c_{\om} \bar d_{\th}.
\]
where we have used the notation $\bar d_{\th}$ defined in \eqref{eq:dth0}. Using \eqref{eq:dth0}, we can simplify $I_1, I_2$ as 
\[
I_1 = 2 \bar c_{\om} d_{\th} , \quad I_2 = 2 c_{\om} \bar d_{\th}.
\]

The $c_{\om}  d_{\th}$ term in $I_2$ and $I_4$ are canceled partially. Using these computations and multiplying both sides of \eqref{eq:dth1} by $d_{\th}$ yields
\beq\label{eq:dth2}
\bal
\f{1}{2}\f{d}{dt} d^2_{\th} &= 2  \bar c_{\om} d_{\th}^2 + c_{\om} \bar d_{\th} d_{\th} -  d_{\th} \int_0^{\infty} \f{ \bar u \th_x }{x^2} dx -  d_{\th}\int_0^{\infty} \f{ \td u \bar {\th}_{x}}{x^2} dx   + d_{\th} \int_0^{\infty} \f{ F_{\th} + N(\th) }{x} dx .
\eal
\eeq
Since $\bar d_{\th}>0$, the term $c_{\om} d_{\th}$ in \eqref{eq:dth2} and \eqref{eq:cw20} have cancellation.

The quadratic parts on the right hand sides in \eqref{eq:cw20}, \eqref{eq:dth2} involve the following terms remained to estimate
\beq\label{eq:ODE_J0}
 J_1 = \la \bar u , \om_x x^{-1} \ra, \  J_2 = \la u , \bar \om_x  x^{-1} \ra , \ J_3 = \la \bar u, \th_x x^{-2} \ra,  \ J_4 =  \la \td{u},  \bar \th_{x} x^{-2} \ra .
 \eeq

We use the idea in Section \ref{sec:model3_nonlocal} to rewrite the integrals of $u$ as the integrals of $\om$ and of $ \td{u} - \f{1}{5} \td{u}_x  x = u_{\D}$ (see \eqref{eq:dth0}). We use the functions $f_i$ defined \eqref{eq:func} to simplify the integrals of $\th_x, \om$. In Appendix \ref{app:ODE_J}, we rewrite $J_i$ as follows 
 \beq\label{eq:ODE_Jdone}
\bal
J_1 + J_2 = \la \om , f_2 \ra +  \la u_{\D} x^{-1}, f_5 \ra, \quad J_3 =  \la \th_x , f_6 \ra,  \quad J_4 = \la u_{\D} x^{-1}, f_9  \ra - \la \om, f_4 \ra.
\eal
\eeq

For some parameters $\lam_2 , \lam_3 > 0$ to be determined, combining \eqref{eq:cw20} and \eqref{eq:dth2}, we yield \[
\bal
\f{1}{2} \f{d}{dt} ( \f{\lam_2 \pi }{2} c_{\om}^2 + \lam_3 d_{\th}^2 )
&=   \f{\pi \lam_2 }{2} (\bar c_{\om} + \bar u_x(0)) c^2_{\om} 
+ \lam_2 c_{\om} (J_1 + J_2)  - \lam_2 c_{\om} d_{\th}  - \lam_2 c_{\om} \la F_{\om } + N(\om), x^{-1} \ra \\
& \quad + 2 \bar c_{\om} \lam_3 d_{\th}^2 + \lam_3 c_{\om} \bar d_{\th} d_{\th}
- \lam_3 d_{\th } J_3 - \lam_3 d_{\th} J_4 + \lam_3 d_{\th} \la F_{\th} + N(\th), x^{-1} \ra .
\eal
\]
Plugging \eqref{eq:ODE_Jdone} in the above ODE, we derive \eqref{eq:ODE}.

\subsubsection{ Derivations of \eqref{eq:ODE_Jdone} in the ODEs}\label{app:ODE_J}

Recall the integrals $J_i$ from \eqref{eq:ODE_J0}. 
We use the idea in Section \ref{sec:model3_nonlocal} to derive the formulas in \eqref{eq:ODE_Jdone}.

Recall $\td{u} = u -u_x(0) x$ from \eqref{eq:ut}. Firstly, we consider $J_2$. Since $\int_0^{\infty}  \bar \om_x dx = 0$, we have 
\[
J_2 = \la u - u_x(0) x , \bar \om_x x^{-1} \ra  = \la \td{u}, \bar \om_x x^{-1} \ra.
\]
We approximate the far field of $\bar \om_x x^{-1}$ by $\f{1}{4} (\f{\bar \om}{x})_x$ and derive
\[
J_2 = \B\la \td{u} , \f{ \bar \om_x}{x} - \f{1}{4}  (\f{\bar \om}{x})_x \B\ra
+ \f{1}{4} \B\la \td{u},  (\f{\bar \om}{x})_x  \B\ra \teq J_{21} + J_{22}. 
\]

Applying integration by parts, \eqref{eq:hil} and \eqref{eq:anti} yields
\[
J_{22} = -\f{1}{4} \la \td{u}_x , \bar \om x^{-1} \ra
= -\f{1}{4} \B\la H \lt( \f{\om}{x} \rt), \bar \om \B\ra
=  \f{1}{4} \B\la \f{\om}{x}, H \bar \om \B\ra = \f{1}{4} \B\la \f{\om}{x}, \bar u_x \B\ra .
\]

In $J_{21}$, the coeifficient 
\[\f{ \bar \om_x}{x} - \f{1}{4}  (\f{\bar \om}{x})_x 
= \f{3}{4} \f{\bar \om_x}{x} + \f{1}{4} \f{\bar \om}{x^2} \] 
decays much faster than $\bar \om_x x^{-1}$ for large $x$. We approximate $\td{u}$ by $\f{1}{5} \td{u}_x x$ 
\[
J_{21} = \B\la \td{u},  \f{3}{4} \f{\bar \om_x}{x} + \f{1}{4} \f{\bar \om}{x^2} \B\ra =
 \B\la \td{u} - \f{1}{5} \td{u}_x x , \f{3}{4} \f{\bar \om_x}{x} + \f{1}{4} \f{\bar \om}{x^2} \B\ra
+  \f{1}{5} \B\la \td{u}_x x ,  \f{3}{4} \f{\bar \om_x}{x} + \f{1}{4} \f{\bar \om}{x^2} \B\ra 
\teq I_1 + I_2.
\]

Using a direct computation and then applying \eqref{eq:hil} and \eqref{eq:anti}, we get
\[
\bal
I_2 &
= \f{1}{5} ( \f{3}{4} \la \td{u}_x , \bar \om_x \ra + \f{1}{4} \la \f{ \td u_x}{x}, \bar \om \ra )
 = \f{1}{5} ( \f{3}{4} \la u_x , \bar \om_x \ra + \f{1}{4} \la H\lt( \f{\om}{x} \rt), \bar \om \ra ) \\
 & = \f{1}{5} ( - \f{3}{4} \la \om, H \bar \om_x \ra- \f{1}{4} \la \f{\om}{x}, H \bar \om \ra ) 
 = -\f{1}{5} \la \om , \f{3}{4} \bar u_{xx} + \f{1}{4} \f{ \bar u_x}{x} \ra ,
\eal
\]
where we have used $\int_0^{\infty} u_x(0) \bar \om_x dx = 0$ in the second identity. Using the notation and function in \eqref{eq:dth0}, \eqref{eq:func}, we can simplify $I_1$ as 
\[
I_1 = \la u_{\D} x^{-1}, f_8 \ra .
\]

Combining the above calculations on $J_{22}, I_1, I_2$, we obtain 
\[
J_2 = I_1 + I_2 + J_{22} = \B\la \om , \f{1}{4} \f{\bar u_x}{x} - \f{1}{5} ( \f{3}{4} \bar u_{xx} + \f{1}{4} \f{\bar u_x}{x}) \B\ra + \B\la \f{ u_{\D}}{x}  , f_8  \B\ra .
\]

For $J_1$ in \eqref{eq:ODE_J0}, using integration by parts, we obtain 
\[
J_1 =  \la \bar u x^{-1},  \om_x \ra 
= - \la \pa_x ( \bar u x^{-1}), \om \ra 
=   \la - \f{\bar u_x}{x} + \f{\bar u}{x^2} , \om \ra.
\]

We can simplify $J_1 + J_2$ using the function $f_2$ in \eqref{eq:func}
\beq\label{eq:ODE_J12}
J_1 + J_2 = \la \om , f_2 \ra +  \la u_{\D} x^{-1}, f_5 \ra.
\eeq

For $J_3$, using $f_6$ defined in \eqref{eq:func}, we get 
\beq\label{eq:ODE_J3}
J_3 =  \la \th_x , \bar u x^{-2} \ra = \la \th_x , f_6 \ra.
\eeq

For $J_4$ in \eqref{eq:ODE_J0}, we use a similar computation to obtain 
\[
\bal
J_4 &= \la \td{u}, \bar \th_{x} x^{-2} \ra 
= \B\la \td{u} ,  \f{ \bar \th_{x}}{ x^2} + \f{3}{5} (\f{\bar \th_x}{x})_x \B\ra
- \f{3}{5} \B\la \td{u},  (\f{\bar \th_x}{x})_x \B\ra  = \la \td{u}, \f{3}{5} \f{ \bar \th_{xx} }{x} + \f{2}{5} \f{\bar \th_x}{x^2} \ra
 - \f{3}{5} \B\la \td{u},  (\f{\bar \th_x}{x})_x \B\ra \\
& = \B\la \td{u} - \f{1}{5} \td{u}_x x , \f{3}{5} \f{ \bar \th_{xx} }{x} + \f{2}{5} \f{\bar \th_x}{x^2} \B\ra + \f{1}{5} \B\la \td{u}_x x ,  \f{3}{5} \f{ \bar \th_{xx} }{x} + \f{2}{5} \f{\bar \th_x}{x^2} \B\ra 
- \f{3}{5} \B\la \td{u},  (\f{\bar \th_x}{x})_x \B\ra 
\teq J_{41} + J_{42} + J_{43} .\\
\eal
 \]

For $J_{41}$, using the notations in \eqref{eq:dth0} and \eqref{eq:func}, we obtain 
\[
J_{41} = \la u_{\D} x^{-1}, f_9  \ra.
\]

For $J_{42}, J_{43}$, using Lemmas \ref{lem:iso} and \ref{lem:anti}, we derive
\[
\bal
J_{42}  &= \f{1}{5} \B\la \td{u}_x, \f{3}{5} \bar \th_{xx} + \f{2}{5} \f{\bar \th_x}{x} \B\ra
= \f{1}{5} \lt( \f{3}{5} \B\la H \om - H\om(0), \bar \th_{xx} \B\ra + \f{2}{5} \B\la H \lt( \f{\om}{x} \rt), \bar \th_x \B\ra  \rt ) \\
&= - \f{1}{5} \lt( \f{3}{5} \la \om, H \bar \th_{xx} \ra + \f{2}{5} \la \f{\om}{x}, H \bar \th_x \ra \rt) , \\
J_{43} & =\f{3}{5} \B\la \td{u}_x , \f{\bar \th_x}{x} \B\ra
= \f{3}{5} \B\la H( \f{\om}{x}) ,\bar \th_x \B\ra 
=  - \f{3}{5} \B\la  \f{\om}{x} , H \bar \th_x \B\ra .
\eal
\]

Combining the above computations and using the notations $\bar u_{\th,x}, f_4$ defined in \eqref{eq:dth0}, \eqref{eq:func}, we yield 
\[
\bal
J_4 &= J_{41} + J_{42} + J_{43} = \la u_{\D} x^{-1}, f_9  \ra
 - \B\la \om , \f{3}{5} \f{\bar u_{\th, x}}{x} 
 + \f{1}{5} \lt( \f{3}{5} \bar u_{\th, xx} + \f{2}{5} \f{\bar u_{\th, x}}{x}  \rt) \B\ra  =  \la u_{\D} x^{-1}, f_9  \ra - \la \om, f_4 \ra.
 \eal
 \]

The formulas in \eqref{eq:ODE_J12}, \eqref{eq:ODE_J3} and the above formula imply \eqref{eq:ODE_Jdone}.

\subsection{Derivations of the commutators in \eqref{eq:commu} }\label{app:commu}

Recall $D_x = x\pa_x$ and the operators in \eqref{eq:linop}. We choose $f= \th_x, g = \om $ in \eqref{eq:linop}. We use the notation $u_x = H \om$. Then $ u = -\Lam^{-1} \om$. 

Firsrly, we compute the commutator related to the transport term. Using  $ (\bar c_l x + \bar u ) \pa_x = (\bar c_l + \f{ \bar u}{x} ) D_x $, for $p = \om$ or $\th_x$, we yield
\beq\label{eq:commu_tran}
\bal
- {[D_x ,  (\bar c_l x + \bar u ) \pa_x ] p} &= - [D_x ,  (\bar c_l  + \f{\bar u}{x} ) D_x ] p=
 - D_x( (\bar c_l + \f{ \bar u}{x} ) D_x p ) + (\bar c_l + \f{ \bar u}{x} ) D_x (D_x p )  \\
&= - D_x( \bar c_l + \f{ \bar u}{x} ) D_x p =  -(\bar u_x - \f{\bar u}{x}) D_x p. 
\eal
\eeq

Next, we compute the velocity corresponding to $D_x \om$. Using Lemma \ref{lem:com}, we get 
\[
H( D_x \om) - H(D_x \om)(0) = x H( \om _x)  = x \pa_x H \om = x u_{xx}. 
\]
Note that $ H(D_x \om)(0)= -\f{1}{\pi} \int_{\R} \om_x dx = 0$. We obtain $D_x u_x = x u_{xx} = H(D_x\om)$. From 
\[
 (x u_x - u)_x =   x u_{xx} = H( D_x \om), \quad (x u_x - u)(0) = 0,
\]
we obtain that $x u_x - x$ is the velocity corresponding to $D_x\om$. Therefore, we have 
\[
H\om = u_x, \quad -\Lam^{-1} \om =u, \quad H (D_x \om)(0)  =0, \quad H( D_x \om) = x u_{xx}, \quad -\Lam^{-1} (D_x \om)  = x u_x - u.
\]

Using these formulas, for $q = \bar \om_x$ or $\bar \th_{xx} $ we obtain 
\beq\label{eq:commu_u}
\bal
& D_x \B( - (- \Lam^{-1} \om - H \om(0) x) q \B)
 - \B(  - (- \Lam^{-1} D_x \om - H D_x \om(0) x) q \B)  \\
 = &D_x ( -  (u -u_x(0) x ) q ) + ( xu_x - u) q \\
 = & - (u -u_x(0) x) D_x q + ( - (x u_x - u_x(0) x) ) q+ ( xu_x - u)   q
= - (u -u_x(0) x ) ( D_x q + q ).
\eal
\eeq

Similarly, we have 
\beq\label{eq:commu_ux}
\bal
&D_x \B( - (H \om - H \om(0) x) q \B)- (  - ( H D_x \om - H D_x \om(0) ) q ) 
=  D_x (- (u_x - u_x(0))q ) + x u_{xx}q   \\
= & - D_x u_x q - (u_x -u_x(0)) D_x q+ x u_{xx} q = - (u_x -u_x(0)) D_x q .
\eal
\eeq

Since $\bar c_{\om} \om, \th_x$ in $\cL_{\om 1}$ \eqref{eq:linop} vanish in the commutator, applying \eqref{eq:commu_tran} with $p = \om$ and \eqref{eq:commu_u} with $q = \bar \om_x$ yields the formula for $[D_x, \cL_{\om1}]$ in \eqref{eq:commu}. Note that 
\[
D_x  (  (2 \bar c_{\om} - \bar u_x) \th_x )
- (2 \bar c_{\om} - \bar u_x) D_x \th_x = - D_x \bar u_x \th_x.
\]

Combining this computation, \eqref{eq:commu_tran} with $p = \th_x$, \eqref{eq:commu_u} with $q = \bar \th_{xx}$ and \eqref{eq:commu_ux} with $q = \bar \th_x$, we derive the formula for $[D_x, \cL_{\th1}]$ in \eqref{eq:commu}.

\subsection{Derivation and computing $C_{opt}$ in Section \ref{sec:cT} } \label{app:Copt}
Recall the inequality \eqref{eq:Copt1}, the functions in \eqref{eq:func2} and the spaces $\Sigma_i$ in \eqref{eq:Copt_space}. We use the argument similar to that in \cite{chen2019finite} to derive and compute $C_{opt}$.

In Section \ref{sec:cT}, we have reduced \eqref{eq:Copt1} to an optimization problem on the finite dimensional space $\Sigma_1 \oplus \Sigma_2 \oplus \Sigma_3$ with $X \in \Sigma_1, Y \in \Sigma_2, Z \in \Sigma_3$. Here, we have a direct sum of spaces since there is no inner product  among $X, Y, Z$. Let $ \{ \ee_1, \ee_2, \ee_3, \ee_4 \}$ be an orthonormal basis (ONB) of $\Sigma_1$ with $\ee_1 = \f{g_1} {||g_1||_2^2}$; 
$ \{ \ee_5, \ee_6, \ee_7 \}$ be that of $\Sigma_2$ with $\ee_2 = \f{g_5} {||g_5||_2^2}$; $ \{ \ee_8, \ee_9 \}$ be that of $\Sigma_3$. Then $\{ \ee_i \}_{i=1}^9$ is an ONB of $\Sigma \teq \Sigma_1 \oplus \Sigma_2 \oplus \Sigma_3$. 

Let $v_i \in \R^9$ be the coordinate of $g_i$ in $\Sigma$ under the basis $\{ \ee_i \}_{i=1}^9$ and $p = (x, y,z) \in \R^4 \times \R^3 \times \R^2$ be that of $X + Y + Z$. The vectors $v_i$ and $p$ are column vectors. By abusing notation, we also use $\la \cdot , \cdot \ra$ to denote the Euclidean inner product in $\R^9$. With these convections, each summand on the left hand side of \eqref{eq:Copt1} is a quadratic form in $p$. For example, we have 
\[
\la X, g_1 \ra \la Y, g_7 \ra = \la p , v_1 \ra \la p, v_7 \ra 
= ( p^T v_1 ) ( v_7^T p) = p^T ( v_1 v_7^T ) p .
\]
Hence, \eqref{eq:Copt1} is equivalent to 
\beq\label{eq:Copt2}
p^T M p \leq C_{opt} p^T D p,
\eeq
where $M$ and $D$ are given by 
\beq\label{eq:Copt_M}
\bal
M &= v_1 v_3^T +  v_1  v_7^T - (\lam_2 - \lam_3 \bar d_{\th}) v_1 v_5^T
+ \lam_2 v_1 v_2^T - \lam_3 v_5 v_6^T + \lam_3 v_5 v_4^T+ \lam_2 v_1 v_8^T - \lam_3 v_5 v_9^T ,\\
D &= Id + s_1 v_1 v_1^T + s_2 v_5 v_5^T .
\eal
\eeq

By definition of $\ee_1, \ee_5$, i.e. $ \ee_1 = \f{g_1}{ ||g_1||_2}, \ee_5 = \f{g_5}{ ||g_5||_2}$, we have $v_1 = ||g_1||_2 E_1, v_5 = || g_5||_2 E_5 $, where $E_i \in \R^9$ is the standard basis of $\R^9$, i.e. the $i$-th coordinate of $E_i$ is $1$ and $0$ otherwise. Therefore, $D$ is a diagonal matrix 
\[
D = \mathrm{ diag }( 1 + s_1 ||g_1||_2^2, \ 1, \ 1, \ 1, \ 1 + s_2 ||g_5||_2^2, \ 1,\ 1,\ 1,\ 1) \in \R^{ 9\times 9}.
\]

Symmetrizing the left hand side of \eqref{eq:Copt2} and using a change of variable $q = D^{1/2} p$, we obtain 
\[
C_{opt} = \lam_{\max}( D^{-1/2} M_{s} D^{-1/2}), \quad M_{s} = \f{1}{2} (M + M^T).
\]

Firstly, $M$ can be written as
\[
\bal
M &=  V_1 V_2^T, \quad\ V_2  = ( v_3, v_7, v_5, v_2, v_6, v_4, v_8, v_9) , \\
V_1 &= ( v_1, v_1, -(\lam_2 - \bar d_{\th} \lam_3) v_1, \lam_2 v_1, -\lam_3 v_5,  \lam_3 v_5, \lam_2 v_1, -\lam_3 v_5 ) . 
\eal
\]

Then $M_{s} = \f{1}{2} ( V_1 V_2^T + V_2 V_1^T) = \f{1}{2} U_1 U_2^T$ with $U_1 = [V_1, V_2], U_2 = [V_2, V_1] \in \R^{ 9 \times 16}$. Using the argument in \cite{chen2019finite}, for any even integer $p \geq 2$,  we obtain
\beq\label{eq:Copt3}
\bal
C_{opt} 
&\leq (\tr |D^{-1/2}M_s D^{-1/2}|^p)^{1/p} 
= 2^{-1} (\tr (D^{-1/2} U_1 U_2^T D^{-1/2} )^p)^{1/p} \\
& = 2^{-1} (\tr (  U_2^T D^{-1} U_1 )^p)^{1/p} .\\
\eal
\eeq
We will explain how to rigorously estimate the bound above in the Supplementary Material \cite{chen2021HLsupp}.

\subsection{Estimate of $\cT$ in Section \ref{sec:est_sum}}\label{app:cT}

For $\lam_2, \lam_3, t_{61} , \kp, r_{c_{\om}} >0$ chosen in \eqref{eq:para3}, Appendix \ref{app:para} and $t_{62}$ determined by these parameters, we define $T_i$ and $s_i$ \beq\label{eq:para1}
\bal
 T_1 &= (-\lam_1 D_{\om}  - A_{\om} \vp^{-1}   - \lam_1 \kp ) \vp  - t_{61} x^{-4} , \quad T_2  = (- D_{\th}  - A_{\th} \psi^{-1}- \kp ) \psi,  \\
T_3 & =  25 t_{61} x^{-4} + t_{62} x^{-4/3}, \quad  s_1 = -\f{\pi}{2}  \lam_2 (\bar c_{\om} + \bar u_x(0)) - r_{c_{\om}} - \f{ \pi\lam_1 e_3 \al_6}{12} - G_c, \\
  s_2 & = - 2 \bar c_{\om} \lam_3 - \kp \lam_3,
\quad 
 \eal
\eeq
We will verify that $T_i >0, s_i > 0$ later. The parameter $r_{c_{\om}}$ is essentially determined by $\kp$. See Appendix \ref{app:para_guide} for the procedure to determine these parameters.
Plugging the above $T_i$ and $s_i$ in \eqref{eq:Copt0}, we can compute the upper bound of $C_{opt}$ in \eqref{eq:Copt0} using \eqref{eq:Copt3} with $p=36$ 
\beq\label{eq:Copt_comp}
C_{opt} \leq 2^{-1} (\tr (  U_2^T D^{-1} U_1 )^p)^{1/p} < 0.9930 < 1,
\eeq
which is verified in \eqref{ver:Copt_comp}, Appendix \ref{app:ver}. Thus from \eqref{eq:Copt0}, we obtain 
\[
 \cT \leq  || \om T_1^{1/2} ||_2^2 + || \th_x T_2^{1/2} ||_2^2  + || \f{ u_{\D}}{x} T_3^{1/2}||_2^2  + s_1 c_{\om}^2 + s_2 d_{\th}^2,
\]
which is exactly \eqref{eq:est_cT}. By definition of $T_1, T_2$, we have 
\[
\bal
 \la ( D_{\th} + A_{\th} \psi^{-1} )  \psi, \th_x^2 \ra + \la T_2, \th_x^2 \ra
 &= - \kp \la   \th_x^2 , \psi\ra  , \\
\la ( \lam_1 D_{\om}  + A_{\om}\vp^{-1} ) \vp, \om^2 \ra + \la T_1 , \om^2 \ra 
&= - \kp \lam_1 \la \om^2, \vp \ra - t_{61} \la \om^2, x^{-4} \ra .
  \eal
\]

Hence, plugging the above estimate on $\cT$ in \eqref{eq:comb1}, we yield 
\beq\label{eq:comb12}
\bal
J &= - \kp || \th_x  \psi^{1/2} ||_2^2
- \kp \lam_1 || \om  \vp^{1/2} ||_2^2  - t_{61} || \om  x^{-2} ||_2^2  + s_1 c_{\om}^2 + s_2 d_{\th}^2  \\
& \quad + || \f{ u_{\D}}{x} T_3^{1/2}||_2^2 - \B( D_u - \f{9}{49} t_{12} - \f{72 \lam_1  }{49}\cdot 10^{-5} \B) || \td{u}_x x^{-2/3}||_2^2  + A(u)
+ G_{c} c_{\om}^2.
\eal
\eeq

It remains to estimate the $u_{\D}$ term. Recall $u_{\D}$ in \eqref{eq:dth0} and $T_3$ in \eqref{eq:para1}. 
A direct calculation yields 
\[
\bal
|| \f{u_{\D}}{x} T_{3}^{1/2}  ||_2^2
= \int_0^{\infty} (\td{u} - \f{1}{5} \td{u}_x x )^2  \cdot 25 t_{61} x^{-6} dx
+ \int_0^{\infty} (\td{u} - \f{1}{5} \td{u}_x x )^2 \cdot  t_{62} x^{-10/3}  dx 
\teq I_1 + I_2
\eal
\]

Using \eqref{eq:IBPu} with $p = 3$ and Lemma \ref{lem:iso}, we get 
\[
I_1 = t_{61} || \td{u}_x x^{-2}||_2^2 = t_{61} || \om x^{-2} ||_2^2.
\]

For $I_2$, using integration by parts and  Lemma \ref{lem:wg} about $\td u$ with $\al = \b = 0$, we get 
\[
\bal
I_2 &= t_{62} \int_0^{\infty} \f{1}{25}  \f{\td{u}^2_x  }{x^{4/3}}
- \f{2}{5} \f{ \td{u} \td{u}_x}{x^{7/3}} + \f {\td{u}^2 }{ x^{10/3}}  dx
= t_{62} \int_0^{\infty} \f{1}{25}  \f{\td{u}^2_x  }{x^{4/3}}
+ \f{1}{5} \td{u}^2 \pa_x x^{-7/3} + \f {\td{u}^2 }{ x^{10/3}}  dx \\
& = t_{62} \int_0^{\infty} \f{1}{25}  \f{\td{u}^2_x  }{x^{4/3}} + (1- \f{7}{15}) \f {\td{u}^2 }{ x^{10/3}}  dx 
\leq t_{62} \int_0^{\infty} \f{\td u^2_x}{x^{4/3}} \lt( \f{1}{25} + \f{8}{15} \cdot \f{36}{49} \rt) dx .
\eal
\]
Combining the estimates of $I_1 , I_2$ yields 
\beq\label{eq:udel}
|| \f{u_{\D}}{x} T_{3}^{1/2}  ||_2^2 \leq t_{61} || \om x^{-2}||_2^2 + 
\B( \f{1}{25} + \f{8}{15} \cdot \f{36}{ 49} \B) t_{62} || \td{u}_x x^{-2/3}||_2^2,
\eeq

We define $t_{62}$ in Appendix \ref{app:para} so that the terms $|| \td u_x x^{-2/3}||_2^2$ in 
\eqref{eq:udel} and \eqref{eq:comb12} are almost canceled. 
We establish \eqref{eq:comb2}, i.e.
\[
J \leq 
- \kp  ||  \th_x  \psi^{1/2} ||_2^2
- \kp \lam_1  || \om  \vp^{1/2} ||_2^2  + (s_1 + G_{c }) c_{\om}^2 + s_2 d_{\th}^2 
- 10^{-6}  || \td{u}_x x^{-2/3}||_2^2  + A(u) .
\]

\section{ Parameters in the estimates}\label{app:para}

\subsection{ Parameters}
Parameters $e_1, e_2, e_3$ introduced in \eqref{eq:ext} are determined by the approximate self-similar profiles 
\beq\label{eq:para_far}
 e_1 = 1.5349, \quad  e_2 = 1.2650, \quad e_3  = 1.3729 .
\eeq

We choose the following parameters for the weights $\psi, \vp$ \eqref{eq:mode},\eqref{eq:wg}
\beq\label{eq:para2}
 \al_1 = 5.3,  \quad \al_2 = 3.3, \quad \al_3 = 0.68, \quad
\al_4 = 12.1, \quad \al_5 = 2.1 , \quad \al_6 = 0.77, 
\eeq
and the following parameters in the linear stability analysis in Section \ref{sec:lin}
\beq\label{eq:para3}
\bal
\lam_1 & = 0.32,  \quad  t_1 = 1.29,   \quad  t_{12} = \f{49}{9} \cdot 0.9 D_u ,  \quad t_2 = 5.5,  \quad t_{22} = 13.5, \quad  t_{31}  = 3.2, \\
   t_{32} &= 0.5,  \quad t_{34} = 2.9,  \quad \tau_1 = 4.7 , \quad t_4  = 3.8  ,    \quad \lam_2   =  2.15, \quad \lam_3 = 0.135 ,  \\
  t_{61} & = 0.16, \quad \kp = 0.03,  \quad r_{c_{\om}} = 0.15. 
   \\
\eal
\eeq

Parameter $\lam_1$ is introduced in \eqref{eq:goal}, \eqref{eq:bad}; $t_{2}, t_{22}$ are introduced in the estimates of $I_n$ \eqref{eq:varcoe_near}, \eqref{eq:est_fast_IR}; $t_{1}, t_{12}$ are introduced in the estimate of $I_f$ \eqref{eq:varcoe_far}, \eqref{eq:Du}; $t_4$ is introduced in the estimate of $I_s$ in \eqref{eq:est_singu}; $( t_{31}, t_{32}), t_{34}, \tau_1$ are introduced in the estimate of $I_{u\om}$ in \eqref{eq:Su1}, \eqref{eq:est_Iuw} and \eqref{eq:Kuw}, respectively; $\lam_2, \lam_3, t_{61}, \kp, r_{c_{\om}}$ are introduced in \eqref{eq:para1} to estimate $\cT$ in \eqref{eq:est_cT}.

The parameter $D_u$ introduced in \eqref{eq:Du}, $t_{62}$ in \eqref{eq:para1} are determined by the above parameters 
\[
\bal
 D_{u} = \f{t_1 \al_3 \lam_1 \al_6}{\sqrt{3}},\quad  t_{62}  = ( D_u - \f{9}{49} t_{12} - \f{72 \lam_1  }{49} \cdot 10^{-5} - 10^{-6}  ) ( \f{1}{25} + \f{8}{15} \cdot \f{36}{49})^{-1}.   \\
        \eal
\]

After we complete the weighted $L^2$ estimate, we choose the following parameters in the weighted $H^1$ estimates and nonlinear stability estimates
\beq\label{eq:para_H1}
\bal
& \kp_2 = 0.024,\quad  t_{71} = 2.8, \quad t_{ 72} = 2,  \quad  t_{81} = 5, \quad  t_{82} = 0.7, \quad t_{91} = 1, \quad t_{92} = 1.2,  \\
&   \quad \g_1 = 0.98, \quad \g_2 = 0.07 , \quad \lam_4 = 0.005, \quad E_* = 2.5 \cdot 10^{-5},\quad a_{H^1} = 0.31.
\eal
\eeq
Parameters $t_{7i}, t_{8i}, t_{9i}$ are introduced in the estimates of $Q_2$ \eqref{eq:Su23}, \eqref{eq:Kuw2}; $\kp_2$ in \eqref{eq:est_H14}; $\g_1, \g_2$ in \eqref{eq:est_H13}; $\lam_4$ in \eqref{eq:para_ineq}. 
Parameter $a_{H1}$ is determined by the above parameters via $A_{\om2}$ \eqref{eq:coer2_cost} and \eqref{eq:est_H15}
\[
a_{H^1} = 0.31 .
\]

\subsection{ Choosing parameters in $\cT$ and determining $\kp$}\label{app:para_guide}

We first choose $r_{ c_{\om}} = \kp \f{\pi}{2 } \lam_2$ with small $\kp =0.001$. The remaining unknown parameters in the linear stability analysis are $\lam_2, \lam_3, t_{61} > 0$. Once $\lam_2, \lam_3, t_{61}$ are chosen, the functions $T_i$ and scalars $s_i$ in \eqref{eq:para1} are determined and then we can compute $C_{opt}$ in \eqref{eq:Copt0} using the argument in Section \eqref{eq:cT} and Appendix \ref{app:Copt}.
 We optimize $\lam_2, \lam_3, t_{61} >0$ subject to the constraints $T_i > 0, s_i > 0$, such that $C_{opt} < 0.98$ and $C_{opt}$ is as small as possible. Then we obtain the approximate values for $\lam_2 ,\lam_3, t_{61}$.

Our goal is to obtain $\kp$ as large as possible. The estimate of $C_{opt}$ depends on all the parameters in \eqref{eq:para2}-\eqref{eq:para3}. We gradually increase $\kp$ until $C_{opt} < 0.98$ is violated. We further refine all the parameters in \eqref{eq:para2}-\eqref{eq:para3} one by one and by modifying them around their approximate values to obtain smaller $C_{opt}$. Then we increase $\kp$ again. Repeating this process several times, we obtain larger $\kp$ and $\kp=0.03$. Finally, we increase $r_{c_{\om}}$ until $C_{opt} < 0.98$ is violated. This allows us to obtain 
a damping term for $c_{\om}^2$ with a larger coefficient in the weighted $L^2$ estimate \eqref{eq:L2}, 
Using this procedure, we determine the parameters in \eqref{eq:para2}, \eqref{eq:para3} and further establish \eqref{eq:L2}.

In our process of determining the parameters, we actually first use the grid point values of the functions and only need to track the constraints, e.g. $T_i>0$, on the grid points instead of every $x \in \R$. 
 After we determine all parameters, we verify the constraints rigorously by using computer-assisted analysis and establish the desired bound $C_{opt} < 0.993<1$ \eqref{eq:Copt_comp}.

\section{Rigorous Verification}\label{app:ver}

This section is a collection of inequalities that will be rigorously verified with the help of computer programs. The methods of computer-assisted verification are introduced and discussed in detail in the Supplementary Material \cite{chen2021HLsupp}. All the numerical computations and quantitative verifications are performed in MATLAB (version 2020a) in double-precision floating-point operations. The MATLAB codes can be found via the link \cite{Matlabcode}. 

\subsection{ Ranges of the parameters}\label{app:range}

Denote by 
\[
 G_1(\lam_1, t_2, t_{22})  \teq  t_2  x^{-4} + \f{t_{22}}{25} x^{-4} + t_2 (\lam_1 \al_5)^2 x^{-2}, 
 \ 
 G_2( t_2, t_{22} ) \teq  \f{1}{4 t_2} (\al_2 x^{-1} + \al_1 x^{-2})^2 + \f{1}{4 t_{22}} (x^{3} \bar \th_{xx}  \psi_n )^2  
\]
the coefficients in \eqref{eq:est_fast}. Applying estimate \eqref{eq:est_fast} on $I_n$, we establish \eqref{eq:goal_para1} with $c=0.01$ if 
\[
\f{1}{\lam_1} G_1(\lam_1, t_2, t_{22}) \vp^{-1} + D_{\om}  \leq -c , \quad 
G_2(t_2, t_{22}) \psi^{-1} + D_{\th} \leq -c, 
\]
where $D_{\om}, D_{\th}$ defined in \eqref{eq:damp} are the coefficients in $D_1, D_2$. To verify the above estimate for $\lam_1 \in [\lam_{1l}, \lam_{1u}] = [0.31, 0.33]$, $t_2 \in  [t_{2l}, t_{2u} ] = [5.0, 5.8], t_{22} \in [t_{22l}, t_{22u}] = [13, 14]$, since $G_1, G_2$ are monotone in $\lam_1, t_2, t_{22}$, it suffices to verify  
\beq\label{eq:ver_f1}
\f{1}{\lam_{1l}}G_1(\lam_{1u}, t_{2u}, t_{22 u}) \vp^{-1} +  D_{\om}  \leq -c , \quad 
G_2(t_{2l}, t_{22 l}) \psi^{-1} + D_{\th} \leq -c. 
\eeq

Similar, in order for $I_f + D_1 + D_2 \leq 
-0.01 ( || \th_x \psi^{1/2}||_2 + \lam_1 || \om \vp^{1/2}||_2^2 )$ with estimate \ref{eq:est_slow}
 on $I_f$ and $\lam_1 \in [\lam_{1l}, \lam_{1u}] = [0.31, 0.33]$, $t_1 \in [t_{1l}, t_{1u}] = [1.2,1.4], t_{12} \in [t_{12l}, t_{12u}] = [0.55, 0.65]$, it suffices to verify 
 \beq\label{eq:ver_f2}
\f{1}{\lam_{1l}} G_3(\lam_{1u}, t_{1u}, t_{12 u}) \vp^{-1} + D_{\om}  \leq -c , \quad 
G_4(t_{1l}, t_{12 l}) \psi^{-1} + D_{\th} \leq - c,
 \eeq
 where 
 \[
 G_3(\lam_1, t_1, t_{12} ) = t_1 \B( \al_3^2 x^{-2} + \f{\al_3 \lam_1 \al_6}{\sqrt{3}} x^{-4/3}  + ( \lam_1 \al_6)^2 x^{-2/3} \B), \ 
 \quad G_4(t_1, t_{12}) = \f{1}{4 t_1} x^{-2/3} + \f{1}{4 t_{12}} (\psi_f \bar \th_{xx} x^{5/3})^2 . 
 \]

In order for $I_s + D_1 + D_2 \leq 
-0.01 ( || \th_x \psi^{1/2}||_2 + \lam_1 || \om \vp^{1/2}||_2^2 )$ with estimate \eqref{eq:est_singu} on $I_s$ and $\lam_1 \in [\lam_{1l}, \lam_{1u}] = [0.31, 0.33]$, $t_4 \in [t_{4l}, t_{4u}] = [3.5, 4.0]$, it suffices to verify 
\beq\label{eq:ver_f3}
\f{1}{\lam_{1l}} G_5(t_{4u}) +  D_{\om}  \leq -c ,
\quad G_6(\lambda_{1u},t_{4l})  + D_{\th} \leq -c, 
\eeq
where
\[G_5(t_4) = t_{4} x^{-3} \vp^{-1}, \quad G_6(\lambda_1, t_4) = \f{ (\lam_1 \al_4)^2}{4 t_4 } x^{-5} \psi^{-1}.\]

\begin{remark}
We do not actually use the above estimates. Yet, they provide a useful guideline to determine the parameters  $t_{ij}$ in the estimates.
\end{remark}

\subsection{ Inequalities on the approximate steady state}

To establish the nonlinear estimates in Sections \ref{sec:lin} and \ref{sec:non}, we have used several inequalities on the approximate steady state and the parameters defined in Appendix \ref{app:para}. These inequalities are summarized below.

In \eqref{eq:damp}, we derive the damping terms in the weighted $L^2$ estimate with coefficients $D_{\th}, D_{\om}$. These coefficients are negative uniformly. That is, for some $c>0$, we have 
\beq\label{ver:damp}
D_{\th}, D_{\om} \leq - c < 0.
\eeq

Recall that we choose the weights $T_i$ and $s_i$ defined in \eqref{eq:para1} and apply the argument in Section \ref{sec:cT} to obtain the sharp estimate of the $\cT$ term defined in \eqref{eq:cT}. This estimate requires that the weights are nonnegative, i.e. 
\beq\label{ver:para1}
\bal
 T_1 &= (-\lam_1 D_{\om}  - A_{\om} \vp^{-1} - \lam_1 \kp ) \vp  - t_{61} x^{-4}  > 0, \\
T_2  &= (- D_{\th}  - A_{\th} \psi^{-1}- \kp ) \psi > 0,  \\
T_3 & =  25 t_{61} x^{-4} + t_{62} x^{-4/3} > 0.
\quad 
 \eal
\eeq
and 
\beq\label{ver:para2}
\bal
 s_1 &= -\f{\pi}{2} \lam_2 (\bar c_{\om} + \bar u_x(0)) - r_{c_{\om}} - \f{ \pi\lam_1 e_3 \al_6}{12} - G_{c } > 0,   \\
s_2  & = - 2 \bar c_{\om} \lam_3 - \kp \lam_3 > 0. 
\quad 
 \eal
\eeq

Using the above $T_i, s_i$ and the argument in Section \ref{app:Copt}, we establish the following estimate for the constant $C_{opt}$ in \eqref{eq:Copt0}
\beq\label{ver:Copt_comp}
C_{opt} \leq 2^{-1} (\tr (  U_2^T D^{-1} U_1 )^p)^{1/p} < 0.9930 < 1.
\eeq
The fact that $C_{opt} < 1$ implies \eqref{eq:est_cT}. 

In the weighted $H^1$ estimates, we have used 
\beq\label{ver:ver_H1_1}
  (x^2 \bar u_{xx} \psi)_x \leq  0.02 \psi 
\eeq
in \eqref{eq:ver_H1_1} to establish \eqref{eq:est_H113}. We have also used 
\beq\label{ver:est_H14}
\bal
D_{\th} + A_{\th} \psi^{-1}  - (\bar u_x - \f{ \bar u}{x}) + B_{\th} \psi^{-1}  & \leq - \kp_2, \\
\lam_1 D_{\om} + A_{\om} \vp^{-1} -\lam_1  (\bar u_x - \f{\bar u}{x}) + B_{\om} \vp^{-1} & \leq -\kp_2 \lam_1 .
\eal
\eeq
and
\beq\label{ver:est_H15}
\quad || A_{\om 2} \vp^{-1} ||_{\infty}\leq a_{H^1} ,
\eeq
originated from \eqref{eq:ver_H1_1} and \eqref{eq:est_H15} to establish \eqref{eq:H1}.

\bibliographystyle{abbrv}

 \bibliography{selfsimilar}

\begin{thebibliography}{10}

\bibitem{castelli2018rigorous}
R.~Castelli, M.~Gameiro, and J.-P. Lessard.
\newblock Rigorous numerics for ill-posed {PDE}s: periodic orbits in the
  {B}oussinesq equation.
\newblock {\em Archive for Rational Mechanics and Analysis}, 228(1):129--157,
  2018.

\bibitem{Cor10}
A.~Castro and D.~C\'ordoba.
\newblock Infinite energy solutions of the surface quasi-geostrophic equation.
\newblock {\em Advances in Mathematics}, 225(4):1820--1829, 2010.

\bibitem{castro2020global}
A.~Castro, D.~C{\'o}rdoba, and J.~G{\'o}mez-Serrano.
\newblock Global smooth solutions for the inviscid sqg equation.
\newblock 2020.

\bibitem{castro2014remarks}
A.~Castro, D.~C{\'o}rdoba, J.~G{\'o}mez-Serrano, and A.~M. Zamora.
\newblock Remarks on geometric properties of {SQG} sharp fronts and
  $\alpha$-patches.
\newblock {\em arXiv preprint arXiv:1401.5376}, 2014.

\bibitem{chen2020slightly}
J.~Chen.
\newblock On the slightly perturbed {D}e {G}regorio model on $ {S}^1$.
\newblock {\em To appear in ARMA, arXiv preprint arXiv:2010.12700}, 2020.

\bibitem{chen2020singularity}
J.~Chen.
\newblock Singularity formation and global well-posedness for the generalized
  {Constantin--Lax--Majda} equation with dissipation.
\newblock {\em Nonlinearity}, 33(5):2502, 2020.

\bibitem{chen2019finite2}
J.~Chen and T.~Y. Hou.
\newblock Finite time blowup of {2D} {Boussinesq} and {3D} {Euler} equations
  with ${C}^{1,\alpha}$ velocity and boundary.
\newblock {\em Communications in Mathematical Physics}, 383(3):1559--1667,
  2021.

\bibitem{Matlabcode}
J.~Chen, T.~Y. Hou, and D.~Huang.
\newblock Matlab codes for computer-aided proofs in the paper ``asymptotically
  self-similar blowup of the {H}ou--{L}uo model for the 3{D} {E}uler
  equations''.
\newblock
  \url{https://www.dropbox.com/sh/qjs6p6d9n3uiq8r/AABCDI-rZeVuTmBxGQuLJbUva?dl=0}.

\bibitem{chen2021HLsupp}
J.~Chen, T.~Y. Hou, and D.~Huang.
\newblock Supplementary materials for the paper for the paper ``asymptotically
  self-similar blowup of the {H}ou--{L}uo model for the 3{D} {E}uler
  equations''.
\newblock {\em arXiv preprint arXiv:2106.05422}.

\bibitem{chen2019finite}
J.~Chen, T.~Y. Hou, and D.~Huang.
\newblock On the finite time blowup of the {D}e {G}regorio model for the 3{D}
  {E}uler equations.
\newblock {\em Communications on Pure and Applied Mathematics},
  74(6):1282--1350, 2021.

\bibitem{choi2014on}
K.~Choi, T.~Hou, A.~Kiselev, G.~Luo, V.~Sverak, and Y.~Yao.
\newblock On the finite-time blowup of a 1{D} model for the 3{D} axisymmetric
  {E}uler equations.
\newblock {\em CPAM}, 70(11):2218--2243, 2017.

\bibitem{choi2015finite}
K.~Choi, A.~Kiselev, and Y.~Yao.
\newblock Finite time blow up for a 1{D} model of 2{D} {B}oussinesq system.
\newblock {\em Comm. Math. Phys.}, 334(3):1667--1679, 2015.

\bibitem{constantin2007euler}
P.~Constantin.
\newblock On the {E}uler equations of incompressible fluids.
\newblock {\em Bulletin of the American Mathematical Society}, 44(4):603--621,
  2007.

\bibitem{constantin1996geometric}
P.~Constantin, C.~Fefferman, and A.~Majda.
\newblock Geometric constraints on potentially singular solutions for the 3{D}
  {E}uler equations.
\newblock {\em Communications in Partial Differential Equations}, 21(3-4),
  1996.

\bibitem{CLM85}
P.~Constantin, P.~D. Lax, and A.~Majda.
\newblock A simple one‐dimensional model for the three‐dimensional
  vorticity equation.
\newblock {\em CPAM}, 38(6):715--724, 1985.

\bibitem{cordoba2005formation}
A.~C{\'o}rdoba, D.~C{\'o}rdoba, and M.~Fontelos.
\newblock Formation of singularities for a transport equation with nonlocal
  velocity.
\newblock {\em Annals of Mathematics}, pages 1377--1389, 2005.

\bibitem{Cor06}
A.~C\'ordoba, D.~C\'ordoba, and M.~A. Fontelos.
\newblock Integral inequalities for the hilbert transform applied to a nonlocal
  transport equation.
\newblock {\em Journal de Math\'ematiques Pures et Appliqu\'ees},
  88(6):529--540, 2006.

\bibitem{cordoba2017note}
D.~C{\'o}rdoba, J.~G{\'o}mez-Serrano, and A.~Zlato{\v{s}}.
\newblock A note on stability shifting for the {M}uskat problem, {II}: From
  stable to unstable and back to stable.
\newblock {\em Analysis \& PDE}, 10(2):367--378, 2017.

\bibitem{DG90}
S.~De~Gregorio.
\newblock On a one-dimensional model for the three-dimensional vorticity
  equation.
\newblock {\em Journal of Statistical Physics}, 59(5-6):1251--1263, 1990.

\bibitem{DG96}
S.~De~Gregorio.
\newblock A partial differential equation arising in a {1D} model for the {3D}
  vorticity equation.
\newblock {\em Mathematical Methods in the Applied Sciences},
  19(15):1233--1255, 1996.

\bibitem{deng2005geometric}
J.~Deng, T.~Hou, and X.~Yu.
\newblock Geometric properties and nonblowup of 3{D} incompressible {E}uler
  flow.
\newblock {\em Communications in Partial Difference Equations},
  30(1-2):225--243, 2005.

\bibitem{duoandikoetxea2001fourier}
J.~Duoandikoetxea and J.~D. Zuazo.
\newblock {\em Fourier analysis}, volume~29.
\newblock American Mathematical Soc., 2001.

\bibitem{elgindi2019finite}
T.~M. Elgindi.
\newblock Finite-time singularity formation for ${C}^{1,\alpha}$ solutions to
  the incompressible {Euler} equations on $\mathbb{R}^3$.
\newblock {\em arXiv:1904.04795}, 2019.

\bibitem{elgindi2019stability}
T.~M. Elgindi, T.-E. Ghoul, and N.~Masmoudi.
\newblock On the stability of self-similar blow-up for ${C}^{1,\alpha}$
  solutions to the incompressible {Euler} equations on $\mathbb{R}^3$.
\newblock {\em arXiv preprint arXiv:1910.14071}, 2019.

\bibitem{Elg19}
T.~M. Elgindi, T.-e. Ghoul, and N.~Masmoudi.
\newblock Stable self-similar blow-up for a family of nonlocal transport
  equations.
\newblock {\em Analysis \& PDE}, 14(3):891--908, 2021.

\bibitem{elgindi2018finite}
T.~M. Elgindi and I.-J. Jeong.
\newblock Finite-time singularity formation for strong solutions to the
  axi-symmetric 3 d euler equations.
\newblock {\em Annals of PDE}, 5(2):1--51, 2019.

\bibitem{Elg17}
T.~M. Elgindi and I.-J. Jeong.
\newblock On the effects of advection and vortex stretching.
\newblock {\em Archive for Rational Mechanics and Analysis}, Oct 2019.

\bibitem{elgindi2017finite}
T.~M. Elgindi and I.-J. Jeong.
\newblock Finite-time singularity formation for strong solutions to the
  boussinesq system.
\newblock {\em Annals of PDE}, 6:1--50, 2020.

\bibitem{enciso2018convexity}
A.~Enciso, J.~G{\'o}mez-Serrano, and B.~Vergara.
\newblock Convexity of {W}hitham's highest cusped wave.
\newblock {\em arXiv preprint arXiv:1810.10935}, 2018.

\bibitem{gabai2003homotopy}
D.~Gabai, G.~R. Meyerhoff, and N.~Thurston.
\newblock Homotopy hyperbolic 3-manifolds are hyperbolic.
\newblock {\em Annals of Mathematics}, 157(2):335--431, 2003.

\bibitem{gibbon2008three}
J.~Gibbon.
\newblock The three-dimensional {E}uler equations: Where do we stand?
\newblock {\em Physica D: Nonlinear Phenomena}, 237(14):1894--1904, 2008.

\bibitem{gomez2019computer}
J.~G{\'o}mez-Serrano.
\newblock Computer-assisted proofs in pde: a survey.
\newblock {\em SeMA Journal}, 76(3):459--484, 2019.

\bibitem{gomez2014turning}
J.~G{\'o}mez-Serrano and R.~Granero-Belinch{\'o}n.
\newblock On turning waves for the inhomogeneous {M}uskat problem: a
  computer-assisted proof.
\newblock {\em Nonlinearity}, 27(6):1471, 2014.

\bibitem{hales2005proof}
T.~Hales.
\newblock A proof of the {K}epler conjecture.
\newblock {\em Annals of Mathematics}, pages 1065--1185, 2005.

\bibitem{hardy1952inequalities}
G.~Hardy, J.~Littlewood, and G.~P{\'o}lya.
\newblock {\em Inequalities}.
\newblock Cambridge university press, 1952.

\bibitem{hoang2018blowup}
V.~Hoang, B.~Orcan-Ekmekci, M.~Radosz, and H.~Yang.
\newblock Blowup with vorticity control for a 2d model of the boussinesq
  equations.
\newblock {\em Journal of Differential Equations}, 264(12):7328--7356, 2018.

\bibitem{hoang2017cusp}
V.~Hoang and M.~Radosz.
\newblock Cusp formation for a nonlocal evolution equation.
\newblock {\em Archive for Rational Mechanics and Analysis}, 224(3):1021--1036,
  2017.

\bibitem{hoang2020singular}
V.~Hoang and M.~Radosz.
\newblock Singular solutions for nonlocal systems of evolution equations with
  vorticity stretching.
\newblock {\em SIAM Journal on Mathematical Analysis}, 52(2):2158--2178, 2020.

\bibitem{hou2009blow}
T.~Hou.
\newblock Blow-up or no blow-up? a unified computational and analytic approach
  to 3{D} incompressible {E}uler and {N}avier-{S}tokes equations.
\newblock {\em Acta Numerica}, 18(1):277--346, 2009.

\bibitem{hou2006dynamic}
T.~Hou and R.~Li.
\newblock Dynamic depletion of vortex stretching and non-blowup of the 3{D}
  incompressible {E}uler equations.
\newblock {\em Journal of Nonlinear Science}, 16(6):639--664, 2006.

\bibitem{hou2015self}
T.~Hou and P.~Liu.
\newblock Self-similar singularity of a 1{D} model for the 3{D} axisymmetric
  {E}uler equations.
\newblock {\em Research in the Mathematical Sciences}, 2(1):1--26, 2015.

\bibitem{Sve19}
H.~Jia, S.~Stewart, and V.~Sverak.
\newblock On the de gregorio modification of the constantin--lax--majda model.
\newblock {\em Archive for Rational Mechanics and Analysis}, 231(2):1269--1304,
  2019.

\bibitem{kenig2006global}
C.~E. Kenig and F.~Merle.
\newblock Global well-posedness, scattering and blow-up for the
  energy-critical, focusing, non-linear {S}chr{\"o}dinger equation in the
  radial case.
\newblock {\em Inventiones mathematicae}, 166(3):645--675, 2006.

\bibitem{kiselev2010regularity}
A.~Kiselev.
\newblock Regularity and blow up for active scalars.
\newblock {\em Mathematical Modelling of Natural Phenomena}, 5(04):225--255,
  2010.

\bibitem{kiselev2018}
A.~Kiselev.
\newblock Small scales and singularity formation in fluid dynamics.
\newblock In {\em Proceedings of the International Congress of Mathematicians},
  volume~3, 2018.

\bibitem{kryz2016}
A.~Kiselev, L.~Ryzhik, Y.~Yao, and A.~Zlatos.
\newblock Finite time singularity for the modified {SQG} patch equation.
\newblock {\em Ann. Math.}, 184:909--948, 2016.

\bibitem{kiselev2013small}
A.~Kiselev and V.~Sverak.
\newblock Small scale creation for solutions of the incompressible two
  dimensional {E}uler equation.
\newblock {\em Annals of Mathematics}, 180:1205--1220, 2014.

\bibitem{kiselev2018finite}
A.~Kiselev and C.~Tan.
\newblock Finite time blow up in the hyperbolic boussinesq system.
\newblock {\em Adv. Math.}, 325:34--55, 2018.

\bibitem{landman1988rate}
M.~Landman, G.~Papanicolaou, C.~Sulem, and P.~Sulem.
\newblock Rate of blowup for solutions of the nonlinear {S}chr{\"o}dinger
  equation at critical dimension.
\newblock {\em Physical Review A}, 38(8):3837, 1988.

\bibitem{lanford2017computer}
O.~E. Lanford.
\newblock A computer-assisted proof of the {F}eigenbaum conjectures.
\newblock In {\em Universality in Chaos}, pages 245--252. Routledge, 2017.

\bibitem{li2008blow}
D.~Li and J.~Rodrigo.
\newblock Blow-up of solutions for a 1d transport equation with nonlocal
  velocity and supercritical dissipation.
\newblock {\em Advances in Mathematics}, 217(6):2563--2568, 2008.

\bibitem{liu2017spatial}
P.~Liu.
\newblock {\em Spatial Profiles in the Singular Solutions of the 3D Euler
  Equations and Simplified Models}.
\newblock PhD thesis, California Institute of Technology, 2017.
\newblock \url{https://resolver.caltech.edu/CaltechTHESIS:09092016-000915850}.

\bibitem{luo2013potentially-2}
G.~Luo and T.~Hou.
\newblock Toward the finite-time blowup of the 3{D} incompressible {E}uler
  equations: a numerical investigation.
\newblock {\em SIAM Multiscale Modeling and Simulation}, 12(4):1722--1776,
  2014.

\bibitem{luo2014potentially}
G.~Luo and T.~Y. Hou.
\newblock Potentially singular solutions of the 3d axisymmetric euler
  equations.
\newblock {\em Proceedings of the National Academy of Sciences},
  111(36):12968--12973, 2014.

\bibitem{lushnikov2020collapse}
P.~M. Lushnikov, D.~A. Silantyev, and M.~Siegel.
\newblock Collapse vs. blow up and global existence in the generalized
  constantin-lax-majda equation.
\newblock {\em arXiv preprint arXiv:2010.01201}, 2020.

\bibitem{majda2002vorticity}
A.~Majda and A.~Bertozzi.
\newblock {\em Vorticity and incompressible flow}, volume~27.
\newblock Cambridge University Press, 2002.

\bibitem{martel2014blow}
Y.~Martel, F.~Merle, and P.~Rapha{\"e}l.
\newblock Blow up for the critical generalized {K}orteweg--de {V}ries equation.
  {I}: Dynamics near the soliton.
\newblock {\em Acta Mathematica}, 212(1):59--140, 2014.

\bibitem{mclaughlin1986focusing}
D.~McLaughlin, G.~Papanicolaou, C.~Sulem, and P.~Sulem.
\newblock Focusing singularity of the cubic {S}chr{\"o}dinger equation.
\newblock {\em Physical Review A}, 34(2):1200, 1986.

\bibitem{merle2005blow}
F.~Merle and P.~Raphael.
\newblock The blow-up dynamic and upper bound on the blow-up rate for critical
  nonlinear {S}chr{\"o}dinger equation.
\newblock {\em Annals of mathematics}, pages 157--222, 2005.

\bibitem{merle1997stability}
F.~Merle and H.~Zaag.
\newblock Stability of the blow-up profile for equations of the type $u_t=
  \delta u + | u|^{ p- 1} u $.
\newblock {\em Duke Math. J}, 86(1):143--195, 1997.

\bibitem{merle2015stability}
F.~Merle and H.~Zaag.
\newblock On the stability of the notion of non-characteristic point and
  blow-up profile for semilinear wave equations.
\newblock {\em Communications in Mathematical Physics}, 333(3):1529--1562,
  2015.

\bibitem{moore2009introduction}
R.~E. Moore, R.~B. Kearfott, and M.~J. Cloud.
\newblock {\em Introduction to interval analysis}, volume 110.
\newblock Siam, 2009.

\bibitem{OSW08}
H.~Okamoto, T.~Sakajo, and M.~Wunsch.
\newblock On a generalization of the constantin–lax–majda equation.
\newblock {\em Nonlinearity}, 21(10):2447--2461, 2008.

\bibitem{rump2010verification}
S.~M. Rump.
\newblock Verification methods: Rigorous results using floating-point
  arithmetic.
\newblock {\em Acta Numerica}, 19:287--449, 2010.

\bibitem{schochet1986explicit}
S.~Schochet.
\newblock Explicit solutions of the viscous model vorticity equation.
\newblock {\em Communications on pure and applied mathematics}, 39(4):531--537,
  1986.

\bibitem{schwartz2009obtuse}
R.~E. Schwartz.
\newblock Obtuse triangular billiards {II}: One hundred degrees worth of
  periodic trajectories.
\newblock {\em Experimental Mathematics}, 18(2):137--171, 2009.

\bibitem{silvestre2014transport}
L.~Silvestre and V.~Vicol.
\newblock On a transport equation with nonlocal drift.
\newblock {\em transactions of the American Mathematical Society},
  368(9):6159--6188, 2016.

\end{thebibliography}

\end{document}


\begin{abstract}
In the main paper \cite{chen2021HL_supp}, we have established the weighted $L^2$ and $H^1$ estimates of the linearized equations of the Hou-Luo (HL) model and obtained the linear stability of the approximate steady state. In the supplementary material, we close the nonlinear energy estimates and establish the nonlinear stability by further estimating the nonlinear terms and the error terms. We will also estimate the nonlinear terms and error terms in the energy estimates for the convergence to the self-similar solution and the uniqueness of the self-similar profile. Moreover, we will discuss the construction of an approximate steady state with an extremely small residual and the procedure to verify several inequalities on the approximate steady state related to the energy estimates with computer assistance. The organization of the Supplementary material will be outlined in Section \ref{sec:organ}.

\end{abstract}

\maketitle

\setcounter{section}{8}

\section{Notations and outlines of the Supplementary material}

In Section \ref{sec:recall}, we first recall the functions, parameters and notations introduced in the main paper \cite{chen2021HL_supp}.
In Section \ref{sec:organ}, we outline the organization in the Supplementary material.

\subsection{Functions and notations}\label{sec:recall}

Recall the following weights in the stability analysis 
\beq\label{eq:wg}
\psi = \psi_n + \psi_f, \quad \vp = \vp_s + \vp_n + \vp_f.
\eeq
where $\psi_n, \psi_f, \vp_s, \vp_n, \vp_f$ are given by 
\beq\label{eq:mode}
\bal
\psi_n &= \f{1}{ \bar \th_x + \f{1}{5} x \bar \th_{xx} + \chi \xi_1 } ( \al_1 x^{-4} + \al_2 x^{-3}),
\quad \psi_f =  \f{1}{\bar \th_x + \f{3}{7} x \bar \th_{xx} + \chi \xi_2} \al_3 x^{-4/3} , \\
\vp_s & = \al_4 x^{-4} , \quad \vp_n  = \al_5 (  \al_1 x^{-3} + \al_2 x^{-2} ) , \quad \vp_f = \al_6 x^{-2/3} , \\
\xi_1 &= e_1 x^{-2/3}  - ( \bar \th_x + \f{1}{5} x \bar \th_{xx}  ) , \quad \xi_2 = 
 e_2 x^{-2/3}  - ( \bar \th_x + \f{3}{7} x \bar \th_{xx}  ) , \quad \xi_3 = -\f{e_3}{3} x^{-4/3} - \bar \om_x.
\eal
\eeq

Due to the symmetry, we restrict the inner product and $L^2$ norm to $\R_+$
\[
\la f, g \ra  \teq \int_0^{\infty} fg dx  , \quad  || f||_2^2  = \int_0^{\infty} f^2 dx .
\]

The nonlinear terms and error terms in the equations linearized around the approximate steady state are given by 
\beq\label{eq:NF}
\bal
F_{\th} &= ( 2 \bar{c}_{\om} - \bar u_x) \bar\th_x - (\bar{c}_l x + \bar{u}) \cdot \bar{\th}_{xx} , \quad
F_{\om} =\bar{\th}_x + \bar{c}_{\om} \bar{\om } -  ( \bar{c}_l x + \bar{u} ) \cdot \bar{\om}_x , \\
N(\th) & = (2 c_{\om} - u_x) \th_x - u \th_{xx},  \quad N(\om)  = c_{\om} \om -  u \om_x , \\
 \eal
 \eeq
 and the associated nonlinear terms and error terms in the weighted $L^2$ estimates are \beq\label{eq:l2h1}
\bal
N_1 & = \la  N(\th)  , \th_x  \psi  \ra, \quad F_1 =\la F_{\th} , \th_x \psi \ra , \quad
N_2 = \la N(\om) , \om  \vp  \ra, \quad F_2 =  \la F_{\om} , \om \vp  \ra . 
\eal
\eeq

The remaining terms in the ODEs, weighted $L^2$ and $H^1$ estimates to be estimated are 
\beq\label{eq:R}
\bal
\cR_{ODE} &\teq  - \lam_2 c_{\om} \la F_{\om} + N(\om), x^{-1} \ra + \lam_3 d_{\th} \la F_{\th} + N(\th), x^{-1} \ra ,\\
\cR_{L^2}  &\teq N_1 + F_1 + \lam_1 N_2 + \lam_1 F_2 + \cR_{ODE} , \\
\cR_{H^1} &= \la D_x N(\th), D_x \th_x \psi \ra + \lam_1 \la D_x N(\om), D_x \om \vp \ra
+ \la D_x F_{\th}, D_x \th_x \psi \ra + \lam_1 \la D_x F_{\om}, D_x \om \vp \ra,
\eal
\eeq
where $d_{\th}, c_{\om}$ are given below. We introduce $\bar d_{\th}, u_{\th}$ and $\td u$ as follows
\beq\label{eq:cwdth}
\bal
c_{\om} &= H \om(0) = u_x(0), \  d_{\th} \teq \la \th_x, x^{-1} \ra,   \
  \bar d_{\th} \teq \la \bar \th_x , x^{-1} \ra,
\quad \bar u_{\th ,x} \teq  H \bar \th_x , \quad \td u = u - u_x(0) x.
\eal
\eeq

The energy in the weighted $L^2$ estimate is 
\beq\label{energy:l2}
E_1^2( \th_x, \om) = \la \th_x^2, \psi \ra + \lam_1 \la \om^2, \psi \ra 
+ \lam_2 \f{\pi}{2} \cdot \f{4}{\pi^2} \la \om, x^{-1}\ra^2  + \lam_3 \la \th_x, x^{-1} \ra^2, 
\eeq
and the energy in the nonlinear estimates is given by 
\beq\label{energy:h1}
\bal
&E^2(\th_x, \om) =  E_1^2(\th_x, \om) +  \lam_4 ( ||  D_x \th_x \psi^{1/2}||_2^2 + \lam_1 || D_x \om \vp^{1/2}||_2^2  )  \\
=&||  \th_x \psi^{1/2}||_2^2 + \lam_1 ||  \om \psi^{1/2}||_2^2
+ \lam_2 \f{\pi}{2} c_{\om}^2 + \lam_3 d_{\th}^2 +  \lam_4 ( ||  D_x \th_x \psi^{1/2}||_2^2 + \lam_1 || D_x \om \vp^{1/2}||_2^2  ) .
\eal
\eeq

Recall that we choose the following parameters for the weights $\psi, \vp$ \eqref{eq:mode},\eqref{eq:wg}
\beq\label{eq:para2}
 \al_1 = 5.3,  \quad \al_2 = 3.3, \quad \al_3 = 0.68, \quad
\al_4 = 12.1, \quad \al_5 = 2.1 , \quad \al_6 = 0.77, 
\eeq
the following parameters in the linear stability analysis 
\beq\label{eq:para3}
\bal
\lam_1 & = 0.32,  \quad  t_1 = 1.29,   \quad  t_{12} = \f{49}{9} \cdot 0.9 D_u ,  \quad t_2 = 5.5,  \quad t_{22} = 13.5, \quad  t_{31}  = 3.2, \\
   t_{32} &= 0.5,  \quad t_{34} = 2.9,  \quad \tau_1 = 4.7 , \quad t_4  = 3.8  ,    \quad \lam_2   =  2.15, \quad \lam_3 = 0.135 ,  \\
  t_{61} & = 0.16, \quad \kp = 0.03,  \quad r_{c_{\om}} = 0.15,
   \\
   \eal
   \eeq
   and the following parameters in the weighted $H^1$ estimates and nonlinear stability estimates
\beq\label{eq:para_H1}
\bal
& \kp_2 = 0.024,\quad  t_{71} = 2.8, \quad t_{ 72} = 2,  \quad  t_{81} = 5, \quad  t_{82} = 0.7, \quad t_{91} = 1, \quad t_{92} = 1.2,  \\
&   \quad \g_1 = 0.98, \quad \g_2 = 0.07 , \quad \lam_4 = 0.005, \quad E_* = 2.5 \cdot 10^{-5}.
\eal
\eeq
  We refer to Appendix C in the main paper \cite{chen2021HL_supp} for more discussions on these parameters.

\subsection{ Organization and outlines of the Supplementary material}\label{sec:organ}

\subsubsection{Outline of Section \ref{sec:non_est}}

In Section 5.2 of the main paper \cite{chen2021HL_supp}, we have established the following estimates 
\[
\bal
\f{1}{2}\f{d}{dt} E^2(\th_x, \om) \leq & -\kp_2 E^2(\th_x, \om) +  \cR_{L^2} + \lam_4 \cR_{H^1} \teq  -\kp_2 E^2(\th_x, \om) + \cR,
 \eal
 \]
 where $\cR_{L^2}, \cR_{H^1}$ defined in \eqref{eq:R} are the nonlinear terms and error terms. Formally, using embedding inequalities, the remainder term can be bounded by 
 \beq\label{eq:outl_boost}
\cR \leq C_1 E^3(\th_x, \om) + \e_1  E(\th_x, \om)
 \eeq
 for some $C_1, \e_1 > 0$ depending on the approximate steady state. The first term on the right hand side bounds the nonlinear terms related to $N(\th)$ and $N(\om)$, and $\e_1$ controls the weighted norm of the residual error $F_{\th}, F_{\om}$. If $C_1 \e_1 < \kp_2^2 / 4 $, using a bootstrap argument yields $E(\th_x, \om)(t) < (\e_1 / C_1)^{1/2}$, which implies the nonlinear stability of the profile. In Section \ref{sec:non_est}, we work out these constants $C_1, \e_1$ and prove 
 \[
\cR \leq 36 E^3(\th_x, \om) +  5.5 \cdot 10^{-7} E(\th_x, \om). 
 \]

Similarly, we estimate the nonlinear terms and error terms in the energy estimates for the convergence to the self-similar solution and the uniqueness of the self-similar profile.

We remark that these estimates are \textit{standard} and do not exploit special structure of the system. Yet, we have to work out the constants in these estimates to verify that $C_1 \e_1 < \kp_2^2 / 4$. The factor $\e_1$ plays a role as a small parameter in these estimates, and we construct an approximate steady state with a small enough residual to close the bootstrap argument. 

\subsubsection{Outline of Section \ref{sec:Computation} }

We discuss the construction of the approximate self-similar profile in Section \ref{sec:Computation}, which is one of the most important steps in the current framework. As we discuss in Section 4 of the main paper \cite{chen2021HL_supp}, an essential step to obtain an approximate steady state with a sufficiently samll residual is to represent the approximate steady state as follows
\beq\label{eq:intro_decomp}
\om(t,x) = \om_b(x) + \om_p(t,x), \quad v(t,x) = v_b(x) + v_p(t,x),
\eeq
where $v = \th_x$, $\om_b, v_b$ are some explicit functions that capture the decay of $\om, \th_x$
 and $\om_p$ and $v_p$ are represented by piecewise polynomials with compact supports.
Once $\om_b, v_b$ are determined, the computation of $\om_p, v_p$ follows an approach similar to that in \cite{chen2019finite_supp}. 

\subsubsection{A new difficulty in constructing an approximate steady state with extremely small residual}\label{sec:new_diff}
Based on the decomposition \eqref{eq:intro_decomp} and the above discussions, the essential new difficult is to construct a smooth and explicit function $\om_b$ which captures the far field behavior of $\om $. Moreover, we need to compute the Hilbert transform of this smooth explicit function with high accuracy. Yet, since the solution has tail behavior as in \eqref{eq:tail_behavior} with a power close to $\bar c_{\om} / \bar c_l$, which is neither an integer nor a simple rational number, an explicit function with all the desired properties is not known to us.
A key observation is 
\[
H (\sgn(x) |x|^{-a}) = - \cot \f{a \pi}{2} |x|^{-a}
\]
for any $a \in (0,1)$, which is proved in the Appendix of the arXiv version of the main paper 
\cite{chen2021HL_supp}. Then we modify it to obtain a smooth function $F_a = \f{x^5}{1 + |x|^{5+a}}\in C^{10, a}$, and choose $\om_b = C F_a$ for some $C$ to approximate $\om$ in the far field. Since $F_a$ is close to $ \sgn(x) |x|^{-a}$ for large $x$, we expect that $H F_a$ can be approximated by $H (\sgn(x) |x|^{-a})$ for large $x$. In fact, we have the following result 
\begin{lem}[Rough statement]\label{lem:ub_rough}
For $k = 0,1,2$ and $x \geq x_0$, we have 
\[
\pa_x^k H F_a  = \pa_x^k ( - \cot \f{a \pi}{2} |x|^{-a}) -  \f{ (-1)^k  (k+1)! \cdot H ( f_a x^2)(0)  }{x^{k+2}} + E_k(x),
\]
where the error $E_k(x)$ and $f_a$ satisfy 
\[
f_a = F_a - \sgn(x) |x|^{-a}, \quad E_k(x ) \leq C(a, k, x_0) x^{ - k - 4}.
\]
\end{lem}
A more precise statement is given in Corollary \ref{cor:ub}, whose proof and related results are one of the main focuses of Section \ref{sec:Fa}. The above result shows that for large $x$, the leading order term of the Hilbert tranform of $F_a$ is given by that of $\sgn(x) |x|^{-a}$ as expected. Moreover, for large $x$, with the next order term with decay $x^{-k-2}$, we can compute $\pa_x^k H F_a$ with a very small error that has a faster decay rate $x^{-k-4}$.

For $x$ that is not too large, computing $H F_a(x)$ numerically is much more accurate because we can use the explicit formula of $F_a$ and its power series expansion, see Lemma \ref{lem:Fa_power}. This computation is discussed in details in Section \ref{sec:compute_ub}. To achieve the high accuracy, the computation of $H F_a$ is a combination of analytic estimates and numerical computation, and thus involves several technical estimates and decompositions of the integrals.

\subsubsection{Outline of Section \ref{sec:Numerical-verification}}

In addition to the numerical computation of the Hilbert transform of $\om_b$ for $x$ not too large, we discuss how to perform rigorous computer-assisted verifications in Section \ref{sec:Numerical-verification}. Many methods follow similar approaches that we have developed in \cite{chen2019finite_supp}, e.g. the estimate of the weighted norm of the error $F_{\om}, F_{\th}$ with singular weight. To reduce the round off error in computing the Hilbert transform, e.g.
\[
u = \f{1}{\pi} \int_0^{\inf} \om(y) \log \B| \f{x-y}{x+y} \B| dy 
\]
for odd $\om$, we use different methods to compute the integral in different regimes where $x/y$ is close to $1$, $x/y$ or $y/x$ is sufficiently small, or otherwise.

\subsubsection{Outline of Section \ref{sec:GQ}}

In Section \ref{sec:GQ}, we follow the approach in \cite{chen2019finite_supp} to establish the error estimate of the integrals computed numerically. Due to different cases and the partition of the integrals in the computation of the Hilbert transform discussed in Section \ref{sec:Numerical-verification}, the error estimates of the integral that we need to estimate in Section \ref{sec:GQ} are more technical. We remark that the errors of these integrals are very small since we use sufficiently high order scheme, and the error estimates can be obtained using basic numerical analysis.

\subsubsection{Outline of Section \ref{sec:Fa}}

In addition to establishing Corollary \ref{cor:ub} similar to Lemma \ref{lem:ub_rough} and the related results,
we estimate high order derivatives of $F_a$, e.g. $\pa_x^{20} F_a$, which are used to estimate the error of the high order quadrature rule. Though $F_a$ is given explicitly, for large $k$, e.g. $k \geq 10$, the explicit formula of $\pa_x^k F_a$ is too lengthy and thus cannot be used to estimate $\pa_x^k F_a$ directly. We use the structure of $ \pa_x^k F_a$ and the power series expansion of $F_a(x)$ for large $x$ to establish such estimates. We remark that the estimates of the power series expansion is \textit{standard}.

\subsubsection{Outline of Section \ref{sec:beyond_LB}} 
In Section \ref{sec:beyond_LB}, we derive the leading order behavior of the approximate steady state for sufficiently large $x$, e.g. $x \geq 10^{10}$. This property enables us to verify several inequalities on the approximate steady state and estimate the residual errors $F_{\om}, F_{\th}$ for sufficiently large $x$, e.g. $x \geq 10^{10}$. We remark that these estimates are \textit{standard}.

In summary, the numerical computation and the rigorous verification in Sections \ref{sec:Computation}-\ref{sec:beyond_LB} follow the routine steps outlined in \cite{chen2019finite_supp}. 
See also Section 3.13 in the main paper \cite{chen2021HL_supp} on the ideas about these steps. The essential new difficulty and additional new steps are due to $F_a$ discussed in Section \ref{sec:new_diff}. Its related computation and estimates are established in Sections \ref{sec:compute_ub} and Section \ref{sec:Fa}. Moreover, to reduce the round off error in computing the Hilbert transform, we need to use more structure of the kernel and its leading order approximation in certain region, instead of treating it as a general explicit function. These lead to a substantial expansion of the supplementary material compared to that of \cite{chen2019finite_supp}.

\section{Nonlinear estimates}\label{sec:non_est}

In Section \ref{sec:Linf}, we establish the pointwise estimate on $u_x, \th_x$ using the weighted $H^1$ energy $E$, which will be used to estimate the nonlinear terms.
In Section \ref{sec:NF}, we estimate the remaining nonlinear terms and error terms in the weighted $L^2$ and $H^1$ estimates, which are used in the main paper \cite{chen2021HL_supp} to close the nonlinear energy estimate. In Section \ref{supp:conv}, we establish the estimate for the convergence to self-similar solution, which follows the argument in \cite{chen2019finite_supp}.
In Section \ref{supp:uni}, we estimate the difference between two solutions and establish the uniqueness of the self-similar profiles. These estimates complete the estimates in Sections 5.2-5.4 in the main paper \cite{chen2021HL_supp}.

In Section \ref{supp:profile}, we prove several properties for the self-similar profiles, which are used to establish the H\"older regularity of the blowup solution in Section 6 of the main paper \cite{chen2021HL_supp}.

\subsection{Estimate of $u_x, \th_x$}\label{sec:Linf}
\begin{lem}\label{lem:comp}
The weights $\vp, \psi$ defined in \eqref{eq:wg} satisfy
\[
\vp_u \teq  \al_4 x^{-4} + p_1 x^{-2}  \leq \vp,  \quad 
 x^{-4/3} + \f{3}{p_2} x^{-2/3} + p_2 x^{-2}   \leq p_3 \vp, \quad 
 x\psi_x < -\f{1}{2} \psi ,  \quad b_3 x^{-2/3}  \leq\psi, 
 \]
where $\al_4$ is defined in \eqref{eq:para2} and 
\[ 
p_1  = \f{7}{4}  \B( \al_5 \al_1  \cdot  (\f{4}{3} \al_6  )^{3/4} \B)^{4/7} \cdot (1-10^{-7})+ \al_5 \al_2 , \quad  p_2 = 6.5, \quad p_3 =0.87 , \quad b_3 = 0.5.
\]
\end{lem}

\begin{proof}
We will verify the last three inequalities directly using the computer-assisted approach outlined in Section \ref{sec:Numerical-verification}. Recall $\vp$ in \eqref{eq:wg}. For the first one, using Young's inequality, we get
\[
\bal
\al_5 \al_1 x^{-3} + \al_6 x^{-2/3} &= 
\al_5 \al_1 x^{-3} + \f{3}{4} \cdot (\f{4}{3} \al_6 x^{-2/3}) 
\geq \f{7}{4} \B( \al_5 \al_1 x^{-3} \cdot  (\f{4}{3} \al_6 x^{-2/3})^{3/4} \B)^{4/7}\\
& = \f{7}{4}  \B( \al_5 \al_1  \cdot  (\f{4}{3} \al_6  )^{3/4} \B)^{4/7} x^{-2} 
> (p_1 - \al_5 \al_2) x^{-2},
\eal
\]
where we have used the definition of $p_1$ to obtained the last inequality. It follows  
\[
\vp =  \al_4 x^{-4} +  \al_5 (  \al_1 x^{-3} + \al_2 x^{-2} ) + \al_6 x^{-2/3} 
\geq  \al_4 x^{-4} + p_1 x^{-2} = \vp_u.
\]
\end{proof}

\begin{lem}\label{lem:ux_inf}
Let $u_x = H \om$. For $x \geq0$, we have 
\[
\bal
&|u_x - u_x(0)| \leq C_u(x) E, \quad  | u - u_x(0) x| \leq  C_u(x) x E, \\
 &C_u(x) \teq  \min( ( \lam_1 \sqrt{\lam_4})^{-1/2} ( x \vp_u)^{-1/2} , \ \B( \f{ (2+\sqrt{3} )  \sqrt{p_3} }{ \lam_1 \sqrt{\al_6 \lam_4 } } \B)^{1/2}   )  < 20,
\eal
\]
where $\lam_1, \al_6, \lam_4 $ and $E$ are defined in \eqref{eq:para2}, \eqref{eq:para_H1} and \eqref{energy:h1}, $p_3$ in Lemma \ref{lem:comp}. 
\end{lem}

\begin{proof}

Denote $\td{u}_x = u_x - u_x(0)$. We first establish the estimate on $\td{u}_x$.

For $x \geq0$, since $(x \vp_u)^{-1}$ is increasing, we yield
\beq\label{eq:ux_inf1}
\bal
\td{u}_x^2(y) &= 2 \int_0^y  \td{u}_x u_{xx} dx
= 2 \int_0^y  \td{u}_x \vp_u^{1/2}  u_{xx} x \vp_u^{1/2} ( x \vp_u(x))^{-1}  dx \\
&\leq 2 (y \vp_u(y) )^{-1} || \td{u}_x \vp_u^{1/2}||_2  || u_{xx} x \vp_u^{1/2}||_2.
\eal
\eeq

Since $H(D_x \om)(0) = - \f{1}{\pi} \int_{\R} \om_x dx = 0$, using Lemma \ref{lem:com}, we get
\beq\label{eq:xuxx}
H(D_x \om)  = H( D_x \om)(0) + x H( \om_x) = 0 + x \pa_x H \om = x u_{xx}. 
\eeq

Using Lemma \ref{lem:comp} and the weighted estimate in Lemma \ref{lem:wg}, we obtain 
\[
\td{u}_x^2 \leq 2 (x \vp_u)^{-1} || \om \vp_u^{1/2}||_2 || D_x \om \vp_u^{1/2}||_2
\leq 2 (x \vp_u)^{-1} || \om \vp^{1/2}||_2 || D_x \om \vp^{1/2}||_2.
\]

Recall $E$ in \eqref{energy:h1}. Using Young's inequality $2ab \leq a^2 + b^2$, we get 
\beq\label{eq:ux_inf2}
E^2 \geq \lam_1 || \om \vp^{1/2}||_2^2 + \lam_1 \lam_4 || D_x \om \vp^{1/2}||_2^2
\geq 2 \lam_1 \lam_4^{1/2} || \om \vp^{1/2}||_2|| D_x \om \vp^{1/2}||_2.
\eeq
Combining the above two estimates obtains the first bound in the minimum in Lemma \ref{lem:ux_inf}.

Using the formulas in \eqref{eq:ux_inf1} and \eqref{eq:xuxx}, we obtain 
\beq\label{eq:u_Linf1}
\td{u}_x^2 \leq 2 || \td{u}_x x^{-2/3}||_2 || u_{xx} x^{2/3}||_2
= 2 || \td{u}_x x^{-2/3}||_2 || H(D_x\om) x^{-1/3}||_2 .
\eeq
Then using Lemma \ref{lem:wg13} and $x^{-2/3} 
\leq \al_6^{-1} \vp$ due to \eqref{eq:mode}, \eqref{eq:wg}, we have another estimate for $\td{u}_x$
\beq\label{eq:u_Linf2}
 \td u_x^2 \leq 2 (2 + \sqrt{3} )|| \td{u}_x x^{-2/3}||_2 || D_x\om x^{-1/3}||_2
\leq \f{ 2 (2 + \sqrt{3} ) }{ \sqrt{\al_6}}  || \td{u}_x x^{-2/3}||_2 || D_x\om \vp^{1/2} ||_2.
\eeq

Using Lemmas \ref{lem:rota} and Young's inequality $2 \sqrt 3 ab \leq \f{1}{p_2} (\sqrt{3}a)^2 + p_2 b^2$, we obtain 
\[
|| \td{u}_x x^{-2/3}||_2^2 
\leq || \om x^{-2/3}||_2^2 + 2 \sqrt{3} \la |\td{u}_x|, \om x^{-4/3} \ra
\leq || \om x^{-2/3}||_2^2 + \f{3}{p_2} || \td{u}_x x^{-1}||_2^2 + p_2 || \om x^{-1/3}||_2^2. 
\]

Using the weighted estimate in Lemma \ref{lem:wg} and then the second inequality in Lemma \ref{lem:comp}, we further obtain 
\beq\label{eq:ux_wg23}
|| \td{u}_x x^{-2/3}||_2^2 \leq \la \om^2 , x^{-4/3} + \f{3}{p_2} x^{-2} + p_2 x^{-2/3} \ra
\leq p_3 \la \om^2, \vp \ra .
\eeq

Combining the above estimates and \eqref{eq:ux_inf2}, we establish 
\[
\td{u}_x^2 
\leq \f{ 2 (2 + \sqrt{3} ) \sqrt{ p_3}  }{ \sqrt{\al_6}}  || \om \vp^{1/2} ||_2 || D_x\om \vp^{1/2} ||_2 \leq
\f{  (2+\sqrt{3} ) \sqrt{p_3} }{ \lam_1 \sqrt{\al_6 \lam_4 } } E^2,
\]
which implies the second bound in the minimum in Lemma \ref{lem:ux_inf}.

Applying the first estimate in Lemma \ref{lem:u} for $\td u_x$ and the fact that $C_u(x)$ is increasing, we get 
\[
| u - u_x(0) x| = |\int_0^x ( u_x(y) - u_x(0)) dy| \leq  \int_0^x 1 d y 
\cdot C_u(x) E = x C_u(x) E.
\]
\end{proof}

Recall $\psi$ in \eqref{eq:wg}. For $\th_x$, we have the following estimates. 

\begin{lem}\label{lem:thx_inf}
For $x \geq 0$, we have 
\[
|\th_x| \leq  ( \lam_4^{-1/2} + 0.5)^{1/2} E \cdot (x \psi)^{-1/2}< \sqrt{15} E \cdot (x \psi)^{-1/2},
\]
where $\lam_4 = 0.005$ given in \eqref{eq:para3} satisfies $\lam_4^{1/2} + 0.5 = 200^{1/2} + 0.5 < 15$.
\end{lem}

\begin{proof}
Consider $y \geq 0$. A direct calculation yields 
\[
\bal
\th_x^2 y \psi(y) =  \int_0^y (\th_x^2 x\psi)_x dx 
= 2\int_0^y \th_x \th_{xx} x\psi dx + \int_0^x \th_x^2 ( \psi + x \psi_x) dx 
\teq I_1 + I_2.
\eal
\]

Recall $E^2$ in \eqref{energy:h1}. Using the inequality on $\psi$ in Lemma \ref{lem:comp}, we get
\[
I_2 \leq 0.5  || \th_x \psi^{1/2}||^2_2 \leq 0.5 E^2.
\]

For $I_1$, using Young's inequality $ a ^2 + \lam_4 b^2 \geq 2 \lam_4^{1/2} ab$, we get
\[
E^2 \geq || \th_x \psi^{1/2}||_2^2 + \lam_4 || D_x \th_x \psi^{1/2}||_2^2
\geq 2 \sqrt{ \lam_4} \int_0^{\infty} |\th_x D_x \th_x \psi | dx  \geq \sqrt{\lam_4} I_2.
\]
Combining the estimates of $I_1, I_2$ and using $\lam_4^{-1/2} + 0.5 < 15$ conclude the proof.
\end{proof}

\subsection{ Estimates of the nonlinear and the error terms}\label{sec:NF}

Recall $N(\th), N(\om), F_{\th}, F_{\om} $ in \eqref{eq:NF}, $N_1, N_2, F_1, F_2$ in \eqref{eq:l2h1}, and the remaining terms $\cR_{ODE}$, $\cR_{L^2}$ and $\cR_{H^1}$ in \eqref{eq:R}. 
Denote 
\beq\label{eq:est_non0}
 \td \psi = \f{ x \psi_x}{\psi}, \quad \td{\vp} = \f{x \vp_x}{\vp}.
\eeq

Using integration by parts, we get 
\beq\label{eq:est_N21}
\bal
N_2 &= \la N(\om), \om \vp \ra 
= \la -u\om_x + c_{\om} \om , \om \vp \ra \\
&= \la \f{1}{2} (u \vp)_x + c_{\om} \vp , \om^2 \ra 
=  \la \f{1}{2} \f{(u \vp)_x}{\vp} + c_{\om}  , \om^2 \vp\ra 
\teq \la Z_2(x), \om^2 \vp \ra. 
\eal
\eeq

Recall $c_{\om} = u_x(0)$ and $\td{u} = u - u_x(0) x$ from \eqref{eq:cwdth}. Using $u_x = \td{u}_x + c_{\om}, u = \td{u} + c_{\om} x$, we derive
\beq\label{eq:est_non1}
\f{ (uf)_x}{f} = \f{ u_x f + uf_x}{f} = u_x + \f{u}{x} \f{xf_x}{f} 
= \td{u}_x +  \f{\td{u}} {x} \f{x f_x}{f} + c_{\om} (1 + \f{x f_x}{f}).
\eeq

From the definition of $E$ in \eqref{energy:h1}, we have 
\beq\label{eq:est_non2}
|c_{\om}| \leq k_c  E,  \quad k_c = (\f{\pi \lam_2}{2})^{-1/2} .
\eeq

Applying \eqref{eq:est_non1} with $f= \vp$, the notation $\td \vp$ in \eqref{eq:est_non0}, Lemma \ref{lem:ux_inf} and \eqref{eq:est_non2}, we get 
\beq\label{eq:est_N22}
Z_2(x)= \f{1}{2} (\td{u}_x + \f{\td u}{x} \td \vp  + c_{\om} (1 + \td \vp))+ c_{\om}
\leq  z_2 E , \quad z_2 \teq  \f{1}{2} || C_u +  C_u |\td \vp|  + k_c |3 + \td{\vp}| \ ||_{L^{\infty}} .
\eeq

Similarly, for $N_1$ in \eqref{eq:l2h1}, we obtain
\beq\label{eq:est_N11}
N_1 = \la - u \th_{xx} + (2 c_{\om} - u_x) \th_x, \th_x \psi \ra
= \la \f{1}{2} \f{ (u \psi)_x}{  \psi} + (2 c_{\om} - u_x), \th_x^2 \psi \ra
\teq \la Z_1(x), \th_x^2 \psi \ra .
\eeq

Applying \eqref{eq:est_non1} with $f= \psi$ and using $2c_{\om} - u_x = c_{\om} - \td{u}_x$, we rewrite $Z_1$ as follows 
\[
Z_1(x) = \f{1}{2} ( \td u_x + \f{\td{u}}{x} \td \psi + c_{\om} (1 + \td{\psi} ) ) + c_{\om} - \td u_x =  \f{1}{2}  (- \td u_x + \f{\td{u}}{x} \td \psi + c_{\om} (3 + \td{\psi} ) ) .
\]
Using Lemma \ref{lem:ux_inf} and \eqref{eq:est_non2}, we obtain 
\beq\label{eq:est_N12}
\bal
Z_1(x) \leq z_1 E, \quad z_1 \teq  \f{1}{2} || C_u + C_u | \td \psi | + k_c |3 + \td{\psi}| \ ||_{L^{\infty}} .
\eal
\eeq
Therefore, we establish
\beq\label{eq:non_N12}
|N_1|  + \lam_1  |N_2|  \leq z_1  \la \th_x^2, \psi \ra E 
+ z_2 \lam_1 \la \om^2, \vp \ra E.
\eeq

Recall the nonlinear term in $\cR_{H^1}$ \eqref{eq:R}. We compute
\beq\label{eq:est_non3}
D_x (u f_x ) = D_x ( \f{u}{x} D_x  f) = \f{u}{x} D_x D_x f + D_x( \f{u}{x}) D_x f
= u \pa_x D_x f + (u_x - \f{u}{x}) D_x f.
\eeq
Recall $N(\om)$ in \eqref{eq:NF}. Applying this formula with $f = \om$ and integration by parts, we obtain 
\[
\bal
N_4 &\teq \la D_x N(\om), D_x \om \vp \ra 
= \la - u \pa_x D_x \om- (u_x - \f{u}{x}) D_x \om + c_{\om} D_x \om, D_x \om \vp \ra \\
& = \la  \f{1}{2}\f{( u\vp)_x}{ \vp} - (u_x - \f{u}{x} ) + c_{\om}, ( D_x\om)^2 \vp \ra
\teq \la Z_4(x), (D_x \om)^2 \vp \ra.
\eal
\]

Using \eqref{eq:est_non1} with $f= \psi$ and $u_x - \f{u}{x} = \td{u}_x - \f{\td u}{x}$, we derive
\[
Z_4(x) = \f{1}{2} (\td{u}_x + \f{\td u}{x} \td \vp  + c_{\om} (1 + \td \vp)) - (\td u_x - \f{\td u}{x}) + c_{\om}
= \f{1}{2} ( -\td u_x+ \f{\td u}{x} (2 + \td \vp) + c_{\om} (3 + \td{\vp})).
\]

Applying Lemma \ref{lem:ux_inf} and \eqref{eq:est_non2}, we yield 
\[
Z_4(x) \leq \f{1}{2} || C_u + C_u |\td \vp + 2| + k_c |3 + \td{\vp}| \ ||_{L^{\infty}} E \teq z_4 E,
\]
which implies 
\beq\label{eq:non_N4}
N_4 \leq z_4 \la ( D_x \om)^2, \vp \ra E.
\eeq

Recall $N(\th)$ in \eqref{eq:NF}. Similarly, applying \eqref{eq:est_non3} with $f = \th_x$, we have 
\[
D_x N(\th) = D_x( -u \th_{xx} + (2 c_{\om} - u_x) \th_x)
= -u \pa_x (D_x \th_x) -(u_x - \f{u}{x}) D_x \th_x +(2 c_{\om} - u_x ) D_x \th_x - (D_x u_x) \th_x.
\]

Applying integration by parts to the transport term, we get 
\[
\la -u \pa_x (D_x \th_x), D_x \th_x \psi \ra = \la \f{1}{2} (u\psi)_x,  (D_x \th_x)^2 \ra
=\la \f{1}{2} \f{ (u\psi)_x}{\psi}  D_x \th_x, D_x \th_x \psi \ra.
\]
Hence, we derive 
\[
\bal
N_3 &\teq \la \pa_x N(\th), D_x \th_x \psi \ra
= \la\f{1}{2} \f{ (u\psi)_x}{ \psi}  D_x \th_x
-(u_x - \f{u}{x}) D_x \th_x +(2 c_{\om} - u_x ) D_x \th_x - (D_x u_x) \th_x,  D_x \th_x \psi \ra \\
& = \B\la \f{1}{2} \f{ (u\psi)_x}{ \psi} -(u_x - \f{u}{x}) + 2c_{\om} - u_x, (D_x \th_x)^2 \psi \B\ra
- \la D_x u_x \th_x , D_x \th_x \psi \ra \teq N_{31} + N_{32}
\eal
\]

Applying \eqref{eq:est_non1} with $f = \psi$, $u_x = \td u_x + c_{\om}$ and $\f{u}{x} = \f{ \td u}{x} + c_{\om}$, we get 
\[
\bal
Z_3(x) & \teq \f{1}{2} \f{ (u\psi)_x}{ \psi} + \B(-(u_x - \f{u}{x}) + 2c_{\om} - u_x \B) \\
&= \f{1}{2} (\td{u}_x + \f{\td u}{x} \td \psi  + c_{\om} (1 + \td \psi)) 
+ \B( -( \td u_x - \f{ \td u}{x} ) +  c_{\om}  -\td u_x  \B) 
 =  \f{ 1}{2} ( -3 \td u_x + \f{\td u}{x} (2 +  \td \psi) + c_{\om} (3 + \td \psi)  ).
\eal
\]

Applying Lemma \ref{lem:ux_inf} and \eqref{eq:est_non2} yields 
\[
Z_3 \leq z_3 E, \quad z_3 \teq  \f{1}{2} || 3 C_u +  C_u |2 + \td \psi| + k_c |3 + \td \psi| \ ||_{L^{\infty}} .
\]
It follows 
\[
N_{31} \leq || Z_3 ||_{\infty} || D_x \th_x \psi^{1/2}||_2^2 
\leq z_3 E || D_x \th_x \psi^{1/2}||_2^2 .
\]

Next, we estimate $N_{32}$. Denote 
\[
S_3  = \al_4 x^{-4} + \al_2 \al_5 x^{-2} + \al_6 (2 + \sqrt{3})^{-2} x^{-2/3}
\]
where $\al_4, \al_2, \al_5, \al_6$ are the same as that in \eqref{eq:mode} and defined in \eqref{eq:para2}. Recall $ D_x u_x = H(D_x \om)$ from \eqref{eq:xuxx}. Using the weighted estimates in Lemma \ref{lem:wg} and Lemma \ref{lem:wg13}, we have 
\beq\label{eq:ux_wgs3}
\bal
|| D_x u_x S_3^{1/2}||_2^2 &= \la ( H D_x\om)^2, ex^{-4} + \al_2 \al_5 x^{-2} \ra 
+ \la ( H D_x\om)^2,  \al_6 (2 + \sqrt{3})^{-2} x^{-2/3} \ra \\
&\leq 
\la ( D_x \om)^2 , ex^{-4} + \al_2 \al_5 x^{-2} + \al_6 x^{-2/3}\ra
\leq \la (D_x \om)^2, \vp \ra = || D_x \om \vp^{1/2}||_2^2,
\eal
\eeq
where we have used the definition of $\vp$ in \eqref{eq:mode}, \eqref{eq:wg} to obtain the last inequality. Applying the Cauchy-Schwarz inequality, Lemma \ref{lem:thx_inf} and the above estimate, we obtain 
\[
\bal
N_{32} &= -\la D_x u_x S_3^{1/2} \th_x S_3^{-1/2} \psi^{1/2}, D_x \th_x \psi^{1/2} \ra  
\leq || \th_x S_3^{-1/2} \psi^{1/2} ||_{L^{\infty}}  || D_x u_x S_3^{1/2}||_2  || D_x \th_x \psi^{1/2} ||_2 \\
&\leq  || (S_3 x)^{-1/2}    ||_{L^{\infty}} \sqrt{15} E \cdot || D_x \om \vp^{1/2}||_2  || D_x \th_x \psi^{1/2} ||_2 .
\eal
\]

Using Young's inequality $ a b \leq \f{1}{2 \lam_1^{1/2}} (a^2 + \lam_1 b^2)$  and combining the estimate on $N_{31}$, we prove 
\beq\label{eq:non_N3}
\bal
N_3 =  N_{31} + N_{32} &\leq \f{ \sqrt{15}}{ 2 \sqrt{\lam_1}}  
|| (S_3 x)^{-1/2}    ||_{L^{\infty}} E \cdot 
\B(  || D_x \th_x \psi^{1/2} ||_2^2  + \lam_1  || D_x \om \vp^{1/2}||_2^2  \B) \\
&+ z_3 || D_x \th_x \psi^{1/2}||_2^2  E. 
\eal
\eeq

Next, we estimate the nonlinear parts in $\cR_{ode}$ \eqref{eq:R}. Recall $N(\om), N(\th)$ in \eqref{eq:NF}. Using integration by parts and $u = \td u + u_x(0) x$, we get 
\beq\label{eq:non_ode21}
I_1 \teq \la - u \om_x, x^{-1} \ra = \la \pa_x(\f{u}{x}), \om \ra
= \la \f{ u_x x - u }{x^2}, \om \ra 
= \la  \f{ \td u_x x - \td u}{x^2} ,\om \ra
= \la \f{4}{7} \f{ \td u_x}{x^{2/3}} + ( \f{3}{7} \f{ \td u_x}{ x^{2/3}} - \f{\td{u}}{x^{5/3}}),\f{\om}{x^{1/3}} \ra.
\eeq

Using \eqref{eq:IBPu} with $p=\f{5}{3}$ and \eqref{eq:ux_wg23} and $x^{-2/3} \leq \al_6^{-1} \vp$ due to \eqref{eq:mode}, \eqref{eq:wg}, we obtain 
\beq\label{eq:non_ode22}
|I_1| \leq (  \f{4}{7} \B| \B| \f{ \td u_x}{x^{2/3}}\B|\B|_2 +
 \B|\B|\f{3}{7} \f{ \td u_x}{ x^{2/3}} - \f{\td{u}}{x^{5/3}}\B|\B|_2
) \B| \B| \f{\om}{x^{1/3}}\B|\B|_2
\leq \B| \B| \f{ \td u_x}{x^{2/3}}\B|\B|_2 \B| \B| \f{\om}{x^{1/3}}\B|\B|_2
\leq \sqrt{ \f{p_3}{\al_6}} ||\om \vp^{1/2}||_2^2.
\eeq

Since $\la \om,x^{-1} \ra = -\f{\pi}{2} c_{\om}$, using the above estimates and \eqref{eq:est_non2}, we obtain 
\beq\label{eq:non_N6}
\bal
N_6 &\teq -\lam_2 c_{\om} \la N(\om), x^{-1} \ra
= -\lam_2 c_{\om} \la - u\om_x + c_{\om} \om, x^{-1}\ra
= - \lam_2 c_{\om} I_1 + \lam_2  \f{\pi}{2} c_{\om}^3 \\
&\leq k_c  E  \cdot ( \lam_2 \sqrt{ \f{p_3}{\al_6}} ||\om \vp^{1/2}||_2^2 + \f{\lam_2 \pi}{2} c^2_{\om} ) 
= k_c \lam_2 \sqrt{ \f{p_3}{\al_6}} E \cdot ||\om \vp^{1/2}||_2^2 + k_c E \cdot \f{\lam_2 \pi}{2} c^2_{\om} .
\eal
\eeq

Recall $d_{\th} = \la \th_x, x^{-1}\ra$ and $ u = \td u + c_{\om} x$ from \eqref{eq:cwdth}. Using integration by parts, we obtain 
\beq\label{eq:non_ode11}
\bal
\la N(\th), x^{-1} \ra &=\la -(u\th_x)_x + 2 c_{\om} \th_x , x^{-1} \ra
= \la  u \th_x, \pa_x x^{-1} \ra + 2 c_{\om} d_{\th}
= \la - u\th_x, x^{-2} \ra + 2 c_{\om} d_{\th}\\
&= \la - (\td{u} + c_{\om} x) \th_x , x^{-2} \ra +2 c_{\om} d_{\th}
= \la -  \td{u}  \th_x , x^{-2} \ra +  c_{\om} d_{\th} \teq I_2 + c_{\om} d_{\th}.
\eal
\eeq

Applying the weighted estimate in Lemma \ref{lem:wg}, \eqref{eq:ux_wg23} and Lemma \ref{lem:comp}, we get
\beq\label{eq:non_ode12}
|I_2 | \leq || \td{u} x^{-5/3}||_2 || \th_x x^{-1/3}||_2
\leq \f{6}{7} || \td{u}_x x^{-2/3}||_2 || \th_x x^{-1/3}||_2
\leq \f{6}{7} \sqrt{ \f{p_3}{b_3}} || \om \vp^{1/2}||_2  || \th_x \psi^{1/2} ||_2.
\eeq

Using Young's inequality $ab \leq \f{1}{2 \sqrt \lam_1} (a^2 + \lam_1 b^2)$, we further obtain 
\[
|I_2| \leq \f{3}{7} \sqrt{ \f{p_3}{\al_6 \lam_1}} (|| \th_x \psi^{1/2} ||_2^2 + \lam_1|| \om \vp^{1/2}||_2^2  ).
\]

By definition \eqref{energy:h1}, we have $d_{\th} \leq \lam_3^{-1/2} E$. Using the above estimate and \eqref{eq:est_non2} , we obtain 
\beq\label{eq:non_N5}
\bal
N_5 & \teq \lam_3 d_{\th} \la N(\th), x^{-1}\ra
= \lam_3 d_{\th} I_2 + \lam_3 c_{\om} d_{\th}^2 \\
&\leq \f{3}{7} \sqrt{ \f{p_3 \lam_3}{\al_6 \lam_1}} E \cdot ( || \th_x \psi^{1/2} ||_2^2 + \lam_1|| \om \vp^{1/2}||_2^2  ) +  k_c E \cdot \lam_3 d_{\th}^2.
\eal
\eeq

Combining \eqref{eq:non_N12}, \eqref{eq:non_N4}, \eqref{eq:non_N3}, \eqref{eq:non_N6} and \eqref{eq:non_N5}, we obtain the estimate on the nonlinear parts in $\cR_{L^2}, \cR_{H^1}$ \eqref{eq:R}
\[
\bal
\cN &= N_1 + \lam_1 N_2 + \lam_4( N_3 + \lam_1 N_4) + N_5 + N_6
\leq q_1 E \cdot  \la \th_x^2 ,\psi \ra
+ q_2 E \cdot  \lam_1 \la \om^2, \vp \ra\\
&+ q_3 E \cdot \lam_4 \la (D_x \th_x)^2, \psi \ra
+ q_4 E \cdot \lam_4 \lam_1 \la (D_x \om )^2, \vp  \ra
+ q_5 E \cdot  \lam_3 d_{\th}^2 + q_6 E \cdot \f{\pi \lam_2}{2} c_{\om}^2,
\eal
\]
where the coefficients $q_i >0$ are given by 
\[
\bal
q_1 &= z_1 +  \f{3}{7} \sqrt{ \f{p_3 \lam_3}{\al_6 \lam_1}} , \quad q_2 = z_2 +  \f{3}{7} \sqrt{ \f{p_3 \lam_3}{\al_6 \lam_1}} + \f{k_c \lam_2}{\lam_1} \sqrt{ \f{ p_3 }{\al_6 } },  \\
q_3 &=  z_3+ \f{ \sqrt{15}}{ 2 \sqrt{\lam_1}}  || (S_3 x)^{-1/2}    ||_{L^{\infty}} \quad
q_4  = z_4 + \f{ \sqrt{15}}{ 2 \sqrt{\lam_1}}  || (S_3 x)^{-1/2}   ||_{L^{\infty}} , \quad q_5 = k_c , \quad q_6 =  k_c.
\eal
\] 

We verify the estimates below using the computer-assisted approach outlined in Section \ref{sec:Numerical-verification}
\beq\label{ver:q}
q_i < 36, \quad i= 1,2,.., 6.
\eeq

Recall the definition of $E$ in \eqref{energy:h1}. As a result, we prove 
\[
|\cN |\leq 36 E^3.
\]

Recall $F_{\th}, F_{\om}$ in \eqref{eq:NF}, $F_1, F_2$ in \eqref{eq:l2h1}, the error terms in $\cR_{ode}, \cR_{L^2}, \cR_{H^1}$ in \eqref{eq:R}, and the energy $E$ in \eqref{energy:h1}.
Using the Cauchy-Schwarz inequality, we obtain 
\[
\bal
\cF =& \la F_{\th}, \th_x \psi \ra 
+\lam_1 \la F_{\om} , \om \vp \ra
+  \lam_4 \la D_x F_{\th}, D_x \th_x \psi \ra 
+ \lam_4 \lam_1 \la D_x F_{\om}, D_x \om \vp \ra
- \lam_2 c_{\om} \la F_{\om}, x^{-1} \ra
+ \lam_3 d_{\th} \la F_{\th}, x^{-1}\ra \\
 \leq & \B( || F_{\th} \psi^{1/2}||_2^2 +
\lam_1 || F_{\om} \vp^{1/2}||_2^2
+ \lam_4 || D_x F_{\th} \psi^{1/2}||_2^2
+ \lam_1 \lam_4 || D_x F_{\om} \vp^{1/2}||_2^2   \\
&+ \f{2 \lam_2}{ \pi } \la F_{\om}, x^{-1}\ra^2 + \lam_3 \la F_{\th}, x^{-1}\ra^2 \B)^{1/2} E
\teq \bar \e \cdot E.
\eal
\]

We verify the following using the computer-assisted approach outlined in Section \ref{sec:Numerical-verification}
\beq\label{ver:error}
\bar \e < 5.5 \cdot 10^{-7}.
\eeq
 Therefore, we prove $\cF  <5.5 \cdot 10^{-7} E$, which along with the estimate on $\cN$ conclude 
\[
\cR = \cR_{L^2} + \cR_{H^1} = \cN + \cF \leq 36 E^3 + 5.5 \cdot 10^{-7} E.
\]

\subsection{Convergence to the self-similar solution}\label{supp:conv}

Recall that we have derived the following estimate for the time derivative of the solutions 
in Section 5.3 in the main paper \cite{chen2021HL_supp}
\beq\label{eq:time_energy1}
\f{1}{2} \f{d}{dt} E_1( (\th_x)_t, \om_t)^2 
\leq -\kp  E_1( (\th_x)_t, \om_t)^2   + \cR_2,
\eeq
where $ E_1, \cR_2$ are given by 
\beq\label{eq:time_ER}
\bal
&E_1( (\th_x)_t, \om_t )= \la (\th_x)_t^2, \psi \ra + \lam_1 \la \om_t^2, \psi \ra 
+ \lam_2 \f{\pi}{2} (\pa_t c_{\om})^2 + \lam_3 (\pa_t d_{\th})^2,\\
&\cR_2 = \la \pa_t N(\th), (\th_x)_t \psi \ra + \lam_1 \la \pa_t N(\om), \om_t \vp \ra 
 - \lam_2 \pa_t c_{\om} \la \pa_t N(\om), x^{-1} \ra + \lam_3  \pa_t d_{\th} \la \pa_t N(\th), x^{-1} \ra. 
 \eal
\eeq

The perturbation $(\th_x, \om)$ enjoys the a-priori estimate $E(\th_x, \om) < E_*$. 

Our goal is to establish the following estimates 
\beq\label{eq:time_energy2}
\f{1}{2} \f{d}{dt} E_1( (\th_x)_t, \om_t)^2 
\leq - 0.02 E_1( (\th_x)_t, \om_t)^2 .
\eeq

We need the following estimate on $u$ using the weighted $L^2$ norm of $\om$.

\begin{lem}\label{lem:u}
Let $u_x = H \om$. Assume that $\om$ is odd. For $x \geq 0$, we have 
\[
 \f{| u - u_{x}(0)  x|}{x} \leq p_4 || \om \vp^{1/2}||_2 + |u_x(0)| , \quad
 p_4 = \B( \f{12 (2 + \sqrt{3}) \sqrt{p_3}}{7 \sqrt{\al_6}} \B)^{1/2}.
\]

\end{lem}

\begin{proof}
Denote $M \teq || \f{ u-u_x(0)}{x}||_{L^{\infty}}$ and $\td u = u-u_x(0) x$. Given $ \om\vp^{1/2} \in L^2$, it is not difficult to apply Lemmas \ref{lem:wg}, \ref{lem:wg13}, \ref{lem:comp} to obtain that $M < +\infty$. For $y \geq 0$, a direct calculation yields 
\[
\bal
\f{\td{u}^2}{y^2} &= \f{2}{y^2} \int_0^y \td{u} \td u_x =\f{2}{y^2} \int_0^y \td{u} u_x dx- \f{2}{y^2} u_x(0) \int_0^y \td u dx 
\leq 2 \int_0^y \f{| \td u | }{x^{5/3}} \f{ |u_x|}{x^{1/3}} dx
+ \f{2}{y^2} M |u_x(0)| \int_0^y x dx   \\
& \leq 2 || \td{u} x^{-5/3}||_2 || u_x x^{-1/3}||_2+ M |u_x(0)| \teq I_1 + I_2.
\eal
\]

Using \eqref{eq:ux_wg23}, Lemmas \ref{lem:wg}, \ref{lem:wg13} and $\al_6 x ^{-2/3} \leq \vp$ (see \eqref{eq:mode}, \eqref{eq:wg}), we get 
\[
|I_1| \leq 2 \cdot\f{6}{7} || \td{u}_x x^{-2/3}||_2 (2 + \sqrt{3}) || \om x^{-1/3}||_2
\leq \f{12}{7} \cdot (2 + \sqrt{3}) \sqrt{ \f{p_3}{\al_6}} || \om \vp^{1/2}||_2^2.
\]

Taking the supremum over $y \geq 0$, we get 
\[
M^2 \leq I_1 + I_2 \leq \f{12}{7} (2 + \sqrt{3}) \sqrt{ \f{p_3}{\al_6}} || \om \vp^{1/2}||_2^2
+ M |u_x(0)| \teq b^2  + c M.
\]

Solving the quadratic inequality, we prove 
\[
M \leq \f{ c + \sqrt{ c^2 + 4 b^2}}{2} \leq b+ c,
\]
which implies the inequality in Lemma \ref{lem:u}.
\end{proof}

Now, we are in a position to estimate $\cR_2$. Recall $E_1$ in \eqref{eq:time_ER} and $E$ in \eqref{energy:h1}. We have the following trivial estimates 
\beq\label{eq:time_EE1}
\bal
&|| \th_x \psi^{1/2}||_2 \leq E,  \quad  || D_x \th_x \psi^{1/2}||_2 \leq \lam_4^{-1/2} E,  \quad
|| \om \vp^{1/2}||_2 \leq  \lam_1^{-1/2} E, \\
&  || D_x \om \vp^{1/2}||_2 \leq (\lam_1 \lam_4)^{-1/2} E, \quad
|| \th_{xt} \psi^{1/2}||_2 \leq E_1, \quad || \om_t \vp^{1/2}||_2 \leq \lam_1^{-1/2} E_1, \\
&  
|c_{\om}| \leq k_c E, \  |\pa_t c_{\om}| \leq k_c E_1 , \quad
|d_{\th}| \leq  \lam_3^{-1/2} E, \quad |\pa_t d_{\th} | \leq \lam_3^{-1/2} E_1,
\eal
\eeq
where $k_c = (\f{\pi \lam_2}{2})^{-1/2}$. See also \eqref{eq:est_non2}.

Using integration by parts, we obtain 
\[
\bal
N_{t1} & \teq \la \pa_t N(\th), (\th_x)_t \psi \ra 
= \la \pa_t ( -u \th_{xx} + (2 c_{\om} -u_x ) \th_x) , \th_{xt} \psi \ra \\
&= \la - u (\th_{xt})_x + (2 c_{\om} -u_x) \th_{xt} , \th_{xt} \psi \ra
+ \la - u_t \th_{xx} - u_{x t} \th_x + 2 \pa_t c_{\om} \th_x, \th_{xt} \psi \ra
\teq I_1 + I_2.
\eal 
\]

The estimate of $I_1$ is the same as that of $N_1$ in \eqref{eq:est_N11} and \eqref{eq:est_N12}. Thus, we get 
\[
|I_1| \leq z_1 E \la \th_{x t}^2, \psi \ra \leq z_1 E E_1^2,
\]
where we have used \eqref{eq:time_EE1} in the last inequality.

Recall $\pa_t c_{\om} = u_{xt}(0), u_t =  \td u_t  + u_{xt}(0)x$ from \eqref{eq:cwdth}. We rewrite
\[
J \teq - u_t \th_{xx} - u_{x t} \th_x + 2 \pa_t c_{\om} \th_x
= \B( - \td u_t \th_{xx}  + \pa_t c_{\om} (\th_x  - x\th_{xx}) \B)-  \td{u}_{xt} \th_x
\teq J_1 + J_2.
\]

Applying Lemma \ref{lem:u} and \eqref{eq:time_EE1}, we get 
\[
\bal
 || J_1 \psi^{1/2}||_2
 &\leq  || \f{ \td{u}_t}{x}||_{L^{\infty}} ||D_x \th_x \psi^{1/2}||_2
 + |\pa_t c_{\om}| ( ||\th_x \psi^{1/2}||_2  + || D_x \th_x \psi^{1/2}||_2) \\
 & \leq (  p_4 || \om_t \vp^{1/2} ||_2 +  |\pa_t c_{\om}|) \lam_4^{-1/2} E
 + |\pa_t c_{\om}| (1 + \lam_4^{-1/2}) E  \\
& = \B( p_4 (\lam_1 \lam_4)^{-1/2} + k_c ( 1 + 2 \lam_4^{-1/2}) \B) E E_1.
      \eal
\]

For $J_2$, applying Lemma \ref{lem:thx_inf}, \eqref{eq:ux_wgs3} with $(D_x u_x, D_x \om)$ replaced by $( u_{x,t}, \om_t)$ yields  
\[
|| J_2 \psi^{1/2}||_2 
\leq ||\th_x S_3^{-1/2} \psi^{1/2}||_{\infty} || \td u_{x t} S_3^{1/2} ||_2
\leq \sqrt{15} || (S_3 x)^{-1/2} ||_{\infty}E \cdot || \om_t \vp^{1/2}||_2 ,
\]
which along with  \eqref{eq:time_EE1} implies 
\[
|| J_2 \psi^{1/2}||_2 \leq  \sqrt{15} || (S_3 x)^{-1/2} ||_{\infty} \lam_1^{-1/2} E E_1.
\]

Combining the above estimates on $I_1, J_1, J_2$ and using \eqref{eq:time_EE1}, we prove 
\[
\bal
N_{t1} &= I_1 + I_2 
\leq I_1 + || J\psi^{1/2} ||_2 ||\th_{xt}\psi^{1/2}||_2
\leq z_1  E E_1^2 + || J \psi^{1/2}||_2 E_1 
\leq p_5 E E_1^2,
\eal
\]
where the constant $p_5$ is given by 
\[
p_5 = z_1 +  p_4 (\lam_1 \lam_4)^{-1/2} + k_c ( 1 + 2 \lam_4^{-1/2}) + \sqrt{15} || (S_3 x)^{-1/2} ||_{\infty} \lam_1^{-1/2}  .
\]

Similarly, we have 
\[
\bal
N_{t2} &\teq 
 \lam_1 \la \pa_t N(\om), \om_t \vp \ra 
= \lam_1 \la \pa_t( -u\om_x + c_{\om} \om), \om_t \vp \ra \\
&= \lam_1 \la - u (\om_t)_x + c_{\om}  \om_t, \om_t \vp \ra
+ \lam_1 \la -u_t \om_x + \pa_t c_{\om} \om, \om_t \vp \ra
\teq I_1 + I_2.
\eal
\]

The estimate of $I_1$ is the same as that of $N_2$ in \eqref{eq:est_N21}, \eqref{eq:est_N22}. Hence, we obtain 
\[
|I_1 | \leq z_2 E \cdot \lam_1 \la \om_t^2, \vp \ra
\leq z_2 E E_1^2,
\]
where the second inequality is due to \eqref{eq:time_EE1}.

The estimate of $I_2$ is similar to that of $J$ in the above. We rewrite
\[
J = -u_t \om_x + \pa_t c_{\om} \om
= - \td u_t \om_x + \pa_t c_{\om}( \om - x \om_x) .
\]

Applying Lemma \ref{lem:u} and \eqref{eq:time_EE1}, we get 
\beq\label{eq:est_Nt2_1}
\bal
|| J \vp^{1/2}||_2
&\leq || \f{\td u_{t} }{x} ||_{\infty} || D_x \om \vp^{1/2}||_2
+ |\pa_t c_{\om}| ( || \om \vp^{1/2}||_2 + || D_x \om \vp^{1/2}||_2) \\
&\leq ( p_4  || \om_t \vp^{1/2}||_2 + | \pa_t c_{\om}| ) (\lam_1 \lam_4)^{-1/2}  E
+ |\pa_t c_{\om}| ( \lam_1^{-1/2} + (\lam_1 \lam_4)^{-1/2}) E \\
& \leq \B(  p_4 (\lam_1^2 \lam_4)^{-1/2} + k_c( \lam_1^{-1/2} +2 (\lam_1 \lam_4)^{-1/2}) \B) E E_1.
\eal
\eeq

Combining the estimates on $I_1, J$ and using \eqref{eq:time_EE1}, we prove 
\beq\label{eq:est_Nt2_2}
N_{t2}  \leq I_1 +\lam_1 || J \vp^{1/2} ||_2 || \om_t \vp^{1/2}||_2
\leq z_2 E E_1^2 +\lam_1^{1/2} || J \vp^{1/2} ||_2 E_1 
\leq p_6 E E_1^2,
\eeq
where $p_6$ is given by 
\beq\label{eq:p6}
p_6 = z_2 + p_4 (\lam_1 \lam_4)^{-1/2} + k_c( 1+2  \lam_4^{-1/2}) .
\eeq

It remains to estimate the nonlinear terms related to the ODE in \eqref{eq:time_ER}. First, we have 
\[
\bal
J& \teq \pa_t \la  N(\om), x^{-1} \ra 
= \pa_t ( \la -u\om_x, x^{-1} \ra + c_{\om} \la \om , x^{-1} \ra )  \\
&=  -\la u_t \om_x, x^{-1} \ra - \la  u \om_{tx}, x^{-1} \ra
- \pa_t  \f{\pi}{2} (c_{\om}^2)  
\teq I_1 + I_2 + I_3.
\eal
\]

For $I_1, I_2$, applying the same estimate as \eqref{eq:non_ode21},\eqref{eq:non_ode22} and then using \eqref{eq:time_EE1}, we obtain 
\[
\bal
|I_1 | + |I_2| 
&\leq || \td u_x x^{-2/3}||_2 || \om_t x^{-1/3}||_2
+ || \td u_{x t} x^{-2/3} ||_2 || \om x^{-1/3}||_2 \\
&\leq 2 \sqrt \f{p_3}{\al_6} || \om_t \vp^{1/2}||_2 || \om \vp^{1/2}||_2
\leq \f{2}{\lam_1} \sqrt \f{p_3}{\al_6}  E E_1.
\eal
\]

Note that $I_3 = -\pi c_{\om} \pa_t c_{\om}$. Using \eqref{eq:time_EE1}, we obtain
\[
\bal
N_{t4} &\teq - \lam_2 \pa_t c_{\om}  \cdot \pa_t \la  N(\om), x^{-1} \ra 
\leq \lam_2 | \pa_t c_{\om}| \cdot ( |I_1 | + |I_2| ) + \lam_2\pi |c_{\om} (\pa_t c_{\om})^2|  \\
& \leq \lam_2 k_c  \f{2}{\lam_1} \sqrt \f{p_3}{\al_6}  E E_1^2 
+ \lam_2 \pi k_c^3  E E_1^2 = p_8 E E_1^2,
\eal
\]
where $p_8$ is given by 
\[
p_8 =\lam_2 k_c \f{2}{\lam_1} \sqrt \f{p_3}{\al_6} + \lam_2 \pi k_c^3
= \lam_2 k_c \f{2}{\lam_1} \sqrt \f{p_3}{\al_6} + 2 k_c.
\]
The last equality is due to $k_c = (\f{\pi \lam_2}{2})^{-1/2}$ and $\pi \lam_2 k_c^2 = 2$.

Next, we estimate the nonlinear term related to $N(\th)$. Using the derivation \eqref{eq:non_ode11}, we get 
\[
J \teq \pa_t  \la N(\th), x^{-1} \ra
= - \la   \td u_t \th_x , x^{-2} \ra
- \la \td u \th_{xt}, x^{-1} \ra 
+ (\pa_t c_{\om } ) d_{\th} + c_{\om} \pa_t d_{\th}
\teq I_1 + I_2 + I_3 + I_4.
\]

For $I_1, I_2$, applying the same estimate as \eqref{eq:non_ode12} and then using \eqref{eq:time_EE1}, we obtain 
\[
|I_1| + |I_2| 
\leq \f{6}{7} \sqrt{ \f{p_3}{b_3}} \B(  || \om_t \vp^{1/2}||_2  || \th_x \psi^{1/2} ||_2 
+  || \om \vp^{1/2}||_2  || \th_{xt} \psi^{1/2} ||_2\B)
\leq \f{6}{7} \sqrt{ \f{p_3}{b_3} }\cdot 2 \lam_1^{-1/2} E E_1. 
\]

The estimates of $I_3, I_4$ are trivial. Using the above estimates and \eqref{eq:time_EE1}, we prove 
\[
\bal
N_{t3} &\teq \lam_3 \pa_t d_{\th} \cdot \pa_t  \la N(\th), x^{-1} \ra
\leq \lam_3 |\pa_t d_{\th}| ( |I_1| + |I_2| ) + \lam_3 |\pa_t c_{\om} | \cdot | \pa_t d_{\th} \cdot d_{\th } |
+ \lam_3 | c_{\om } | \cdot  ( \pa_t d_{\th})^2 \\
& \leq \lam_3 \lam_3^{-1/2} E_1 \cdot \f{6}{7} \sqrt{ \f{p_3}{b_3} } \cdot 2 \lam_1^{-1/2} E E_1
 +  \lam_3 k_c E_1 \cdot \lam_3^{-1 } E E_1
+\lam_3 k_c E \cdot   \lam_3^{-1} E_1^2 = p_7 E E_1^2,\\
\eal
\]
where $p_7$ is given by 
\[
p_7 =   \f{12}{7} \sqrt \f{p_3 \lam_3}{b_3 \lam_1} + 2k_c.
\]

Recall $\cR_2$ in \eqref{eq:time_ER}. Combining the estimates on $N_{ti}$, we prove 
\beq\label{eq:est_cR2}
\bal
\cR_2 = N_{t1} + N_{t2} + N_{t3} + N_{t4}
&\leq (p_5 + p_6 + p_7 + p_8) E E_1^2.
 \eal
\eeq

We verify the following estimate using the computer-assisted approach outlined in Section \ref{sec:Numerical-verification}
\beq\label{ver:p}
p_5 + p_6 + p_7 + p_8 < 300 , \eeq

Then using the a-priori estimate $E < E_* = 2.5 \cdot 10^{-5}$, we prove 
\beq\label{eq:est_cR2_2}
\cR_2 \leq 300 E_* E_1^2 < 0.01 E_1^2.
\eeq

Plugging the above estimates in \eqref{eq:time_energy1} and recall $\kp = 0.03$, we prove 
\[
\f{1}{2} \f{d}{dt} E_1( \th_{xt}, \om_t)^2 \leq - \kp E_1^2 + 0.01 E_1^2
\leq -0.02 E_1^2,
\]
which is \eqref{eq:time_energy2}.

\subsection{Uniqueness of the self-similar profile}\label{supp:uni}

Recall in Section 5.4 in the main paper \cite{chen2021HL_supp} that we have obtained the energy estimate of the difference between two solutions $(\om_1, \th_1)$ and $(\om_2, \th_2)$
\beq\label{eq:EE_uni_1}
\f{1}{2} \f{d}{dt} E_1(  \d \th_x, \d \om )^2 \leq - \kp  E_1( \d \th_x, \d \om)^2 + \cR_3,
\eeq
where the energy notation $E_1$ is defined in \eqref{energy:l2}, and $\d \th_x, \d \om, \cR_3$ are given by 
\beq\label{eq:dN0}
\bal
&\d \om  \teq  \om_1 - \om_2, \quad \d \th \teq  \th_{1} - \th_{2}, \quad 
  \d N_{\th} = N(\th_1) - N(\th_2), \quad \d N_{\om} =  N(\om_1) - N(\om_2),  \\
  &\cR_3 = \la \d N_{\th}, \d \th_x \psi \ra 
+ \lam_1 \la \d N_{\om}, \d \om \vp \ra 
- \lam_2 c_{\om}( \d \om) \cdot \la \d N_{\om}, x^{-1} \ra
+  \lam_3 d_{\th}( \d \th_x) \cdot \la \d N_{\eta}, x^{-1} \ra.
\eal
\eeq

We show that $\cR_3$ enjoys the same estimate of $\cR_2$ in \eqref{eq:est_cR2}. We estimate a typical term $  \lam_1 \la \d N_{\om}, \d \om \vp \ra $. Using  the elementary identity
\[
a_1 b_1 - a_2 b_2 
= a_1( b_1 - b_2) + (a_1 - a_2) b_2
= a_1 \cdot \d b + \d a \cdot b_2, \quad \d a = a_1 - a_2, \quad \d b = b_1 - b_2,
\]
 we rewrite $\d N_{\om}$ as follows 
\[
\bal
\d N(\om) &= ( c_{\om}(\om_1) \om_1 - u(\om_1) \om_{1,x} )
- ( c_{\om}(\om_2) \om_2 - u(\om_2) \om_{2,x} ) \\
&= ( c_{\om}(\om_1) \d \om - u(\om_1) \d \om_{x} )
+ c_{\om}( \d \om) \om_2 - u( \d \om) \om_{2,x}.
\eal
\]

Thus, we derive 
\[
\bal
N_{\d2} &\teq 
 \lam_1 \la \d N_{\om}, \d \om \vp \ra 
= \lam_1 \la c_{\om}(\om_1) \d \om - u(\om_1) \d \om_{x}  , \d \om \vp \ra
+ \lam_1  \la c_{\om}( \d \om) \om_2 - u( \d \om) \om_{2,x} , \d \om \vp \ra
\teq I_1 + I_2.
\eal
\]

The estimate of $I_1$ is the same as that of $N_2$ in \eqref{eq:est_N21}, \eqref{eq:est_N22}. Hence, we obtain 
\[
|I_1 | \leq z_2 E( \th_{x, 1}, \om_1) \cdot \lam_1  || \d \om \vp^{1/2}||_2^2
\leq z_2 E_* E_1^2,
\]
where the second inequality is due to the definition of $E_1( \d \th_x, \d \om)$ \eqref{energy:l2} and the a-priori estimate $E(\th_{x,1}, \om_1) < E_*$. For $I_2$, we rewrite
\[
J = - u(\d \om) \om_{2, x} +  c_{\om}(\d \om) \om_2
= - \td u(\d \om) \om_x +  c_{\om}(\d \om) (  \om_2 - x \om_{2,x}) .
\]

Applying Lemma \ref{lem:u} and \eqref{eq:time_EE1}, we get 
\[
\bal
|| J \vp^{1/2}||_2
&\leq \B|\B| \f{\td u_x(\d \om)}{x} \B|\B|_{\infty} || D_x \om_2 \vp^{1/2}||_2
+ | c_{\om}(\d \om) |  \cdot ( || \om_2 \vp^{1/2}||_2 + || D_x \om_2 \vp^{1/2}||_2). \\
\eal
\]

The remaining steps are exactly the same as those in \eqref{eq:est_Nt2_1} and \eqref{eq:est_Nt2_2}. We obtain 
\[
N_{\d2} \leq p_6 E_* E_1^2( \d \th_x,  \d \om ),
\]
where $p_6$ is defined in \eqref{eq:p6}. 

The estimates of other terms in $\cR_3$ in \eqref{eq:dN0} follow the same estimates of $N_{t1}, N_{t3}, N_{t4}$ in Section \ref{supp:conv}. Thus, we yield the same estimates for $\cR_3$ as \eqref{eq:est_cR2}-\eqref{eq:est_cR2_2}, 
\[
\cR_3 \leq (p_5 + p_6 + p_7 + p_8) E_* E_1^2( \d \th_x, \d \om)
< 0.002 E_1^2( \d \th_x, \d \om).
\]

Plugging the above estimate of $\cR_3$ in \eqref{eq:EE_uni_1}, we prove 
\[
\f{1}{2} \f{d}{dt} E_1(  \d \th_x, \d \om )^2 \leq - \kp  E_1^2( \d \th_x, \d \om) + 
0.01 E_1^2( \d \th_x, \d \om)  \leq -0.02 E_1^2( \d \th_x, \d \om) .
\]

\subsection{ Estimates on the self-similar profiles }\label{supp:profile}
\

\textbf{ Notations } In this Section, we use the notation $A \les B$ if there exists some finite constant $C > 0$, such that $A \leq C B $. The constant $C$ \textit{can} depend on the norms of the approximate steady state $( \bar \th, \bar \om)$ and the self-similar profile $(\th_{\inf}, \om_{\inf})$ constructed in the main paper \cite{chen2021HL_supp}, e.g. $|| \th_x||_{\inf}, || \bar \th_x||_{\inf}$, as long as these norms are finite. 
In the main paper \cite{chen2021HL_supp}, we use the following estimates on the self-similar profile
\beq\label{eq:profile_est}
  | c_{l,\inf} x + u_{\infty}(x) | \geq 0.3 |x| ,
 \quad |u_{\infty}(x) | \les |x|^{5/6},  
   \quad u_{\inf, x} \in L^{\inf},  \quad \th_{\infty}(1) \neq 0, \quad \th_{x, \inf} \in L^{\inf},
 \eeq
 To prove these estimates, we need the following properties about the approximate steady state 
\beq\label{supp:ver1}
 |\bar c_l x + \bar u| \geq 0.4 x,  \quad \bar u_x \in L^{\inf},  \quad \bar \th_x \in C^1 \cap L^{\inf}, \quad \bar \th_{xx}(0) \neq 0, \quad \bar \om \in L^2( x^{-2/3}),
\eeq
whose verifications are deferred to Section \ref{supp:profile2}.

Next, we prove the estimates in \eqref{eq:profile_est} in order. Due to symmetry, it suffices to consider the case $x \geq 0$. Denote $\td \th = \th_{\inf} - \bar \th, \td \om = \om_{\inf} - \bar \om $. Recall the a-priori estimate on the profile 
\beq\label{eq:profile_est2}
E(\td \th_x, \td \om) =  E(  \th_{\inf, x} -\bar \th_x ,  \om_{\inf} - \bar \om ) \leq E_* = 2.5 \cdot 10^{-5}.
\eeq

\paragraph{\bf{First estimate}} Recall the normalization condition $c_{l,\inf} =\bar c_l, \td c_{\om, \inf}- \bar c_{\om}= u_x( \td \om)$. We get  \[
 | c_{l,\inf} x +  u_{\infty}| 
 \geq | \bar c_{l } x + \bar  u | - |  u_{\inf} - \bar u| 
  = | \bar c_{l } x + \bar  u | - | u(\td \om)  | 
  \geq | \bar c_{l } x + \bar  u | - | u(\td \om) - H \td \om (0) x |
  - | H \td \om(0) x |
\]

Using \eqref{supp:ver1}, Lemma \ref{lem:ux_inf}, $||C_u(x)||_{\inf} < 20$, and $ E(\td \th_x, \td \om) > | H \td \om(0)|$, we yield 
\[
 | c_{l,\inf} x +  u_{\infty}|  
 \geq ( 0.4 - || C_u(x)||_{\inf} E( \td \th_x, \td \om) -  E( \td \th_x, \td \om) ) |x|
 \geq ( 0.4 - 20  E_*- E_* ) |x| \geq 0.3 |x|. 
\]

\vspace{0.1in}
\paragraph{\bf{Second estimate}}
Using $ || \td \om x^{-1/3}||_2 \les E( \td \th_x, \td \om) \les 1 $ and $ \bar \om \in L^2(x^{-2/3})$, we get $ \om_{\inf} \in L^2(x^{-2/3})$. Firstly, using a direct estimate, we obtain 
\[
\B|\f{u_{\inf}}{x}\B| \leq \int_x^{\inf} \B| \pa_y \f{u_{\inf}(y)}{y} \B| dy
\les \int_x^{\inf} \B| \f{u_{\inf,x}}{y} \B| + \B| \f{u_{\inf}}{y^2} \B| dx 
\les ( || u_{\inf,x} x^{-1/3}||_2  + || u_{\inf} x^{-4/3}||_2 ) ( \int_x^{\inf} y^{-4/3} dy  )^{1/2}.
\]

Applying the Hardy inequality, Lemma \ref{lem:wg13} and $\om_{\inf} \in L^2(x^{-2/3} )$, we yield 
\[
\B|\f{u_{\inf}}{x}\B| \les || u_{\inf,x} x^{-1/3}||_2 |x|^{-1/6} 
\les || \om_{\inf} x^{-1/3} ||_2 |x^{-1/6} | \les |x|^{-1/6}.
\]

\vspace{0.1in}
\paragraph{\bf{Third estimate}}

From \eqref{supp:ver1}, we get  $\bar u_x \in L^{\inf}$. It suffices to prove $u_{\inf, x} - \bar u_x = H \td \om \in L^{\inf}$. Using Lemma \ref{lem:ux_inf},  we get $H  \td \om - H \td \om(0) \in L^{\inf}$. Since $|H \td \om (0)| <
E(\td \th_x, \td \om) < 1$, we prove $H \td \om \in L^{\inf}$.

\vspace{0.1in}

\paragraph{\bf{Fourth estimate}}
From \eqref{supp:ver1}, we get that $\bar \th_x $ vanishes linearly near $x = 0$. Since  $\th_{x,\inf} - \bar \th_x \in L^2(\psi)$ with the singular weight $\psi$, we must have a nontrivial profile 
$
 \th \not \equiv 0 .$ In the main paper \cite{chen2021HL_supp}, we derive
\[
\th_{\inf} = \th_{\inf}(1) \exp( \int_1^x \f{c_{\th,\inf}} {c_{l, \inf } x + u_{\inf} })
\]

Since $\th_{\inf}$ is nontrivial, we must have 
\[
 \th_{\infty}(1) \neq 0,  \quad \th_{\inf}(x) \neq 0  \mathrm{ \ for \ all \ } x \neq 0. 
 \]

\vspace{0.1in}
\paragraph{\bf{Fifth estimate }}
Since $|x \psi| \geq C |x|^{1/3} $ for large $x$ and $|x \psi| \geq C |x|^{-3}$ for small $x$
with some constant $C$ depending on the approximate steady state $(\bar \th, \bar \om)$, we get $ |x \psi| \geq C$.  Using Lemma \ref{lem:thx_inf} and \eqref{supp:ver1}, we yield 
\[
| \th_{\inf, x} | \leq || \bar \th_x ||_{\inf}
+ || \td \th_x ||_{\inf}
\les 1  + E(\td \th_x, \td \om) \les 1.
\]

\subsubsection{ Regularity of the approximate steady states}\label{supp:profile2}

We verify the properties of the approximate steady state in \eqref{supp:ver1}, i.e.
\beq\label{supp:ver12}
\bar u_x \in L^{\inf},  \quad \bar \th_x \in C^1 \cap L^{\inf}, \quad \bar \om \in L^2( x^{-2/3}),
\eeq
and 
\beq\label{supp:ver13}
 |\bar c_l x + \bar u| \geq 0.4 x,  \quad  |\bar \th_{xx}(0)| > 0. 
\eeq

Recall that the profiles enjoy the decomposition $\bar v = \bar \th_x = v_p + v_b, \bar \om = \om_b + \om_p$ \eqref{eq:solution_decomposition}, where $\om_p, v_p$ are piecewise polynomials with degree $\geq 3$ and compact support, and $\om_b, v_b$ admit the analytic forms \eqref{eq:explicit_part}. We yield $ \bar \th_x, \bar \om \in C^1 \cap L^{\inf} $. Since $\om \in L^{\inf}$, $\om_b$ has a decay rate $ |x|^{-a_{\om}}, a_{\om} \geq \f{1}{4}$ and $\om_p$ has compact support, we get $ \om = \om_p + \om_b \in L^2(x^{-2/3})$.
Recall the notation $\td u_x = u_x - u_x(0) $. Using the estimates \eqref{eq:u_Linf1}, \eqref{eq:u_Linf2} and \eqref{eq:ux_wg23}, we obtain 
\[
| \bar u_x - \bar u_x(0)|^2
\les || D_x \om \cdot x^{-1/3} ||_2 || \om (x^{-2} + x^{-4/3}  + x^{-2/3})^{1/2}||_2,
\]
where $D_x = x \pa_x$. Since $\om ,\om_x \in L^{\inf}$, $\om_b , x \om_{b,x}$ have a decay rate at least $|x|^{-1/4}$ and $\om_p$ has compact support, the above right hand side is finite. The term $\bar u_x(0)$ 
\[
\bar u_x(0) = -\f{2}{\pi} \int_0^{\inf} \f{\om(y)}{y} dy
\]
is finite since $\om \in L^{\inf}$ and $\om$ has a decay rate at least $|x|^{-1/4}$. We yield $\bar u_x \in L^{\inf}$.
Thus, we verify the properties in \eqref{supp:ver12}.

For the inequalities in \eqref{supp:ver13} on the approximate steady state, we will verify them rigorously using the numerical methods to be introduced in Sections \ref{sec:Computation} and 
\ref{sec:Numerical-verification}, in particular the method introduced in Section \ref{sec:maxmin_estimate}.

 \section{Construction of the approximate self-similar profile}\label{sec:Computation}

In this section, we describe how we numerically construct the approximate steady state $\bar{\om},\bar{\theta}_x$ and the corresponding rescaling parameters $\bar{c}_l,\bar{c}_{\om}$ of the dynamic rescaling equations of the HL model:
\beq\label{eq:HLdyn}
\bal
\om_t + (c_l x + u  ) \om_x &= c_{\om} \om + \th_x , \\
(\th_x)_t + (c_l x + u) \th_{xx} &= (2c_{\om} - u_x) \th_x, \\ 
u_x &= H \om,
\eal
\eeq
subject to the normalization conditions 
\beq\label{eq:normal0}
c_l = 2 \f{\th_{xx}(0) }{\om_x(0) }, \quad c_{\om} = \f{1}{2} c_l + u_x(0). 
\eeq
Our construction procedure basically follows from the approach in \cite[Section 4]{chen2019finite_supp} with adaption to the HL model. The main idea is to solve the equations \eqref{eq:HLdyn} numerically for long enough time so that the solution converges to some approximate steady state. We use the residuals to measure the accuracy of the numerical solution. When the maximum residual drops to some prescribed threshold, we stop the computation and use the numerical solution as the approximate steady state. For notational simplicity, we will write $v = \theta_x$ in what follows. 

\subsection{Setup of the solution} For the sake of numerical accuracy, we choose to decompose the solution $\om,v$ into an explicit part and a perturbation part: 
\beq\label{eq:solution_decomposition}
\om(t,x) = \om_b(x) + \om_p(t,x), \quad v(t,x) = v_b(x) + v_p(t,x).
\eeq
The explicit parts $\om_b, v_b$ are fixed and can be expressed by explicit functions, while the perturbation parts $\om_p, v_p$ will be computed numerically. Solving the equations for $\om,v$ now transfer to solving the equations for the perturbations $\om_p, v_p$. 

\subsubsection{Construction of the explicit parts} We construct the explicit parts of $\om$ and $v$ as
\beq\label{eq:explicit_part}
\om_b(x) = \frac{b_\om x^5}{1 + |x|^{5 + a_\om}} + \frac{s_\om x}{1 + (r_\om x)^2}\quad \text{and}\quad v_b(x) = \frac{b_v x^5}{1 + |x|^{5 + a_v}} + \frac{s_v x}{1 + (r_v x)^2}, \quad x\geq 0, 
\eeq
for some constants $a_\om,b_\om, r_\om, s_\om, a_v, b_v, r_v, s_v$. The values of $\om_b$ and $v_b$ for $x<0$ are obtained using the odd symmetry. In the following discussions, we will always assume that $x\geq 0$. The parameters in \eqref{eq:explicit_part} are chosen under one principle: the explicit parts should approximate the steady state $\bar{\om}, \bar{v}$ well enough for both $x$ near $0$ and $x$ sufficiently large. In such a way, the smallness of the perturbation parts $\om_p, v_p$ in both the near field and the far field can help avoid large round-off errors in our computation.

The second terms in both expressions in \eqref{eq:explicit_part} aim to capture the behavior of $\om$ and $v$ near $x=0$, especially the slops at $x=0$. In particular, we choose the parameters in the second terms as
\[r_\om = 1.3468, \quad s_\om = 0.6734, \quad r_v = 1.0101, \quad s_v = 1.0101.\]
We remark that the choice of these parameters is not essential; they only serve to normalize the approximate steady state.

The first terms in the expressions in \eqref{eq:explicit_part} serve to capture the tail behavior of $\om$ and $v$ in the far field. The choice of the corresponding parameters is trickier and is based on the following method.
We first need to do some preliminary computation on solving the equations \eqref{eq:HLdyn} to get a priori information of the far field decay rates (in $x$) of the approximate steady state. As we have argued in Section 4 of the Main Paper, the steady state solutions of $\om$ and $v$ should behave like
\beq\label{eq:tail_behavior}
\bar{\om}(x) \sim \bar{b}_\om x^{\bar{c}_\om/\bar{c}_l}, \quad \bar{v}(x) \sim \bar{b}_v x^{2\bar{c}_\om/\bar{c}_l}
\eeq
for $x$ large enough and some scaling constants $\bar{b}_\om, \bar{b}_v$. In other words, we need to numerically obtain approximations of the ratio $\bar{c}_\om/\bar{c}_l$ and the coefficients $\bar{b}_\om, \bar{b}_v$, and then we choose the parameters $a_\om, b_\om, a_v, b_v$ accordingly. In particular, we choose 
\beq\label{eq:aw}
a_\om = \frac{1.00043212}{3}, \quad b_\om = 1.37954, \quad a_v = 2 a_\om, \quad b_v = 1.7711,
\eeq
according to our preliminary computations. We remark that the values of $a_{\om}, a_v$ only depend on the ratio $c_\om/c_l$, which is intrinsic to the HL model; while the values of $b_\om, b_v$ are scalable and depend on the choice of $r_\om, s_\om, r_v, s_v$.

\subsubsection{Representation of the perturbation parts} On the one hand, the perturbations parts $\om_p,v_p$ will be computed numerically; on the other hand, they need to be well defined for all $x$ so that we can use them for computer-aided analysis. Out of these considerations, we choose to represent $\om_p, v_p$ by piecewise quintic polynomials (fifth-degree polynomials) on a compact support $[-L,L]$ for some sufficiently large $L$. Beyond this domain, we set $\om_p = v_p = 0$. In our computation, we choose $L = 10^{16}$. Again, due to the odd symmetry of $\om_p, v_p$, we only need to consider the computation for $x\geq 0$.

First of all, we construct an adaptive mesh $\{x_i\}_{i=0}^n$ on the domain $[0,L]$ with $x_{i-1}< x_i$, $x_0 = 0$ and $x_n = L$. The adaptive grid points are chosen such that either $x_i - x_{i-1} \leq 10^{-3}$ or $(x_i - x_{i-1})/x_i \leq 0.05$.

Given the grid points, the functions $\om_p, v_p$ are represented by piecewise quintic polynomials on the mesh intervals. For example, on each interval $[x_{i-1},x_i]$, the solution $\om(x)$ has the expression
\[\om_p(x) = p_i(x) := \sum_{l=0}^5c_{i,l}\left(\frac{x-x_{i-1}}{x_i-x_{i-1}}\right)^l,\quad x\in[x_{i-1},x_i].\]
For each $i=1,\dots,n$, the coefficients $c_{i,l}, l=0,\dots,5$ are uniquely determined by the six values 
\[(\om_p)_j = p_i(x_j),\quad (\partial_x\om_p)_j = p_i'(x_j),\quad (\partial_x^2\om_p)_j = p_i''(x_j),\quad j = i-1,i.\]
Here $\{(\om_p)_i\}_{i=0}^n$, $\{(\partial_x\om_p)_i\}_{i=0}^n$ and $\{(\partial_x^2\om_p)_i\}_{i=0}^n$ are the values, the first derivatives, and the second derivatives of $\om_p$, respectively, at the grid points. We require that $(\om_p)_n = (\partial_x\om_p)_n = (\partial_x^2\om_p)_n = 0$ due to the compact support condition. Moreover, we impose the condition that the third derivative of the piecewise polynomial $\om_p(x)$ is continuous on $[0,L]$. As a result, the second derivatives $\{(\partial_x^2\om_p)_i\}_{i=0}^n$ are obtained from $\{(\om_p)_i\}_{i=0}^n$, $\{(\partial_x\om_p)_i\}_{i=0}^n$ by solving a linear system corresponding to the standard quintic spline interpolation on the grid points $\{x_i\}_{i=0}^n$ subject to the boundary conditions
\[(\partial_x^2\om_p)_0 = 0\ (\text{odd symmetry})\quad\text{and}\quad (\partial_x^2\om_p)_n = 0\ (\text{compact support}).\] 
Therefore, the piecewise quintic polynomial $\om$ is completely represented by the values $\{(\om_p)_i\}_{i=0}^n$, $\{(\partial_x\om_p)_i\}_{i=0}^n$. Similarly, the piecewise quintic polynomial $v_p$ is completely represented by the values $\{(v_p)_i\}_{i=0}^n$, $\{(\partial_xv_p)_i\}_{i=0}^n$. In this way, the numerical solution $\om_p,v_p$ are well defined on the support $[0,L]$ and are in the class $C^3\cap H^4$.

\subsubsection{Discussions on the decomposition} We remark that our decomposition \ref{eq:solution_decomposition} of the approximate steady state is not used directly in our analysis; it is designed to help us numerically find an approximate steady state $\bar{\om}, \bar{v}$ that has a sufficiently small residual.

As we have argued, the real steady state solutions should have the tail behavior as in \eqref{eq:tail_behavior}. Therefore, a nice approximate steady state must also follow these power laws up to a high precision in the far field; otherwise the residuals would not be small enough (relatively) for large $x$. This would require the approximate steady state to be functions defined on the whole real line. However, standard numerical discretization methods, such as the piecewise polynomial interpolation, can only work on a finite domain. If we use compactly supported functions to approximate the steady state, the residual cannot be small enough for our analysis in the far field; unless we use a uniformly dense mesh over a super large domain (up to $x = 10^{30}$), which is way beyond our current computational capacity. 

To alleviate this difficulty, we choose to numerically compute the profiles of the approximate steady state only in a relatively smaller domain $[-L, L]$; for the tail behavior \eqref{eq:tail_behavior} beyond this computational domain, we approximate them with some explicit functions. This is why we introduce the decomposition \eqref{eq:solution_decomposition} for the solution. The first terms in both expressions in \eqref{eq:explicit_part} with the form $F_a = x^5/(1+|x|^{5+a})$ have a decay rate $x^{-a}$ for large $x$ and serve to capture the leading asymptotic behavior of the solution beyond $[-L, L]$.  We choose the power $x^5$ so that $F_a$ is smooth enough near $x=0$. In particular, we have $F_a \in C^{10, a}$.

Moreover, the high accuracy in computing the Hilbert transform of $\om$ is critical in obtaining a good approximate steady state. Though we can compute $H(\om_p)$ for the perturbation part $\om_p$ accurately using its piecewise quintic polynomial representation, it stills introduces non-negligible round-off errors if $\om_p$ is large, especially in the far field. 
This difficulty cannot be overcome by simply using a finer mesh. To overcome this difficulty, we need to construct a smooth and explicit function which captures the far field behavior of $\om$ and its Hilbert transform can be computed accurately. The construction of such function and the computation of its Hilbert transform are the essential new difficulties in constructing the approximate steady state. We refer to Section \ref{sec:new_diff} for more discussion.

\subsection{Computing the Hilbert transform}\label{sec:compute_u} We will use the values of $u$ and its derivatives in both our computation and analysis. Correspondingly, $u$ can also be decomposed as $u = u_b + u_p$, such that $\partial_x u_b= H(\om_b)$ and $\partial_x u_p= H(\om_p)$. The values of $u_b$ and its derivatives can be calculated or estimated semi-analytically using the explicit expression of $\om_b$ in \eqref{eq:explicit_part}.

Moreover, with the explicit piecewise representation of $\om$, we can accurately evaluate $u_p$ and its derivatives at any point $x$ using the following exact formulas:
\begin{equation}\label{eq:compute_u}
	\begin{split}
	&u_p(x) = \f{1}{\pi} \int_0^L \omega_p(y)\log\left|\frac{x-y}{x+y}\right|dy,\quad \partial_x u_p(x) = \f{1}{\pi}\int_0^L \frac{2y}{x^2-y^2}\omega_p(y)dy,\\
	&\partial_x^2 u_p(x) = \f{1}{\pi}\int_0^L \frac{2x}{x^2-y^2}\partial_x \omega_p(y)dy,\quad \partial_x^3 u_p(x) = \f{1}{\pi}\int_0^L \frac{2y}{x^2-y^2}\partial_x^2\omega_p(y)dy.
	\end{split}
\end{equation}
Note that we do not need to use $\partial_x^3 u_p$ to compute the approximate steady state; we will need it in the numerical verification procedure. The accurate calculation of these integrals will be discussed in Section \ref{sec:Numerical-verification}.

\subsection{Solving the dynamic rescaling equations} 
\begin{subequations}\label{eq:ds_numerical}
Regarding the preceding setting, to obtain an approximate steady state $\bar{\om}$ means to obtain the grid point values $\{(\om_p)_i\}_{i=0}^n$ and $\{(\partial_x\om_p)_i\}_{i=0}^n$ for the perturbation part (since the explicit part is given and fixed). We thus choose to numerically solve the dynamic rescaling equations for both $\om$ and its derivate $\om_x$:
\begin{align}
(\om_p)_t + (c_l x + u  ) \om_x &= c_{\om} \om + v, \label{eq:ds_w}\\
(\partial_x \om_p)_t + (c_l x + u  ) \om_{xx} &= - (c_l + u_x)\om_x + c_{\om}\om_x + v_x. \label{eq:ds_wx}
\end{align}
We have used that $\om_t = (\om_p)_t$. That is to say, we update the first derivatives $\{(\partial_x\om_p)_i\}_{i=0}^n$ by numerically solving the equation \eqref{eq:ds_wx} in time, instead of obtaining them directly from the grid point values $\{(\om_p)_i\}_{i=0}^n$. Likewise, we solve the dynamic rescaling equations for both $v=\theta_x$ and $v_x = \theta_{xx}$:
\begin{align}
(v_p)_t + (c_l x + u  ) v_x &= (2c_\om - u_x) v, \label{eq:ds_v}\\
(\partial_x v_p)_t + (c_l x + u  ) v_{xx} &=  (2c_\om - c_l - 2u_x) v_x -u_{xx}v. \label{eq:ds_vx}
\end{align}
\end{subequations}
We will denote by $\om^k,v^k$ the solutions at the discrete time steps $t_k, k\geq 0$, and by $(\om_p)^k, (v_p)^k$ the corresponding perturbation parts. That is, $\om^k = \om_b + (\om_p)^k$ and $v^k = v_b + (v_p)^k$. The superscript $k$ in this section stands for the time discretization. The over algorithm for numerically solving the equations \eqref{eq:ds_numerical} is described as follows:
\begin{enumerate}
	\item The initial data $\{(\om_p)_i^0\},\{(\partial_x \om_p)_i^0\},\{(v_p)_i^0\},\{(\partial_x v_p)_i^0\}$ for the perturbation parts are all chosen to be $0$.
	\item The second derivatives $\{(\partial_x^2\om_p)_i^k\}$ (or $\{(\partial_x^2v_p)_i^k\}$) are obtained from $\{(\om_p)_i^k\},\{(\partial_x \om_p)_i^k\}$ (or $\{(v_p)_i^k\},\{(\partial_x v_p)_i^k\}$) using the standard quintic spline interpolation subject to the boundary conditions $(\partial_x^2\om_p)_0^k = (\partial_x^2\om_p)_n^k = 0$ (or $(\partial_x^2v_p)_0^k = (\partial_x^2v_p)_n^k = 0$).
	\item The grid point values of $u^k$ are given by $u_i^k = (u_b)_i + (u_p)_i^k$, where $(u_p)_i^k$ is computed using the formula in \eqref{eq:compute_u}. The grid point values of $u_x^k$ and $u_{xx}^k$ are computed similarly.
	\item The rescaling parameters $c_l^k,c_\om^k$ are computed according to the normalization conditions \eqref{eq:normal0}: 
	\[c_l^k = 2\frac{v_{x,0}^k}{w_{x,0}^k}\quad \text{and}\quad c_\om^k = \frac{c_l^k}{2} + u_{x,0}^k.\]
	Note that that these normalization conditions ensure that 
	\[\om_{x,0}^k = \om_{x,0}^0 = (\partial_x \om_b)_0 \quad \text{and}\quad v_{x,0}^k = v_{x,0}^0 = (\partial_x v_b)_0 \quad \text{for all $k\geq 0$}, \]
	which further implies that 
	\[(\partial_x \om_p)_0^k = (\partial_x v_p)_0^k = 0 \quad \text{for all $k\geq 0$}.\]
	\item Due to the odd symmetry and compact support condition of $\om_p,v_p$, we fix $(\om_p)_0^k = (v_p)_0^k = (\om_p)_n^k = (\partial_x \om_p)_n^k = (v_p)_n^k = (\partial_x v_p)_n^k = 0$. By the design of the normalization conditions \eqref{eq:normal0}, we also fix $(\partial_x \om_p)_0^k = (\partial_x v_p)_0^k = 0$. The other grid point values $(\om_p)_i^k,(\partial_x \om_p)_i^k,(v_p)_i^k,(\partial_x v_p)_i^k, i = 1,\dots,n-1$ are marched in time according to the equations \eqref{eq:ds_numerical} using a $4$th order Runge--Kutta scheme with adaptive time stepping. The discrete time step size $\Delta t_k = t_{k+1}-t_k$ is given by 
	\[\Delta t_k = \frac{1}{10}\min_{1\leq i\leq n-1}\frac{\min\{x_i-x_{i-1}, x_{i+1}-x_i\}}{c_l^k x_i + u_i^k},\]
	respecting the CFL stability condition.
\end{enumerate}

The accuracy of the approximate steady state solution $\om^k,v^k$ is measured by the magnitude of the residuals:
\begin{equation}\label{eq:residual}
\begin{split}
F_\om^k &:= -(c_l^kx + u^k)w_x^k + c_{\om}^k \om^k + v^k,\\
F_v^k &:= - (c_l^k x + u^k) v_x^k + (2c_\om^k - u_x^k) v^k.
\end{split}
\end{equation}
Note that $F_v$ has the same meaning as $F_\th$ defined in \eqref{eq:NF}; we will consistently use $F_v$ for $F_\th$ in what follows. The computation is stopped when the grid point values of residuals satisfy
\[\max_{0\leq i\leq n} \{|F_{\om,i}^k|,|F_{v,i}^k|\} \leq 10^{-10}.\] 
We then use $\bar{\om} = \om^k, \bar{v} = v^k$ as the approximate steady state. The corresponding scaling parameters are
\[
\bar{c}_l = 3 \quad \text{and}\quad \bar{c}_\om = - 1.00043212. 
\] 

We remark that we observe at least fourth order convergence in space and in time for the numerical method described above. However, we do not actually need to do convergence study (by refining the discretization) for our scheme, as we can measure the accuracy of our approximate self-similar profile {\it a posteriori}. The criterion for a good approximate self-similar profile is that it is piecewise smooth and has a small residual in the weighted energy norm that is used in our following analysis.

\subsection{Computer-assisted proof}
All the numerical computations and quantitative verifications are performed by MATLAB (version 2020a) in the double-precision floating-point operation. The MATLAB codes can be found via the link \cite{Matlabcode_supp}. To make sure that our computer-assisted proof is rigorous in the verification stage, we adopt the standard method of interval arithmetic (see \cite{rump2010verification_supp,moore2009introduction_supp}). In particular, we use the MATLAB toolbox INTLAB (version 11 \cite{Ru99a_supp}) for the interval computations. Every single real number $p$ in MatLab is represented by an interval $[p_l,p_r]$ that contains $p$, where $p_l,p_r$ are precise floating-point numbers of 16 digits. Every computation of real number summation, multiplication or division is performed using the interval arithmetic, and the outcome $P$ is hence represented by the resulting interval $[P_l,P_r]$ that strictly contains $P$. We can thus obtain a rigorous upper bound on $|P|$ by rounding up $\max\{|P_l|,|P_r|\}$ to 2 significant digits (or 4 when necessary). More detailed discussions on rigorous, computer-assisted verifications of quantities and inequalities are presented in Appendix \ref{sec:Numerical-verification}.

\section{Rigorous computer-assisted verifications}\label{sec:Numerical-verification}
As we can see in the preceding sections, our stability analysis for the dynamic rescaling equations relies on the validity of some inequalities that depend on the approximate steady state $\bar{\om},\bar{v}$. We have discussed in Section \ref{sec:Computation} how we construct the approximate self-similar profiles $\bar{\om},\bar{v}$ as a combination of the explicit parts ($\om_b, v_b$) that have analytic expressions and the perturbation parts ($\om_p,v_p$) that are represented by piecewise quintic polynomials. In this section, we will describe some strategies for rigorous numerical verifications of the claimed inequalities, based on the explicit expressions of the approximate steady state.

\subsection{The overall strategy} Note that our approximate steady state $\bar{\om},\bar{v}$ are defined on the whole real line $\mathbb{R}$ since the the explicit parts $\om_b, v_b$ are so defined. To verify an inequality involving $\bar{\om}_b$ usually requires to verify a point-wise inequality on $\mathbb{R}$ or to compute an integral on $\mathbb{R}$. Using the odd symmetry of $\om_b$ (or even symmetry of other functions), we only need to do so on the half line $\mathbb{R}_+$, which, however, is still an infinite domain. In order to accomplish the verifications using computer-aided methods, we first need to decompose the infinite domain into two parts:
\[\mathbb{R}_+ = [0,L_B] \cup (L_B, \infty),\]
where $L_B$ is a sufficiently large number. In particular, we want $L_B \geq L$ so that $[-L_B, L_B]$ contains the support of the perturbation part $\om_p$. In practice, we choose $L_B = 100L = 10^{18}$.

For the verification work on the finite domain $[0,L_B]$, we will use mesh-based methods. We first put an adaptive dense mesh $0= x_0<x_1< \cdots < x_N = L_B$ over the domain $[0,L_B]$. Note that this mesh for verification is much finer than the mesh we use for computing the approximate steady state in Section \ref{sec:Computation}. Then we evaluate the functions we need on the grid points and obtain accurate information of these functions in each small interval $[x_i,x_{i-1}]$ using estimates. For example, we can accurately estimate the maximum of a function $f$ on $[0,L_B]$ by checking over its maximum on each interval $[x_i,x_{i-1}]$, which can be bounded above using the grid point values $f_{i-1}, f_i$ and the estimates on the derivative $f_x$. We will further discuss the idea of mesh-based estimates in the upcoming sections. 

For the verification work on the infinite domain $(L_B, \infty)$, we will use rigorous analysis based on the decay properties of the $\bar{\om},\bar{v}$ and other involved functions.

\subsection{Interval arithmetic} As mentioned at the end of Section \ref{sec:Computation}, in the verification procedure, every single real number $p$ is represented by an interval $[p_l,p_r]$ that contains $p$, where $p_l,p_r$ are some floating-point numbers. In particular, if $p$ is evaluated directly (e.g. $p = \pi$), then $p_l = \lfloor p\rfloor_f$ and $p_r = \lceil p\rceil_f$. Here and below, $\lfloor\cdot\rfloor_f$ and $\lceil\cdot\rceil_f$ stand for the rounding down and rounding up to the nearest floating-point value, respectively. If $p$ is the result of some floating-point operation, then $p_l,p_r$ are obtained using the standard interval arithmetic. For example, if $p = a+b$, then $p_l = a_l + b_l$ and $p_r = a_r + b_r$. All the following computer-assisted estimates are based on the rigorous interval arithmetic. 

\subsection{Accurate grid values of the Hilbert transform.} \label{sec:compute_up}
As mentioned in Section \ref{sec:Computation}, the values of $u_p$ and its derivatives can be computed accurately with rigorous error bound from the integral formulas \eqref{eq:compute_u} using the piecewise polynomial representation of $\om_p$. To explain the algorithm, we will need the following notation: given any $x\geq0$, if $x\in[0,L)$, then find interval $[x_{j-1},x_j)$ such that $x\in [x_{j-1},x_j)$ and define $h(x) = x_j - x_{j-1}$; If $x\geq L$, then define $h(x) = x_n - x_{n-1}$. Here in this subsection, the mesh $\{x_i\}_{i=0}^n$ refers to the computational mesh we use in Section \ref{sec:Computation} for computing the approximate steady state.

As an example, we explain the computation of point values of $u_p$. The idea of accurately computing the point values of $\partial_x u_p, \partial_x^2 u_p$ and $\partial_x^3u_p$ are totally similar. Given any $x\geq0$, we decompose the integration for $u_p(x)$ as 
\[u_p(x) = \f{1}{\pi} \int_0^L \omega_p(y)\log\left|\frac{x-y}{x+y}\right|dy = \f{1}{\pi}\sum_{i=1}^n\int_{x_{i-1}}^{x_i}\omega_p(y)\log\left|\frac{x-y}{x+y}\right|dy.\]
Consider the contribution to the above integral from a single interval $[x_{i-1},x_i]$:
\[\int_{x_{i-1}}^{x_i}\omega_p(y)\log\left|\frac{x-y}{x+y}\right|dy = \int_{x_{i-1}}^{x_i}\omega_p(y)\log|x-y|dy - \int_{x_{i-1}}^{x_i}\omega_p(y)\log|x+y|dy := I_i + J_i.\]
The method of computing this integral depends on the position of $x$ relative to the interval $[x_{i-1},x_i]$.
\begin{enumerate}
\item If $x=0$, then $u_p(x) = 0$ due to odd symmetry.
\item If $|x - x_{i-1}| \leq m\cdot h(x)$ for some constant $m$, we compute the part $I_i$ explicitly using the quintic polynomial representation of $\om_p$ in this interval. Similarly, if $|x + x_{i-1}| \leq m\cdot h(x)$, we compute the part $J_i$ explicitly using the quintic polynomial representation of $\om_p$. Note that the case $x\in[x_{i-1},x_i]$ is included in this branch, and the kernel $\log|(x-y)|$ is integrable on this interval.
\item If $|x - x_{i-1}| > m\cdot h(x)$, directly using the explicit integration formula may cause large round-off error. Instead, we adopt an $8$-point Legendre-Gauss quadrature rule over the interval $[x_{i-1},x_i]$ to compute the part $I_i$. Similarly, if $|x + x_{i-1}| > m\cdot h(x)$, then we compute the part $J_i$ using an $8$-point Legendre-Gauss quadrature rule.
\item In particular, if $x_i/x \leq \delta$ or $x/x_{i-1}\leq \delta$ for some small constant $\delta$, we will instead compute the integral $I_i+J_i$ together using an approximation of the kernel $\log|(x-y)/(x+y)|$. This will help avoid large round-off errors in evaluating the kernel. Define $r(y,x) = y/x$ if $x_i/x \leq \delta$ or $r(y,x) = x/y$ if $x/x_{i-1}\leq \delta$ (so $r$ is always small than $10^{-6}$). We approximate the kernel by a Taylor expansion:
\[\log\left|\frac{x-y}{x+y}\right| = \log\left(\frac{1-r}{1+r}\right) = - 2r - \frac{2}{3}r^3 + O(r^5).\]
Then we approximate the integration over $[x_{i-1}, x_{i}]$ as 
\[\int_{x_{i-1}}^{x_i}\om_p(y)\log\left|\frac{x-y}{x+y}\right|dy \approx \int_{x_{i-1}}^{x_i}\om_p(y)\left(- 2r(y,x) - \frac{2}{3}r(y,x)^3\right)dy,\]
and the right-hand side is then approximately computed using an $8$-point Legendre-Gauss quadrature rule.
\end{enumerate}
In our computation, we choose $m=10$ and $\delta = 10^{-6}$.

The idea of accurately computing the point values of $\partial_x u_p, \partial_x^2 u_p$ and $\partial_x^3u_p$ are totally similar, yet with the following exceptions:
\begin{itemize}

\item For the computation of $\partial_x^2u_p$, we approximate the integral kernel $2x/(x^2-y^2)$ (see the formula in \eqref{eq:compute_u}) by a Taylor expansion only when $x_i/x \leq \delta$: 
\[\frac{2x}{x^2-y^2} = \frac{2}{x}\left(1 + \left(\frac{y}{x}\right)^2 + \left(\frac{y}{x}\right)^4 + O\left(\left(\frac{y}{x}\right)^6\right)\right),\quad y\in [x_{i-1},x_i], \ x_i/x\leq \delta.\]
We also apply  case (1) above to obtain $\partial_x^2u_p(0) = 0$ since it is an odd function of $x$.

\item For the computation of $\partial_xu_p$ and $\partial_x^3u_p$, we \textit{only} use the methods in (2) and (3) above, since $\partial_xu_p, \partial_x^3u_p$ are even functions and the evaluation of the corresponding integral kernels does not induce large round-off errors metioned 
in case (4) above.

\end{itemize}

There are three types of errors in this computation. The first type is the round-off error, the second type is the approximation error of the Gaussian quadrature introduced in Steps (3) and (4), and the third type is the truncation error in the Taylor expansion. With the estimates of the derivatives of $\om_p$, we can establish error estimates of the Gaussian quadrature; see Section \ref{sec:GQ}. The truncation error in the Taylor expansion is estimated in Section \ref{sec:Taylor}. These error bounds will be taken into account in the interval representations of $\bar{u}_p$, $\bar{u}_{p,x}$, $\bar{u}_{p,xx}$ and $\bar{u}_{p,xxx}$. For example, each $\bar u(x)$ will be represented by $[\lfloor\bar u(x)-\epsilon\rfloor_f,\lceil\bar u(x)+\epsilon\rceil_f]$ in any computation using the interval arithmetic, where $\epsilon$ is the error bound for $u(x)$. We remark that we will need the values of $\bar u_p(x)$ at finitely many points only. The same arguments apply to $\bar u_{o,x}(x)$, $\bar u_{p,xx}(x)$ and $\bar u_{p,xxx}$ as well.

\subsection{Computing the Hilbert transform of $\om_b$} \label{sec:compute_ub}

Computing the Hilbert transform of $\om_b$ accurately is more complicated since $\om_b$ captures the far field behavior of the solution, which has very slow decay. This can cause large round off error if we directly use method similar to the one in the previous Section. Recall $\om_b$ from \eqref{eq:explicit_part}. Note that 
\[
H (\f{x}{1 + x^2} ) = \f{-1}{1 + x^2}. 
\]
Using a rescaling argument, we obtain the Hilbert transform of $ \frac{s_\om x}{1 + (r_\om x)^2} $. Hence, we only need to compute $H F_a$ for $\f{x^5}{1 +x^{5+a}}$ and $a = a_{\om}$ \eqref{eq:aw} close to $\f{1}{3}$. To obtain accurate value of $H F_a$, we use a combination of analytic estimates and numerical computation similar to the computation of $H \om_p$. For $x$ in the far field $ x \geq M_1 $, intermediate region $ M_2 \leq x < M_1$ and the near field $x < M_2$, we use different methods. We choose $M_1 = 10^5, M_2 = 4$.

\subsubsection{Far field $x \geq M_1$} 
Denote 
\[
C_a = - \cot \f{\pi a}{2} \f{1}{1-a}, \quad C_u = H( (F_a - \sgn(x) |x|^{-a}) x^2)(0).
\]
For large $x$, since $F_a$ is close to $\sgn(x) |x|^{-a}$, we use 
\[
H  \sgn(x) |x|^{-a}  = C_a (1 -a) |x|^{-a}
\]
to approximate $H F_a$. We focus on $x \geq 0$. We have 
\[
u(F_a)(x) \approx C_a x^{1-a} + \f{C_u}{x} ,  \quad
\pa_x^{k+1} u( F_a)(x) = \pa_x^k H F_a  \approx  C_a \pa_x^k x^{1-a} -  \f{ (-1)^k  (k+1)! C_u}{x^{k+2}}
\]  
for $k =0,1,2$. The error is estimated in Corollary \ref{cor:ub}. 

\subsubsection{   $  x < M_1$}

For $0\leq x < M_1$, we use the computation of $u(F_a)$ as an example. Firstly, we partition the integrals as follows 
\[
u(F_a) = \f{1}{\pi} \int_0^L F_a(y) \log \B| \f{x-y}{x+y} \B| dy 
+ \f{1}{\pi } \int_L^{\inf} F_a(y) \log \B| \f{x-y}{x+y} \B| dy  \teq I + II.
\]

For part $II$, since $y \geq L$ is much larger than $x$, we use a Taylor expansion argument similar to the Step (4) in Section \ref{sec:compute_up} and compute
\[
II \approx - \f{2}{\pi a} L^{-a} x.
\]
For $u_x, u_{xx}, u_{xxx}$, the leading order contributions of integrals in the region $y \geq L$ are obtained similarly and the error estimates are established in Lemma \ref{lem:HF_integ_far}.

The computation of $I$ is similar to that in Section \ref{sec:compute_up}. We partition the integrals as follows
\[
I = \f{1}{\pi} \sum_{i=1}^n  ( \int_{ x_{i-1} }^{x_i} F_a(y) \log |x-y|  d y 
- \int_{x_{i-1} }^{ x_i} F_a( y) \log |x+y| dy ).
\]

For index $i$ that falls in Steps (1), (3), (4) in Section \ref{sec:compute_up}, i.e. $[x_{i-1}, x_i]$ away from $x$, we use the same method as in Section \ref{sec:compute_up} to compute the integral. The only difference is that $F_a$ is an explicit function that can be evaluated directly rather than a quintic spline. 

For index $i$ that falls in Step (2). We focus on $i$ with $|x- x_{i-1}| \leq m h(x)$. The treatment for $i$ with $|x+ x_{i-1}|\leq m h(x)$ is similar. For $i$ with $|x- x_{i-1}| \leq m h(x)$, the region $[x_{i-1}, x_i]$ is close to the singularity. The idea is to approximate $F_a$ by some function $ f$, so that the main part in the integral can be obtained analytically and the error part is small.

We consider two approximations based on the location of $x$.

\paragraph{\bf{Case 1: $ M_2 \leq x < M_1$}}
For $x > M_2 = 4$, $F_a = \f{x^5}{1 + x^{5+a}} \approx x^{-a}$. Since $a_{\om}$ \eqref{eq:aw} is close to $\f{1}{3}$, we approximate $F_a - x^{-1/3}$ on $[2, L ] $ using a quintic spline $S_{a,1}$. Then we get 
\[
\bal
\int_{x_{i-1}}^{x_i} F_a \log |x-y| dy  & = 
\int_{x_{i-1}}^{x_i} S_a  \log |x-y| dy    
 + \int_{x_{i-1}}^{x_i} y^{-1/3} \log |x-y| dy \\
& \quad + \int_{x_{i-1}}^{x_i} 
(F_a - S_a - y^{-1/3}) \log |x-y| dy 
 \teq Q_1 + Q_2 + Q_3. \\
\eal
\]

The computation of $Q_1$ follows from the method in Step (2) in Section \ref{sec:compute_up} and is done analytically. For $Q_2$, using a change of variable $y = z^3$, we yield 
\[
Q_2 = \int_{x_{i-1}^{1/3}}^{ x_i^{1/3}} 3 z  \log | x - z^3| dz ,
\]
which can also be obtained analytically. For example, the formula of the indefinite integral can be obtained using \textit{Mathematica}. We treat $Q_3$ as the error part and defer the estimate to Section \ref{sec:appr_err}. 

\vspace{0.1in}

\paragraph{\bf{Case 2: $ 0 \leq x <  M_2$ }}

In the near field $x < M_2$, since we have much finer mesh, we directly approximate $F_a$ on $[0, L]$ using a quintic spline. We get
\[
\int_{x_{i-1}}^{x_i} F_a \log |x-y| dy   = 
\int_{x_{i-1}}^{x_i} S_a  \log |x-y| dy  + \int_{x_{i-1}}^{x_i} 
(F_a - S_a ) \log |x-y| dy \teq Q_1  + Q_3. \\
\]
The computation of $Q_1$ is similar to the one described above, and $Q_3$ is treated as error and is estimated in Section \ref{sec:appr_err}.

\subsubsection{ Estimate the error part $Q_3$} \label{sec:appr_err}

It remains to estimate the error part $Q_3$. For fix $x$, we combine the estimate of the $Q_3$ integral for all $i$ that falls in Step (2) in Section \ref{sec:compute_up}, i.e. $|x- x_{i-1}| \leq m h(x)$ or $|x + x_{i-1}| \leq m h(x)$. Denote by $P_a = F_a - S_a - y^{-1/3}$ in the case of $ M_2 \leq x < M_1$ and $P_a = F_a - S_a$ for $x \leq M_2$. We extend $P_a$ naturally  to an odd function. Then the total contribution from $P_a$ is given by 
\[
 III =  \sum_{ i \in J_1 } \int_{x_i}^{x_{i+1}} P_a \log |x-y| dy 
  - \sum_{ i \in J_3} \int_{x_i}^{ x_{i+1}} P_a \log |x+y| dy, 
\]
where we have changed the dummy variable $i$ to $i+1$ and 
\[
J_1 = \{ i : |x- x_i| \leq m h(x) \}, \quad J_3 = \{ i : | x + x_i| \leq m h(x) \}. 
\]
It is not difficult to obtain that the indices in $J_1$ are consecutive, and those in $J_3$ are consecutive. Moreover, we have that
\[
SG(x) \teq \cup_{i \in J_1} [x_i, x_{i+1}] \cup 
\cup_{i \in J_3}  (-[x_i, x_{i+1}] )
\]
forms an interval, which can be interpreted as the singular region near $x$. Thus we get
\[
III = \int_{SG(x)} P_a \log |x-y| dy.
\]

This integral and similar integrals for $u_x, u_{xx}, u_{xxx}$ are estimated using Lemma \ref{lem:uxinf2}.

These error bounds will be taken into account in the interval representations of $u(F_a)$, 
$ u(F_a)_x, u(F_a)_{xx}, u(F_a)_{xxx}$, which are similar to those in Section \ref{sec:compute_up}. We refer to Section \ref{sec:compute_up} for more discussion.

\subsection{Rigorous estimates maximum/minimum}\label{sec:maxmin_estimate}
To supplement our analysis in the Main Paper, we need to numerically verify some point-wise inequalities that involve some explicit functions and the steady state solution we have constructed. For such purpose, we will need to rigorously estimate the maximum or minimum of some function $g(x)$ on $[0,L_B]$. In other words, we want to obtain $c_1,c_2$ such that $c_1\leq g(x)\leq c_2$ for $x\in [0,L_B]$, and we want these inequalities to be as sharp as possible. A straightforward way to do so is by constructing two sequences of values $g^{up}=\{g^{up}_i\}_{i=1}^N,g^{low}=\{g^{low}_i\}_{i=1}^N$ such that 
\[g^{up}_i \geq \max_{x\in[x_{i-1},x_i]}g(x)\quad\text{and}\quad g^{low}_i \leq \min_{x\in[x_{i-1},x_i]}g(x).\]
Then we have
\[\min_{i} g^{low}_i \leq g(x) \leq \max_{i} g^{up}_i, \quad x\in [0,L_B].\]
In most cases, we will construct $g^{up}$ and $g^{low}$ from the grid point values of $g$ and an estimate of its first derivative. Let $f^{\max}=\{f^{\max}_i\}_{i=1}^N$ denote the sequence such that $f^{\max}_i = \max\{|f^{up}_i|,|f^{low}_i|\}$. Then if we already have $g_x^{\max}$, we can construct $g^{up}$ and $g^{low}$ as 
\begin{equation}\label{eq:max_estimate}
\begin{split}
g^{up}_i &= \max\{g(x_{i-1}),g(x_i)\} + \frac{(x_i - x_{i-1})}{2}\cdot(g_x^{\max})_i,\\
g^{low}_i &= \min\{g(x_{i-1}),g(x_i)\} - \frac{(x_i - x_{i-1})}{2}\cdot(g_x^{\max})_i.
\end{split}
\end{equation}

We can use this method to construct the piecewise upper bounds and lower bounds for many functions we need. For example, the perturbation part $\om_p$ in our approximate steady state $\bar{\om}$ is constructed to be piecewise quintic polynomial using the standard quintic spline interpolation. Since $\partial_x^5\om_p$ is piecewise constant, we have $(\partial_x^5\om_p)^{up}$ and $(\partial_x^5\om_p)^{low}$ for free from the grid point values of $\om_p,\partial_x\om_p$. Then we can construct $(\partial_x^k\om_p)^{up/low}, k= 0,\dots,4$ recursively using \ref{eq:max_estimate}.

In our verification procedure, we also need to obtain the $up/low$ sequences for $u$ and its derivatives, and we will start from constructing $(\partial_x^3u)^{up/low}$. We defer the details of how we achieve this to Appendix \ref{sec:ux3_estimate}. 

Note that for some explicit functions, we can construct the associated sequences of their piecewise upper bounds and lower bounds more explicitly. For example, for a monotone function $g$, the sequences $g^{up}$ and $g^{low}$ are simply given by the grid point values.

Moreover, we can construct the piecewise upper bounds and lower bounds for more complicated functions. For instance, if we have $f_a^{up/low}$ and $f_b^{up/low}$ for two functions, then we can construct $g^{up/low}$ for $g = f_af_b$ using standard interval arithmetic. In this way, we can estimate the integral of all the functions we need in our computer-aided arguments. 

Sometimes we need to handle the ratio between two functions, which may introduce a removable singularity. Directly applying interval arithmetic to the ratio near a removable singularity can lead to large errors. We hence need to treat this issue carefully. For example, let us explain how to reasonably construct $g^{up/low}$ for a $g(x) = f(x)/x$ such that $f(x)$ has continuous first derivatives and $f(0)=0$. Suppose that we already have $f^{up/low}$ and $f_x^{up/low}$. Then for some small number $\varepsilon>0$, we let 
\[g^{up}_i = \max\left\{\frac{f^{up}_i}{x_{i-1}},\frac{f^{low}_i}{x_{i-1}},\frac{f^{up}_i}{x_i},\frac{f^{low}_i}{x_i}\right\}\quad \text{for each index $i$ such that $x_{i-1}\geq \varepsilon$}.\]
Otherwise, for $x\in [0,\varepsilon)$, we have
\[g(x) = \frac{f(x)}{x} = f_x(\xi(x))\quad \text{for some $\xi(x)\in[0,x)\subset[0,\varepsilon)$}.\] 
Then we choose $g^{up}_i = \max_{x\in[0,\varepsilon]} f_x^{up}$ for every index $i$ such that $x_i\leq \varepsilon$. The parameter $\varepsilon$ needs to be chosen carefully. On the one hand, $\varepsilon$ should be small enough so that the bound $f_x(\xi(x))\leq \max_{\tilde{x}\in[L_B-\varepsilon,L_B]}|f_x(\tilde{x})|$ is sharp for $x\in[L_B-\varepsilon,L_B]$. On the other, the ratio $\varepsilon/h$ must be large enough so that $f^{up}_i/x_{i-1},f^{low}_i/x_{i-1},f^{up}_i/x_i $ and $f^{low}_i/x_i$ are close to each other for $x_{i-1}\geq \varepsilon$. Other types of removable singularities can be handled in a similar way.

\subsection{Rigorous estimates of integrals}\label{sec:integral_estimate} In many of our discussions, we also need to rigorously estimate the integral of some function $g(x)$ on $\mathbb{R}_+$. In particular, we want to obtain $c_1,c_2$ such that $c_1\leq \int_0^{+\infty}g(x)dx\leq c_2$. To do so, we first decompose the integral into two parts:
\[\int_0^{+\infty}g(x)dx = \int_0^{L_B}g(x)dx + \int_{L_B}^{+\infty}g(x)dx.\]
for a sufficiently large number $L_B$. As stated earlier, we choose $L_B$ to be $10^{18}$, which is large enough such that the second part of the integral above can be regarded as an error. 

For the first part, namely the integral of $g$ on $[0,L_B]$, we can also bound it rigorously using the sequences $g^{up}=\{g^{up}_i\}_{i=1}^N$ and $g^{low}=\{g^{low}_i\}_{i=1}^N$ that are constructed in the previous subsection. Evidently, we can bound 
\[
\sum_{i=1}^N (x_i - x_{i-1})g^{low}_i \leq \int_0^{L_B}g(x)dx\leq \sum_{i=1}^N (x_i - x_{i-1})g^{up}_i.
\]

For the second part of the integral, we need to make use of the asymptotic estimate of $g$ from its explicit expression. That is, we bound $|g(x)|\leq a x^{-b}$ for $x\geq L_B$ and for some constant $a>0,b>1$, using the formula of $g$. Then we estimate the tail part of the integral as 
\[\left|\int_{L_B}^{+\infty}g(x)dx\right|\leq \int_{L_B}^{+\infty}|g(x)|dx\leq \frac{a}{b-1}L_B^{-b+1}. \]

Combining these estimates, we obtain 
\[
\sum_{i=1}^N (x_i - x_{i-1})g^{low}_i - \frac{a}{b-1}L_B^{-b+1}\leq \int_0^{+\infty}g(x)dx\leq \sum_{i=1}^N (x_i - x_{i-1})g^{up}_i + \frac{a}{b-1}L_B^{-b+1}.
\]
This will provide the lower and upper bounds that are accurate enough for most of the integrals we need to estimate.

\subsection{Accurate estimates of certain integrals} Our rigorous estimate for integrals in the preceding subsection is only first order accurate. However, this method is not accurate enough (in the relative sense) if the target integral is supposed to be a very small number.

Recall that in Section \ref{sec:NF}, we need to bound some weighted norms of the residuals:
\beq\label{eq:error_est}
 \int_0^\infty  F_{\om}^2 \vp dx ,\quad \int_0^\infty  (x\partial_x F_{\om})^2 \vp dx, \quad \int_0^\infty F_v^2 \psi dx, \quad \text{and}\quad  \int_0^\infty (x\partial_x F_v)^2 \psi dx.
\eeq
The residual functions $F_\om$ and $F_v$ are defined in \eqref{eq:NF}. Note that $F_v$ has the same meaning as $F_\th$ defined in \eqref{eq:NF}; using $F_v$ is more accurate since it is the residual of the equation of $v = \th_x$. We will consistently use $F_v$ for $F_\th$ in what follows. To simplify calculations, we first bound the weight functions by some linear combinations of inverse powers of $x$:
\[
\vp \leq \b_1 x^{-4} + \b_2 x^{-2} + \b_3 x^{-2/3}, 
\quad \psi \leq \b_4 x^{-5} + \b_5 x^{-4} + \b_6 x^{-2} + \b_7 x^{-2/3}.
\]
Then we reduce the estimates of \eqref{eq:error_est} to the estimates of the $L^2(x^{-\alpha })$ norm
\beq\label{eq:error_general}
\int_0^\infty F^2 x^{-\alpha} dx
\eeq
for $\alpha = 5, 4, 2$ or $\f{2}{3}$ and $F = F_\om, F_v, x\partial_xF_\om$ or $x\partial_x F_v$. Note that $F_w, F_v$ are degenerate as $x^3$ around $x=0$ and decay faster than $x^{-1/3}$ in the far field, thus the above integral is well defined for our choices of $\alpha$.

The integral in \eqref{eq:error_general} is small since the residual functions $F_w, F_v$ are smooth and small and decay fast in the far field. To obtain finer estimate of this integral, we first decompose it into two parts
\[\int_0^\infty F^2 x^{-\alpha} dx = \int_0^{L_B} F^2 x^{-\alpha} dx + \int_{L_B}^\infty F^2 x^{-\alpha} dx.\]
For the first part, we will estimate it using mesh-based composite quadrature rules on $[0,L_B]$ with rigorous error bounds. For the second part, we will estimate it analytically using the decay properties of $F$.

Here we explain some methods and error bounds we use for the first part. A typical composite quadrature method that we use is the Trapezoidal rule:
\[T(f,x_0,x_n) \teq \sum_{i=1}^n \frac{x_i-x_{i-1}}{2}(f(x_{i-1}) + f(x_i)), \]
where $ x_0 < x_1 < \cdots < x_n $ are some grid points over the interval $[x_0,x_n]$. For the Trapezoidal rule, the error estimates rely on the following standard error formulas:
\beq\label{eq:Trapezoidal_error}
\bal
\int_a^b f(x) dx - \f{b-a}{2} (f(a) + f(b)) &= - \int_a^b  f_x(x) ( x - \f{a+b}{2}) dx \\
&= \int_a^b f_{xx}(x) \f{ (x-a)(x-b)}{2} dx.
\eal
\eeq
For example, when we do not have access to the second derivative of the integrand $f$, we shall use the first order estimate 
\beq\label{eq:first_order}
\left|\int_{x_0}^{x_n} f dx - T(f,x_0,x_n)\right| \leq \frac{1}{4}\sum_{i=1}^n (f_x)_i^{\max}(x_i-x_{i-1})^2. 
\eeq
The definition of $(f_x)_i^{\max}$ is given in Section \ref{sec:maxmin_estimate}. When we have access to estimates for the second derivative, we can use more delicate estimates.

The following lemma allows us to estimate the weighted $L^2$ norm of a function $f$ with weight functions that is singular at $x=0$.

\begin{lem}\label{lem:Trule1}
For $k > 1$, let $f$ be such that $f(x) = O(x^{k/2+1})$ around $x=0$. Then given a mesh $0= x_0< x_1 < \cdots < x_n=M$, we have
\[
\int_0^M \f{f^2}{x^k} d x - T_h( \f{f^2}{x^k}, 0, M)  dx 
\leq \f{h^2}{k^2 - 1}  \int_0^M \f{ f_{xx}^2}{x^{k-2}} dx
\leq \f{h^2}{k^2 - 1} \sum_{i=1}^n (x_i-x_{i-1})\left(\f{ f_{xx}^2}{x^{k-2}}\right)_i^{max},   
\]
where $h = \max_{i}|x_i-x_{i-1}|$.
\end{lem}

The proof is deferred to Appendix \ref{supp:integral_err}. We will use the above Lemma with $k = 2,4$ and $M \approx 0.1$. We remark that we only require the upper bound of the error rather than the bound of the absolute value of the error, since we will use 
\[
0 \leq \int_0^M \f{f^2}{x^k} dx \leq T_h( \f{f^2}{x^k}, 0, M)  + Error 
\]
later to verify that the integral of $f^2$ is small.

When the weight function is more singular such that the bound in Lemma \ref{lem:Trule1} is not well defined, we use a weaker estimate as follows, whose proof is deferred to Appendix \ref{supp:integral_err}.

\begin{lem}\label{lem:Trule2} For $k \geq 0$, let $f$ be such that $f(x) = O(x^{k/2+1})$ around $x=0$. Then given a mesh $0= x_0< x_1 < \cdots < x_n=M$, we have
\[
\left|\int_0^M \f{f^2}{x^{k+1}} d x - T_h( \f{f^2}{x^{k+1}}, 0, M)  dx\right|
\leq \f{h}{k+1}  \int_0^M \f{ f_x^2}{x^k} dx
\leq \f{h}{k+1} \sum_{i=1}^n (x_i-x_{i-1})\left(\f{f_x^2}{x^k}\right)_i^{max},
\]
where $h = \max_{i}|x_i-x_{i-1}|$.
\end{lem}

The next lemma is also useful for estimating the integral $\|f\|_{L^2(x_0, x_n)}$ when $f$ is sufficiently smooth and the integral $\|f\|_{L^2(x_0, x_n)}$ itself is small. It relies on a different composite quadrature rule based on the piecewise linear interpolation.

\begin{lem}\label{lem:int_high}
For any $x \in [a, b]$, suppose that $\hat f(x)$ is the linear interpolation of $f(x)$
\[
\hat f(x) = f(a) + \f{ f(b) - f(a)}{b- a} (x-a).
\]
Then the error $e(x) = f(x) - \hat f(x)$ satisfies 
\beq\label{eq:int_err2}
| e(x ) | \leq \f{ (b-x)(x-a)}{ \sqrt{3 (b-a)} } ( \int_a^b f_{xx}^2(x)  dx )^{1/2}.
\eeq
In particular, we have the following error estimate for a composite quadrature rule on an interval $[x_0,x_n]$ with grid points $x_0 < x_1 < ... < x_n$:
\begin{align*}
\| f \|_{L^2(x_0,x_n)}
&\leq \| \hat f \|_{L^2(x_0,x_n)} + \| e \|_{L^2(x_0,x_n)}\\
&\leq \| \hat f \|_{L^2(x_0,x_n)} +
\B( \sum_{i=1}^N \f{ (x_i - x_{i-1})^5 }{90} (f^2_{xx})_i^{\max} dx \B)^{1/2},
\end{align*}
where $\hat f$ is the piecewise linear interpolation of $f(x)$ with $\hat f(x_i) = f(x_i)$, and $e(x) = f(x) - \hat f(x)$.
\end{lem}

We defer the proof to Appendix \ref{supp:integral_err}.

The concrete use of the lemmas in this subsection will appear in our Matlab code for computer-aided verification. Here we give one example to illustrate how we use these lemmas to bound the integral
\[\int_0^\infty \frac{F_\om^2}{x^4}dx = \int_0^M \frac{F_\om^2}{x^4}dx + \int_M^{L_B} \frac{F_\om^2}{x^4}dx + \int_{L_B}^\infty \frac{F_\om^2}{x^4}dx := I_1 + I_2 + I_3.\]
One can see that we have further decomposed the integral on $[0,L_{B}]$ into two parts, $I_1$ and $I_2$. In the first interval $[0,M]$, estimating the second derivative of the integrand $F_\om^2/x^4$ is not numerically stable since it has a removable singularity at $x=0$. We thus use Lemma \ref{lem:Trule1} with $f = F_\om$ to estimate the integral $I_1$. Moreover, we will choose $M$ to be a small number so that the mesh sizes $(x_i - x_{i-1})$ in the domain $[0,M]$ are uniformly small. We thus can rely on the smallness of $h = \max_{x_i<= M}|x_i-x_{i-1}|$. 

In the second interval $[M,L_{B}]$, the integrand $F_\om^2/x^4$ is sufficiently smooth and does not have any singularity. Thus we can use Lemma \ref{lem:int_high} with $f = F_\om/x^2$ to estimate the integral $I_2$. The reason that we do not use Lemma \ref{lem:Trule1} for this part is because the mesh size $(x_i - x_{i-1})$ can be very large for $x_i$ close to $L_B$ and hence using a uniform mesh size bound $h$ is not acceptable. However, the smallness of the error bound in Lemma \ref{lem:int_high} relies on the fast decay of $f_{xx}$ in the far field.

For the last part $I_3$, we will use the theoretical decay property of $F_\om$ to directly estimate this integral, owing to the fact that $L_B$ is large enough. The decay estimates of the error $F_{\om}, F_v$ are established in Section \ref{sec:err_decay}.

\subsection{Rigorous estimates of $C_{opt}$.} 
We have argued in Appendix B.5 of the Main Paper that the constant $C_{opt}$ can be bounded as
\beq\label{eq:Copt_bound}
C_{opt}\leq 2^{-1} (\tr (  U_2^T D^{-1} U_1 )^p)^{1/p}.
\eeq
Here we will follow the notations in Appendix B.5 of the Main Paper. To evaluate this bound requires to compute the entries of the matrix $U_2^T D^{-1} U_1$, each of which is of the form $c v_i^T D^{-1} v_j$ for some $1 \leq i, j \leq 9$ and some scalar $c$. Since $v_i$ is the coordinate
of $g_i$ in $\Sigma$ under the orthonormal basis $\{ e_i\}_{i=1}^9$, we have the following isometry 
\[
\la v_i, v_j \ra_{l^2(\R^9)} = \la g_i , g_j \ra_{L^2}, \quad 1 \leq i, j \leq 9.
\]

We discuss the calculation of these entries in three cases. In the first case, $v_i, v_j$ are from different spaces $\Sigma_k, \Sigma_l, k \neq l$. Without loss of generality, we assume $k=1, l \neq 1$. Then from $v_i \in \Sigma_1, v_j \in \Sigma_l$, we have $i \leq 4, j \geq 5$. Since $v_i$ is supported in the first four entries, $v_j$ is supported in the last five entries and $D^{-1}$ is a diagonal matrix, we get $v_i^T D^{-1} v_j = 0$.

Now, suppose that $v_i, v_j$ are from the same space $\Sigma_k$. Without loss of generality, we assume $k=1$. Then we have $i, j \leq 4$. In the second case, we consider at least one of $i, j$ is $1$ and assume $i=1$ without loss of generality.
Since $v_1$ is only supported in the first entry, we get $v_1^T D^{-1} = D_{11}^{-1} v_1^T$. Since $v_i$ is the coordinate of $g_i$ under the ONB $\{ \ee_i\}$, we get 
\[
v_1^T D^{-1} v_j = D_{11}^{-1} v_1^T v_j = D_{11}^{-1} \la g_1, g_j \ra.
\]

In the third case, $2 \leq i, j\leq 4$. Since $v_i, v_j \in \Sigma_1$, $v_i, v_j$ are supported in the first four entries. From $D^{-1}(1:4, 1:4) = \mathrm{ diag} ( d_{11}^{-1}, 1, 1, 1 ), d_{11} = 1 + s_1 ||g_1||_2^2$, we derive
\[
v_i^T D^{-1} v_j = v_i^T v_j + (d_{11}^{-1} - 1) v_{i, 1} v_{j, 1}
=  v_i^T v_j + (d_{11}^{-1} - 1) \la v_i, \ee_1 \ra \la v_j , \ee_1 \ra.
\]
Since $ v_1 = ||g_1 || \ee_1$ and $\la v_k, v_l \ra = \la g_k, g_l\ra$ for $1\leq k,l\leq 4$, we further obtain 
\[
v_i^T D^{-1} v_j = \la v_i, v_j \ra+  (d_{11}^{-1} - 1) \f{\la v_i, v_1 \ra \la v_j , v_1 \ra }{||g_1||_2^2} 
= \la g_i, g_j \ra 
+  (d_{11}^{-1} - 1) \f{ \la g_i, g_1 \ra \la g_j , g_1 \ra}{ ||g_1||_2^2}.
\]

In summary, to compute the entries of $U_2^T D^{-1} U_1$, we only need to compute some of the inner products among $\{g_i\}_{i=1}^9$ (we do not need to compute the coordinate vectors $v_i$ explicitly). This can be done using the method introduced in Section \ref{sec:integral_estimate}. Therefore, each entry of $U_2^T D^{-1} U_1$ is represented by a pair of numbers that bound the entry rigorously from above and below. Once we have the interval estimate of $U_2^T D^{-1} U_1$, we can compute an upper bound of $\tr[(U_2^T D^{-1} U_1)^p]$ stably and rigorously by interval arithmetic, which then gives a rigorous bound on $C_{opt}$ via the inequality \eqref{eq:Copt_bound}. In particular, we choose $p=36$ in our computation, and we can rigorously verify that $C_{opt}<0.9930$.

\begin{remark} Let $M$ be a matrix whose entries are represented by intervals. That is, each $M_{ij}$ is represented by a known interval $[M_{ij}^{low}, M_{ij}^{up}]$ that contains the actual (but unknown) value of $M_{ij}$. Computing the interval representation of the quantity $(\tr[M^p])^{1/p}$ requires interval arithmetics. However, when $p$ is large, performing exact interval arithmetics will induce large deviation due to the large number of matrix multiplication (interval arithmetic is not friendly to matrix multiplication). To avoid such issue, we first invoke the triangle inequality for the Schatten $p$-norm of matrices:
\[(\tr[M^p])^{1/p} \leq (\tr[\Delta M^p])^{1/p} + (\tr[\bar{M}^p])^{1/p}.\]
Here the matrix $\bar{M}$ is the matrix such that 
\[\bar{M}_{ij}\quad \text{is represented by}\quad \left[ \lfloor M_{ij}^{mean}\rfloor_f,\ \lceil M_{ij}^{mean}\rceil_f\right]\quad \text{with}\quad M_{ij}^{mean} = \frac{M_{ij}^{low} + M_{ij}^{up}}{2},\]
and the matrix $\Delta M$ is given by $\Delta M = M - \bar{M}$ via interval arithmetics. Then we compute the quantities $(\tr[\Delta M^p])^{1/p}$ and $(\tr[\bar{M}^p])^{1/p}$ and use them to bound $C_{opt}$. Intuitively, $\bar{M}$ is the mean part of $M$ with much smaller fluctuation for each entry, while $\Delta M$ is the error part that accounts for most of the fluctuation. Empirically, computing $(\tr[\Delta M^p])^{1/p}$ and $(\tr[\bar{M}^p])^{1/p}$ with interval arithmetic will have mush smaller deviation than computing $(\tr[M^p])^{1/p}$.
\end{remark}

\section{Error estimates of the Gaussian quadrature}\label{sec:GQ}

In Section \ref{sec:compute_u}, we discussed how to compute the quantities $\bar{u}, \bar{u}_x, \bar{u}_{xx}$ and different types of errors in the computation. In this section, we estimate the error in the computation of integrals using the numerical Gaussian quadrature. We follow the ideas in \cite{chen2019finite_supp}. Recall the $K-$th order numerical Gaussian quadrature rule 
\[
NGQ_K(f, a, b) = \f{b-a}{2} \sum_{i=1}^K \bar A_j f( \f{b-a}{2} \bar z_i + \f{b+a}{2}),
\]
where $\bar A_j, \bar z_i$ are the numerical approximations of the weights and nodes in the $K-$th order Gaussian quadrature. 
Note that the exact values of these weights and nodes $A_i, z_i$ are related to the roots of some high order polynomial, for which we do not have closed-form formulas. In the numerical realization of the Gaussian quadrature, we need to use their approximations $\bar A_j, \bar z_j $ and $NGQ_K(f; a; b)$ instead. Denote 
\[
\e_k \teq \B| \sum_{j=1}^K \bar A_j \bar z_j^{2k} - \f{2}{2k+1}  \B|, \quad  c_k \teq \sum_{j=1}^K \bar A_j \bar z_j^{2k}.
\]

We remark that if the weights and nodes $\bar A_j, \bar z_j$ are exact, $\e_k$ is zero. For the numerical values $\bar A_j, \bar z_j$, $\e_k$ is very small and of order $10^{-16}$. We have the following error estimate for the numerical Gaussian quadrature. 

\begin{lem}\label{lem:GQ_err}
Suppose that $ K_1 \in \mathbb{Z}_+, K_1 \leq K $ and $F \in C^{2 K_1}$. Denote $t_0 = \f{1}{2}$.
We have 
\[
\bal
\B|\int_0^1 F(t) dt - NQG(F, 0, 1) \B|
&\leq \sum_{1\leq k \leq K_1-1 }\f{ |\pa_t^{2k} F( t_0) | }{ (2k)! \cdot 2^{2k+1}} \e_k  \\
 & \qquad + \max_{t \in (0,1)} \f{ |\pa_t^{2K_1} F_(t) | }{ (2 K_1) !} \cdot \B( \f{1}{ (2 K_1 + 1) 2^{2 K_1}}  
+ \f{c_{K_1}}{ 2^{ 2 K_1 + 1}} \B)
\teq E_{GQ}(F, K_1).
\eal
\]
\end{lem}

The case $K_1 = K$ is proved in the Supplementary material of \cite[Section 13]{chen2019finite_supp}, which can be found in its arXiv version \cite{chen2019finite_supp}.
The proof of other cases $K_1=1, 2, .., K$ is exactly the same. 

In view of the above Lemma, in order to control the error, it remains to control the derivatives of $F$.

\subsection{ Estimates of $\pa_t^{2k} F$}

In this Seciton, we apply Lemma \ref{lem:GQ_err} to estimate the error in computing the velocity. Denote 
\[ 
X_j = [x_j, x_{j+1}] , \quad h_j = x_{j+1} - x_j.
\]

Recall the definition of $h(x)$ from the begining of Section \ref{sec:compute_up}
\[
h(x) =  h_j  \quad \forall x \in [x_j, x_{j+1} ), \quad h(x) = h_n \ \forall x \geq L. 
\]

We fix $x$ in the following discussion. Denote $J_0 = \{0, 1,2, ..., n-1\}$. We define the region of the integral, where the integral is computed using the Taylor expansion or the Gaussian quadrature according to the discussion in Section \ref{sec:compute_up}. We define the index $J_{TL}$ associated to the Taylor expansion part
\beq\label{eq:Taylor}
\bal
J_{TL}(x; u)  &\teq \{ j :  x_{j+1}  \leq \d x \mathrm{ \ or \ }  x \leq  \d x_j \} ,   \quad  J_{TL}(x; u_{xx})  \teq  \{ j : x_{j+1}  \leq \d x  \} ,  \\
 J_{TL}(x ; u_x) &=  J_{TL}(x; u_{xxx} ) = \emptyset,  \quad \d = 10^{-6}. 
\eal
\eeq
The empty sets $J_{TL}(x ; u_x) ,  J_{TL}(x; u_{xxx} )$ mean that we do not apply the Taylor expansion argument when we compute $u_x, u_{xxx}$.

We define the index associated to the singular region $J_1, J_3$ and the nonsingular region $J_2, J_4$ as follows
\beq\label{eq:J1234}
\bal
J_1(x) \teq \{ j : |x - x_j| \leq m h(x) \}
, \quad J_2(x) \teq J_0 \bsh ( J_1(x) \cup J_{TL} ) \\
J_3(x) \teq \{ j : |x + x_j| \leq m h(x) \}, \quad J_4(x) \teq J_0 \bsh ( J_3(x)  \cup J_{TL} ).\\
\eal
\eeq

Suppose that $x \in [x_i, x_{i+1} )$. We get $ h(x) = h_i$. We define the subsets and supersets of the above sets as follows, which only depend on the interval $i$ rather than $x$
\beq\label{eq:super_JTL}
\bal
&J_1^l  \teq \{ j :\max( |x_i - x_j|, |x_{i+1} - x_j|  ) \leq m h_i \}, 
\quad J_1^u \teq \{ j :\min( |x_i - x_j|, |x_{i+1} - x_j|  ) \leq m h_i \},  \\
&J_3^l  \teq \{ j :\max( |x_i + x_j|, |x_{i+1} + x_j|  ) \leq m h_i \}, 
\quad J_3^u \teq \{ j :\min( |x_i + x_j|, |x_{i+1} + x_j|  ) \leq m h_i \},  \\
&J^l_{TL}(u)  \teq \{ j : x_{j+1 } \leq \d x_i \mathrm{ \ or \ } x_{i+1} \leq \d x_j \}, \quad 
J^u_{TL}(u_{xx})  \teq \{ j : x_{j+1 } \leq \d x_{i+1} \mathrm{ \ or \ } x_i \leq \d x_j \} .
\eal
\eeq
Since $m = 10$ and $m h_i > h_i = x_{i+1} - x_i$, it is not difficult to verify that 
\[
J_i^l \subset J_i \subset J_i^u , \quad J^l_{TL}(x; u) \subset J_{TL}(x ; u) \subset J_{TL}^u(x; u)
\]
for $i=1,3$. Similarly, we define the subset and superset of $J_{TL}(x ; u_{xx})$. Thus, we can define the superset of $J_2, J_4$ as follows 
\beq\label{eq:super_J24}
J_2^u \teq J_0 \bsh J_1^l  , \quad J_4^u \teq J_0 \bsh J_3^l ,
\quad J_2 \subset J_2^u, \quad J_4 \subset J_4^u.
\eeq

We introduce these supersets and subsets so that we can estimate the error of $ \pa_x^k u(x)$ for all $x \in X_j$ simultaneous.

Recall that in $[x_j, x_{j+1}] , j \in J_1 \cup J_3$, we use analytic integration to compute the integrals \eqref{eq:compute_u}. In $[x_j, x_{j+1}], j \in  J_2 \cup J_4$, we use the numerical Gaussian quadrature to compute them. See Section \ref{sec:compute_up}. We need to estimate the error. Suppose that $f$ is some odd function. We focus on how to estimate the error in the computation of $u_x = H f$. Other terms $u(f), u_{xx}(f), u_{xxx}(f)$ can be estimated similarly. The integrals related to $J_2, J_4$ are given by 
\beq\label{eq:err_step3}
I = \f{1}{\pi} \sum_{ j \in J_2} \int_{x_j}^{x_{j+1}} \f{ f(y)}{ x- y} dy
- \f{1}{\pi} \sum_{j \in J_4} \int_{x_j}^{x_{j+1}} \f{ f(y)}{ x+ y}  dy  
\teq I_1 + I_2. 
\eeq

We focus on $I_1$ and fix $x$. Using a change of variable $y = x_j + t h_j$, we yield
\[
I_1 = \f{1}{\pi} \sum_{j \in J_2} \int_0^1 \f{ f( x_j + t h_j)}{ x - x_j -  t h_j } h_j dt 
= \f{1}{\pi} \sum_{j \in J_2} \int_0^1 \f{ f( x_j + t h_j)}{ \f{ x - x_j}{h_j} -  t  }  dt .
\]
Denote 
\beq\label{eq:def_dj}
\bal
z_j^-(x) & \teq \f{ x - x_j}{ h_j}, \ z_j^+(x) \teq \f{ x+ x_j}{h_j} , \\ 
d_-(t.j) &\teq z_j^-(x) - t, \  d_-^l(j)  \teq \min_{t \in [0,1]} |d_-(t, j)|,  \
d_-^u(j) \teq \max_{t \in [0,1]} |d_-(t, j)|,   \\
d_+(t.j) &\teq z_j^+(x) + t, \  d^l_+(j)  \teq \min_{t \in [0,1]} |d_+(t, j)|,  \
d_+^u(j) \teq \max_{t \in [0,1]} |d_+(t, j)|,   \\
F_j  & \teq \f{ f( x_j + t h_j)}{ z_j^-(x) -  t  }, \
 G_j \teq \f{ f( x_j + t h_j)}{ z_j^+(x) +  t  }, \ I_{1,j} \teq \int_0^1 F_j(t) dt,
\ I_{2, j} \teq \int_0^1 G_j(t) dt,
 \eal
\eeq
where $l, u$ means \textit{low} and $\textit{up}$ and we have dropped $x$ in the $d(t,j), d(j)$ terms to simplify the notations. By definition, we have 
\[
I_1 = \f{1}{\pi} \sum_{j \in J_2} I_{1, j} ,  \quad
I_2 = \f{1}{\pi} \sum_{j \in J_4} I_{2, j}.
\]

Next, we apply Lemma \ref{lem:GQ_err} to estimate the error of the numerical Gaussian quadrature. We need to estimate $\pa_t^{2k} F_j$. Note that for $j \in J_2$, the interval $[x_j, x_{j+1}]$ is away from $x$ and we have $d^l_-(j) > 0$. For any $k \in Z_+$, a direct calculation yields
\[
\bal
\pa_t^k  d_-(t, j)^{-1} &= k! \cdot  (z_j^-(x) - t)^{-k-1} = k! \cdot d_-(t, j)^{-k-1}, \quad 
| \pa_t^k  d_-(t, j)^{-1}| \leq k! \cdot d_-^l(j)^{-k-1}, \\
| \log |d_-(t,j)| | &\leq \max( | \log d_-^u(j)|, |\log d_-^l(j)| ),
\quad | \pa_t^{k} \log |d_-(t,j)|  | \leq (k-1)! d_-^l(j)^{-k}.
\eal
\] 
We remark that the quantity $d_-^u(j)$ is only used to bound the logarithm term.

Using the Leibniz rule, for $ t \in [0,1]$, we derive 
\beq\label{eq:LeibF}
\bal
|\pa_t^{2 k} F_j(t)| &= \B| \sum_{0\leq l \leq 2k} \binom{2 k}{l}  \pa_t^l ( f (x_j + t h_j))\cdot (2 k- l)!  \cdot d_-(t,j)^{- (2k-l + 1)} \B|  \\
& \leq 
\sum_{0\leq l \leq 2k} \binom{2k}{l}
h_j^l || \pa_x^l f ||_{L^{\inf} (X_j)}
\cdot (2 k- l)!  \cdot  d_-^l(j)^{-(2k-l+1)} ,
\eal
\eeq
where we have used $ x_j + t h_j \in [x_j, x_{j+1}] = X_j$.

Similarly, the functions related to the integrals of $u, u_{xx}, u_{xxx}$ are $F^{\al}, F^{\b}, F^{\g}$ given below respectively
\[
F_j^{\al}(t) \teq h_j f(x_j + t h_j) \log |z_j^-(x) - t|,
\quad F_j^{\b}(t) \teq  \f{ (f_x)(x_j + t h_j) }{ z_j^-(x) -t}, 
\quad F_j^{\g}(t) \teq \f{ (f_{xx})(x_j + t h_j) }{ z_j^-(x) -t}.
\]

The derivatives bound for these functions are given by 
\beq\label{eq:LeibF2}
\bal
|\pa_t^{2k} F_j^{\al}(t)|
&\leq \sum_{ l < 2k} 
\binom{2k}{ l}
h_j^{l+1} || \pa_x^l f ||_{L^{\inf} (X_j)}
\cdot (2 k- l-1)!  \cdot  d_-^l(j)^{-(2k-l  )}  \\
& \quad + h_j^{2k+1} || \pa_x^{2k} f ||_{L^{\inf}(X_j)}
\max( | \log d_-^u(j)|, |\log d_-^l(j)| ),   \\
|\pa_t^{2k} F_j^{\b}(t)|
&\leq \sum_{ l \leq 2k} 
\binom{2k}{l}
h_j^{l} || \pa_x^{l +1} f ||_{L^{\inf} (X_j)}
\cdot (2 k- l)!  \cdot  d_-^l(j)^{-(2k-l + 1)} , \\
|\pa_t^{2k} F_j^{\b}(t)|
&\leq \sum_{ l \leq 2k} 
\binom{2k}{l}
h_j^{l } || \pa_x^{l +2} f ||_{L^{\inf} (X_j)}
\cdot (2 k- l )!  \cdot  d_-^l(j)^{-(2k-l + 1)} . \\
\eal
\eeq

The estimate for the $G_j$ part, or related to the kernel $\log(x+y), \f{1}{x+y}$ are estimated similarly.

To simplify the above summation, we denote 
\[
Pow( a_l, a_u, k) = \max( a_l^{k}, a_u^{k} ) \quad \forall k \leq -1, \quad  Pow( a_l, a_u, 0)  = \max( |\log |a_l||, |\log|a_u| |) ,
\]
and introduce 
\beq\label{eq:err_WD1}
WD( l, i, k ) 
\teq \sum_{j \in J_2^u} || \pa_x^l f||_{L^{\inf}(X_j)} h_j^i   \cdot 
Pow( d_-^l(j), d_-^u(j), -k)
+ \sum_{j \in J_4^u} || \pa_x^l f||_{L^{\inf}(X_j)} h_j^i  \cdot
Pow( d_+^l(j), d_+^u(j), -k). 
\eeq

Using the above notation and \eqref{eq:LeibF}, we obtain 
\beq\label{eq:err_DFbd}
  \sum_{ j \in J_2^u} | \pa_t^{2k} F_j(t)|
+ \sum_{ j \in J_4^u} | \pa_t^{2k} G_j(t)|
\leq \sum_{ l \leq 2k} \binom{2k}{l}   (2k-l)! \cdot WD(l, l, 2k -l + 1)
\teq DF_{bd}(k ;  u_x). 
\eeq

Applying Lemma \ref{lem:GQ_err} and the above estimates, we obtain that the error in computing the integral $I$ \eqref{eq:err_step3} in $H f(x)$ for $x \in [x_i, x_{i+1}]$, or equivalently the error in Step (3) in Section \ref{sec:compute_up}, is given by 
\beq\label{eq:LeibF3}
\bal
Error & \leq \f{1}{\pi} \sum_{j \in J_2} E_{GQ}( F_j, K_1)
+  \f{1}{\pi} \sum_{j \in J_4} E_{GQ}( G_j, K_1)
 \leq 
\f{1}{\pi} \sum_{j \in J_2^u} E_{GQ}( F_j, K_1)
+  \f{1}{\pi} \sum_{j \in J^u_4} E_{GQ}( G_j, K_1), \\
&  \leq
\sum_{1\leq k \leq K_1-1 }\f{ \e_k }{ (2k)! \cdot 2^{2k+1}}  
DF_{bd}(k ; u_x)
+  \f{ 1 }{ (2 K_1) !} \cdot \B( \f{1}{ (2 K_1 + 1) 2^{2 K_1}}  
+ \f{c_{K_1}}{ 2^{ 2 K_1 + 1}} \B) DF_{bd}( K_1; u_x),
\eal
\eeq
where $1\leq K_1 \leq K = 8$ and we have used \eqref{eq:super_J24} to obtain the second inequality. Note that the upper bound depends on the interval $i$, where $x$ lies in. Thus, we can apply the above estimate to obtain the error for all $x \in [x_i,x_{i+1}]$.

The error estimates for $u, u_{xx}, u_{xxx}$ are completely similar.

\subsubsection{ Improved estimates for the error related to $u_x$ and $u_{xxx}$ and $j \in J_2 \cap J_4$}\label{sec:ux_improv}

In practice, the above error estimate is not sharp enough for the integrals in $u_x$ and $u_{xxx}$ with large $x$. For $j \in J_2 \cap J_4$, we combine the error estimate of $F_j$ and $G_j$ 
\[
\int_0^1 F_j(t) - G_j(t) dt = \int_0^1 f( x_j + t h_j ) ( \f{1}{ z_j^-(x) - t  } - \f{1}{ z_j^+(x) + t}) dt.
\]
Note that the above quantity appears in \eqref{eq:err_step3}. For $x$ much larger than $x_j$, we have $ z_j^-(x) - t = \f{x - x_j}{h_j} - t \approx \f{x+ x_j}{h_j} + t$ and there is a cancellation in the kernel. We explore this cancellation to improve the error estimate. For even $ k \geq 1$ , we get 
\[
S = \pa_t^k ( \f{1}{ z_j^-(x) - t  } - \f{1}{ z_j^+(x) + t})
= k! ( \f{1}{ (z_j^-(x) - t)^{k+1}} -  \f{1}{ (z_j^+(x) + t )^{k+1} }). 
\]

Denote $A = z_j^-(x) - t$ and $B =z_j^+(x) + t  $. By definition \eqref{eq:def_dj} we have $|A| \leq B$. Using 
 \[
 |A^{k+1} - B^{k+1}|
=| (A-B)(A^k + A^{k-1} B + ... + B^k)| \leq  (k+1) |A-B|  |B|^k ,\]
we obtain 
\[
| \f{1}{A^{k+1}} - \f{1}{B^{k+1}} |  = \f{ |B^{k+1}  - A^{k+1} |}{ |A^{k+1} B^{k+1}|}
\leq \f{ (k+1)|A-B| B^k}{ |A^{k+1} B^{k+1}|}
\leq \f{ (k+1) |A-B|}{ | A^{k+1} B | } .
\]
From \eqref{eq:def_dj}, we have $B- A = \f{2x_j}{h_j} + 2 t \leq 2( \f{x_j}{h_j} + 1)$ and $|A| \geq d^l_-(j), |B| \geq d^l_+(j)$. It follows 
\[
|S| \leq  k! \min( 2( \f{x_j}{h_j} + 1) \f{k+1}{ d_-^l(j)^{k+1} d_+^l(j) }, d_-^l(j)^{-k-1} + d_+^l(j)^{-k-1}  ),
\]
where we have applied the trivial bound, e.g. $ |(z_j^- - t)^{-k-1}| \leq d_-^l(j)^{-k-1}$, to obtain the second bound in the minimum. We will further estimate the $d(j)$ in Section \ref{sec:est_dj}.

Suppose that $x \in [x_i, x_{i+1}]= X_i$. In general, the set $J_2 \cap J_4$ may be different for different $x \in X_i$. In order to perform the improved estimate simultaneously for $x \in X_i$, 
we find a subset of $J_2 \cap J_4$ that only depends on $i$. Note that for $u_x, u_{xxx}$, we do not apply Taylor expansion and have $J_{TL}(x; u_x) = J_{TL}(x; u_{xxx} ) = \emptyset $, see \eqref{eq:Taylor}. Thus, the sets in \eqref{eq:Taylor}-\eqref{eq:super_J24} satisfy 
\[
J_5^l \teq  J_0 \bsh (J_1^u \cup J_3^u)
\subset J_0 \bsh (J_1 \cup J_3)
= (J_0 \bsh J_1 ) \cap (J_0 \bsh J_3) = J_2 \cap J_4. 
\]
When we estimate the error for $u_x, u_{xxx}$, we apply the above improved estimate for the index $j \in J_5^l \subset J_2 \cap J_4$, and then obtain an estimate on 
\[
\sum_{ j \in J_5^l} | \pa_t^{2k} F_j - \pa_t^{2k} G_j|
+ \sum_{ j \in J_2^u \bsh J_5^l} |\pa_t^{2k} F_j |
+   \sum_{ j \in J_4^u \bsh J_5^l} |\pa_t^{2k} F_j |,
\]
which is similar to \eqref{eq:err_DFbd}. In this case, we modify the definition of $WD$ in \eqref{eq:err_WD1} correspondingly.

\subsubsection{Two types of functions $f$}

We apply the Gaussian quadrature to compute the Hilbert transform of two types of $f$:  $f=  \f{x^5}{1 +x^{5+a}}$ with some $a \in (0,1)$ and quintic spline on $[0, L]$. 

For $f$ being a quintic spline, the derivative bound of $f$ in each interval $X_j$ will be estimated using the method discussed in Section \ref{sec:maxmin_estimate}. We remark that in the estimates \eqref{eq:LeibF}, \eqref{eq:LeibF2}, we only require the derivatives bound of $f$ in the interval $X_j$. Since $f$ is piecewise fifth order polynomial, we have $\pa_x^k f = 0$ in $X_j$ for $k \geq 6$.

For $f = \f{x^5}{1 + x^{5+a}}$ with $a \in (0,1)$, we use the estimate and the method established in Section \ref{sec:deri_wb} to estimate $|| \pa_x^k f||_{L^{\inf}(X_j)}$.

\subsubsection{Estimate of $d(t, j)$ }\label{sec:est_dj}

Suppose that $x \in [x_i, x_{i+1}]$. The mesh satisfies that $ i-1, i, i+1 \notin J_2$, i.e. the nonsingular region is at least two intervals away from $[x_i, x_{i+1}]$. Thus, for $j \in J_2$, we have either $j \geq i+2$ or $j \leq i-2$. Moreover, we have $|x - x_j | \geq m h(x) = m h_i$.
If $j \geq i+2$, then we have $x_j > x_{i+1} \geq x$ and
\[
\bal
d_-^l(j) &= \min_{t \in[0,1]} |d_-(t, j )| = \min_{t\in[0,1]} | \f{x - x_j}{h_j} - t | 
= \min_{t\in[0,1]}\f{x_j - x}{h_j} + t \geq \max( \f{x_j - x_{i+1}}{ h_j}, \f{ m h_i}{h_j}  ), \\
d_-^u(j) & = \max_{t \in [0,1]} |d_-(t, j) | =\max_{t \in [0,1]} \f{x_j - x}{h_j} + t \leq  \f{x_j -x_i}{h_j} + 1, 
\eal
\]
where we have used $ x_j - x = |x-x_j| \geq m h_i$ to get the second bound in the maximum.

If $j \leq i-2$, we get $x \geq x_{j+1} = x_j + h_j \geq x_j + t h_j$. Hence, we get 
\[
\bal
d_-^l(j) &=  \min_{t\in[0,1]} |d_-(t, j)| = \min_{t\in[0,1]} \f{ x - x_j}{h_j} - t \geq 
\f{ x - x_j}{h_j} - 1
\geq \max( \f{x_i - x_j}{ h_j} - 1,  \f{m h_i}{h_j}- 1 ), \\
d_-^u(j) &  = \max_{t \in [0,1]}  |d_-(t, j)|   = \max_{t \in [0,1]} \f{ x - x_j}{h_j} - t
\leq \f{x_{i+1} - x_j}{h_j} ,
\eal
\]
where we have used $x -x_j = |x-x_j| \geq m h_i$. Therefore, we obtain the estimate for $d_-^l(j)$ and $d_-^u(j)$.

Similarly, we have the following estimate for $ z_j^+(x ) + t, x \geq 0$ and $ j \in J_4(x)$
\beq\label{eq:dj_plus}
\bal
d_+^l(j) &=\min_{t \in [0,1]} | z_j^+(x) + t| = \min_{t \in [0,1]} ( \f{x + x_j}{h_j} + t ) 
\geq \max( \f{ m h_i}{h_j} , \f{x_i + x_j}{h_j}  ) ,\\
d_+^u(j) & =  \max_{t \in [0,1]} | z_j^+(x) + t|    =
\max_{t \in [0,1]} ( \f{x + x_j}{h_j} + t ) 
\leq  \f{ x_{i+1} + x_j}{ h_j} + 1.
\eal
\eeq

With the above estimate on $\pa_t^{2k} F_j$ and the upper and lower bound of $d_{\pm}(t, j)$, 
we estimate $\pa_t^{2k} F_j$ using \eqref{eq:LeibF} and further estimate the integral error using Lemma \ref{lem:GQ_err}.

\subsection{ Error estimate for the integrals using Taylor expansion}\label{sec:Taylor}
Recall in Section \ref{sec:compute_up} that in the computation of $u$, for $ j \in J_{TL}(x; u)$, we use 
\[
\int_{x_j}^{x_{j+1}} f(y) \log \B| \f{x-y}{x+y} \B| dy =  \int_{x_j}^{x_{j+1}} f(y) ( -2 r(y, x) - \f{2}{3} r(y,x)^3 + \cR(x, y)) dy  \teq I + Error,
\]
where $r(y,x) = \f{y}{x}$ if $x_{j+1} \leq x \d$ and $r(y, x) = \f{x}{y}$ if $ x \leq \d x_j$, 
$\cR$ is the remainder in the expansion of the kernel and the error is the integral associated to $\cR$. Note that we drop $\pi^{-1}$ in the integral to simplify the discussion. This factor will be included in our computation.

For $I$, using a change of variable $ y = x_j + t h_j$, we get 
\[
I = \int_0^1 f(x_j + t h_j) ( - 2 r( y , x ) - \f{2}{3} r( y, x)^3  ) h_j d t
= \int_0^1 F_j(t) dt , 
 \]
 where 
 \[
F_j(t) = h_j f(x_j + t h_j) ( - 2 r( y , x ) - \f{2}{3} r( y , x)^3  )  , \quad y = x_j + t h_j.
 \]

To compute $I$, we use the numerical Gaussian quadrature. Similar to the previous discussion, to estimate the error, we need the estimate of $\pa_t^{2k} F_j$. 

Suppose that $x \in [x_i, x_{i+1}]$. If $ x_{j+ 1} \leq \d x$, we get $r(y, x) = \f{y}{x}$. Using the Leibniz rule, we yield 
\beq\label{eq:TL_dFj1}
\bal
&| \pa_t^{2k} F_j | 
\leq \sum_{l \leq 2k} \binom{2k}{l}   h_j \B|  \pa_t^l f( x_j + t h_j) \pa_t^{2k-l}  
( 2 \f{y}{x} + \f{2}{3} \f{y^3}{x^3}) \B| \\
\leq & \sum_{l \leq 2k} \binom{2k}{l}  h_j^{2k+1} || \pa_x^l f||_{L^{\inf}(X_j)}
(  2 x_i^{-1} D(x_{j+1}, 1, 2k -l) + \f{2}{3} x_i^{-3} D(x_{ j+1}, 3, 2k-l)),
\eal
\eeq
where the function $D(y,k,l)$ is defined as
\[
D(y, k,l ) = \pa_y^l y^k =  \one_{l \leq k} \f{k!}{ (k-l)!} y^{k-l} 
\]
and we have used the fact that $D(y,k,l)$ is an increasing function of $y$ to get $\pa_y^l y^k \leq D(x_{j+1}, k, l)$.

If $ x \leq \d x_j$, we get $r(y, x) = \f{x}{y}$. Similarly, we have 
\beq\label{eq:TL_dFj2}
\bal
&| \pa_t^{2k} F_j | 
\leq \sum_{l \leq 2k} \binom{2k}{l}   h_j \B|  \pa_t^l f( x_j + t h_j) \pa_t^{2k-l}  
( 2 \f{x}{y} + \f{2}{3} \f{x^3}{y^3}) \B| \\
\leq & \sum_{l \leq 2k} \binom{2k}{l}  h_j^{2k+1} || \pa_x^l f||_{L^{\inf}(X_j)}
(  2 x_{i+1}   (2k-l)! x_j^{- (2k-l+1)}
+ \f{2}{3} x_{i+1}^3 
 \f{ (2k-l + 2)!}{2} x_j^{-(2k-l+3)} ,
\eal
\eeq
where we have used $y = x_j + t h_j$ and 
\[
| \pa_y^{2k-l} y^{-1}| = (2k-l)! y^{-(2k-l+1)} 
\leq  (2k-l)! x_j^{- (2k-l+1)}, \quad
 | \pa_y^{2k-l} y^{-3}| \leq \f{ (2k-l + 2)!}{2} x_j^{-(2k-l+3)} . 
\]

\subsubsection{The remainder term $\cR$}

Using the expansion 
\[
\log \B| \f{x-y}{x+y} \B|
 = -2 \sum_{k\geq 1} \f{1}{2k-1} r(y, x)^{2k-1}
 = - 2 r(y, x) - \f{2}{3} r(y, x)^3 + \cR,
\]
and $|r(y, x)| \leq \d =10^{-6}$, we get 
\[
| \cR | \leq \f{2}{5} \sum_{k\geq 3} r(y,x)^{2k-1}
= \f{2}{5} \f{r(y,x)^5}{ 1 - r(y,x)^2}
\leq \f{2}{5} \f{r(y,x)^5}{ 1 - \d^2} 
\leq 0.41 \min( \d^5, r(y,x)^5),
\]
where we have used $\f{2}{5 (1-\d^2)} = 0.4 (1-10^{-12})^{-1} \leq 0.41$.

If $ j \in J_{TL}(x; u)$ and $ x_{j+1} \leq \d x$, we further use $r(y,x) = \f{y}{x}
\leq \f{x_{j+1}}{x_i}$ to obtain  
\[
| Error | = \int_{x_j}^{x_{j+1}} |\cR |  |f(y)| dy 
\leq 0.41 || f||_{L^{\inf}(X_j)} h_j   \min( \f{x^5_{j+1}}{ x^5_i} , \d^5 ).
\]
Similarly, if $x \leq \d x_j$, we get $r(y, x) = \f{x}{y}$ and 
\[
| Error | = \int_{x_j}^{x_{j+1}} |\cR |  |f(y)| dy 
\leq 0.41 || f||_{L^{\inf}(X_j)} h_j   \min( \f{x^5_{i+1} }{ x^5_j} , \d^5 ).
\]

\subsubsection{Estimate related to $u_{xx}$ }
For $u_{xx}$, the integral in $[x_j, x_{j+1}]$ with $j \in J_{TL}(x; u_{xx})$ is given by 
\[
\f{1}{\pi} \int_{x_j}^{x_{j+1}} f_x(y) 
( \f{1}{x-y} + \f{1}{x+y} ) dy 
= \f{1}{\pi} \int_{x_j}^{x_{j+1}} f_x(y) ( \f{2}{x} (1 + \f{y^2}{x^2} + \f{y^4}{x^4}) + \cR(x,y)) dy \teq I + Error,
\]
where $\cR$ is the remainder in the expansion of the kernel. Note that by definition \eqref{eq:Taylor}, for each $j \in J_{TL}(x; u_{xx})$, we have $x_{j+1} \leq \d x$. We apply a similar argument to compute the main term $I$ and estimate the integral error. To estimate the remainder term, we use 
\[
\B| \f{1}{x-y} + \f{1}{x+y} -  \f{2}{x} (1 + \f{y^2}{x^2} + \f{y^4}{x^4}) \B|
=  \f{2}{x} \sum_{k \geq 3} \f{y^{2k}}{ x^{2k}}  
\leq \f{2}{x} 1.01 \cdot \f{ y^6}{x^6}
\leq \f{2.02}{x_i}  \min( \f{x_{j+1}^6}{x_i^6}, \d^6)
\]
for $y \leq x_{j+1} \leq \d x$ and $x \in [x_i, x_{i+1}]$ to obtain 
\[
|Error | = \int_{x_j}^{x_{j+1}} |\cR(x, y)  f_x(y) | dy
 \leq  \f{2.02}{\pi} \f{h_j}{x_i}  \min( \f{x_{j+1}^6}{x_i^6}, \d^6) || f_x||_{L^{\inf}(X_j)}.
\]

\subsubsection{ Superset of $J_{TL}(x)$}

In order to apply the above error estimate simultaneously for $x \in [x_i, x_{i+1}]$, we sum the above error estimates over a superset of $J_{TL}$. In particular, we introduce the following weighted function similar to $WD$ in \eqref{eq:err_WD1}
\[
TL(u, i, k) =  \sum_{j \in J^u_TL(u)} || \pa_t^i (h_j  f(x_j + t h_j) ||_{L^{\inf}[0,1]} 
|| \pa_t^{k} (2 r(x,y) + 2 r^3(x,y) /3 ) ||_{L^{\inf}[0,1]} ,
\]
where $\pa_t^{k} (2 r(x,y) + 2 r^3(x,y) /3 )$ is estimated in \eqref{eq:TL_dFj1}, \eqref{eq:TL_dFj2}. Similarly, we define $TL(u_{xx}, i, k)$.
The idea and the remaining estimates are similar to those in \eqref{eq:LeibF3}. We refer more discussion to the paragraphs before Section \ref{sec:ux_improv}.

\section{Estimates of the analytic approximate profile }\label{sec:Fa}

In this section, we estimate the derivatives of the analytic approximate profile and its Hilbert transform.

\subsection{ Estimate the derivatives of the analytic approximate profile}\label{sec:deri_wb}
Recall that we use 
\[
F_a = \f{x^5}{ 1 + |x|^{5 +a}}
\]
for $a = a_{\om} \in (0,1) $ or $a = a_v \in (0,1)$ to approximate the far field of the approximate steady state. See \eqref{eq:solution_decomposition}, \eqref{eq:explicit_part} and \eqref{eq:aw}. In this Section, we estimate $\pa_x^i F_a$ for $ 0 \leq i \leq 2 K +2$ so that we can apply Lemma \ref{lem:GQ_err} for error estimates. Though we can comptute these derivatives explicitly using symbolic calculation, the formulas are very complicated and lengthy, which makes it very difficult to estimate them. In the following Lemma, we explore the structure of the derivatives of $F_a$ that enables us to estimate them effectively.

\begin{lem}\label{lem:deri_Fa}
Denote $g(x) = \f{1}{1+x}$.  For $k \in \mathbb{Z}$ and $k \geq 0$, we have 
\beq\label{eq:Fa_recur0}
\pa_x^k F_a = \sum_{ j = 0}^k  (\pa^j g )(x^{5 + a}) \cdot x^{ j a + 5 - k} P_{j,k}(x) ,
\eeq
where $P_{j,k}, 0 \leq j \leq k$ is a polynomial in $x$ satisfying 
\beq\label{eq:Fa_recur}
P_{j,k} = x^5 (5 + a) P_{j-1, k-1} + ( j a + 6 - k) P_{j, k-1} + x \pa_x P_{j, k-1}
\eeq
for $ 0\leq j \leq k , k \geq 1$, and we adopt the convention $P_{j,k}=0$ if $j <0$ or $j > k$. Moreover, $P_{0,k}$ satisfies 
\[
P_{0, k} = \f{5 !}{ (5-k)!} \mathrm{ \ for \  } k \leq 5 , \quad P_{0,k } = 0 , \mathrm{ \ for \  } k > 5,
\]
and $P_{j,k}(x)$ is a multiple of $x^5$ for $ 1 \leq j \leq k$, i.e. $P_{j,k}(x) = x^5 Q_{j,k}(x)$ for some polynomial $Q_{j,k}$.
\end{lem}

The advantage of the above expression is that each summand has a simple form 
\[
I_j = (\pa^j g )(x^{5 + a}) \cdot x^{ j a + 5 - k} P_{j,k}(x) 
=(-1)^j \f{1}{  (1 + x^{5 +a})^{j+1}} x^{ja + 5 - k} ( \sum_{i=0}^n a_i x^i  )
\]
for some $a_i \in \R$. For $x> 0$, the polynomial can be bounded by 
\[
|P_{j,k}(x)| \leq \sum_{i=0}^n |a_i| x^i \teq |P|_{j,k}(x),
\]
which is monotone increasing in $x$. In MatLab, we can compute the coefficients of $P_{j, k}$ easily using the recursive formula \eqref{eq:Fa_recur}. Taking absolute values on these coefficients, we assemble the polynomials $|P|_{j,k}$. Then $I_j$ can be bounded by the product of three monotone functions 
\[
|I_j| \leq \f{1}{ (1 + x^{5 + a})^{j+1}} x^{ja + 5 - k} |P|_{j,k}(x).
\]
Thus for $ x \in [x_i, x_{i+1}]$, we get 
\beq\label{eq:deri_Fa1}
|I_j| \leq \f{1}{ (1 + x_i^{5 + a})^{j+1}} \max( x_i^{ja + 5 - k} , x_{i+1}^{ja + 5-k} )
|P|_{j,k}(x_{i+1}). 
\eeq
For $ 1 \leq j \leq k$, since $P_{j,k}$ and $|P|_{j,k}$ are multiples of $x^5$, we obtain that 
$x^{j a + 5 - k} |P|_{j,k}$ is monotone increasing for $k \leq 10$ and the following refined estimate 
\[
|I_j| \leq \f{1}{ (1 + x_i^{5 + a})^{j+1}} x_{i+1}^{ja + 5-k} 
|P|_{j,k}(x_{i+1}) , \quad \forall x \in [x_i, x_{i+1}].
\]
For $k \leq 10$, the above estimate is useful for $x$ near $0$ since it does not blowup at $0$, while the upper bound in \eqref{eq:deri_Fa1} blows up if $i=0, x_i =0$ and $ja + 5 < k$. It is consistent with the regularity $F_a \in C^{10, a}$.

To obtain a sharp estimate of $\pa_x^k F_a$ in each interval $X_j = [x_j, x_{j+1}]$, we first use the above method to obtain the bound for $\pa_t^k F_a $ with $k = 2K+3$ and $2K + 4$, where $ K =8$ is the order of the Gaussian quadrature. Then we apply a simple interpolation result in Lemma \ref{lem:lin_interp} to obtain
\beq\label{eq:dF_interp}
\bal
|| \pa_{x}^{k } F_a||_{L^{\inf}(X_j)} = \max_{y \in [x_j, x_{j+1}]}| \pa_x^k F_a(y) | &\leq \max( |\pa_x^k F_a(x_j)|, | \pa_x^k F_a(x_{j+1}) ) \\
& \quad + \f{ (x_{j+1} - x_j)^2}{8} || \pa_{x}^{k+2} F_a||_{L^{\inf}(X_j)} 
\eal
\eeq
for each $k \leq 2K+2$ recursively.

\begin{proof}

Due to symmetry, we focus on $x \geq 0$. We prove Lemma \ref{lem:deri_Fa} by induction. For $k = 0$, we get 
\[
F_a = \f{x^5}{1 + x^{5 + a}} = g(x^{5 + a}) \cdot x^5 \cdot P_{0,0}(x), \quad P_{0,0}(x) = 1.
\]

Suppose that the statement in the Lemma holds for the cases $\leq k-1$ with $k\geq1$. We shall prove the case for $k$. Using the induction hypothesis for \eqref{eq:Fa_recur0} at $k-1$, we compute $ \pa_x^k F_a$
\[
\bal
\pa_x^k F_a &= \pa_x ( \pa_x^{k-1} F_a) = \pa_x \sum_{j=0}^{k-1} ( \pa^j g )(x^{5 + a}) \cdot x^{ja + 6- k} P_{j,k-1}(x) \\
&= \sum_{j=0}^{k-1} \B( (\pa^{j+1} g) (x^{5 + a}) \cdot (5 + a)x^{4 + a} \cdot x^{ja + 6- k} P_{j,k-1}(x) \\
& \quad + ( \pa^j g )(x^{5 + a}) \cdot ( (ja+ 6 - k)  x^{ja + 5-k } P_{j,k-1}(x)
+ x^{ja + 6 - k } \pa_x P_{j,k-1}(x) )  \B) \\
& = \sum_{j=0}^{k-1} \B( (\pa^{j+1} g)(x^{5+a}) \cdot x^{ (j+1)a + 5- k}  \cdot (5+a)x^5 P_{j,k-1} \\
& \qquad + (\pa^j g)(x^{5+a}) \cdot x^{ ja + 5 -k}  (  (ja + 6 -k) P_{j, k-1}
+ x \pa_x P_{j,k-1}) \B) \\
& = \sum_{j=0}^k (\pa^j g)(x^{5+a}) \cdot x^{ ja + 5 -k} 
( (5 + a) x^5 P_{j-1, k-1} 
+   (ja + 6 -k) P_{j, k-1}+ x \pa_x P_{j,k-1})
), 
\eal
\]
where we have used the convention $P_{j,k} = 0$ if $j<0$ or $j > k$. It follows the recursive formula \eqref{eq:Fa_recur} and the expression \eqref{eq:Fa_recur0}. 
The induction hypothesis that $P_{j,k-1}, 0 \leq j \leq k-1$ is a polynomial and \eqref{eq:Fa_recur} imply that $P_{j,k}$ is a polynomial. 

Taking $j = 0$ in \eqref{eq:Fa_recur}, we get 
\[
P_{0,k} = (6-k) P_{0, k-1}. 
\]
Since $P_{0,0}=1$, we yield the desired formula for $P_{0,k}$ in the Lemma. For $ j \geq 1$, using the induction hypothesis, we obtain that each term on the right hand side of \eqref{eq:Fa_recur} is a multiple of $x^5$ (note $0 = 0 \cdot x^5$). It follows that $P_{j,k}$ is a multiple of $x^5$.

The induction argument implies that the Lemma holds for all $k, 0\leq j\leq k$. 
\end{proof}

\subsubsection{ Estimate of derivatives for large $x$}

The above estimates do not work for large $x$ in MatLab due to overflow. In fact, for $x \geq 10^5$, the denominator in $\pa^j g = \f{1}{(1 + x^{5+a})^j }$ for $j \geq 16$ exceeds $ 10^{5 \cdot 16 \cdot 5} = 10^{400}$, which is out of the range of the double precision and treated as $\infty$ in MatLab. Thus, we develop another estimates using power series for $x \geq 10$. 

For $x \geq 10$, $F_a$ with $a\in(0,1)$ enjoys the following expansion 
\[
F_a = \f{x^5}{1 + x^{5+a}} = \f{x^{-a}}{1 + x^{-5-a}}
= \sum_{i \geq 0} (-1)^i x^{ - (5 +a)i - a}.
\]

Since $ x \geq 10$, the above series converges absolutely. Moreover, it is easy to obtain 
\beq\label{eq:Fa_power}
\pa_x^k F_a =  \sum_{i \geq 0} (-1)^{i+k} C_{i,k} x^{ - (5 +a)i - a - k}, \quad C_{i,k} = \prod_{j=0}^{k-1}( (5+a) i + a +j) ,
\eeq
for any $ k\in \mathbb{Z}, k \geq 0$, and the series also converges absolutely. 

We have the following simple estimate for $C_{i,k}$.

\begin{lem}\label{lem:Fa_power_coe}
For $x \geq 9, i \geq 10, k \leq 20, 0 < a < 1$, we have $C_{i,k} \leq x^{ 2i + k}$.
\end{lem}

\begin{proof}
Clearly, we have $C_{i,k} \leq ( (5+a) i + a + k-1)^k$. It suffices to prove 
$( (5+a) i + a + k-1)^k \leq 9^{2 i + k} $, which is equivalent to that
\[
G(i, k) =  \f{2 i + k}{k} \log 9 -\log( (5+a) i + a +k-1) 
\]
is nonnegative. For $i \geq 10, k \leq 20 \leq 2i, a<1$, note that 
\[
\pa_1 G(i,k) = \f{2}{k} \log 9 - \f{5+a}{  (5+a) i + a +k-1) } 
> \f{2}{k} \log 9 - \f{6}{ 3k} >0. 
\]

Thus we obtain
\[
\bal
G(i,k) &\geq G(10, k)  = \f{20 +k}{k }\log 9 - \log( 10 \cdot (5 + a) + a + k-1)  \\
&> 2 \cdot \log 9 - \log( 60 + 20) = \log 81 - \log 80  >0.
\eal
\]
This proves the desired estimate. 
\end{proof}

\begin{lem}\label{lem:Fa_power}
Let $x_0 \geq 9$ and $0\leq k \leq 20$. For $x \geq x_0$, we have 
\beq\label{eq:Fa_power2}
\pa_x^k F_a = \sum_{i=0}^9 (-1)^{i+k} C_{i,k} x^{-(5+a)i- a -k} + Err_{a,k}(x),
\quad |Err_{a,k}(x)| \leq \f{x^{-30 -a}}{1- x^{-3}}.
\eeq
As a result, we yield 
\[
|\pa_x^k ( F_a - x^{-\al} )|
\leq C(x_0, a, k) x^{-5-2a - k} , \quad C(x_0, a, k) =  \sum_{i=1}^9 C_{i, k} x_0^{ - (5+a) (i-1) } +\f{ x_0^{- 25 + a + k }}{ 1-  x_0^{-3}},
\]
where $C_{i,k}$ is defined in \eqref{eq:Fa_power}.
\end{lem}

\begin{proof}
Denote by $Err_{a,k}(x)$ the tail part of \eqref{eq:Fa_power}
\[
 Err_{a,k}(x) \teq \sum_{i\geq 10} (-1)^{i+k} C_{i,k} x^{-(5+a)i- a -k}.
\]
For $x \geq 9, k \leq 20$, applying Lemma \ref{lem:Fa_power_coe}, we obtain 
\[
|Err_{a,k}(x)| \leq \sum_{i\geq 10} |C_{i,k}| x^{-(5+a)i- a -k}
\leq \sum_{i\geq 10} x^{ 2i +k}  x^{-(5+a)i- a -k}
= \sum_{i\geq 10} x^{- 3i -a}
= \f{x^{-30 -a}}{1- x^{-3}},
\]
which is the estimate in \eqref{eq:Fa_power2}. Note that $\pa_x^k x^{-a}$ is the leading order term in \eqref{eq:Fa_power2}. Using $x \geq x_0 \geq 9$, the triangle inequality and the above estimates, we yield 
\[
|\pa_x^k( F_a-x^{-a})|
\leq \sum_{i\geq 1}^9 C_{i,k} x^{- (5 +a)i -a - k} + Err_{a,k}(x)
\leq x^{- 5-2a - k} \B( \sum_{i=1}^9 C_{i,k} x_0^{- (5 +a)(i-1)}  
+ \f{x_0^{-25 + a +k}}{ 1 - x_0^{-3}} \B),
\]
where we have used $-25  + a + k < 0$ to yield $x^{-25 + a + k} \leq x_0^{-25 + a + k}$ in the last inequality.
\end{proof}

Using the same argument, we obtain the following estimate for the other profile $\f{x}{1+x^2}$. 
Note that the associated velocity $u_g$ with $u_g(0) = 0, \pa_x u_g = H G =- \f{1}{1 + x^2}$ is given by $u_g = -\arctan(x)$. 

\begin{lem}\label{lem:g_power}
Let $ G = \f{x}{1 + x^2}$ and $u_g = - \arctan x$. For $ x \geq 10^{10}$ and $k=0,1,2,3,4$, we have 
\[
|\pa_x^k G| \leq (1 + \e_2) k ! x^{-k-1}, \quad
| \pa_x^{k+1} u_g|  \leq  (1+\e_2) (k+1) ! x^{-2 - k}, \quad \e_2 = 0.02 .
\]
For $u_g$, we have the trivial bound 
\[
|u_g (x) | \leq \f{\pi}{2}.
\]
\end{lem}

We use the rescaled version of $G$ to construct the approximate profile. In particular, suppose that $ G_{A, B} = A\cdot G(B x) $ and $u_{g, A,B}$ be the associated velocity. For $B x \geq 10^{10}$, we have 
\beq\label{eq:g_power_rescale}
|\pa_x^k G_{A, B}(x)| \leq A B^k  | (\pa_x^k G) (Bx) |,
\quad  |\pa_x^k u_{g, A, B} (x)| \leq A B^{k-1}  |(\pa_k u_g) (Bx)|.
\eeq

\subsection{ Leading order approximations and error estimates for $HF_a$ for large $x$ }\label{sec:ub_lead}

In this Section, we derive the leading order approximation of the Hilbert transform of $F_a$  and estimate the error for large $x$. Denote 
\beq\label{eq:Fa}
F_a = \f{x^5}{1 + |x|^{5 + a}}, \quad \bar F_a = \sgn(x) |x|^{-a}, \quad f_a = F_a - \bar F_a
= - \sgn(x) \f{ |x|^{-a} }{ 1 + |x|^{5+a}}.
\eeq

Since $|f_a| = |F_a - \bar F_a| \les x^{-5 -2a}$ is small, we expect that the Hilbert transform of $F_a$ can be approximated by $H \bar F_a$ for large $x$. 
We have the following result for the leading order part of $ \pa_x^k H F_a$ for large $x$. 

\begin{lem}\label{lem:ub}
For $k = 0,1,2$, we have 
\[
x^k \pa_x^k H F_a = x^k \pa_x^k H (\bar F_a ) - \f{1}{x^2} H (x^{k+2} \pa_x^k f_a )(0)
+ \f{ H(x^{k+4} \pa_x^k f_a ) - H (x^{k+4} \pa_x^k f_a)(0)}{x^4}.
\]
\end{lem}

We remark that $f_a$ has decay rates $|x|^{-5-2a}$ for large $x$ and $ |f_a x^2| \les |x|^{-3-2a}, |f_a x^4| \les |x|^{-1-2a}$ for large $x$. For large $x$, the first, the second and the third term on the right hand side have decay rates $|x|^{-a}, x^{-2}, x^{-4}$, respectively. Thus the third term can be regarded as the error term. To estimate the error term, we control the pointwise values of $H( x^{k+4} \pa_x^k f_a)$. 

\begin{lem}\label{lem:ub_err}
Suppse that $g$ is odd with $g \in L^1, g_x \in L^{\inf}$. We have 
\[
|H g(x)| \leq \f{2}{\pi |x|} || g||_{L^1} + \f{2}{\pi} \sup_{|y| \geq \f{|x|}{2}} |y g_x(y)|.
\] 
\end{lem}

We will apply the above estimate with $g = x^{k+4} \pa_x^k f_a$ to estimate the error term in Lemma \ref{lem:ub}. We use the above estimate rather than the standard Sobolev embedding to estimate $Hg$, e.g. $ || H g||_{L^{\inf}} \les || H g||_{L^2}^{1/2} || \pa_x H g||_{L^2}^{1/2}
= || g||_{L^2}^{1/2} || g_x||_{L^2}^{1/2}$, since the upper bound decays for large $x$ and it is easier to establish the bounds for the $|| g||_{L^1}, \sup_{|y| \geq \f{|x|}{2}} |y g_x(y)|$ rigorously. 

A consequence of the above two lemmas is the following estimate for $H F_a$. 

\begin{cor}\label{cor:ub}
For $k = 0,1,2$ and $x \geq x_0$, we have 
\beq\label{eq:cor_ub1}
\pa_x^k H F_a  = \pa_x^k H \bar F_a -  \f{ (-1)^k  (k+1)! C_u}{x^{k+2}} + E_k(x), \quad C_u = H ( f_a x^2)(0),
\eeq
where the error term satisfies
\beq\label{eq:cor_ub2}
\bal
E_k(x) & \leq C_{u,err}(k) x^{-k-4},  \\ 
 C_{u,err}(k) & = \f{ (3+k)!}{ 3\pi} (\f{1}{4-a} + \f{1}{1 + 2a}) + \f{2}{\pi x_0} || x^{k+4} \pa_x^k f_a||_{L^1} + \f{2}{\pi} \sup_{ y \geq x_0 / 2} 
| y \pa_x( x^{k+4} \pa_x^k f_a ) |.
\eal
\eeq
Moreover, suppose that $u_a$ is the associated velocity such that $u_{a,x} = H F_a$ and $u_a(0)$. We have 
\[
u_a = - \cot \f{\pi a}{2} \f{x^{1-a}}{1-a} + \f{C_u}{x} + E, \quad E \leq \f{C_{u,err}(0)}{3 x^3}.
\]

\end{cor}

We will first prove the Corollary and then prove Lemmas \ref{lem:ub} and Lemma \ref{lem:ub_err}. 

\begin{proof}[Proof of Corollary \ref{cor:ub}]

We first estimate $|H( x^{k+4} \pa_x^k f_a)(0)|$ that appears in Lemma \ref{lem:ub}.
From \eqref{eq:Fa}, we get $ |x^k \pa_x^k f_a | \les \min( |x|^{-a}, x^{-5-2a})$. Using integration by parts, we yield 
\[
\bal
|H (x^{k+4} \pa_x^k f_a)(0)| &= \f{1}{\pi} \B| \int_{\R} x^{k+3} \pa_x^k f_a(x) dx \B|
= \f{1}{\pi} \B| \int_{\R} \pa_x^k (x^{k+3}) f_a(x) dx \B|  \\
&= \f{(k+3)!}{6 \pi} \B|\int_{\R} x^3 f_a(x) dx \B|
\leq \f{(k+3)!}{6 \pi} \cdot 2 \B|\int_0^{\inf} \f{x^{3-a}}{1 + x^{5+a}} dx \B| \\
&\leq \f{(k+3)!}{3 \pi} ( \int_0^1 x^{3-a} dx + \int_1^{\inf} x^{-2-2a} dx)
= \f{(k+3)!}{3 \pi} ( \f{1}{4-a} + \f{1}{1 + 2a}).
\eal
\] 

Similarly, using integration by parts, we have 
\[
\bal
H (x^{k+2} \pa_x^k f_a)(0)
&= -\f{2}{\pi} \int_0^{\pi} x^{k+1} \pa_x^k f_a(x) d x
= -\f{2}{\pi } (-1)^k \int_0^{\inf} \pa_x^k( x^{k+1}) f_a(x) dx \\
&= -\f{2}{\pi} (-1)^k  (k+1)! \int_0^{\inf} x f_a(x) dx 
=  (-1)^k(k+1)! H( x^2 f_a)(0).
\eal
\]

Plugging the above estimates and the estimate in Lemma \ref{lem:ub_err} with $g = x^{k+4} \pa_x^k f_a$ into Lemma \ref{lem:ub}, we prove \eqref{eq:cor_ub1} and \eqref{eq:cor_ub2}. 

Next, we establish the estimate for $u_a$. Since $u_a$ is odd, it suffices to consider $x>0$. Denote by $u(f) = \int_0^x Hf(y) dy$ the velocity associated to $f$. Firstly, we show that $\lim_{|x|\to \inf} u(f_a)(x) = 0$. Since $f_a$ is odd, using a change of variable $y = xz$, we yield 
\[
u(f_a)(x) = \f{1}{\pi} \int_0^{\inf} f_a(y) \log \B| \f{x- y}{x+y} \B| dy 
=  \f{1}{\pi} \int_0^{\inf} f_a(xz) x  \log \B| \f{1- z}{1+z} \B| dz .
\]
From \eqref{eq:Fa}, we get $|f_a(x) | \les |x|^{-3/2} $, which implies 
\[
|f_a(x z) x | \les (xz)^{-3/2} x \leq x^{-1/2} z^{-3/2} .
\]

In particular, for $z>0$ and $x \geq 1$, we obtain 
\[
\lim_{x \to \inf} f_a(x z) x  = 0, \quad |f_a(x z) x | \leq z^{-3/2}.
\]

Since $ \log \B| \f{1- z}{1+z} \B| z^{-3/2} \in L^1(0,\inf)$ (it is of order $z^{-1/2}$ near $z=0$), using the Dominated Convergence Theorem, we yield 
\[
\lim_{x \to \infty} u(f_a)(x) = 0.
\]

Using \eqref{eq:cor_ub1} with $k = 0$, we have 
\[
H f_a = H F_a - H \bar F_a = - \f{C_u}{x^2} + E_0(x) .
\]

Since $ u(f_a)(\infty) = 0$, integrating the above formula from $x$ to $\infty$, we obtain 
\[
u(f_a) = -\int_x^{\inf} H f_a(y) dy
= - \int_x^{\inf} - \f{C_u}{y^2} + E_0(y) dy 
= \f{C_u}{x} - \int_x^{\inf} E_0(y) dy.
\]

Note that $H \bar F_a =-  \cot \f{\pi a}{2} |x|^{-a}, u(\bar F_a) =- \sgn(x) \cot \f{\pi a}{2} \f{ |x|^{1-a} }{1-a} $. Using the above identity and \eqref{eq:cor_ub2} with $k=0$, we establish 
\[
\bal
\B| u(F_a) + \sgn(x) \cot \f{\pi a}{2} \f{ |x|^{1-a} }{1-a} - \f{C_u}{x}\B|
&= \B| u(f_a) - \f{C_u}{x}\B|
\leq \int_x^{\inf} | E_0(y) | dy  \\
&\leq C_{u,err}(0) \int_x^{\inf} y^{-4} dy 
= \f{C_{u,err}(0)}{3} x^{-3 } ,
\eal
\]
which is the desired estimate.
\end{proof}

\subsubsection{Proof of Lemmas \ref{lem:ub} and Lemma \ref{lem:ub_err}}

To prove Lemma \ref{lem:ub}, we need the following identity.

\begin{lem}\label{lem:ub1}
Suppose that $f$ is in the Schwartz space and is odd. We have 
\[
H(f x^4) = H (fx^4)(0) +  H(fx^2)(0) \cdot x^2 + x^4 H f(x)
\]
\end{lem}

\begin{proof}
Since $f$ is odd, for $k=1,3$, we obtain 
\[
H( f x^k )(0) = -\f{1}{\pi} \int_{\R} f x^{k-1} dx = 0 .
\]

Hence, applying Lemma \ref{lem:com}, we get 
\[
\bal
H (fx^4) &= H (fx^4)(0) + x H (fx^3), \quad H (fx^3) = x H (f x^2), \\
 H (fx^2) &= H (fx^2)(0) + x H (fx) , \quad  H (fx) = x H f.
\eal
\]

It follows
\[
\bal
H ( fx^4) &= H (f x^4)(0) + x H( fx^3) = H (fx^4)(0) + x^2 H ( fx^2)\\
 &= H (fx^4)(0) + x^2 H ( fx^2)(0) + x^3 H(fx) = 
  H (fx^4)(0) + x^2 H ( fx^2)(0) + x^4 H f.
  \eal
\]
\end{proof}

Now, we prove Lemma \ref{lem:ub}.

\begin{proof}[Proof of Lemma \ref{lem:ub}]

Firstly, for Schwartz function $g$, note that $H (x g_x)(0) = - \f{1}{\pi}\int_{\R} g_x dx = 0 $. Using Lemma \ref{lem:com}, we obtain 
\[
H( x\pa_x g) = x H (\pa_x g) = x \pa_x H g.
\]
Thus $D_x\teq x\pa_x$ commutes with $H$. By definition \eqref{eq:Fa}, for $k \in \mathbb{Z}$ and $k \geq 0$, we yield $ |D_x^k f_a | \les  \min( |x|^{-a}, x^{-5})$ for some $a \in (0,1)$ and thus  $ D_x^k f_a \in L^p$ for some $p \in (1,\infty)$. For a nonnegative integer $k$, using a standard approximation argument, we obtain 
\[
H ( D_x^k f_a) = D_x^k H f_a,
\]
which along with $F_a = \bar F_a + f_a$ implies 
\[
D_x^k H F_a = D_x^k  H \bar F_a + D_x^k H f_a = D_x^k  H \bar F_a +  H ( D_x^k f_a).
\]

Since the Hilbert transform is linear and $x^k \pa_x^k$ can be written as a linear combinations of $D_x^i, i=0,1,.., k$, we derive 
\[
x^k \pa_x^k H F_a = x^k \pa_x^kH\bar F_a + H( x^k \pa_x^k f_a).
\]

Applying Lemma \ref{lem:ub1} with $f = x^k\pa_x^k f_a, k \in \mathbb{Z}$ and $k\geq0$, we yield
\[
\bal
H(x^k \pa_x^k f_a) &= \f{ H(f x^4) - H(fx^4)(0) - H(fx^2)(0) \cdot x^2}{x^4}  \\
&= -\f{1}{x^2} H( x^{k+2} \pa_x^k f_a) + \f{1}{x^4} ( H( x^{k+4} \pa_x^k f_a) -  H( x^{k+4} \pa_x^k f_a)(0) ).
\eal
\]
We conclude the proof of Lemma \ref{lem:ub}.\\
\end{proof}

Next, we prove Lemma \ref{lem:ub_err}.
\begin{proof}[Proof of Lemma \ref{lem:ub_err}]
Without loss of generality, we consider $x > 0$. We decompose the integral as follows 
\[
Hg(x) = \f{1}{\pi} \int_{x /2}^{3x/2} \f{g(y)}{x-y} dy 
+ \f{1}{\pi}\int_{ |x-y|> x/2} \f{g(y)}{x-y} dy \teq I + II.
\]
For $I$, using a change of variable $y = x + z$, we yield 
\[
\bal
|I| &= \B|\f{1}{\pi}\int_{-x/2}^{x/2} \f{g(x+z)}{-z} dz \B|
= \f{1}{\pi} \B|\int_0^{x/2} \f{g(x-z) - g(x+z)}{z} dz \B| \\
&\leq \f{1}{\pi} \sup_{y \geq x/2} |g_x(y)| \int_0^{z/2} 2 dz
= \f{1}{\pi} \sup_{y \geq x/2} |g_x(y)| x 
\leq \f{2}{\pi} \sup_{ y \geq x/2} | y g_x(y)|.
\eal
\]
The estimate of $II$ is trivial and is given below:
\[
|II| \leq \f{1}{\pi}\cdot \f{1}{ x/ 2} ||g||_{L^1} \leq \f{2}{\pi |x|} || g||_{L^1}.
\]
Combining the estimates of $I$ and $II$ yields the desired result. 
\end{proof}

\subsubsection{ Computation of $H(f_a x^2)(0)$}

Recall from \eqref{eq:Fa} that
\[
F_a = \f{x^5}{1 + |x|^{5+a}}, \quad  f_a = - \sgn(x) \f{ |x|^{-a}}{ 1 + |x|^{5+a}} = F_a - \sgn(x) |x|^{-a} .
\]
Since $f_a$ is not smooth near $x=0$, to compute $H(f_a x^2)$ with small error, we cannot apply a standard quadrature rule directly. For $ 0 < A < 1 < B$, we decompose the integral $H(f_a x^2)(0)$ as follows 
\[
\bal
 H (f_a x^2)(0)
&= - \f{2}{\pi} \int_0^{\inf} x f_a(x) dx 
= - \f{2}{\pi} ( \int_0^A x F_a(x) dx + \int_A^B  x F_a(x) dx 
- \int_0^B x \cdot x^{-a} + \int_B^{\inf} x f_a (x) dx )  \\
&\teq -\f{2}{\pi} (I_1 + I_4 - I_2 + I_3).
\eal
\]
We choose $A \approx 0.1, B \approx 10$. 
The idea to compute $I_1$ and $I_4$ is similar to performing a power series expansion. For $I_1$, we get 
\[
I_1 = \int_0^A \f{x^6}{1 + x^{5+a}} dx 
=  \int_0^A \f{ x^6(1 - x^{20 + 4a})}{1 + x^{5+a}} + \f{x^{26 + 4a}}{1 + x^{5+a}} dx
\teq I_{1,1} + I_{1,2}.
\]
Using the identity $\f{1 - y^{2k}}{1 + y} = \sum_{0\leq i \leq 2k-1} (-1)^i y^i$ and a direct calculation yields 
\[
I_{1,1} =  \int_0^A x^6 \B(\sum_{i=0}^3 (-1)^i x^{ (5+a)i} \B)  dx
= \sum_{i=0}^3 (-1)^i \f{ A^{7 + (5+a)i}}{7 + (5+a)i}.
\]
The estimate for the error term $I_{1,2}$ is trivial and is given below
\[
|I_{1,2} |\leq \int_0^A x^{26 + 4a} dx = \f{ A^{27 + 4a}}{27 + 4a}.
\]

The estimate for $I_3$ is similar. We first perform a decomposition
\[
I_3 = - \int_B^{\inf} \f{x^{1-a}}{1 + x^{5+a}} dx 
=- \int_B^{\inf} \f{ x^{-4-2a}}{1 + x^{-5-a}} dx
= - \int_B^{\inf}  \f{ x^{-4-2a} (1 - x^{-20-4a})}{1 + x^{-5-a}} 
- \f{ x^{-24 -6a}}{ 1 + x^{-5-a}} dx \teq I_{3,1}  + I_{3,2}.
\]
Then we get 
\[
\bal
I_{3,1} & = -\int_B^{\inf} \sum_{i=0}^3 x^{-4-2a} (-1)^i x^{-(5+a) i} dx
= \sum_{i=0}^3 (-1)^{i+1} \f{ B^{-(5+a)i -3 - 2a}}{ 3 + 2a + (5+a)i }  ,\\
|I_{3,2} |& \leq \int_B^{\inf} x^{-24 - 6a} dx 
= \f{ B^{-23- 6a}}{  23 + 6a}. 
\eal
\]

For $I_2$, using a direct calculation yields 
\[
I_2 = \int_0^B x^{1-a} dx = \f{B^{2-a}}{2-a}.
\]

For $I_4$, we partition $[A, B]$ into $ y_i = A + ih, i=0, 1, .., n_1$ with $n_1 h = B - A$ and use the Simpson's rule  to obtain 
\[
I_4 = \int_A^B x F_a(x) dx 
= \f{h}{6} \sum_{i=1}^{n_1} ( y_{i-1} F_a(y_{i-1} ) 
+ y_i F_a(y_i) + 4 ( y_{i}-\f{h}{2} ) F_a( y_i - \f{h}{2}) ) + I_{4,err},
\]
where the error is bounded by 
\[
I_{4,err} \leq \f{h^5}{2880} \sum_{i=1}^{n_1} || \pa_x^4 (x F_a)||_{L^{\inf}[y_{i-1}, y_i]}.
\]
Since $F_a$ is given explicit, we use Mathematica to derive the explicit formula of $\pa_x^4 (xF_a)$. Then using the triangle inequality, we further bound it by 
\[
 |\pa_x^4 (x F_a) | \leq \f{P(x)}{ (1 + x^{5+a})^5} 
\]
for some increasing function $P(x)$, which is given in the MatLab and Mathematica codes. Then we can obtain the bound for  $ |\pa_x^4 (x F_a) |$ in each interval $[y_{i-1}, y_i]$ using the endpoint of $P(x)$.

\subsubsection{ Rigorous estimate of some norms related to $\pa_x^k f_a$}

To apply Corollary \ref{cor:ub} with $g = x^{k+4} \pa_x^k f_a, k =0,1,2$, we need to estimate $||g||_{L^1}$ and $\sup_{y \geq x_0/2} | y g_x(y)|$. 
Firstly, note that 
\[
 y g_x(y) = (k+4) y^{k+4} \pa_x^k f_a + x^{k+5} \pa_x^{k+1} f_a.  
\]

For $y_0 \geq 10$ and $y \geq y_0$, applying the estimate in Lemma \ref{lem:Fa_power} to $\pa_x^k f_a = \pa_x^k (F_a - x^{-a})$, we get 
\[
 | y g_x(y)| 
 \leq C(y_0, a, k) (k+4) y^{-1- 2a} + C(y_0, a, k+1) y^{-1-2a}.
\]
Choosing $y_0 = \f{M_1}{2}, M_1 = 10^5$, we yield 
\[
\sup_{y \geq M_1/ 2}  | y g_x(y)| 
\leq \B( C(M_1 /2, a, k) (k+4) + C(M_1/2, a, k+1) \B) (M_1/ 2)^{-1-2a}.
\]

To estimate the $L^1$ norm of $g$, we first decompose the integrals as follows 
\[
\bal
|| g||_{L^1} &= 2 \int_0^{\inf} |x^{k+4} \pa_x^k f_a |dx  \leq   2 \int_0^{x_i} |x^{k+4} \pa_x^k F_a| +| x^{k+4} \pa_x^k x^{-a}| dx  \\
 & \qquad + 2\int_{x_i}^{x_j} |x^{k+4} \pa_x^k f_a| dx + 2 \int_{x_j}^{\inf} |x^{k+4} \pa_x^k f_a| dx
 \teq I_{1,1} + I_{1,2} + I_2 + I_3,
 \eal
\]
 where $x_i, x_j$ are some grid points with  $ 0\leq i < j \leq n$.

The motivation for the above decomposition is the following. Near $x=0$, $x^{-\al}$ is singular and we do not have a sharp estimate for $\pa_x^k (F_a - x^{-a})$. Thus, we estimate the contribution from $F_a$ and $x^{-a}$ separately. 

We establish the derivative bounds for $F_a$ and $f_a$ in each intervals $[x_l, x_{l+1}]$ 
recursively using the method similar to \eqref{eq:dF_interp} in MatLab. Then the estimates of $I_{1,1}$ and $I_{2}$ are simple 
\[
I_{1,1} + I_2
\leq 2 \sum_{l=1}^{i} x_l^{k+4}( x_l - x_{l-1}) \max_{y \in [x_{l-1}, x_l]} | \pa_x^k F_a(y)|
+ 2 \sum_{l=i+1}^{j} x_l^{k+4}( x_l - x_{l-1}) \max_{y \in [x_{l-1}, x_l]} | \pa_x^k f_a(y)|.
\]

Since $|x^{k+4}\pa_x^k x^{-a}| =x^{4-a} \prod_{l=0}^{k-1} (a + l) $, we yield 
\[
I_{1,2} \leq 2 \B( \prod_{l=0}^{k-1}  (a + l) \B) \int_0^{x_i} x^{4-a} dx 
= 2 \B(\prod_{l=0}^{k-1} (a + l) \B)\f{ x_i^{5-a} }{5-a}.
\]

Using Lemma \ref{lem:Fa_power}, we get 
\[
I_3 \leq 2 \int_{x_j}^{\inf} x^{k+4} C(x_j, a, k) x^{-5-2a -k} dx 
= 2 C(x_j, a,k) \f{x_j^{-2a}}{2a}.
\]
Since we only use $k =0,1,2$, by definition of $C(x_j, a, k)$, the above upper bound is decreasing in $x_j$.

\subsection{ Error estimates for $H F_a$ with $x$ of $O(1)$}\label{sec:ub_errO1}

Recall in Section \ref{sec:compute_ub} that for $x \leq M_1$, we decompose the integral in $H F_a$ into the part $x \leq L$ and $x > L$. For $x \leq M_1$, we have the following estimates for the contribution from $ |y| \geq L$.

\begin{lem}\label{lem:HF_integ_far}
Let $K_i(x, y)$ be the kernel associated to the derivatives of the velocity \eqref{eq:compute_u}, \eqref{eq:ker}.
For $F_a = \f{x^5}{1 + x^{5+a}}$ and $0\leq x \leq M_1$, we have 
\[
\bal
 \f{1}{\pi}\int_L^{\inf} F_a(y) K_1(x, y) dy  &= - \f{2}{\pi a} L^{-a} x + E(u,1)(x),  \\
\f{1}{\pi}\int_L^{\inf} F_a(y) K_2(x, y) dy &= - \f{2}{\pi a} L^{-a}  + E(u,2)(x) , \\
\f{1}{\pi}\int_L^{\inf} \pa_yF_a(y) K_3(x, y) dy &=  \f{2a}{\pi (2 + a)} L^{-a-2}  x + E(u,3)(x) \\
\f{1}{\pi}\int_L^{\inf} \pa_y^2F_a(y) K_2(x, y) dy &= - \f{2a(a+1)}{\pi (2 + a)} L^{-a-2}   + E(u,4)(x)
\eal
\]
where the error terms satisfy 
\[
\bal
|E(u, 1)(x)| &\leq \f{2(1+\e_2)}{\pi} \B(  C( L, a, 0)   \f{L^{-5-2a}}{5 + 2a}  x
+ \f{ L^{-2-a}}{2 + a} \f{x^3}{3}  \B) ,\\
|E(u, 2)(x)| &\leq \f{2(1+\e_2)}{\pi} \B(  C( L, a, 0)   \f{L^{-5-2a}}{5 + 2a} 
+ \f{ L^{-2-a}}{2 + a} x^2  \B) ,\\
|E(u,3)(x)| &\leq \f{2(1+\e_2)}{\pi} \B(  C( L, a, 1)   \f{L^{-7-2a}}{7 + 2a} x
+ \f{a L^{-4-a}}{4 + a} x^3  \B) ,\\
|E(u,4)(x)| &\leq \f{2(1+\e_2)}{\pi}   \B( C( L, a, 2)   \f{L^{-7-2a}}{7 + 2a} 
+ a(a+1)    \f{L^{-4-a}}{4 + a}  x^2 \B)  ,
\eal
\]
 $\e_2 = 0.02$, and $C(L, a, k)$ is defined in Lemma \ref{lem:Fa_power}.
\end{lem}

\begin{proof}
We establish the last estimate related to $\pa_y^2 F_a$. Other estimates can be proved similarly. For large $y$, from Lemma \ref{lem:Fa_power} and Lemma \ref{lem:ker_xsm}, we know that $\pa_y^2 F_a \approx \pa_y^2 y^{-a}, K_2(x, y) \approx - \f{2}{y}$. Recall $f_a = F_a - x^{-a}$ from \eqref{eq:Fa}. We decompose the integral as follows 
\[
\bal
I_4(x)  \teq \f{1}{\pi}\int_L^{\inf} \pa_y^2F_a(y) K_2(x, y) dy 
&= \f{1}{\pi} \int_L^{\inf} \pa_y^2  f_a  K_2(x, y) dy 
+ \f{1}{\pi} \int_L^{\inf} \pa_y^2  y^{-a}  ( K_2(x, y) + \f{2}{y}) dy \\
& \qquad - \f{1}{\pi} \int_L^{\inf} \pa_y^2  y^{-a}  \f{2}{y} dy 
\teq J_1 + J_2 + J_3.
\eal
\]

The main term in $I_4(x)$ is given by $J_3$. Using a direct calculation, we yield 
\[
J_3 = -\f{2}{\pi} a (a+1) \int_L^{\inf}  y^{-a-3} dy = -\f{2}{\pi} a (a+1) \f{L^{-a-2}}{a+2} .
\]
Using Lemma \ref{lem:Fa_power} and Lemma \ref{lem:ker_xsm}, we estimate $J_1, J_2$ as follows 
\[
\bal
|J_1| &\leq \f{1}{\pi}  2(1 + \e_2) C( L, a, 2)\int_L^{\inf} y^{-5-2a - 2} y^{-1} dy
= \f{1}{\pi}  2(1 + \e_2)C( L, a, 2)   \f{L^{-7-2a}}{7 + 2a} ,\\
|J_2| &\leq \f{a(a+1)}{\pi} 2 (1 + \e_2 )  \int_L^{\inf} y^{-a-2} \f{x^2}{y^3} dy
\leq \f{a(a+1)}{\pi} 2 (1 + \e_2)  x^2 \f{L^{-4-a}}{4 + a},
\eal
\]
where $\e_2 = 0.02$. \end{proof}

\section{ Estimate of solution and error beyond the computation domain}\label{sec:beyond_LB}

In this section, we estimate the approximate steady and the error $(x \pa_x)^i F_{\om}, (x\pa_x)^i F_v, i=0,1$ for $x\geq L_B$. 

Recall that the approximate steady state can be decomposed into $\om = \om_b + \om_p,
\th_x = v = v_b + v_p$, where the perturbation part $\om_p, v_p$ is supported in $[-L, L], L < L_B$. Thus $\om_p(x) =0, v_p(x) = 0$ for $x> L_B$.

\subsection{Estimate of the derivatives of $H \om_p$ for $x \geq L_B$}

Since $\om_p$ is supported in $[-L, L]$ and is odd, symmetrizing the kernel, we get 
\beq\label{eq:up_far1}
u_p = \f{1}{\pi}\int_0^{L} \om_p(y) \log \B|\f{y-x}{y+x} \B| dy, \quad
\pa_x^{k+1} u_p = \f{1}{\pi} \int_0^{L} \pa_x^k \om_p(y) ( \f{1}{x-y} - (-1)^k \f{1}{x+y}) dy.
\eeq

For $x \geq L_B$, it is away from the support of $\om_p$ and we have $\f{y}{x} \leq \f{L}{L_B}$ for $y\geq 0$ and $ y \in \mathrm{supp}(\om_p)$. To estimate $u_p, \pa_x^{k+1} u_p$, we expand the kernel. For $ 0\leq \f{y}{x} \leq \f{1}{100}$, we have 
\beq\label{eq:ker_TL1}
\bal
K_1(x, y) &\teq \log \B|\f{x-y}{x+y}\B|
= \log \f{1 - \f{y}{x}}{1 + \f{y}{x}}
= -2 \sum_{i\geq 1} \f{1}{2i-1} ( \f{y}{x})^{2i-1} \\
K_2(x, y) &\teq \f{1}{x-y} - \f{1}{x+y} 
= \f{1}{x} ( \f{1}{1 - \f{y}{x}} - \f{1}{1 + \f{y}{x} } )
= \f{1}{x} \sum_{i\geq 0} (\f{y}{x})^i ( 1 - (-1)^i  )
= \f{2}{x} \sum_{i \geq 1} (\f{y}{x})^{2i-1}, \\ 
K_3(x, y)& \teq \f{1}{x-y} + \f{1}{x+y} 
= \f{1}{x} ( \f{1}{1 - \f{y}{x}} + \f{1}{1 + \f{y}{x} } )
= \f{1}{x} \sum_{i\geq 0} (\f{y}{x})^i ( 1 + (-1)^i  )
= \f{2}{x} \sum_{i \geq 0} (\f{y}{x})^{2i} . 
\eal
\eeq

Denote by $m_i$ the moments of $\om_p$ 
\[
m_i =  \f{2}{\pi} \int_0^L \om_p(y) y d y  , \quad m_3 =   \f{2}{\pi} \int_0^L \om_p(y) y^3 d y ,
\quad m_5 = \f{2}{\pi} \int_0^L |\om_p y^5|  dy.
\]

We will use the following simple Lemma to estimate the remainder term in the power series.

\begin{lem}\label{lem:power_coe2}
For $k \geq 3,i = 1,2$ and $0\leq \f{y}{x} \leq \f{1}{100}$, we get 
\[\prod_{j=0}^{i-1} (2k+j) \leq \B(  \prod_{j=0}^{i-1}( 6 +j)  \B) (\f{y}{x})^{3-k}.\]
\end{lem}

\begin{proof}
The case of $k=3$ is trivial. Clearly, we have $\prod_{j=0}^{i-1} \f{2k + j}{6 +j }  \leq (\f{2k}{6})^i = (\f{k}{3})^i$. For $k \geq 4, i=1,2, \f{x}{y} \geq 100$, using Bernoulli's inequality, we get 
\[
(\f{x}{y})^{ \f{k-3}{i}} \geq 10^{k-3}= (1 + 9)^{k-3}
\geq 1 + 9(k-3) = 9k - 26 \geq k  > \f{k}{3}.
\]
It follows $(\f{x}{y})^{ k-3} \geq \prod_{j=0}^{i-1} \f{2k + j}{6 +j }$. Rearranging the inequality yields the desired result.
\end{proof}

Denote $\e_2 = 0.02$. For $x \geq 100 y$, using \eqref{eq:ker_TL1}, we get 
\[
\bal
K_1(x, y) & = - 2 \f{y}{x} - 2 \f{y^3}{3 x^3} + E_1(x, y) ,  \quad 
K_2(x, y) = \f{2 y}{x^2} + \f{2y^3}{x^4}+ E_2(x, y),  \\
\pa_y K_3(x, y) & = \f{2}{x} \sum_{i\geq 1}  \f{2i}{x} (\f{y}{x})^{2i-1}
= \f{4 y}{x^3} + \f{8 y^3}{x^5} + E_3(x, y), \\
\pa_y^2 K_2(x, y) & = \f{2}{x} \sum_{i\geq 2} \f{ (2i-1)(2i-2)}{x^2}(\f{y}{x})^{2i-3}=
\f{2}{x} \sum_{i\geq 1} \f{ (2i+1)(2i)}{x^2}(\f{y}{x})^{2i-1}
= \f{12 y}{x^4} + \f{40 y^3}{x^6} + E_4(x, y).
\eal
\]
We focus on the $\pa_y^2 K_2(x, y)$ term. A direct estimate yields 
\[
\bal
 |E_4(x, y)|
 &\leq \f{2}{x} \sum_{i \geq 3} \f{(2i+1)(2i)}{x^2}  (\f{y}{x} )^{2i-1}
 \leq \f{2}{x} \sum_{i \geq 3} \f{1}{x^2} 42 (\f{y}{x})^{3-i} (\f{y}{x} )^{2i-1}
 = \f{84}{x^3} \sum_{i\geq 3} (\f{y}{x})^{i+2} \\
 & =  \f{84 }{x^3} \f{ (y/x)^5}{1 - y/x}  \leq \f{84(1 + \e_2) y^5 }{x^8},
\eal
\]
and we have used $ \f{1}{1 - y/x} \leq \f{1}{1 - 0.01} \leq 1+ \e_2$ and Lemma \ref{lem:power_coe2} to estimate the coefficient $2i (2i+1)$. Similarly, for $E_1, E_2, E_3$, we get 
\[
\bal
 |E_1(x, y) | \leq \f{2 (1 + \e_2)}{5} (\f{y}{x})^5 , \quad 
|E_2(x, y)| \leq 2 (1 + \e_2) \f{y^5}{x^6} ,  \quad  |E_3(x, y)|  \leq \f{12 ( 1 + \e_2)y^5}{x^7} . \\
\eal
\]

Now, plugging the above estimates of the kernels into \eqref{eq:up_far1} and then using integration by parts, we yield 
\beq\label{eq:up_far}
\bal
u_p &= \f{1}{\pi}\int_0^L \om_p(y) K_1(x, y) dy 
=  \f{1}{\pi}\int_0^L \om_p(y) (- 2 \f{y}{x} - 2 \f{y^3}{3 x^3} + E_1(x, y) ) dy 
= - \f{m_1}{x} - \f{m_3}{3 x^3} + E(u_p,1)(x) ,\\
\pa_x u_{p} &= \f{1}{\pi} \int_0^L \om_p(y) K_1(x, y) dy 
= \f{1}{\pi} \int_0^L \om_p(y) (
\f{2 y}{x^2} + \f{2y^3}{x^4}+ E_2(x, y) ) dy
= \f{m_1}{x^2} + \f{m_3}{x^4} + E(u_p, 2)(x) , \\
\pa_x^2 u_{p} &= \f{1}{\pi} \int_0^L \partial_x \om_p(y) K_3(x, y) dy 
= \f{1}{\pi} \int \om_p(y) (-\pa_y K_3(x, y) ) dy \\
&= - \f{1}{\pi} \int \om_p(y) (  \f{4 y}{x^3} + \f{8 y^3}{x^5} + E_3(x, y)) dy
= - \f{2 m_1}{x^3} - \f{4m_3}{x^5} + E(u_p, 3)(x) , \\
\pa_x^3 u_{p} &= \f{1}{\pi} \int_0^L \partial_x^2\om_p(y) K_2(x, y) dy 
= \f{1}{\pi} \int_0^L \om_p(y) (\pa^2_y K_2(x, y) ) dy \\
&=  \f{1}{\pi} \int \om_p(y) ( \f{12 y}{x^4} + \f{40 y^3}{x^6} + E_4(x, y) ) dy
=  \f{6 m_1}{x^4} + \f{20 m_3}{x^6} + E(u_p, 4)(x) , \\
\eal
\eeq
where we have used $\om_p(0)=  \om_p(L) = \om_{p, x}(L) = 0$ and $K_2(x,0) = 0$ when we applied integration by parts in the derivations of $\pa_x^2 u_p, \pa_x^3 u_p$.  The error terms $E(u_p,i)$ satisfy 
\beq\label{eq:up_far_err}
\bal
|E(u_p, i)(x)| &\leq \f{1}{\pi} \int_0^L |\om_p(y) E_i(x,y)| dy 
\leq \f{1}{\pi} 2(1+\e_2) C_i   \int_0^L |\om_p(y) \f{y^5}{x^{4+i}}| dy \\
&\leq C_i( 1 + \e_2)  m_5 x^{-4-i}.
\eal
\eeq
and $C_1 = \f{1}{5}, C_2 = 1, C_3 = 6, C_4 = 42$.

We remark that the leading order terms of $\pa_x^k u_p$ enjoy the differential relation, and the upper bound for the error term also enjoys a similar relation :
$ |\pa_x C_i(1 + \e_2) m_5 x^{-4-i}| = |C_{i+1} (1 + \e_2) m_5 x^{-4-i-1}|$.

\subsection{ Estimate of the error for $x \geq L_{B}$}\label{sec:err_decay} 

We estimate the residual errors $(x\pa_x)^i F_{\om}, (x \pa_x)^i F_v$ for the approximate steady state \eqref{eq:solution_decomposition}. Due to symmetry, we focus on $x \geq 0$. Recall that the approximate steady state can be decomposed as 
\beq\label{eq:solu_decomp}
\om = \om_p + \om_b = \om_p  + \f{b_{\om} x^{5}}{ 1 + x^{5 + a_{\om}}}
+ \f{s_{\om} x}{1 + (r_{\om} x)^2} ,
\quad v =v_p +  \frac{b_v x^5}{1 + x^{5 + a_v}} + \frac{s_v x}{1 + (r_v x)^2},
\eeq
for $x \geq 0$, where $\om_p $ and $v_p$ are quintic splines.

We first derive the leading order behavior for the profile $\om, v$ and $u$ in the far field $x \geq L_B$. Since the quintic spline $\om_p, v_p $ is supported in $[-L, L]$ and $L_{B} > L$, we get $\om_p = v_p = 0$ for $x \geq L_B$. Clearly, the leading order parts of $\om, v$ are given by 
\beq\label{eq:decomp_wv}
\bal
\om &= b_{\om} x^{-a_{\om}} + \B( b_{\om}  \f{x^{5}}{ 1 + x^{5 + a_{\om}}} - 
 b_{\om} x^{-a_{\om}} + \f{s_{\om} x}{1 + (r_{\om} x)^2}
  \B)
  \teq \om_{M} + \om_{R},  \\
v & = b_v x^{-2 a_v} + 
\B(
\frac{b_v x^5}{1 + x^{5 + a_v}}  - b_v x^{-a_v}+ \frac{s_v x}{1 + (r_v x)^2}
\B) \teq v_{M} + v_R. 
\eal
\eeq
For $x \geq L_{B} \geq 100 L$, the remainder terms $\om_R, v_R$ are estimated using Lemma \ref{lem:Fa_power} and the rescaled version of Lemma \ref{lem:g_power}, i.e. \eqref{eq:g_power_rescale}. Notice that we have 
\[
\f{s_{\om} x}{1 + (r_{\om} x)^2} = \f{s_{\om}}{ r_{\om}} G( r_{\om} x)
 \teq G_{ s_{\om} / r_{\om}, r_{\om}}(x)
, \quad
\f{s_{v} x}{1 + (r_{v} x)^2} = \f{s_v }{ r_v } G( r_{v } x), \quad G = \f{x}{1 + x^2}.
\]

For the velocity, we have 
\[
u = u_p + u_b = u_p + b_{\om} u( F_{a_{\om}})
+  u_{g, s_{\om}/ r_{\om}, r_{\om} },
\]
where $u_p, b_{\om} u( F_{a_{\om}}),  u_{g, s_{\om}/ r_{\om}, r_{\om} }$ are the velocities associated to different parts in $\om$ \eqref{eq:solu_decomp}, respectively. The slowest decay part of $u$ is from $b_{\om} u(F_{a_{\om}})$ and we further decompose $u$ as follows 
\beq\label{eq:decomp_up}
u = u_{M} + (  b_{\om} u(\bar F_{a_{\om}} )- u_M + u_p 
+  u_{g, s_{\om}/ r_{\om}, r_{\om} }  )
\teq u_M + u_R,  \quad
u_M = - b_{\om} \cot \f{\pi a_{\om}}{2} \f{ x^{1- a_{\om}}}{1 - a_{\om}} .
\eeq
We remark that $- \sgn(x) \cot \f{\pi a_{\om}}{2} \f{ |x|^{1- a_{\om}}}{1 - a_{\om}}$ is the velocity associated to $ \bar F_{a_{\om}}  = \sgn(x) |x|^{- a_{\om}}$. See Corollary \ref{cor:ub}. The first term in $u_R$, i.e. $b_{\om} u(\bar F_{a_{\om}} )- u_M$, is estimated using Corollary \ref{cor:ub}, the second term $u_p$ using \eqref{eq:up_far} and \eqref{eq:up_far_err}, and the third term 
$ u_{g, s_{\om}/ r_{\om}, r_{\om} } $ using the rescaled version of Lemma \ref{lem:g_power}, i.e. \eqref{eq:g_power_rescale}.

With the above preparations, we are in a position to estimate $ (x\pa_x)^i F_{\om}, (x\pa_x)^i F_v$. We seek to establish the following bounds 
\beq\label{eq:Err_decay}
\bal
| (x \pa_x)^i F_{\om} | &\leq C_{F_{\om}}(i, 1) x^{-a_{\om}}
+  C_{F_{\om}}(i, 2) x^{-2 a_{\om}},    \\
| (x \pa_x)^i F_{v} | &\leq C_{F_{v}}(i, 1) x^{-2 a_{\om}}
+  C_{F_{v}}(i, 2) x^{-3 a_{\om}} ,
\eal
\eeq
for some constants $C_{F_{\om}}(i, j), C_{F_v}(i,j)$ to be determined. The asymptotic behavior is characterized by the slowest decaying parts in the approximate steady state $\om_M, v_M, u_M$. Recall 
\[
F_{\om} = - (\bar c_l x + u ) \om_x + \bar c_{\om} \om + v , 
F_v = (- \bar c_l x + u ) v_x + (2 \bar c_{\om} - u_x) v. 
\]
We focus on $F_{\om}$ and $x \geq 0$. Consider the following decomposition 
\[
\bal
 - \bar c_l x \om_x + \bar c_{\om} \om
 &= ( - \bar c_l x \om_{M,x} + \bar c_{\om} \om_{M})
 + (- \bar c_l x \om_{R,x} + \bar c_{\om} \om_{R} ) \teq I_M + I_R , \\
- u \om_x + v 
&= ( - u_M \om_{M, x} + v_M  ) 
+ (- u_M \om_{R, x} + v_R - u_R \om_x  )
\teq II_M + II_R.
 \eal
\]

Denote $C_{a_{\om}}= - \cot \f{\pi a_{\om}}{2} \f{1}{1- a_{\om}}$. From \eqref{eq:decomp_up}, we get $u_M = b_{\om} C_{a_{\om}} x^{1-a_{\om}}$. Using the definitions \eqref{eq:decomp_wv} and \eqref{eq:decomp_up}, we obtain 
\[
\bal
(x \pa_x)^i I_M &=b_{\om} (x\pa_x)^i 
( \bar c_l a_{\om} x^{-a_{\om}} + \bar c_{\om} x^{-a_{\om}})
= b_{\om} ( \bar c_l a_{\om}  + \bar c_{\om}  )  (x\pa_x)^i x^{-a_{\om}} \\
&= b_{\om} ( \bar c_l a_{\om}  + \bar c_{\om}  )  (-a_{\om})^i x^{-a_{\om}}  , \\
(x \pa_x)^i II_M &
=(x \pa_x)^i
 ( - C_{a_{\om}} b_{\om} x^{1 - a_{\om}} \cdot b_{\om} (-a_{\om}) x^{-1 - a_{\om} }
+ b_v x^{-2 a_{\om}}) \\
& = ( C_{a_{\om}} b_{\om}^2 a_{\om} + b_v) 
 (x \pa_x)^i x^{-2 a_{\om}}
 = ( C_{a_{\om}} b_{\om}^2 a_{\om} + b_v)  (-2 a_{\om})^{i} x^{-2 a_{\om}}. 
\eal
\]

For the remainder term $I_R, II_R$, we use the Leibniz rule and the triangle inequality to estimate them directly. For example, we have 
\[
| (x\pa_x) ( u_M \om_{R, x} ) |
\leq  | (x\pa_x) u_M|  \cdot  |\om_{R, x} | 
+ |u_M| \cdot  | (x\pa_x) \om_{R, x}|,
\]
and then apply the formula of $u_{M}$ and Lemma \ref{lem:Fa_power} to estimate them. In particular, we can obtain the bound 
\[
|(x\pa_x)^i I_R |+
|(x\pa_x)^i II_R |  \leq Err_i(x) = \sum_{j=1}^{m_i} c_{i,j} x^{ -\al_j}
\]
with power $\al_i \geq 1 $ and coefficients $c_{i,j}> 0$. 
To see this, for example, we use
\[
| x \pa_x ( b_{\om} F_{a_{\om}} - b_{\om} x^{-a_{\om}})| \leq 
|b_{\om}| x \cdot C(L_{B}, a_{\om}, 1 ) x^{-5 -2 a_{\om} - 1}
\leq |b_{\om}| C(L_{B}, a_{\om}, 1 ) x^{-5-2 a_{\om}}
\]
from Lemma \ref{lem:Fa_power} to bound the first term in $ x\pa_x \om_R$.

Since each $\al_i > 2 a_{\om}$, for $x \geq L_B$, we further get 
\[
Err_i(x) \leq \sum_{j=1}^{m_i} c_{i,j} x^{-2 a_{\om}} L_B^{ - ( \al_{j}  - 2 a_{\om} ) }
=  L_B^{2 a_{\om}} x^{-2 a_{\om}} Err_i( L_B).
\]

Thus, in practice, we only need to estimate the error at $x = L_B$ to get $Err_i(L_B)$ rather than tracking the coefficients $c_{i,j}$. Combining the above estimates, we obtain 
\[
\bal
|(x\pa_x)^i  F_{\om} |
&\leq 
|b_{\om} ( \bar c_l a_{\om}  + \bar c_{\om}  )  (a_{\om})^i| x^{- a_{\om}}
+ ( | ( C_{a_{\om}} b_{\om}^2 a_{\om} + b_v)  (2 a_{\om})^{i} | 
+ L_B^{2 a_{\om}} Err_i( L_B))  x^{-2 a_{\om}} \\
&\teq C_{F_{\om}}(i, 1) x^{- a_{\om}}
+ C_{F_{\om}}(i, 2) x^{- 2a_{\om}}
\eal
\]
for $i = 1,2$. This gives the constants in \eqref{eq:Err_decay}. Since we choose the factor $a_{\om}$ sufficiently close to $ - \f{\bar c_{\om}}{\bar c_l}$, $\bar c_l a_{\om}  + \bar c_{\om} $ 
is of order $10^{-8}$. The coefficient $C_{F_{\om}}(i, 2)$ is of order $1$. Due to the smallness 
in $C_{F_{\om}}(i, 1)$ and faster decay $x^{-2a_{\om}}$, the error is very small in the far field.

\appendix
\setcounter{section}{4}
\section{Estimates of the Hilbert transform and functional inequalities}

\subsection{Properties and estimates of the Hilbert transform}
\begin{lem}\label{lem:com}
Suppose that $f$ is in the Schwartz space. We have 
\[
H (fx) - H(fx)(0) = x Hf(x).
\]

\end{lem}

The following Lemma is proved in the main paper \cite{chen2021HL_supp}.
\begin{lem}\label{lem:wg}
Assume that $\om$ is odd and $ \om \in L^2( x^{-4} + x^{-2/3})$. Let $u_x = H \om$. For any $\al, \b, \g \geq 0$, we have 
\[
\bal
|| (u_x - u_x(0) ) ( \al x^{-4} + \b x^{-2})^{1/2} ||_2^2 &= || \om ( \al x^{-4} + \b x^{-2})^{1/2}||_2^2 \\
\B|\B| (u - u_x(0)  x ) \B( \f{\al}{x^6} + \f{\b}{x^4}  +  \f{\g}{x^{10/3}} \B)^{1/2} \B|\B|_2^2 & \leq \B|\B| \om \B( \f{4\al}{25 x^4}  + \f{ 4\b}{9 x^2} \B)^{1/2} \B|\B|_2^2  +  \f{36\g}{49}  || (u_x - u_x(0)) x^{-2/3}||_2^2. 
\eal
\]
\end{lem}

The following weighted estimate was established in \cite{Cor06_supp}.

\begin{lem}\label{lem:wg13}
For $x^{-1/3} f \in L^2$, we have 
\[
|| x^{-1/3} H f ||_2 \leq  \cot \f{\pi}{12} || x^{-1/3}f ||_2 = (2 + \sqrt{3})|| x^{-1/3}f ||_2 .
\]
\end{lem}

The following Lemma is proved in the arXiv version of the main paper \cite{chen2021HL_supp}.
\begin{lem}\label{lem:rota}
Assume that $\om \in L^2( |x|^{-4/3} + |x|^{-2/3})$ is odd and $u_x = H \om$. We have 
\[
\bal
\int_{\R} \f{ (u_x (x) - u_x(0))^2} {  |x|^{4/3} } =  \int_{\R} \B( \f{ w^2} { | x|^{4/3} }  +
2  \sqrt{3} \cdot \mathrm{sgn}(x) \f{  \om  ( u_x(x) - u_x(0) ) }{|x|^{4/3} } \B) dx .
\eal
\]
\end{lem}

The following Lemma is proved in the main paper \cite{chen2021HL_supp}.
\begin{lem}\label{lem:IBPu}
For odd $u$ with $u_x x^{- (p-1)} \in L^2, p \in [1,3]$, we have
\beq\label{eq:IBPu}
 || (\td{u} - \f{1}{ 2p  - 1} \td{u}_x x ) x^{- p}||_2^2 = \f{1}{(2p-1)^2}  \int_{\R_+} \f{\td{u}^2_x}{ x^{ 2p -2}}  dx.
\eeq
\end{lem}

\begin{lem}\label{lem:uxinf2}
Suppose that $f, f_x \in L^{\inf}([a, b]) $ and $x \in (a, b)$. We have 
\[
\bal
&|\int_a^b \f{ f(y)}{x-y} dy |
\leq || f||_{L^{\inf}[a,b]} ( \log|a-x|  + \log |b-x| - 2 \log \d )+ 2 \d || f_x||_{L^{\inf}[a,b]}, \\
&| \int_a^b \log |x- y| f(y) dy |
\leq || f||_{L^{\inf}} ( \int_0^{x-a} |\log y | dy + \int_0^{b-x} |\log y| dy ),
\eal
\]
where $ \d = \min( |x-a|, |x-b|,  \f{ || f||_{L^{\inf} } }{ || f_x||_{L^{\inf}[a,b]} } )$. The integral of $|\log y|$ is further given by 
\[
\int_0^z |\log y| dy = -\int_0^{z_0}  \log y dy + \int_{z_0}^z \log y dy 
= g(z) - 2 g(z_0)
\]
for $z >0$, where $z_0 = \min(z, 1)$ and $g(y) = \int_0^y \log y dy = y \log y - y$.

\end{lem}

\begin{proof}
Denote $A =  || f||_{L^{\inf}[a,b] },  B =  || f_x||_{L^{\inf}[a,b]}$, $\d_0 = \min( |x-a|, |x-b|)$. For some $ s \in (0, \d_0]$ to be determined, we rewrite the integral as follows 
\[
\int_a^b \f{ f(y)}{x-y} dy
= \int_{a-x}^{b-x } \f{ f( x+ z)}{-z} d z
= \int_{-s}^{s} \f{ f(x+z)}{-z} dz
+ \int_{ [a-x, b-x] \bsh [-s, s]} \f{ f( x+ z)}{-z} d z \teq I + II.
\]
We estimate $I, II$ as follows 
\[
\bal
|I | &= \B| \int_{-s }^{s} \f{ f(x+z) - f(x)}{-z} dz \B|
\leq B \int_{-s}^{s} 1 dz = 2 s B, \\
|II| &\leq A 
\int_{ [a-x, b-x] \bsh [-s, s]} \f{1}{|z|} d z 
= A ( \log |a-x| + \log |b-x| - 2 \log s).
\eal
\]
Optimizing the above estimates by choosing $s = \d$ concludes the proof of the first inequality. The proof of the second inequality is trivial and can be obtained as follows
\[
| \int_a^b \log |x- y| f(y) dy |
\leq || f||_{L^{\inf}}  \int_a^b | \log |x-y| | dy 
= || f||_{L^{\inf}} ( \int_0^{x-a} |\log y | dy + \int_0^{b-x} |\log y| dy ).
\]

The integral formula for $|\log y|$ is straightforward and we omit the proof.
\end{proof}

\subsection{Estimate of $\partial_x^3u$} \label{sec:ux3_estimate}

In our numerical verification procedure, we need to obtain estimates of $u$ and its derivatives on each small interval $[x_{i-1},x_i]$ for the mesh $\{x_i\}_{i=1}^N$ over domain $[0,L_B]$. In particular, we need to obtain $(\partial_x^ku)^{up/low}$ for $k = 0,1,2,3$ (see the definition of $f^{up/low}$ for a function $f$ in Section \ref{sec:maxmin_estimate}). We will start from obtaining $(\partial_x^3u)^{up/low}$, and then we use the method introduced in Section \ref{sec:maxmin_estimate} to obtain $(\partial_x^ku)^{up/low}$ for $k = 0,1,2$, given that we have computed the grid point values of these functions as in Sections \ref{sec:compute_up} and \ref{sec:compute_ub}.

To obtain $u_{xxx}^{up/low}$ means to estimate the maximum and minimum of $\partial_x^3u$ on each mesh interval $[x_{i-1},x_i]$. We will achieve this by calculating
\begin{align*}
\partial_x^3u(x) = \partial_x^3u(x_{i-1}) + \int_{x_{i-1}}^x \partial_x^4u(y) dy = \partial_x^3u(x_i) + \int_x^{x_i} \partial_x^4u(y) dy, \quad x\in[x_{i-1}, x_i].
\end{align*}
Thus we have the estimate
\[\max\{\partial_x^3u(x_{i-1}),\partial_x^3u(x_i)\} - I_i \leq \partial_x^3u(x)\leq \min\{\partial_x^3u(x_{i-1}),\partial_x^3u(x_i)\} + I_i, \quad x\in[x_{i-1}, x_i],\]
where 
\[I_i \teq \int_{x_{i-1}}^{x_i} |\partial_x^4u(y)| dy.\]
We then need to bound $I_i$ for each mesh interval $[x_{i-1}, x_i]$. 

For a relatively small $x_i$ such that the interval size $h_i = x_i - x_{i-1}$ is sufficiently small, we bound $I_i$ as 
\[
I_i \leq \sqrt{h_i} \left(\int_{x_{i-1}}^{x_i} |\partial_x^4u(y)|^2 dy\right)^{1/2} \leq \sqrt{h_i} \|\partial_x^4u\|_2= \sqrt{h_i}\|\partial_x^3\om\|_2\teq I_{i,1}.
\]
The first equality above uses the fact that $\partial_x^4u$ is odd in $x$, and the second equality uses the isometry $\|\pa_x^4 u \|_2 = \| \pa_x^3 \om \|_2$ that relies on the relation $\pa_x^4 u = H (\pa_x^3 \om)$ and the fact that $\om \in C^{2, 1}$. We can compute $\|\partial_x^3\om\|_2$ sufficiently accurately using the method in Section \ref{sec:integral_estimate}.

For a larger $x_i$, the interval size $h_i = x_i - x_{i-1}$ might be too large that the estimate described above is too crude for our analysis. Instead, we need to make use of the decay property of $\partial_x^4u$ in the far field. To do so, we first apply \ref{lem:wg13} to obtain
\begin{align*}
 \| x^{3-1/3} \partial_x^4u \|_2 = \| x^{3-1/3} H(\partial_x^3\om) \|_2 = \| x^{-1/3} H(x^3\partial_x^3\om) \|_2 \leq (2 + \sqrt 3) \| x^{3-1/3} \partial_x^3\om \|_2.
\end{align*}
The first equality above uses $\partial_x^4 u = H(\partial_x^3\om)$; the second equality uses Lemma \ref{lem:com} repeatedly and the fact that 
\[H(x^k\partial_x^3\om)(0) = -\frac{1}{\pi}\int_{\mathbb{R}} x^{k-1}\partial_x^3\om(x)dx = 0 \quad \text{for $k = 1,2,3$}.\]
Then we can bound $I_i$ as follows
\begin{align*}
I_i &\leq \left(\int_{x_{i-1}}^{x_i}y^{-16/3}dy\right)^{1/2}\left(\int_{x_{i-1}}^{x_i} y^{16/3}|\partial_x^4u(y)|^2 dy\right)^{1/2} \\
&\leq \left(\frac{x_{i-1}^{-13/3}-x_i^{-13/3}}{13/3}\right)^{1/2} \| x^{3-1/3} \partial_x^4u \|_2 \\
&\leq  (2 + \sqrt 3)\left(\frac{x_{i-1}^{-13/3}-x_i^{-13/3}}{13/3}\right)^{1/2} \| x^{3-1/3} \partial_x^3\om \|_2\teq I_{i,2}.
\end{align*}
The smallness of $I_{i,2}$ relies on the decay of $x^{-13/3}$ in the far field. Again, we can compute $\|x^{3-1/3} \partial_x^3\om\|_2$ accurately using the method in Section \ref{sec:integral_estimate}.

In summary, we can bound $I_i$ as 
\[I_i \leq \min\{I_{i,1}\,,\, I_{i,2}\},\]
where $I_{i,1}$ and $I_{i,2}$ can both be estimated rigorously given the information of $\om$. Then we construct $(\partial_x^3u)^{up/low}$ as
\[(\partial_x^3u)_i^{up} = \min\{\partial_x^3u(x_{i-1}),\partial_x^3u(x_i)\} + I_i,\quad (\partial_x^3u)_i^{low} = \max\{\partial_x^3u(x_{i-1}),\partial_x^3u(x_i)\} - I_i\]
for $i=1, \cdots, N$ over the whole mesh.

\subsection{ Error estimates of the integral}\label{supp:integral_err}

In this section, we prove Lemmas \ref{lem:Trule1}, \ref{lem:Trule2} and \ref{lem:int_high}.

\begin{proof}[Proof of Lemma \ref{lem:Trule1}]

Denote by $P$ a piecewise quadratic polynomial
\[
P(x) \teq \f{(x- x_i)(x-  x_{i+1} )}{2} , \ \forall x \in [x_i, x_{i+1}], \ i=0, 1,2.., n-1.
\]
Clearly, we have 
\beq\label{eq:trap_P}
P(x) \leq 0, \quad |P(x)| \leq \f{ (x_{i+1}-x_i)^2}{8} \leq \f{h^2}{8}.
\eeq

Using \eqref{eq:Trapezoidal_error}, we get 
\[
\bal
\D_1 \teq &\int_0^M  \f{F^2}{x^k} dx  - T_h \lt(  \f{F^2}{x^k} ,0, M  \rt)  = \int_0^{M}
 \lt(\f{F^2}{x^k}\rt)_{xx}  P(x) dx ,\\
\lt(\f{F^2}{x^k}\rt)_{xx} =& \f{2F_{xx} F + 2F^2_x}{x^k} - 4 k \f{F_x F}{x^{k+1}} + k(k+1) \f{F^2}{x^{k+2}}  \\
=& 2 \lt(\f{F_{xx}}{x^{k/2-1}} - (k-1) \f{F_x}{x^{k/2}} \rt) \f{F}{x^{k/2+1}} +  2 \f{F_x^2}{x^k} - (2k+2) \f{F_x F}{x^{k+1}} + k(k+1) \f{F^2}{x^{k+2}} .
\eal
\]
For $k>1$, we have $ 2 a^2 + k(k+1) b^2 \geq 2 \sqrt{2 k (k+1)} |ab| \geq (2k+2) |ab|$. Thus, we yield 
\[
\lt(\f{F^2}{x^k}\rt)_{xx} \geq 2 \lt(\f{F_{xx}}{x^{k/2-1}} - (k-1) \f{F_x}{x^{k/2}} \rt) \f{F}{x^{k/2+1}} \teq 2 J  \f{F}{x^{k/2+1}} ,
\]
where $J =\f{F_{xx}}{x^{k/2-1}} - (k-1) \f{F_x}{x^{k/2}} $. Using  \eqref{eq:trap_P} and the Cauchy-Schwarz inequality, we obtain 
\[
\D_1 \leq \int_0^M 
2 J \f{F}{x^{k/2+1}} P(x) dx 
\leq \f{h^2}{4} \int_0^M \B| J \f{F}{x^{k/2+1}} \B| dx 
\leq \f{h^2}{4} (\int_0^M J^2 dx )^{1/2}
( \int_0^M \f{F^2}{x^{k + 2}} dx )^{1/2}
\teq I_1^{1/2} I_2^{1/2}.
\]

For $I_1$, expanding the square 
\[
I_1 = \int_0^M J^2 d x =  \int_0^M \f{F_{xx}^2}{x^{k-2}}
- 2(k-1) \f{F_{xx} F_x}{ x^{k - 1}}
+ (k-1)^2 \f{ F_x^2}{x^k} dx \teq I_{1,1} + I_{1,2} + I_{1,3},
\]
and then using integration by parts and $k > 1$, 
\[
\bal
I_{1,2} &= -(k-1) \int_0^M  \pa_x (F_x^2) x^{- (k-1)} dx
= (k-1) \int_0^M F_x^2 \pa_x x^{-(k-1)} dx
- (k-1) F_x^2 x^{-(k-1)} \B|_0^M  \\
&\leq (k-1) \int_0^M F_x^2 \pa_x x^{-(k-1)} dx
= - (k-1)^2 \int_0^M F_x^2 x^{-k} dx = - I_{1,3},
\eal 
\]
 we derive 
 \[
I_1 \leq I_{1,1}.
 \]

For $I_2$, applying the Hardy-type inequality in Lemma \ref{lem:hardy}, we yield 
\[
I_2  \leq \f{4}{ (k+1)^2} \int_0^M \f{F_x^2}{x^k} dx
\leq \f{16}{ (k+1)^2(k-1)^2} \int_0^M \f{F_{xx}^2}{x^{k-2}} dx.
\]
Combining the estimate of $I_1, I_2$, we prove the first inequality in the Lemma. The second inequality is trivial.
\end{proof}

Next, we prove Lemma \ref{lem:Trule2}.

\begin{proof}[Proof of Lemma \ref{lem:Trule2}]
Let $I$ denote the left-hand side of the desired inequality. Using the first line in \eqref{eq:Trapezoidal_error}, we have
\begin{align*}
I &\leq \frac{h}{2}\int_0^M\left| \partial_x\left(\frac{f^2}{x^{k+1}}\right)\right|dx = \frac{h}{2}\int_0^M\left| \frac{2f_xf}{x^{k+1}} - \frac{(k+1)f^2}{x^{k+2}}\right|dx  \\
&\leq \frac{h}{2}\left(\int_0^M\left(\frac{2f_x}{x^{k/2}} - \frac{(k+1)f}{x^{k/2+1}}\right)^2dx\right)^{1/2}\left(\int_0^M\frac{f^2}{x^{k+2}}dx\right)^{1/2} \teq \frac{h}{2}I_1^{1/2} I_2^{1/2}.
\end{align*}
Using integration by parts, we obtain 
\begin{align*}
I_1 &= \int_0^M\left(\frac{4f_x^2}{x^k} - \frac{4(k+1)f_xf}{x^{k+1}} + \frac{(k+1)^2f^2}{x^{k+2}}\right)dx\\
&= \int_0^M\left(\frac{4f_x^2}{x^k} - \frac{2(k+1)^2f^2}{x^{k+2}} + \frac{(k+1)^2f^2}{x^{k+2}}\right)dx - \frac{2(k+1)f^2}{x^{k+1}}\Big|_0^M\\
&\leq \int_0^M\left(\frac{4f_x^2}{x^k} - \frac{(k+1)^2f^2}{x^{k+2}}\right)dx.
\end{align*}
Therefore, we get
\[I_1^{1/2} I_2^{1/2} \leq \frac{1}{k+1}\frac{I_1 + (k+1)^2I_2}{2} \leq \frac{2}{k+1} \int_0^M\frac{f_x^2}{x^k}dx,\]
which then yields the desired inequality. 
\end{proof}

Next, we prove Lemma \ref{lem:int_high}.

\begin{proof}[Proof of Lemma \ref{lem:int_high}]
We only need to prove  \eqref{eq:int_err2}. A direct calculation yields 
\[
\bal
e(x)  & = \int_a^x f^{\prime}( y) dy - \f{x-a}{b-a} \int_a^b f^{\prime}(y) dy
= \f{b-x}{b-a} \int_a^x f^{\prime}(y) dy - \f{x-a}{b-a} \int_x^b f^{\prime}(y) dy \\
& = \f{b-x}{b-a} \int_a^x ( f^{\prime}(y) - f^{\prime}(x) ) dy - \f{x-a}{b-a} \int_x^b ( f^{\prime}(y) - f^{\prime}(x) )dy \\
& = - \f{b-x}{b-a} \int_a^x \int_y^x f^{\prime \prime}(z) d z  dy - \f{x-a}{b-a}  \int_x^b \int_x^y f^{\prime \prime}(z)  dz  dy \\
&= - \f{b-x}{b-a} \int_a^x (z-a)  f^{\prime \prime}(z) d z  - \f{x-a}{b-a}  \int_x^b (b-z) f^{\prime \prime}(z)  dz .
\eal
\]
Using the Cauchy--Schwarz inequality, we get 
\[
|e(x) | \leq  \B(\int_a^b f_{xx}^2(x) dx \B)^{1/2} 
\B( \f{ (b-x)^2}{(b-a)^2} \int_a^x  (z-a)^2 dz 
+ \f{ (x-a)^2 }{ (b-a)^2 } \int_x^b (b-z)^2 dz   \B)^{1/2}.
\]
Evaluating the integral, we prove \eqref{eq:int_err2}.
\end{proof}

\subsection{Other Lemmas}

\begin{lem}[Hardy-type inequality]\label{lem:hardy}
Suppose that $f x^{-p/2}, f_x x^{1-p/2} \in L^2(0,M)$. For $p > 1$ and $M > 0$, we have 
\[
\int_0^{ M} \f{ f^2 }{x^{ p}} dx \leq \f{4}{(p-1)^2}  \int_0^{ M}  \f{f_x^2 }{x^{p-2}} dx .
\]
\end{lem}

\begin{proof}
The proof follows by expanding the square 
\[
0 \leq \int_0^M  (  \f{f_x}{ x^{ p/2 - 1} }  - \f{(p-1)}{2} \f{ f }{x^{p/2}} )^2 dx ,
\]
and then applying integration by parts to the cross term 
\[
\bal
 &- (p-1) \int_0^M \f{f_x f}{x^{p-1}} dx
 = - \f{p-1}{2} \int_0^M \f{ \pa_x f^2} {x^{p-1}} dx  \\
 =& - \f{ (p-1)^2}{2} \int_0^M \f{f^2}{x^{p}} dx
 - \f{p-1}{2} \f{f^2}{x^{p-1}} \B|_0^M \leq   - \f{ (p-1)^2}{2} \int_0^M \f{f^2}{x^{p}} dx.
 \eal
\]
\end{proof}

\begin{lem}[Linear interpolation]\label{lem:lin_interp}
Assume that $f \in C[a, b] \cap C^2[a,b]$. We have 
\[
|| f||_{L^{\inf}} \leq \max( |f(a)|, |f(b)|) + \f{(b-a)^2}{8} || f_{xx}||_{\inf}.
\]

\end{lem}

\begin{proof}
The standard linear (or Lagrangian) interpolation theory states that 
\[
f(x) = f(a) + \f{f(b)-f(a)}{b-a}(x-a)+ \f{ (x-a)(x-b) }{2} f^{\prime \prime}(\xi(x))
\teq \hat f(x)  + \f{ (x-a)(x-b) }{2} f^{\prime \prime}(\xi(x))
\]
for some $\xi(x) \in [a,b]$. Since $ |\hat f| \leq \max( |f(a)| , |f(b)|)$ and $|(x-a)(x-b)| \leq \f{(b-a)^2}{4}$ for any $x \in [a, b]$, using the triangle inequality, we prove the desired result. 
\end{proof}

Recall that the kernels associated to $ d_x^k u$ are the following 
\beq\label{eq:ker}
K_1(x, y) = \log \B|\f{x-y}{x+y}\B|, \quad K_2(x, y) = \f{1}{x-y} - \f{1}{x+y} , 
\quad 
K_3(x, y)=  \f{1}{x-y} + \f{1}{x+y} .
\eeq

We have the following simple estimates when $ \f{x}{y} $ is small.
\begin{lem}\label{lem:ker_xsm}
Suppose that $ \f{x}{y} < \f{1}{1000}$ and $\e_2 = 0.02$. Then we have 
\[
\bal
|K_1 | &\leq 2 (1 + \e_2) \f{x}{y} ,\quad | K_1 + \f{2 x}{y} | \leq \f{2(1 + \e_2)}{3} \f{x^3}{y^3},
\quad  |K_2| \leq 2 (1 + \e_2) y^{-1}, \quad 
|K_2 + \f{2}{y}| \leq 2 (1 + \e_2) \f{x^2}{y^3},  \\
|K_3| & \leq 2 (1 + \e_2)\f{x}{y^2}, 
\quad |K_3 + \f{2x}{y^2}| \leq 2 (1 + \e_2)\f{x^3}{y^4} .
\eal
\]

\end{lem}

\begin{proof}
For $ y \geq 1000 x$, similar to \eqref{eq:ker_TL1}, we get 
\[
K_1(x, y) = - 2 \sum_{i \geq 1} \f{1}{2i-1} (\f{x}{y})^{2i-1}.
\]
It follows 
\[
\bal
|K_1(x, y)|& \leq 2 \sum_{i\geq 1}  (\f{x}{y})^{2i-1} 
= \f{2x}{y} \f{1}{ 1 - (x/y)^2} \leq 2 (1 +\e_2) \f{x}{y}, \\
 | K_1(x, y) + \f{2x}{y}|
&\leq \f{2}{3} \sum_{i \geq 2}  (\f{x}{y})^{2i-1} 
= \f{2x^3}{3 y^3} \f{1}{1 - (x/y)^2} < \f{2}{3} (1 + \e_2) \f{x^3}{y^3}.
\eal
\]

For $K_2$, a direct calcultion yields 
\[
\bal
|K_2 | & = | \f{2y}{x^2 - y^2} | = \f{2}{y} \f{1}{1 - (x/y)^2} \leq \f{2 (1+\e_2)}{y}, \\
| K_2 + \f{2}{y}| &= | \f{2y}{x^2 - y^2} + \f{2}{y}| = \f{2 x^2}{ y(y^2 - x^2)} 
= \f{2x^2}{y^3}\f{1}{1 - (x/y)^2} \leq 2(1+\e_2) \f{x^2}{y^3}. \\
\eal
\]
The proof of $K_3$ is completely similar.
\end{proof}

\bibliographystyle{abbrv}